\documentclass[10.99pt]{article}
\usepackage{epic,latexsym,amssymb}
\usepackage{color}
\usepackage{tikz}
\usepackage{comment}
\usepackage{lineno}
\usepackage{pstricks,pst-node,pst-plot}
\usetikzlibrary{patterns}
\usepackage{cite}
\usepackage{geometry,graphicx}

\usetikzlibrary{calc}

\usepackage{amsmath}

\newtheorem{theorem}{Theorem}
\newtheorem{lemma}{Lemma}

\newtheorem{fact}{Fact}
\newtheorem{claim}{Claim}

\newtheorem{proposition}{Proposition}
\newtheorem{conjecture}{Conjecture}

\newtheorem{subclaim}{Claim}[claim]
\newtheorem{subsubclaim}{Claim}[subclaim]

\newcommand{\QED}{$\Box$}
\newcommand{\cp}{\mathbin{\Box}}
\newcommand{\cB}{\mathcal{B}}

\newcommand{\smallqed}{{\tiny ($\Box$)}}
\newcommand{\w}{{\rm w}}

\newcommand{\cO}{{\cal O}}

\newcommand{\bb}[1]{{\color{red} #1}}

\newcommand{\tc}{{\rm tc}}

\newcommand{\proof}{\noindent\textbf{Proof. }}

\let\oldenumerate\enumerate
\renewcommand{\enumerate}{
  \oldenumerate
  \setlength{\itemsep}{0pt}
  \setlength{\parskip}{0pt}
  \setlength{\parsep}{0pt}
}

\begin{document}


\title{A proof of the $\frac{3}{8}$-conjecture for \\ independent domination in cubic graphs}

\author{Bo\v{s}tjan Bre\v{s}ar$^{a,b}$, Tanja Dravec$^{a,b}$, \, and \, Michael A. Henning$^{c}$\\ \\
$^a$ Faculty of Natural Sciences and Mathematics \\
University of Maribor, Slovenia \\
$^b$ Institute of Mathematics, Physics and Mechanics,
Ljubljana, Slovenia\\
\small \tt Email: bostjan.bresar@um.si \\
\small \tt Email: tanja.dravec@um.si \\
\\
$^{c}$ Department of Mathematics and Applied Mathematics \\
University of Johannesburg, South Africa\\
\small \tt Email: mahenning@uj.ac.za
}

\date{}
\maketitle

\begin{abstract}
A set $S$ of vertices in a graph $G$ is a dominating set of $G$ if every vertex not in $S$ is adjacent to a vertex in~$S$. An independent dominating set in $G$ is a dominating set of $G$ with the additional property that it is an independent set. The domination number, $\gamma(G)$, and the independent domination number, $i(G)$, are the minimum cardinalities among all dominating sets and independent dominating sets in $G$, respectively. By definition, $\gamma(G) \le i(G)$ for all graphs $G$. Let $G$ be a connected cubic graph of order~$n$. In 1996 Reed [Combin.\ Probab.\ Comput.\ 5 (1996), 277--295] proved a breakthrough result that $\gamma(G) \le \frac{3}{8}n$. We prove the stronger result that if $G$ is different from $K_{3,3}$ and the $5$-prism $C_5 \, \Box \, K_2$, then $i(G) \le \frac{3}{8}n$. This proves a known conjecture. The bound is tight in the sense that there are infinite families of connected cubic graphs that achieve equality in this bound.
\end{abstract}


\section{Introduction}

A set $S$ of vertices in a graph $G$ is a \emph{dominating set} if every vertex in $V(G) \setminus S$ is adjacent to a vertex in~$S$. The \emph{domination number} of $G$, denoted by $\gamma(G)$, is the minimum cardinality of a dominating set.  A set is \emph{independent} in $G$ if no two vertices in the set are adjacent. An \emph{independent dominating set}, abbreviated ID-set, of $G$ is a set that is both dominating and independent in $G$. Equivalently, an ID-set of $G$ is a maximal independent set in $G$. The \emph{independent domination number} of $G$, denoted by~$i(G)$, is the minimum cardinality of an ID-set. An ID-set of cardinality $i(G)$ in $G$ is called an $i$-set of $G$. Independent dominating sets have been studied extensively in the literature (see, for  example,~\cite{GoHe-23,GoLy-12,KnSkTe-21,KuRa-23,Ly-14,OWest,SuWa-99}). A thorough treatise on dominating sets can be found in the so-called ``domination books''~\cite{HaHeHe-20,HaHeHe-21,HaHeHe-23}. A survey on independent domination in graphs can be found in~\cite{GoHe-13}.

For graph theory notation and terminology, we generally follow~\cite{HaHeHe-23}.  Specifically, let $G$ be a graph with vertex set $V(G)$ and edge set $E(G)$, and of order $n(G) = |V(G)|$ and size $m(G) = |E(G)|$. Two adjacent vertices in $G$ are \emph{neighbors}. The \emph{open neighborhood} of a vertex $v$ in $G$ is $N_G(v) = \{u \in V \, \colon \, uv \in E\}$ and the \emph{closed neighborhood of $v$} is $N_G[v] = \{v\} \cup N_G(v)$. We denote the degree of $v$ in $G$ by $\deg_G(v)$, and so $\deg_G(v) = |N_G(v)|$. An \emph{isolated vertex} in $G$ is a vertex of degree~$0$ in $G$. For $k \ge 1$ an integer, we let $[k]$ denote the set $\{1,\ldots,k\}$ and we let $[k]_0 = [k] \cup \{0\} = \{0,1,\ldots,k\}$.

A \emph{cycle} on $n$ vertices is denoted by $C_n$ and a \emph{path} on $n$ vertices by $P_n$. The  \emph{complete graph} on $n$ vertices is denoted by $K_n$, while the \emph{complete bipartite graph} with one partite set of size~$n$ and the other of size~$m$ is denoted by $K_{n,m}$. A \emph{cubic graph} is graph in which every vertex has degree~$3$, while a \emph{subcubic graph} is a graph with maximum degree at most~$3$. If $G$ is a subcubic graph, then we let $V_j(G)$ denote the set of vertices of degree~$j$ in $G$ and we let $n_j(G) = |V_j(G)|$ denote the number of vertices of degree~$j$ in $G$ for $j \in [3]_0$.

For subsets $X$ and $Y$ of vertices of $G$, we denote the set of edges with one end in $X$ and the other end in $Y$ by $[X,Y]$. The \emph{girth} of $G$, denoted $g(G)$, is the length of a shortest cycle in $G$. A component of a graph $G$ isomorphic to a graph $F$ we call an $F$-\emph{component} of $G$. For a set $S \subseteq V(G)$, the subgraph induced by $S$ is denoted by $G[S]$. Further, the subgraph of $G$ obtained from $G$ by deleting all vertices in $S$ and all edges incident with vertices in $S$ is denoted by $G - S$; that is, $G - S = G[V(G)\setminus S]$.

\section{Motivation and known results}
\label{S:known}

In 1996 Reed~\cite{Re-96} proved the following breakthrough result establishing a best possible upper bound on the domination number of a cubic graph. 

\begin{theorem}{\rm (\hskip-0.5pt \cite{Re-96})}
\label{thm:Reed}
If $G$ is a cubic graph of order~$n$, then $\gamma(G) \le \frac{3}{8}n$.
\end{theorem}

The importance of Reed's paper is reflected in over 100 citations one can find on the MathSciNet. Subsequently, determining best possible bounds on the independent domination number of cubic graphs attracted much interest (see, for example,~\cite{AkAk-22,BrHe-19,ChChKwPa-23,ChChPa-23,ChKim-23,DoHeMoSo-12,HHLS-12,HeLoRa-14}). In 1999 Lam, Shiu, and Sun~\cite{LSS-99} established the following upper bound on the independent domination number of a connected cubic graph.


\begin{theorem}{\rm (\hskip-0.5pt \cite{LSS-99})}
\label{thm:Lam}
If $G \not\cong K_{3,3}$ is a connected, cubic graph of order~$n$, then $i(G)   \le \frac{2}{5}n$.
\end{theorem}

\begin{figure}[htb]
\begin{center}
\begin{tikzpicture}
[thick,scale=0.75,
	vertex/.style={circle,draw,inner sep=0pt,minimum size=1.5mm},
	vertexlabel/.style={circle,draw=white,inner sep=0pt,minimum size=0mm}]
\def \irad {.7}
\def \orad {1.5}
\def \ydist {1}
\coordinate (A) at (0,0);
\draw (A)
{	
	node [fill=white] [vertex] at (-1.5,-\ydist) (x1){}
	node [fill=white][vertex] at (0,-\ydist) (x2){}
	node [fill=white][vertex] at (1.5,-\ydist) (x3){}
	node [fill=white][vertex] at (-.75,\ydist) (y1){}
	node [fill=white][vertex] at (.75,\ydist) (y2){}
	(y1)--(x1)--(y2)--(x2)--(y1)--(x3)--(y2)
	(A) +(-0.25,-2) node [rectangle, draw=white!0, fill=white!100]{\small (a) $K_{2,3}$}
};
\coordinate (B) at (6,0);

\foreach \i in {0,1,2,3,4}
{
	\path (B) node [fill=white][vertex] at +(90 - \i*360/5:\irad) (x\i){};
	\path (B) node [fill=white][vertex] at +(90 - \i*360/5:\orad) (y\i){};
	\draw (x\i)--(y\i);
}
	\draw (x0)--(x1)--(x2)--(x3)--(x4)--(x0);
	\draw (y0)--(y1)--(y2)--(y3)--(y4)--(y0);

\draw (B)
{	
	+(0,-2) node [rectangle, draw=white!0, fill=white!100]{\small (b) $C_5 \, \Box \, K_2$}
};
\end{tikzpicture}\end{center}
\vskip -0.6 cm \caption{The graphs $K_{2,3}$ and $C_5 \, \Box \, K_2$.} \label{f:K33}
\end{figure}
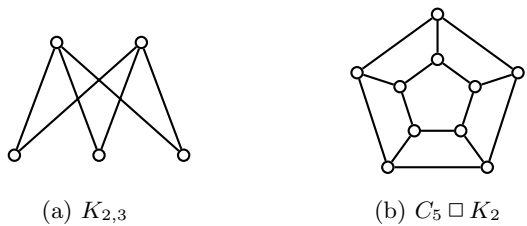

The $\frac{2}{5}$-upper bound on the independent domination number in Theorem~\ref{thm:Lam} is best possible in the sense that there exists at least one graph that achieves equality in the bound, namely the $5$-prism, $C_5 \cp K_2$, which is the Cartesian product of a $5$-cycle with a copy of $K_2$. The graphs $K_{2,3}$ and $C_5 \cp K_2$ are shown in Figure~\ref{f:K33}(a) and \ref{f:K33}(b), respectively.

However despite the pleasing result given in Theorem~\ref{thm:Lam}, prior to this paper it remained an open problem to characterize the graphs achieving equality in Theorem~\ref{thm:Lam}. In their 2013 survey paper on independent domination in graphs, Goddard and Henning~\cite{GoHe-13} conjectured that the $5$-prism, $C_5 \cp K_2$, is the only extremal graph achieving equality in the $\frac{2}{5}$-upper bound on the independent domination number $i(G)$ established by Lam, Shiu, and Sun.

The authors in~\cite{GoHe-13} conjectured that Reed's $\frac{3}{8}$-upper bound on the domination number, $\gamma(G)$, also holds for the independent domination, $i(G)$, in a connected, cubic graph, with the exception of the two graphs $K_{3,3}$ and $C_5 \, \Box \, K_2$.

\begin{conjecture}{\rm (\hskip-0.5pt \cite{GoHe-13})}
\label{conj:Idom-cubic}
If $G \notin \{K_{3,3}, C_5 \, \Box \, K_2\}$ is a connected cubic graph of order~$n$, then $i(G) \le \frac{3}{8}n$.
\end{conjecture}

In 2015 Dorbec et al.~\cite{DoHeMoSo-12} gave support to Conjecture~\ref{conj:Idom-cubic}
by proving that the conjecture holds if the cubic graph $G$ contains no subgraph isomorphic to $K_{2,3}$. Recall that in a subcubic graph $G$, we let $n_j(G)$ denote the number of vertices of degree~$j$ in $G$ for $j \in [3]_0$.

\begin{theorem}{\rm (\hskip-0.5pt \cite{DoHeMoSo-12})}
\label{thm:Dorbec1}
If $G$ is a subcubic graph that does not have a subgraph isomorphic to $K_{2,3}$ and which has no $(C_5 \cp K_2)$-component, then
\[
8 i(G) \le 8  n_0(G) + 5  n_1(G) + 4  n_2(G) + 3  n_3(G).
\]
\end{theorem}

As a consequence of Theorem~\ref{thm:Dorbec1}, we have the main result in~\cite{DoHeMoSo-12}, noting that if $G$ is a cubic graph of order~$n$, then $n_0(G) = n_1(G) = n_2(G)=0$ and $n = n_3(G)$.

\begin{theorem}{\rm (\hskip-0.5pt \cite{DoHeMoSo-12})}
\label{thm:Dorbec2}
If $G \ne C_5 \cp K_2$ is a connected cubic graph of order $n$ that does not have a subgraph isomorphic to $K_{2,3}$, then $i(G) \le \frac{3}{8}n$.
\end{theorem}

\section{Main result}

We give a complete proof of Conjecture~\ref{conj:Idom-cubic} by adding forbidden subgraphs isomorphic to $K_{2,3}$ in the statement of Theorem~\ref{thm:Dorbec2} back into the mix. We will prove the following result.

\begin{theorem}
\label{thm:main}
If $G \notin \{K_{3,3}, C_5 \, \Box \, K_2\}$ is a connected cubic graph of order~$n$, then $i(G) \le \frac{3}{8}n$.
\end{theorem}

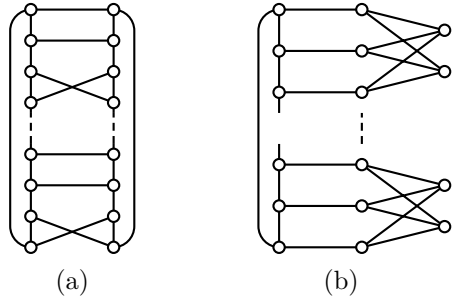
\begin{figure}[htb]
\begin{center}
\begin{tikzpicture}[thick,scale=.55]%
\tikzstyle{every node}=[circle, draw, fill=black!0, inner sep=0pt,
minimum width=.15cm]
  \draw [densely dashed] (0,0) {
    +(0,2.5) -- +(0,3.25) +(2,2.5) -- +(2,3.25)};
  \draw [rounded corners=8pt] (0,0) {
    +(0,0) -- +(-.5,0) -- +(-.5,5.75) -- +(0,5.75)
    +(2,0) -- +(2.5,0) -- +(2.5,5.75) -- +(2,5.75)};
  \draw (0,0) {
    +(0,0) -- +(2,.75) +(0,.75) -- +(2,0)
    +(0,1.5) -- +(2,1.5) +(0,2.25) -- +(2,2.25)};
  \draw (0,3.5) {
    +(0,0) -- +(2,.75) +(0,.75) -- +(2,0)
    +(0,1.5) -- +(2,1.5) +(0,2.25) -- +(2,2.25)};
  \draw (0,0) {
    +(0,0.00) node{} -- +(0,0.75) node{} -- +(0,1.5) node{} -- +(0,2.25) node{} -- +(0,2.5)
    +(0,3.25) -- +(0,3.5) node{} -- +(0,4.25) node{} -- +(0,5.0) node{} -- +(0,5.75) node{}};
  \draw (2,0) {
    +(0,0.00) node{} -- +(0,0.75) node{} -- +(0,1.5) node{} -- +(0,2.25) node{} -- +(0,2.5)
    +(0,3.25) -- +(0,3.5) node{} -- +(0,4.25) node{} -- +(0,5.0) node{} -- +(0,5.75) node{}};
  \node[rectangle, draw=white!0, fill=white!100] at (1,-0.85) {(a)};
  \draw [densely dashed] (8,0) {
    +(0,2.5) -- +(0,3.25)};
  \draw [rounded corners=8pt] (6,0) {
    +(0,0) -- +(-.5,0) -- +(-.5,5.75) -- +(0,5.75)};
  \draw (6,0) {
    +(2,2) -- +(4,.5)        +(2,2) -- +(4,1.5)
    +(2,1) -- +(4,.5)        +(2,1) -- +(4,1.5)
    +(2,0) -- +(4,.5) node{} +(2,0) -- +(4,1.5) node{}};
  \draw (6,3.75) {
    +(2,2) -- +(4,.5)        +(2,2) -- +(4,1.5)
    +(2,1) -- +(4,.5)        +(2,1) -- +(4,1.5)
    +(2,0) -- +(4,.5) node{} +(2,0) -- +(4,1.5) node{}};
  \draw (6,0) {
    +(0,0) -- +(0,2.5) +(0,3.25) -- +(0,5.75)};
  \draw (6,0) {
    +(0,0) node{} -- +(2,0) node{}
    +(0,1) node{} -- +(2,1) node{}
    +(0,2) node{} -- +(2,2) node{}
    +(0,3.75) node{} -- +(2,3.75) node{}
    +(0,4.75) node{} -- +(2,4.75) node{}
    +(0,5.75) node{} -- +(2,5.75) node{}};
  \node[rectangle, draw=white!0, fill=white!100] at (7.5,-0.85) {(b)};
\end{tikzpicture}
\caption{Two infinite families of cubic graphs $G$ of order~$n$ satisfying $i(G) = \frac{3}{8}n$}
\label{fig:cubic}
\end{center}
\end{figure}

As a consequence of Theorem~\ref{thm:main}, the $5$-prism, $C_5 \cp K_2$, is the unique extremal graph in Theorem~\ref{thm:Lam}. The $\frac{3}{8}$-upper bound on the independent domination number given in Theorem~\ref{thm:main} is tight. Two infinite families of connected cubic graphs $G$ satisfying $i(G) = \frac{3}{8}n$ are constructed in~\cite{GoHe-13}. These families are also defined in~\cite[Chapter~6]{HaHeHe-23}, so we do not redefine them here. However, graphs in these families are illustrated in Figure~\ref{fig:cubic}(a) and~(b), respectively. Akbari et al.~\cite{AkAk-22} constructed additional families (which we do not define here) of connected cubic graphs $G$ satisfying $i(G) = \frac{3}{8}n$.

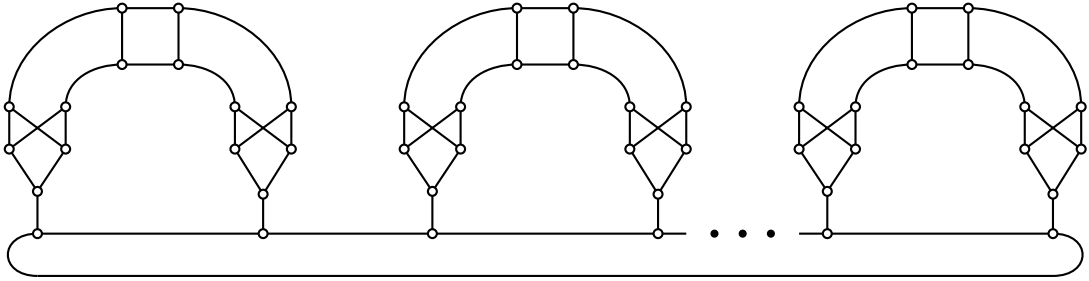
\begin{figure}[htb]
\begin{center}
\begin{tikzpicture}[scale=.75,style=thick,x=1cm,y=1cm]
\def\vr{2.25pt}
\def\vrn{1.5pt}
\path (0,1.5) coordinate (a2);
\path (0,2.25) coordinate (a4);
\path (0.5,0) coordinate (a);
\path (0.5,-0.75) coordinate (x);
\path (0.5,0.75) coordinate (a1);
\path (1,1.5) coordinate (a3);
\path (1,2.25) coordinate (a5);
\path (2,3) coordinate (b1);
\path (2,4) coordinate (b2);
\path (3,3) coordinate (b3);
\path (3,4) coordinate (b4);
\path (4,2.25) coordinate (c4);
\path (5,2.25) coordinate (c5);
\path (5,1.5) coordinate (c3);
\path (4,1.5) coordinate (c2);
\path (4.5,0.7) coordinate (c1);
\path (4.5,0) coordinate (c);
\draw (a)--(a1)--(a2)--(a4)--(a3)--(a5)--(a2);
\draw (a1)--(a3);
\draw (b1)--(b2)--(b4)--(b3)--(b1);
\draw (c)--(c1)--(c2)--(c4)--(c3)--(c5)--(c2);
\draw (c1)--(c3);
\draw (a4) to[out=90,in=180, distance=1cm] (b2);
\draw (a5) to[out=90,in=180, distance=0.5cm] (b1);
\draw (b3) to[out=0,in=90, distance=0.5cm] (c4);
\draw (b4) to[out=0,in=90, distance=1cm] (c5);
\path (7,1.5) coordinate (d2);
\path (7,2.25) coordinate (d4);
\path (7.5,0) coordinate (d);
\path (7.5,0.75) coordinate (d1);
\path (8,1.5) coordinate (d3);
\path (8,2.25) coordinate (d5);
\path (9,3) coordinate (e1);
\path (9,4) coordinate (e2);
\path (10,3) coordinate (e3);
\path (10,4) coordinate (e4);
\path (11,2.25) coordinate (f4);
\path (12,2.25) coordinate (f5);
\path (12,1.5) coordinate (f3);
\path (11,1.5) coordinate (f2);
\path (11.5,0.7) coordinate (f1);
\path (11.5,0) coordinate (f);
\draw (d)--(d1)--(d2)--(d4)--(d3)--(d5)--(d2);
\draw (d1)--(d3);
\draw (e1)--(e2)--(e4)--(e3)--(e1);
\draw (f)--(f1)--(f2)--(f4)--(f3)--(f5)--(f2);
\draw (f1)--(f3);
\draw (d4) to[out=90,in=180, distance=1cm] (e2);
\draw (d5) to[out=90,in=180, distance=0.5cm] (e1);
\draw (e3) to[out=0,in=90, distance=0.5cm] (f4);
\draw (e4) to[out=0,in=90, distance=1cm] (f5);
\path (14,1.5) coordinate (g2);
\path (14,2.25) coordinate (g4);
\path (14.5,0) coordinate (g);
\path (14.5,0.75) coordinate (g1);
\path (15,1.5) coordinate (g3);
\path (15,2.25) coordinate (g5);
\path (16,3) coordinate (h1);
\path (16,4) coordinate (h2);
\path (17,3) coordinate (h3);
\path (17,4) coordinate (h4);
\path (18,2.25) coordinate (i4);
\path (19,2.25) coordinate (i5);
\path (19,1.5) coordinate (i3);
\path (18,1.5) coordinate (i2);
\path (18.5,0.7) coordinate (i1);
\path (18.5,0) coordinate (i);
\path (18.5,-0.75) coordinate (y);
\path (12,0) coordinate (z1);
\path (12.5,0) coordinate (n1);
\path (13,0) coordinate (n2);
\path (13.5,0) coordinate (n3);
\path (14,0) coordinate (z2);
\draw (g)--(g1)--(g2)--(g4)--(g3)--(g5)--(g2);
\draw (g1)--(g3);
\draw (h1)--(h2)--(h4)--(h3)--(h1);
\draw (i)--(i1)--(i2)--(i4)--(i3)--(i5)--(i2);
\draw (i1)--(i3);
\draw (g4) to[out=90,in=180, distance=1cm] (h2);
\draw (g5) to[out=90,in=180, distance=0.5cm] (h1);
\draw (h3) to[out=0,in=90, distance=0.5cm] (i4);
\draw (h4) to[out=0,in=90, distance=1cm] (i5);
\draw (a)--(c)--(d)--(f)--(z1);
\draw (i)--(z2);
\draw (x)--(y);
\draw (a) to[out=180,in=180, distance=0.7cm] (x);
\draw (i) to[out=0,in=0, distance=0.7cm] (y);
\draw (a) [fill=white] circle (\vr);
\draw (a1) [fill=white] circle (\vr);
\draw (a2) [fill=white] circle (\vr);
\draw (a3) [fill=white] circle (\vr);
\draw (a4) [fill=white] circle (\vr);
\draw (a5) [fill=white] circle (\vr);
\draw (b1) [fill=white] circle (\vr);
\draw (b2) [fill=white] circle (\vr);
\draw (b3) [fill=white] circle (\vr);
\draw (b4) [fill=white] circle (\vr);
\draw (c) [fill=white] circle (\vr);
\draw (c1) [fill=white] circle (\vr);
\draw (c2) [fill=white] circle (\vr);
\draw (c3) [fill=white] circle (\vr);
\draw (c4) [fill=white] circle (\vr);
\draw (c5) [fill=white] circle (\vr);
\draw (d) [fill=white] circle (\vr);
\draw (d1) [fill=white] circle (\vr);
\draw (d2) [fill=white] circle (\vr);
\draw (d3) [fill=white] circle (\vr);
\draw (d4) [fill=white] circle (\vr);
\draw (d5) [fill=white] circle (\vr);
\draw (e1) [fill=white] circle (\vr);
\draw (e2) [fill=white] circle (\vr);
\draw (e3) [fill=white] circle (\vr);
\draw (e4) [fill=white] circle (\vr);
\draw (f) [fill=white] circle (\vr);
\draw (f1) [fill=white] circle (\vr);
\draw (f2) [fill=white] circle (\vr);
\draw (f3) [fill=white] circle (\vr);
\draw (f4) [fill=white] circle (\vr);
\draw (f5) [fill=white] circle (\vr);
\draw (g) [fill=white] circle (\vr);
\draw (g1) [fill=white] circle (\vr);
\draw (g2) [fill=white] circle (\vr);
\draw (g3) [fill=white] circle (\vr);
\draw (g4) [fill=white] circle (\vr);
\draw (g5) [fill=white] circle (\vr);
\draw (h1) [fill=white] circle (\vr);
\draw (h2) [fill=white] circle (\vr);
\draw (h3) [fill=white] circle (\vr);
\draw (h4) [fill=white] circle (\vr);
\draw (i) [fill=white] circle (\vr);
\draw (i1) [fill=white] circle (\vr);
\draw (i2) [fill=white] circle (\vr);
\draw (i3) [fill=white] circle (\vr);
\draw (i4) [fill=white] circle (\vr);
\draw (i5) [fill=white] circle (\vr);
\draw (n1) [fill=black] circle (\vrn);
\draw (n2) [fill=black] circle (\vrn);
\draw (n3) [fill=black] circle (\vrn);
\end{tikzpicture}
\caption{An infinite family of cubic graphs $G$ of order~$n$ satisfying $i(G) = \frac{3}{8}n$}
\label{f:new-family}
\end{center}
\end{figure}

Indeed, the family of extremal graphs achieving the $\frac{3}{8}$-upper bound in Theorem~\ref{thm:main} is rich and varied. As a further example of extremal graphs, the graphs $G$ of order~$n$ in the infinite family illustrated in Figure~\ref{f:new-family} also satisfy $i(G) = \frac{3}{8}n$.
We remark that the cubic graphs $G$ in the family illustrated in Figure~\ref{fig:cubic}(a) do not contain a subgraph isomorphic to $K_{2,3}$ and satisfy $i(G) = \frac{3}{8}n$, while the cubic graphs $G$ in the family illustrated in Figure~\ref{fig:cubic}(b) and Figure~\ref{f:new-family} do contain a subgraph isomorphic to $K_{2,3}$ and satisfy $i(G) = \frac{3}{8}n$.


We proceed as follows. In Section~\ref{S:overview} we give an overview of our proof to give the main ideas and direction of the proof. In Section~\ref{S:family}, we formally define the (infinite) special families of subgraphs we encounter in our proof, and we prove desirable properties of the families. Thereafter in Section~\ref{S:main-proof}, we present a detailed proof of our key result, namely Theorem~\ref{thm:main2}.

\section{Overview of proof of Theorem~\ref{thm:main}}
\label{S:overview}

Before presenting a proof of Theorem~\ref{thm:main}, we give a sketch proof to provide the reader an overview of the detailed proof that follows and the direction of the arguments.
In order to prove Theorem~\ref{thm:main}, we relax the $3$-regularity condition to a bounded maximum degree~$3$ condition to make the inductive hypothesis easier to handle, that is, we relax the requirement that $G \notin \{K_{3,3}, C_5 \, \Box \, K_2\}$ is a connected cubic graph to the requirement that $G$ is a subcubic graph that contains no $K_{3,3}$-component and no $(C_5 \, \Box \, K_2)$-component. We define next what we have coined the vertex weight and the structural weight of such graphs.

\subsection{The vertex weight}

For a subcubic graph $G$, we assign weights~$8$,~$5$,~$4$ and~$3$ to the vertices of $G$ of degrees~$0$,~$1$,~$2$ and~$3$, respectively. For each vertex $v$ of $G$, we denote the \emph{weight} of $v$ in $G$ by $\w_G(v)$ (see Table~1).
\[
\begin{array}{||c|c|c|c|c||} \hline
\deg_G(v) & 0 & 1 & 2 & 3 \\ \hline
\w_G(v) & 8 & 5 & 4 & 3 \\ \hline
\end{array}
\]
\begin{center}
\textbf{Table~1.} The weight $\w_G(v)$ of a vertex $v$ in $G$
\end{center}

For a subset $X$ of vertices in $G$, we define the \emph{vertex weight} of $X$ in $G$ as the sum of the weights in $G$ of vertices in $X$; that is,
\[
\w_G(X) = \sum_{v \in X} \w_G(v).
\]

We define the \emph{vertex weight} of $G$, denoted by $\w(G)$, as the sum of the weights of vertices in $G$, that is,
\[
\w(G) = 8n_0(G) + 5n_1(G) + 4n_2(G) + 3n_3(G).
\]

\subsection{The structural weight}

We would like to show that $8i(G) \le \w(G)$ since in the special case when $G$ is a cubic graph of order~$n$, we note that $n_0(G) = n_1(G) = n_2(G) = 0$ and $n = n_3(G)$, and so the weight of $G$ simplifies to~$\w(G) = 3n$, from which we would infer that $i(G) \le \frac{3}{8}n$, as desired.

However, relaxing the cubic condition results in the so-called ``bad'' family of subcubic graphs. We construct such a family $\cB$ of subcubic graphs in Section~\ref{S:familyB} with the property that if $G \in \cB$, then $8i(G) = \w(G) + 2$. We will refer to graphs $G$ in this family $\cB$ as \emph{bad graphs} since they satisfy $8i(G) > \w(G)$. We denote by $b(G)$ the number of components of $G$ that belong to the family $\cB$.

Moreover relaxing the cubic condition results in so-called ``troublesome'' configurations of subcubic graphs.  In Section~\ref{S:trouble-subgraph}, we define these troublesome configurations, which we call $T$-\emph{configurations}. We distinguish two types of troublesome configurations, namely $T_1$-\emph{configurations} and $T_2$-\emph{configurations}. A $T_i$-configuration is a troublesome configuration that is joined to vertices not in the subgraph by exactly $i$ edges. We denote by $\tc(G)$ the number of vertex disjoint troublesome configurations in $G$. We define the $\Theta$-\emph{weight} of $G$, which we call the \emph{structural weight} of $G$ since it is determined by structural properties of $G$, by
\[
\Theta(G) = 2\tc(G) + 2b(G).
\]

\subsection{The total weight}

We define the \emph{total weight} of $G$, which we denote by $\Omega(G)$, as the sum of the vertex weight and the structural weight of a subcubic graph $G$, and so
\[
\Omega(G) = \w(G) + \Theta(G).
\]

We call the total weight of $G$ the $\Omega$-\emph{weight of $G$.} Our key result, given in Theorem~\ref{thm:main2}, gives a tight upper bound on the independent domination number of a subcubic graph in terms of its $\Omega$-weight.

\begin{theorem}
\label{thm:main2}
If $G$ is a subcubic graph of order~$n$ that contains no $K_{3,3}$-component and no $(C_5 \, \Box \, K_2)$-component, then $8i(G) \le \Omega(G)$.
\end{theorem}

In order to prove our main result, namely Theorem~\ref{thm:main}, we need the above stronger statement given in Theorem~\ref{thm:main2}. Each graph in the family $\cB$ of subcubic graphs contains vertices of degree~$1$ or~$2$. Moreover, every so-called ``troublesome configuration'' contains a vertex of degree~$2$ in $G$. From these properties, a cubic graph $G$ has no component that belongs to the family $\cB$ and contains no troublesome configuration. We therefore infer that if $G$ is a cubic graph, then $\tc(G) = b(G) = 0$, implying that its structural weight is zero, that is, $\Theta(G) = 0$. Further if $G$ is a cubic graph of order~$n$, then $n_0(G) = n_1(G) = n_2(G) = 0$ and $n = n_3(G)$, and so its vertex weight is~$3n$, that is, $\w(G) = 3n$. Hence in this case when $G$ is a cubic graph, the inequality in the statement of Theorem~\ref{thm:main2} simplifies to $8i(G) \le 3n$. Our main result, namely Theorem~\ref{thm:main}, therefore follows readily from Theorem~\ref{thm:main2}.

\section{Special families of graphs and subgraphs}
\label{S:family}

In this section, we define the family $\cB$ of so-called `bad subgraphs' and we define what we have coined `troublesome configurations' of a subcubic graph $G$. In order to prove desirable properties of subcubic graphs, we first define `exit edges' associated with a set $X$ of vertices of a graph $G$.

\subsection{Exit edges of a vertex set}
\label{S:exit-edge}

Let $X$ be a proper subset of vertices in a subcubic graph $G$ and let $V = V(G)$. We call an edge an \emph{$X$-exit edge} in $G$ if it joins a vertex in $X$ to a vertex in $V \setminus X$. We denote the number of $X$-exit edges in $G$ by $\xi_G(X)$. Those vertices in $G - X$ that are incident in $G$ with at least one $X$-exit edge have smaller degree in $G - X$ than in $G$, and therefore have a larger weight in $G-X$ than in $G$. We refer to the sum of these weight increases as the $\w$-\emph{cost of removing} $X$ from $G$ and denote the weight increase by $\Phi_G(X)$. Thus,
\[
\Phi_G(X) =  \sum_{v \in V \setminus X} \big( \w_{G-X}(v) - \w_{G}(v) \big).
\]

We prove next the fact that the $\w$-cost of removing a set $X$ of vertices from $G$ is equal to the number of $X$-exit edges in $G$ plus twice the number of isolated vertices in $G-X$. Consequently, if removing $X$ yields no isolated vertices, then the cost of removing $X$ from $G$ will be precisely the number of $X$-exit edges in $G$. For a vertex $v \in V \setminus X$, we denote by $\Phi_G(v)$ the number of $X$-edges in $G$ that are incident with the vertex~$v$.

\begin{fact}
\label{fact:F1}
If $X \subset V(G)$, then $\Phi_G(X) = \xi_G(X) + 2 n_0(G-X)$.
\end{fact}
\proof
We consider the contribution of a vertex $v \in V \setminus X$ to the cost $\Phi_G(X)$. If the vertex $v$ is incident with no $X$-exit edge in $G$, then $\w_{G-X}(v) = \w_{G}(v)$, and the contribution of~$v$ to the sum $\Phi_G(X)$ is zero. Hence, we may assume that $v$ is incident with at least one $X$-exit edge in $G$, that is, $\Phi_G(v) \ge 1$. If $v$ is not isolated in $G-X$, then $\w_{G-X}(v) - \w_{G}(v)$ is precisely the number of exit edges of $X$ incident with $v$ in $G$, namely~$\Phi_G(v)$. If $v$ is isolated in $G-X$, then $\w_{G-X}(v) - \w_{G}(v)$ is the number of exit edges of $X$ incident with $v$ in $G$ plus an additional~$2$ noting that in this case $\w_{G-X}(v) = 8$ and~$\w_{G}(v) = 6 - \Phi_G(v)$.~\QED

\subsection{The bad family $\cB$}
\label{S:familyB}

\begin{figure}[htb]
\begin{center}
\begin{tikzpicture}[scale=.75,style=thick,x=1cm,y=1cm]
\def\vr{2.25pt}
\path (0,1) coordinate (v1);
\path (1,0) coordinate (v2);
\path (1,1) coordinate (v3);
\path (1,2) coordinate (v4);
\path (2,1) coordinate (v5);
\path (-1,1) coordinate (v6);
\draw (v1)--(v2)--(v5)--(v4)--(v1);
\draw (v2)--(v3)--(v4);
\draw (v1)--(v6);
%
\draw (v1) [fill=white] circle (\vr);
\draw (v2) [fill=white] circle (\vr);
\draw (v3) [fill=white] circle (\vr);
\draw (v4) [fill=white] circle (\vr);
\draw (v5) [fill=white] circle (\vr);
\draw (v6) [fill=white] circle (\vr);
%
\draw[anchor = east] (v6) node {{\small $r$}};
\end{tikzpicture}
\caption{The base graph $B_1$}
\label{fig:B1}
\end{center}
\end{figure}
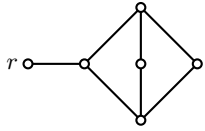

In this section, we define the so-called ``bad family'' $\cB$ of graphs. Let $B_1$ be the graph obtained from $K_{2,3}$ by adding a new vertex $r$ and adding an edge joining $r$ and a vertex of degree~$2$ in the copy of $K_{2,3}$. The graph $B_1$, illustrated in Figure~\ref{fig:B1}, we call the \emph{base graph} used in our construction of graphs in the family $\cB$.

Let $\cB$ be the family of subcubic graphs that contain the base graph $B_1$ and is closed under the operation $\cO_1$ listed below that extends a graph $G' \in \cB$ to a new graph $G \in \cB$. As remarked earlier, we refer to graphs in the family $\cB$ as \emph{bad graphs}.

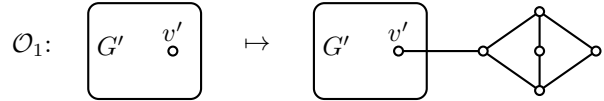
\begin{figure}[htb]
\begin{center}
\begin{tikzpicture}[scale=.75,style=thick,x=1cm,y=1cm]
\def\vr{2.25pt}
\path (1.5,0.65) coordinate (v);
\draw (v) [fill=white] circle (\vr);
\draw [rounded corners] (0,-.2) rectangle (2,1.5);
\draw (.4,0.75) node {$G'$};
\draw (-1,.75) node {$\cO_1$:};
\draw[anchor = south] (v) node {$v'$};
\draw (3,0.75) node {$\mapsto$};
\path (5.5,0.65) coordinate (g);
\path (7,0.65) coordinate (a);
\path (8,-0.05) coordinate (b);
\path (8,0.65) coordinate (c);
\path (8,1.35) coordinate (d);
\path (9,0.65) coordinate (e);
%
\draw (g)--(a)--(b)--(e)--(d)--(a);
\draw (b)--(c)--(d);
\draw (g) [fill=white] circle (\vr);
\draw (a) [fill=white] circle (\vr);
\draw (b) [fill=white] circle (\vr);
\draw (c) [fill=white] circle (\vr);
\draw (d) [fill=white] circle (\vr);
\draw (e) [fill=white] circle (\vr);
\draw [rounded corners] (4,-.2) rectangle (6,1.5);
\draw (4.4,0.75) node {$G'$};
\draw[anchor = south] (g) node {$v'$};
\end{tikzpicture}
\caption{The operation $\cO_1$}
\label{f:cO1}
\end{center}
\end{figure}

\begin{itemize}
\item \emph{Operation} $\cO_1$: The graph $G$ is obtained from $G' \in \cB$ by selecting a vertex $v'$ of degree at most~$2$ in $G'$, adding a vertex disjoint copy of $K_{2,3}$, and adding an edge joining $v'$ and a vertex of degree~$2$ in the added copy of $K_{2,3}$. See Figure~\ref{f:cO1}.
\end{itemize}

Thus, if $G \in \cB$, then either $G = B_1$ or $G$ is constructed from the base graph $B_1$ by $k - 1$ applications of Operation~$\cO_1$ where $k \ge 2$. Thus, either $G = B_1$ or $k \ge 2$ and $G = B_k$ where $B_1,B_2,\ldots,B_k$ is a sequence of graphs in the family $\cB$ and the graph $B_{i+1}$ is obtained from the graph $B_{i}$ by applying Operation~$\cO_1$ for all $i \in [k-1]$. We call $B_{i+1}$ a \emph{bad}-\emph{extension} of the graph $B_{i}$ for $i \in [k-1]$.

If $G \in \cB$, then we define the \emph{root} of the graph $G$ as the vertex of $G$ that does not belong to any copy of $K_{2,3}$, equivalently, the root of $G$ is the (unique) vertex that belongs to no cycle. We note that the root has degree~$1$,~$2$ or~$3$ in the graph $G$. We define
$\cB_i = \{B \in \cB \colon \mbox{the root of $B$ has degree~$i$ in $B$} \}$ for $i \in [3]$,
and so
\[
\cB = \bigcup_{i=1}^3 \cB_i.
\]

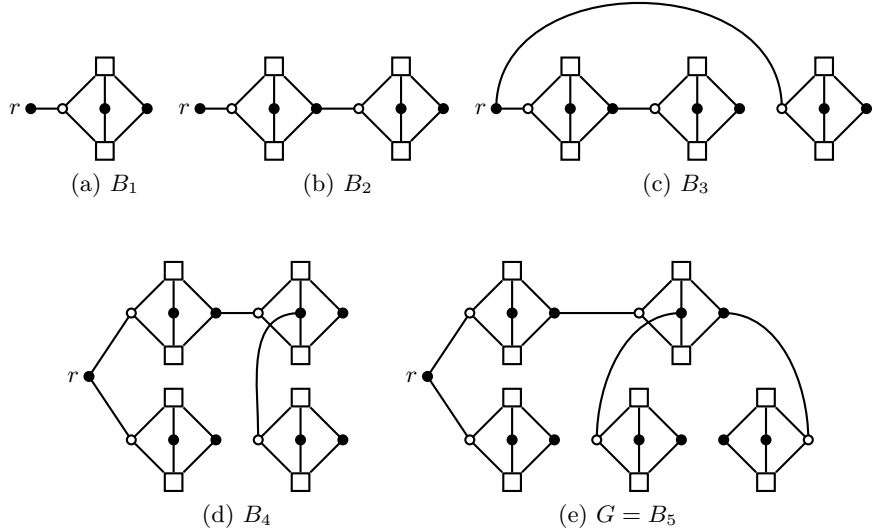
\begin{figure}[htb]
\begin{center}
\begin{tikzpicture}[scale=.75,style=thick,x=0.75cm,y=0.75cm]
\def\vr{2.25pt}
\path (0,1) coordinate (v1);
\node[rectangle,draw] (v2) at (1,0) {};
\path (1,1) coordinate (v3);
\node[rectangle,draw] (v4) at (1,2) {};
\path (2,1) coordinate (v5);
\path (-0.75,1) coordinate (v6);
\draw (v1)--(v2)--(v5)--(v4)--(v1);
\draw (v2)--(v3)--(v4);
\draw (v1)--(v6);
%
\draw (v1) [fill=white] circle (\vr);
\draw (v3) [fill=black] circle (\vr);
\draw (v5) [fill=black] circle (\vr);
\draw (v6) [fill=black] circle (\vr);
\draw (1,-0.8) node {{\small (a) $B_1$}};
\draw[anchor = east] (v6) node {{\small $r$}};
\path (3.25,1) coordinate (v6);
\path (4,1) coordinate (v1);
\node[rectangle,draw] (v2) at (5,0) {};
\path (5,1) coordinate (v3);
\node[rectangle,draw] (v4) at (5,2) {};
\path (6,1) coordinate (v5);
\path (7,1) coordinate (u1);
\node[rectangle,draw] (u2) at (8,0) {};
\path (8,1) coordinate (u3);
\node[rectangle,draw] (u4) at (8,2) {};
\path (9,1) coordinate (u5);
\draw (v1)--(v2)--(v5)--(v4)--(v1);
\draw (v2)--(v3)--(v4);
\draw (u1)--(u2)--(u5)--(u4)--(u1);
\draw (u2)--(u3)--(u4);
\draw (v5)--(u1);
\draw (v1)--(v6);
%
\draw (v1) [fill=white] circle (\vr);
\draw (v3) [fill=black] circle (\vr);
\draw (v5) [fill=black] circle (\vr);
\draw (v6) [fill=black] circle (\vr);
\draw (u1) [fill=white] circle (\vr);
\draw (u3) [fill=black] circle (\vr);
\draw (u5) [fill=black] circle (\vr);
\draw (6.5,-0.8) node {{\small (b) $B_2$}};
\draw[anchor = east] (v6) node {{\small $r$}};
\path (10.25,1) coordinate (v6);
\path (11,1) coordinate (v1);
\node[rectangle,draw] (v2) at (12,0) {};
\path (12,1) coordinate (v3);
\node[rectangle,draw] (v4) at (12,2) {};
\path (13,1) coordinate (v5);
\path (14,1) coordinate (u1);
\node[rectangle,draw] (u2) at (15,0) {};
\path (15,1) coordinate (u3);
\node[rectangle,draw] (u4) at (15,2) {};
\path (16,1) coordinate (u5);
\path (17,1) coordinate (w1);
\node[rectangle,draw] (w2) at (18,0) {};
\path (18,1) coordinate (w3);
\node[rectangle,draw] (w4) at (18,2) {};
\path (19,1) coordinate (w5);
%
\draw (v1)--(v2)--(v5)--(v4)--(v1);
\draw (v2)--(v3)--(v4);
\draw (u1)--(u2)--(u5)--(u4)--(u1);
\draw (u2)--(u3)--(u4);
\draw (w1)--(w2)--(w5)--(w4)--(w1);
\draw (w2)--(w3)--(w4);
\draw (v5)--(u1);
\draw (v1)--(v6);
%
\draw (v6) to[out=90,in=90, distance=2.5cm] (w1);
%
%
\draw (v1) [fill=white] circle (\vr);
\draw (v3) [fill=black] circle (\vr);
\draw (v5) [fill=black] circle (\vr);
\draw (v6) [fill=black] circle (\vr);
\draw (u1) [fill=white] circle (\vr);
\draw (u3) [fill=black] circle (\vr);
\draw (u5) [fill=black] circle (\vr);
\draw (w1) [fill=white] circle (\vr);
\draw (w3) [fill=black] circle (\vr);
\draw (w5) [fill=black] circle (\vr);
\draw (14.5,-0.8) node {{\small (c) $B_3$}};
\draw[anchor = east] (v6) node {{\small $r$}};
\end{tikzpicture}
\end{center}

\vskip 0.2 cm

\begin{center}
\begin{tikzpicture}[scale=.75,style=thick,x=0.75cm,y=0.75cm]
\def\vr{2.25pt}
\path (6,-0.5) coordinate (v6);
\path (7,1) coordinate (v1);
\node[rectangle,draw] (v2) at (8,0) {};
\path (8,1) coordinate (v3);
\node[rectangle,draw] (v4) at (8,2) {};
\path (9,1) coordinate (v5);
\path (10,1) coordinate (u1);
\node[rectangle,draw] (u2) at (11,0) {};
\path (11,1) coordinate (u3);
\node[rectangle,draw] (u4) at (11,2) {};
\path (12,1) coordinate (u5);
\path (7,-2) coordinate (w1);
\node[rectangle,draw] (w2) at (8,-3) {};
\path (8,-2) coordinate (w3);
\node[rectangle,draw] (w4) at (8,-1) {};
\path (9,-2) coordinate (w5);
\path (10,-2) coordinate (x1);
\node[rectangle,draw] (x2) at (11,-3) {};
\path (11,-2) coordinate (x3);
\node[rectangle,draw] (x4) at (11,-1) {};
\path (12,-2) coordinate (x5);
\draw (v1)--(v2)--(v5)--(v4)--(v1);
\draw (v2)--(v3)--(v4);
\draw (u1)--(u2)--(u5)--(u4)--(u1);
\draw (u2)--(u3)--(u4);
\draw (w1)--(w2)--(w5)--(w4)--(w1);
\draw (w2)--(w3)--(w4);
\draw (x1)--(x2)--(x5)--(x4)--(x1);
\draw (x2)--(x3)--(x4);
\draw (v5)--(u1);
\draw (v1)--(v6);
\draw (w1)--(v6);
\draw (x1) to[out=90,in=180, distance=1cm] (u3);
%
\draw (v1) [fill=white] circle (\vr);
\draw (v3) [fill=black] circle (\vr);
\draw (v5) [fill=black] circle (\vr);
\draw (v6) [fill=black] circle (\vr);
\draw (u1) [fill=white] circle (\vr);
\draw (u3) [fill=black] circle (\vr);
\draw (u5) [fill=black] circle (\vr);
\draw (w1) [fill=white] circle (\vr);
\draw (w3) [fill=black] circle (\vr);
\draw (w5) [fill=black] circle (\vr);
\draw (x1) [fill=white] circle (\vr);
\draw (x3) [fill=black] circle (\vr);
\draw (x5) [fill=black] circle (\vr);
\draw (9.5,-3.8) node {{\small (d) $B_4$}};
\draw[anchor = east] (v6) node {{\small $r$}};
\path (14,-0.5) coordinate (v6);
\path (15,1) coordinate (v1);
\node[rectangle,draw] (v2) at (16,0) {};
\path (16,1) coordinate (v3);
\node[rectangle,draw] (v4) at (16,2) {};
\path (17,1) coordinate (v5);
\path (19,1) coordinate (u1);
\node[rectangle,draw] (u2) at (20,0) {};
\path (20,1) coordinate (u3);
\node[rectangle,draw] (u4) at (20,2) {};
\path (21,1) coordinate (u5);
\path (15,-2) coordinate (w1);
\node[rectangle,draw] (w2) at (16,-3) {};
\path (16,-2) coordinate (w3);
\node[rectangle,draw] (w4) at (16,-1) {};
\path (17,-2) coordinate (w5);
\path (18,-2) coordinate (x1);
\node[rectangle,draw] (x2) at (19,-3) {};
\path (19,-2) coordinate (x3);
\node[rectangle,draw] (x4) at (19,-1) {};
\path (20,-2) coordinate (x5);
\path (21,-2) coordinate (y1);
\node[rectangle,draw] (y2) at (22,-3) {};
\path (22,-2) coordinate (y3);
\node[rectangle,draw] (y4) at (22,-1) {};
\path (23,-2) coordinate (y5);
\draw (v1)--(v2)--(v5)--(v4)--(v1);
\draw (v2)--(v3)--(v4);
\draw (u1)--(u2)--(u5)--(u4)--(u1);
\draw (u2)--(u3)--(u4);
\draw (w1)--(w2)--(w5)--(w4)--(w1);
\draw (w2)--(w3)--(w4);
\draw (x1)--(x2)--(x5)--(x4)--(x1);
\draw (x2)--(x3)--(x4);
\draw (y1)--(y2)--(y5)--(y4)--(y1);
\draw (y2)--(y3)--(y4);
\draw (v5)--(u1);
\draw (v1)--(v6);
\draw (w1)--(v6);
\draw (x1) to[out=90,in=180, distance=1cm] (u3);
\draw (u5) to[out=0,in=90, distance=1cm] (y5);
%
\draw (v1) [fill=white] circle (\vr);
\draw (v3) [fill=black] circle (\vr);
\draw (v5) [fill=black] circle (\vr);
\draw (v6) [fill=black] circle (\vr);
\draw (u1) [fill=white] circle (\vr);
\draw (u3) [fill=black] circle (\vr);
\draw (u5) [fill=black] circle (\vr);
\draw (w1) [fill=white] circle (\vr);
\draw (w3) [fill=black] circle (\vr);
\draw (w5) [fill=black] circle (\vr);
\draw (x1) [fill=white] circle (\vr);
\draw (x3) [fill=black] circle (\vr);
\draw (x5) [fill=black] circle (\vr);
\draw (y1) [fill=black] circle (\vr);
\draw (y3) [fill=black] circle (\vr);
\draw (y5) [fill=white] circle (\vr);
\draw (18.5,-3.8) node {{\small (e) $G = B_5$}};
\draw[anchor = east] (v6) node {{\small $r$}};
\end{tikzpicture}

\caption{A graph $G$ in the family $\cB_2$ with root vertex~$r$}
\label{fig:B0B5}
\end{center}
\end{figure}

For example, the graph $G$ illustrated in Figure~\ref{fig:B0B5}(e) belongs to the family~$\cB$ since it can be constructed from a sequence $B_1,B_2,B_3,B_4,B_5$ of graphs in $\cB$ shown in Figure~\ref{fig:B0B5}(a)--\ref{fig:B0B5}(e), where $B_{i+1}$ a bad-extension of $B_i$ for $i \in [4]$ and where $G = B_5$ (in this case, $G = B_k$ where $k = 5$). In this example, the root is the vertex labelled~$r$ in Figure~\ref{fig:B0B5}, and has degree~$2$ in the resulting graph $G = B_5$, and so the graph $G$ belongs to the subfamily~$\cB_2$ of $\cB$.

We define the \emph{canonical ID}-\emph{set} of a bad graph $G$ (that belongs to the family $\cB$) as follows. Initially if $G = B_1$, we define the canonical ID-set to consist of the root and the two vertices of degree~$2$ in $B_1$ (as illustrated by the shaded vertices in Figure~\ref{fig:B0B5}(a)). Suppose that $k \ge 2$ and let $B_1,B_2,\ldots,B_k$ be the sequence of graphs in the family $\cB$ used to construct $G \in \cB$, and so $G = B_k$ and $B_{i+1}$ a bad-extension of the graph $B_{i}$ for $i \in [k-1]$. When we apply Operation~$\cO_1$ to extend the graph $B_{i}$ to the graph $B_{i+1}$ we add to the current canonical ID-set in $B_{i}$ the two new vertices of degree~$2$ in $B_{i+1}$ that belong to the added copy of $K_{2,3}$ for all $i \in [k-1]$. We note that by construction, every vertex of degree~$2$ in $G$ belongs to the canonical ID-set of $G$. For example, the canonical ID-set of each of the graphs $B_1, B_2, B_3, B_4, B_5$ used to construct the graph $G = B_5$ is given by the shaded vertices in Figures~\ref{fig:B0B5}(a)--\ref{fig:B0B5}(e).

A vertex of degree~$3$ in $G$ we call a \emph{large vertex}. Each copy of $K_{2,3}$ in $G \in \cB$ contains exactly two vertices that belong to the canonical ID-set of $G$, and these two vertices have two common large neighbors that belong to the same copy of $K_{2,3}$. We refer to these two large vertices as \emph{non}-\emph{canonical vertices} of $G$. Thus each copy of $K_{2,3}$ in $G \in \cB$ contains two non-canonical vertices, namely the two vertices of degree~$3$ in that copy of $K_{2,3}$ that are open twins, where two vertices are open twins if they have the same open neighborhood. We define the \emph{non}-\emph{canonical independent set} of $G \in \cB$ to be the set consisting of all non-canonical vertices of~$G$. Thus if $G \in \cB$ is constructed from the base graph $B_1$ by $k - 1$ applications of Operation~$\cO_1$ where $k \ge 1$, then $G$ contains $k$ vertex disjoint copies of $K_{2,3}$ and the $2k$ non-canonical vertices of $G$ (two from each copy of $K_{2,3}$ in $G$) form the non-canonical independent set of $G$. For example, the non-canonical independent set of each of the graphs $B_1, B_2, B_3, B_4, B_5$ used to construct the graph $G = B_5$ is given by the vertices represented by squares in Figures~\ref{fig:B0B5}(a)--\ref{fig:B0B5}(e).

We define the set of vertices in $G \in \cB$ that belong to neither the canonical ID-set of $G$ nor the non-canonical independent set of $G$ as the \emph{core independent set} in $G$, and we refer to the vertices in the core independent set as \emph{core vertices}. Thus if $G \in \cB$ is constructed from the base graph $B_1$ by $k - 1$ applications of Operation~$\cO_1$ where $k \ge 1$, then $G$ contains $k$ core vertices, one from each of the $k$ vertex disjoint copies of $K_{2,3}$. For example, the core independent set of each of the graphs $B_1, B_2, B_3, B_4, B_5$ used to construct the graph $G = B_5$ is given by the vertices represented by white circled vertices in Figures~\ref{fig:B0B5}(a)--\ref{fig:B0B5}(e).

We shall need the following lemmas.

\begin{figure}[htb]

\begin{center}
\begin{tikzpicture}[scale=.75,style=thick,x=1cm,y=1cm]
\def\vr{2.25pt}
\path (5.5,0.75) coordinate (g);
\path (7,0.75) coordinate (a);
\path (8,0.0) coordinate (b);
\path (8,0.75) coordinate (c);
\path (8,1.5) coordinate (d);
\path (9,0.75) coordinate (e);
%
\draw (g)--(a)--(b)--(e)--(d)--(a);
\draw (b)--(c)--(d);
\draw (g) [fill=white] circle (\vr);
\draw (a) [fill=white] circle (\vr);
\draw (b) [fill=white] circle (\vr);
\draw (c) [fill=white] circle (\vr);
\draw (d) [fill=white] circle (\vr);
\draw (e) [fill=white] circle (\vr);
\draw [rounded corners] (3,-.2) rectangle (6,1.5);
\draw (4.4,0.75) node {$G'$};
\draw[anchor = north] (a) node {$v_3$};
\draw[anchor = north] (g) node {$v'$};
\draw[anchor = south] (d) node {$u_2$};
\draw[anchor = north] (b) node {$u_1$};
\draw[anchor = west] (c) node {$v_2$};
\draw[anchor = west] (e) node {$v_1$};
\end{tikzpicture}
\caption{The graphs $G'$ and $G$ in the proof of Lemma~\ref{lem:lem1}}
\label{fig:lem1}
\end{center}
\end{figure}
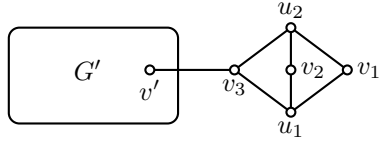

\begin{lemma}
\label{lem:lem1}
If a subcubic graph $G$ is obtained from a subcubic isolate-free graph $G'$ by adding a vertex disjoint copy of $K_{2,3}$, and adding an edge joining a vertex $v'$ of degree at most~$2$ in $G'$ and a vertex of degree~$2$ in the added copy of $K_{2,3}$, then the following properties hold. \\ [-22pt]
\begin{enumerate}
\item[{\rm (a)}] $i(G) = i(G') + 2$.
\item[{\rm (b)}] $\w(G) = \w(G') + 16$.
\end{enumerate}
\end{lemma}
\proof Let $G$, $G'$ and $v'$ be as defined in the statement of the lemma. Let $F$ be the copy of $K_{2,3}$ added to $G'$ when constructing the graph $G$, where $v_1,v_2$ and $v_3$ are the three small vertices (of degree~$2$) in $F$ and $u_1$ and $u_2$ are the two large vertices (of degree~$3$) in $F$. Further, let $v_3$ be the vertex of $F$ adjacent to $v'$ in $G$. Let $X = V(F)$. The construction of $G$ is illustrated in Figure~\ref{fig:lem1}.

(a) Every $i$-set of $G'$ can be extended to an ID-set of $G$ by adding to it the vertices $u_1$ and $u_2$, and so $i(G) \le i(G') + 2$. Conversely, let $I$ be an $i$-set of $G$. If $\{u_1,u_2\} \subset I$, then we let $I' = I \setminus \{u_1,u_2\}$. If $\{u_1,u_2\} \not\subset I$, then $\{v_1,v_2\} \subset I$. Moreover if in this case $v_3 \notin I$, then $v' \in I$ and we let $I' = I \setminus \{v_1,v_2\}$. Suppose that $\{v_1,v_2,v_3\} \subseteq I$. In this case, if a neighbor of $v'$ different from $v_3$ belongs to the set $I$, then the set $(I \setminus \{v_1,v_2,v_3\}) \cup \{u_1,u_2\}$ is an ID-set of $G$ of cardinality~$|I|-1$, contradicting the minimality of the set $I$. Hence, the vertex $v_3$ is the only neighbor of $v'$ that belongs to the set $I$. We now let $I' = (I \setminus \{v_1,v_2,v_3\}) \cup \{v'\}$. In all the above cases, the set $I'$ is an ID-set of $G'$ of cardinality~$|I| - 2$, and so $i(G') \le |I| - 2 = i(G) - 2$. As observed earlier, $i(G) \le i(G') + 2$. Consequently, $i(G) = i(G') + 2$. This proves part~(a).

\medskip
(b) By supposition, the graph $G$ is a subcubic graph and the graph $G'$ is isolate-free. We infer that the vertex $v'$ has degree~$1$ or~$2$ in $G'$ and therefore has degree~$2$ or~$3$ in $G$, respectively. Thus adopting our earlier notation, the edge $v'v_3$ is the only $X$-exit edge in $G$. Thus, the number of $X$-exit edges is $\xi_G(X) = 1$. By Fact~\ref{fact:F1}, the cost of removing $X$ from $G$ is $\Phi_G(X) = \xi_G(X) = 1$ noting that $G' = G - X$ is isolate-free. Thus, $\w(G) = \w(G') + \w_G(X) - \Phi_G(X) = \w(G') + 17 - 1 = \w(G') + 16$, which proves part~(b).~\QED

\medskip
Suppose that $G \in \cB$, and so $G$ is constructed from the base graph $B_1$ by $k - 1$ applications of Operation~$\cO_1$ where $k \ge 1$. We note that if $k = 1$, then $G = B_1$. Further, we note that $i(B_1) = 3$ and $\w(B_1) = 22$. Hence by repeated applications of Lemma~\ref{lem:lem1}, we infer the following properties of graphs in the family~$\cB$.

\begin{proposition}
\label{prop:cB}
If $G \in \cB$ is constructed from the base graph $B_1$ by $k - 1$ applications of Operation~$\cO_1$ where $k \ge 1$, then the following properties hold. \\ [-24pt]
\begin{enumerate}
\item[{\rm (a)}] $i(G) = 2k + 1$, $\w(G) = 16k + 6$, and  $8i(G) = \w(G) + 2$.
\item[{\rm (b)}] The canonical ID-set of $G$ is an $i$-set of $G$.
\item[{\rm (c)}] The non-canonical independent set of $G$, together with the root vertex of $G$, is an $i$-set of $G$.
\item[{\rm (d)}] If $v$ is the root vertex of $G$ or a vertex of degree~$2$ in $G$, then $i(G - v) = i(G) - 1$.
\item[{\rm (e)}] The core independent set of $G$ is an ID-set of the graph obtained by removing all vertices of degree~$2$ from $G$.
\item[{\rm (f)}] The graph $G$ contains at most one vertex of degree~$1$, and no two adjacent vertices of $G$ both have degree~$2$ in $G$.
\end{enumerate}
\end{proposition}

\subsection{Troublesome configurations}
\label{S:trouble-subgraph}

In this section, we define what we have coined a `troublesome configuration' of a subcubic graph $G$. We define a \emph{troublesome configuration} of $G$, abbreviated $T$-\emph{configuration of $G$}, as a subgraph of $G$ obtained from a bad graph $B \in \cB_1$ with root vertex~$v_1$ (of degree~$1$ in $B$) by applying the following process. \\ [-20pt]
\begin{enumerate}
\item[$\bullet$] adding a new vertex~$v_2$ to~$B$,
\item[$\bullet$] adding an edge joining $v_2$ to the vertex $v_1$ and to a vertex $w_2$ of degree~$2$ in $B$ different from~$v_1$,
\item[$\bullet$] adding an edge joining $v_2$ to a vertex not in $V(B) \cup \{v_2\}$, and
\item[$\bullet$] adding at most one edge joining $v_1$ to a vertex not in $V(B) \cup \{v_2\}$.
\end{enumerate}

\begin{figure}[htb]
\begin{center}
\begin{tikzpicture}[scale=.75,style=thick,x=0.75cm,y=0.75cm]
\def\vr{2.25pt}
\def\vrn{1.25pt}
\path (-3,3) coordinate (z0);
\path (4,3) coordinate (z1);
\path (4,6) coordinate (z2);
\path (0.5,6) coordinate (z3);
\path (0.5,7) coordinate (z4);
\path (-3,7) coordinate (z5);
\path (-2,4.5) coordinate (b1);
\path (-1,3.5) coordinate (b2);
\path (-1,4.5) coordinate (b3);
\path (-1,5.5) coordinate (b4);
\path (0,4.5) coordinate (b5);
\path (-2,6.5) coordinate (e1);
\path (-2,7.5) coordinate (e11);
\path (1,4.5) coordinate (f1);
\path (2,3.5) coordinate (f2);
\path (2,4.5) coordinate (f3);
\path (2,5.5) coordinate (f4);
\path (3,4.5) coordinate (f5);
%
\path (3,6.5) coordinate (g5);
\path (3,7.75) coordinate (g51);
\draw (b1)--(e1);
\draw (e1)--(g5);
\draw (b1)--(b2)--(b3)--(b4)--(b5)--(b2);
\draw (b1)--(b4);
%
%
\draw (f5)--(g5)--(g51);
%
\draw (f1)--(f2)--(f3)--(f4)--(f5)--(f2);
\draw (f1)--(f4);
\draw (b5)--(f1);
%
\draw [style=dashed,rounded corners] (-4,2.75) rectangle (4.35,7.25);
\draw[densely dashed] (z0)--(z1)--(z2)--(z3)--(z4)--(z5)--(z0);
\draw (e1) [fill=black] circle (\vr);
%
\draw (b1) [fill=white] circle (\vr);
\draw (b2) [fill=white] circle (\vr);
\draw (b3) [fill=black] circle (\vr);
\draw (b4) [fill=white] circle (\vr);
\draw (b5) [fill=black] circle (\vr);
%
%
\draw (g5) [fill=white] circle (\vr);
\draw (f1) [fill=white] circle (\vr);
\draw (f2) [fill=white] circle (\vr);
\draw (f3) [fill=black] circle (\vr);
\draw (f4) [fill=white] circle (\vr);
\draw (f5) [fill=black] circle (\vr);
\draw[anchor = east] (e1) node {{\small $v_1$}};
\draw[anchor = west] (g5) node {{\small $v_2$}};
\draw[anchor = west] (f5) node {{\small $w_2$}};
\draw (-2.5,4) node {{\small $B$}};
\draw (-3.6,5) node {{\small $T$}};
\draw (0.5,1.75) node {{\small (a) A $T_{v_1,v_2}^1$-configuration}};
\path (8,3) coordinate (z0);
\path (15,3) coordinate (z1);
\path (15,6) coordinate (z2);
\path (11.5,6) coordinate (z3);
\path (11.5,7) coordinate (z4);
\path (8,7) coordinate (z5);
\path (9,4.5) coordinate (b1);
\path (10,3.5) coordinate (b2);
\path (10,4.5) coordinate (b3);
\path (10,5.5) coordinate (b4);
\path (11,4.5) coordinate (b5);
\path (9,6.5) coordinate (e1);
\path (9,7.75) coordinate (e11);
\path (12,4.5) coordinate (f1);
\path (13,3.5) coordinate (f2);
\path (13,4.5) coordinate (f3);
\path (13,5.5) coordinate (f4);
\path (14,4.5) coordinate (f5);
%
\path (14,6.5) coordinate (g5);
\path (14,7.75) coordinate (g51);
\draw (b1)--(e1)--(e11);
\draw (e1)--(g5);
\draw (b1)--(b2)--(b3)--(b4)--(b5)--(b2);
\draw (b1)--(b4);
%
%
\draw (f5)--(g5)--(g51);
%
\draw (f1)--(f2)--(f3)--(f4)--(f5)--(f2);
\draw (f1)--(f4);
\draw (b5)--(f1);
%
%
\draw [style=dashed,rounded corners] (7,2.75) rectangle (15.35,7.25);
\draw[densely dashed] (z0)--(z1)--(z2)--(z3)--(z4)--(z5)--(z0);
\draw (e1) [fill=white] circle (\vr);
%
\draw (b1) [fill=white] circle (\vr);
\draw (b2) [fill=black] circle (\vr);
\draw (b3) [fill=white] circle (\vr);
\draw (b4) [fill=black] circle (\vr);
\draw (b5) [fill=white] circle (\vr);
%
%
\draw (g5) [fill=black] circle (\vr);
\draw (f1) [fill=white] circle (\vr);
\draw (f2) [fill=black] circle (\vr);
\draw (f3) [fill=white] circle (\vr);
\draw (f4) [fill=black] circle (\vr);
\draw (f5) [fill=white] circle (\vr);
\draw[anchor = east] (e1) node {{\small $v_1$}};
\draw[anchor = west] (g5) node {{\small $v_2$}};
\draw[anchor = west] (f5) node {{\small $w_2$}};
\draw (8.5,4) node {{\small $B$}};
\draw (7.3,5) node {{\small $T$}};
\draw (11.5,1.75) node {{\small (b) A $T_{v_1,v_2}^2$-configuration}};
\end{tikzpicture}
\vskip -0.25 cm
\caption{Type-$1$ and type-$2$ troublesome configurations in $G$}
\label{fig:type1-type2}
\end{center}
\end{figure}
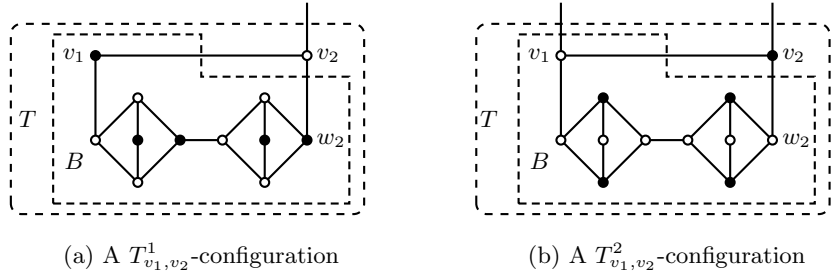

If $v_1$ has degree~$2$ in $G$, then we call the troublesome configuration a \emph{type-$1$ troublesome configuration} of $G$, abbreviated $T_{v_1,v_2}^1$-\emph{configuration of $G$}. A $T_{v_1,v_2}^1$-configuration is therefore joined to the rest of the graph by exactly one edge, namely an edge that joins $v_2$ to one vertex not in the subgraph. If $v_1$ has degree~$3$ in $G$, then we call the troublesome configuration a \emph{type-$2$ troublesome configuration} of $G$, abbreviated $T_{v_1,v_2}^2$-\emph{configuration of $G$}. A $T_{v_1,v_2}^2$-configuration is therefore joined to the rest of the graph by exactly two edges, namely an edge that joins each of $v_1$ and $v_2$ to a vertex not in the subgraph. In a troublesome configuration, we call the vertices $v_1$ and $v_2$ the \emph{link vertices} and we call all other vertices the \emph{non-link vertices} of the configuration. Thus the degree of every non-link vertex of a troublesome configuration is the same as its degree in $G$. A $T_{v_1,v_2}^1$-configuration is illustrated in Figure~\ref{fig:type1-type2}(a), and a $T_{v_1,v_2}^2$-configuration is illustrated in Figure~\ref{fig:type1-type2}(b).

The shaded vertices in Figure~\ref{fig:type1-type2}(a) indicate an $i$-set of $T_{v_1,v_2}^1$ and this set is the canonical ID-set of $B$. Moreover the shaded vertices in Figure~\ref{fig:type1-type2}(b) indicate an $i$-set of $T_{v_1,v_2}^2$ that contains the vertex~$v_2$ and the non-canonical vertices of $B$. More generally, we have the following properties of troublesome configurations of~$G$. These properties follows readily from properties of bad graphs (that belong to the family $\cB$) established in Section~\ref{S:familyB}.

\begin{proposition}
\label{prop:troublesome-subgraph}
If $T = T_{v_1,v_2}^1$ or $T = T_{v_1,v_2}^2$ is a $T_1$- or $T_2$-configuration of $G$ with link vertices $v_1$ and $v_2$ that is obtained from a bad graph $B \in \cB$ with root vertex $v_1$, then the following properties hold. \\ [-24pt]
\begin{enumerate}
\item[{\rm (a)}] $i(T) = i(B)$.
\item[{\rm (b)}] The canonical ID-set of $B$ is an $i$-set of $T$.
\item[{\rm (c)}] There exists an $i$-set of $T$ that contains the link vertex~$v_2$.
\item[{\rm (d)}] If $S$ is an arbitrary set of vertices of degree~$2$ in $T$ that contains neither $v_1$ nor $v_2$, then $i(T - S) = i(T) - |S|$. Furthermore, there exists an $i$-set of $T - S$ that contains neither $v_1$ nor $v_2$, and there exists an $i$-set $I_j$ of $T - S$ that contains $v_j$ for $j \in [2]$.
\item[{\rm (f)}] Every vertex in $T$ has degree at least~$2$, and no two adjacent vertices of $T$ both have degree~$2$.
\end{enumerate}
\end{proposition}

\subsection{Examples of graphs achieving equality in Theorem~\ref{thm:main2}}

In this section we present a small sample of graphs achieving equality in the upper bound of Theorem~\ref{thm:main2}.

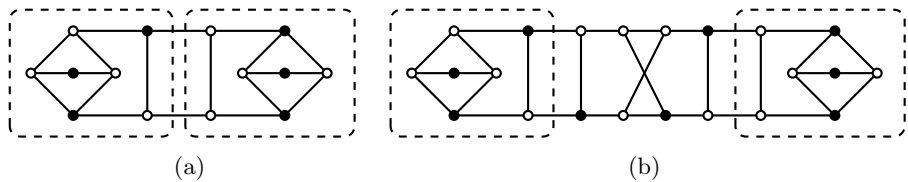
\begin{figure}[htb]
\begin{center}
\begin{tikzpicture}[scale=.75,style=thick,x=0.75cm,y=0.75cm]
\def\vr{2.25pt}
\path (-2,1) coordinate (v1);
\path (-1,0) coordinate (v2);
\path (-1,1) coordinate (v3);
\path (-1,2) coordinate (v4);
\path (0,1) coordinate (v5);
\path (0.75,0) coordinate (w1);
\path (0.75,2) coordinate (w2);
\path (2.25,0) coordinate (z1);
\path (2.25,2) coordinate (z2);
\path (3,1) coordinate (u1);
\path (4,0) coordinate (u2);
\path (4,1) coordinate (u3);
\path (4,2) coordinate (u4);
\path (5,1) coordinate (u5);
\draw (v1)--(v2)--(v5)--(v4)--(v1);
\draw (v1)--(v3)--(v5);
\draw (u1)--(u2)--(u5)--(u4)--(u1);
\draw (u1)--(u3)--(u5);
\draw (v2)--(w1)--(z1)--(u2);
\draw (v4)--(w2)--(z2)--(u4);
\draw (w1)--(w2);
\draw (z1)--(z2);
\draw [style=dashed,rounded corners] (-2.5,-0.5) rectangle (1.35,2.5);
\draw [style=dashed,rounded corners] (1.65,-0.5) rectangle (5.65,2.5);
%
\draw (v1) [fill=white] circle (\vr);
\draw (v2) [fill=black] circle (\vr);
\draw (v3) [fill=black] circle (\vr);
\draw (v4) [fill=white] circle (\vr);
\draw (v5) [fill=white] circle (\vr);
\draw (w1) [fill=white] circle (\vr);
\draw (w2) [fill=black] circle (\vr);
\draw (z1) [fill=white] circle (\vr);
\draw (z2) [fill=white] circle (\vr);
\draw (u1) [fill=white] circle (\vr);
\draw (u2) [fill=black] circle (\vr);
\draw (u3) [fill=black] circle (\vr);
\draw (u4) [fill=black] circle (\vr);
\draw (u5) [fill=white] circle (\vr);
\draw (1.75,-1.25) node {{\small (a)}};
\path (7,1) coordinate (v1);
\path (8,0) coordinate (v2);
\path (8,1) coordinate (v3);
\path (8,2) coordinate (v4);
\path (9,1) coordinate (v5);
\path (9.75,0) coordinate (w1);
\path (9.75,2) coordinate (w2);
\path (11,0) coordinate (z1);
\path (11,2) coordinate (z2);
\path (12,0) coordinate (x1);
\path (13,0) coordinate (x2);
\path (14,0) coordinate (x3);
\path (15.25,0) coordinate (x4);
\path (12,2) coordinate (y1);
\path (13,2) coordinate (y2);
\path (14,2) coordinate (y3);
\path (15.25,2) coordinate (y4);
\path (16,1) coordinate (u1);
\path (17,0) coordinate (u2);
\path (17,1) coordinate (u3);
\path (17,2) coordinate (u4);
\path (18,1) coordinate (u5);
\draw [style=dashed,rounded corners] (6.5,-0.5) rectangle (10.35,2.5);
\draw [style=dashed,rounded corners] (14.65,-0.5) rectangle (18.65,2.5);
%
\draw (v1)--(v2)--(v5)--(v4)--(v1);
\draw (v1)--(v3)--(v5);
\draw (u1)--(u2)--(u5)--(u4)--(u1);
\draw (u1)--(u3)--(u5);
\draw (v2)--(w1)--(z1)--(x1)--(x2)--(x3)--(x4)--(u2);
\draw (v4)--(w2)--(z2)--(y1)--(y2)--(y3)--(y4)--(u4);
\draw (w1)--(w2);
\draw (x1)--(y2);
\draw (x2)--(y1);
\draw (x3)--(y3);
\draw (x4)--(y4);
\draw (z1)--(z2);
\draw (v1) [fill=white] circle (\vr);
\draw (v2) [fill=black] circle (\vr);
\draw (v3) [fill=black] circle (\vr);
\draw (v4) [fill=white] circle (\vr);
\draw (v5) [fill=white] circle (\vr);
\draw (w1) [fill=white] circle (\vr);
\draw (w2) [fill=black] circle (\vr);
\draw (z1) [fill=black] circle (\vr);
\draw (z2) [fill=white] circle (\vr);
\draw (x1) [fill=white] circle (\vr);
\draw (x2) [fill=black] circle (\vr);
\draw (x3) [fill=white] circle (\vr);
\draw (x4) [fill=white] circle (\vr);
\draw (y1) [fill=white] circle (\vr);
\draw (y2) [fill=white] circle (\vr);
\draw (y3) [fill=black] circle (\vr);
\draw (y4) [fill=white] circle (\vr);
\draw (u1) [fill=white] circle (\vr);
\draw (u2) [fill=black] circle (\vr);
\draw (u3) [fill=black] circle (\vr);
\draw (u4) [fill=black] circle (\vr);
\draw (u5) [fill=white] circle (\vr);
\draw (12.5,-1.25) node {{\small (b)}};
\end{tikzpicture}
\vskip -0.25 cm
\caption{Examples of graphs achieving equality in the upper bound of Theorem~\ref{thm:main2}}
\label{fig:example-1}
\end{center}
\end{figure}

\noindent
\textbf{Example~1.} If $G$ is the graph shown in Figure~\ref{fig:example-1}(a), then $i(G) = 6$ and the shaded vertices indicate an $i$-set in~$G$. In this example, $\w(G) = 12 \times 3 + 2 \times 4 = 44$ and $\Theta(G) = 2\tc(G) = 2 \times 2 = 4$ (where the two vertex disjoint $T_2$-configurations are indicated by the subgraphs in the dashed boxes), and so $8i(G) = 48 = 44 + 4 = \w(G) + \Theta(G) = \Omega(G)$.

Moreover if $G$ is the graph shown in Figure~\ref{fig:example-1}(b), then $i(G) = 9$ and the shaded vertices indicate an $i$-set in~$G$. In this example, $\w(G) = 20 \times 3 + 2 \times 4 = 68$ and $\Theta(G) = 2\tc(G) = 2 \times 2 = 4$ (where the two vertex disjoint $T_2$-configurations are indicated by the subgraphs in the dashed boxes), and so $8i(G) = 72 = 44 + 4 = \w(G) + \Theta(G) = \Omega(G)$.

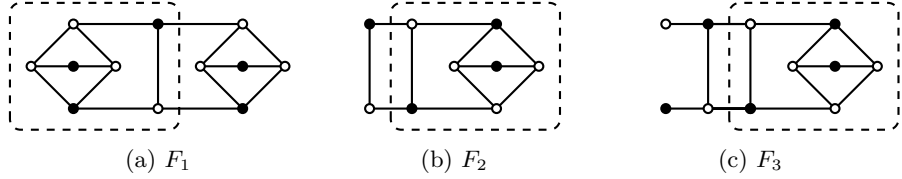
\begin{figure}[ht]
\begin{center}
\begin{tikzpicture}[scale=.75,style=thick,x=0.75cm,y=0.75cm]
\def\vr{2.25pt}
\path (0,1) coordinate (v1);
\path (1,0) coordinate (v2);
\path (1,1) coordinate (v3);
\path (1,2) coordinate (v4);
\path (2,1) coordinate (v5);
\path (3,0) coordinate (w1);
\path (3,2) coordinate (w2);
\path (4,1) coordinate (u1);
\path (5,0) coordinate (u2);
\path (5,1) coordinate (u3);
\path (5,2) coordinate (u4);
\path (6,1) coordinate (u5);
\draw (v1)--(v2)--(v5)--(v4)--(v1);
\draw (v1)--(v3)--(v5);
\draw (u1)--(u2)--(u5)--(u4)--(u1);
\draw (u1)--(u3)--(u5);
\draw (v2)--(w1)--(u2);
\draw (v4)--(w2)--(u4);
\draw (w1)--(w2);
\draw [style=dashed,rounded corners] (-0.5,-0.5) rectangle (3.5,2.5);
\draw (v1) [fill=white] circle (\vr);
\draw (v2) [fill=black] circle (\vr);
\draw (v3) [fill=black] circle (\vr);
\draw (v4) [fill=white] circle (\vr);
\draw (v5) [fill=white] circle (\vr);
\draw (w1) [fill=white] circle (\vr);
\draw (w2) [fill=black] circle (\vr);
\draw (u1) [fill=white] circle (\vr);
\draw (u2) [fill=black] circle (\vr);
\draw (u3) [fill=black] circle (\vr);
\draw (u4) [fill=white] circle (\vr);
\draw (u5) [fill=white] circle (\vr);
\draw (3,-1.25) node {{\small (a) $F_1$}};
\path (8,0) coordinate (v1);
\path (8,2) coordinate (v2);
\path (9,0) coordinate (w1);
\path (9,2) coordinate (w2);
\path (10,1) coordinate (u1);
\path (11,0) coordinate (u2);
\path (11,1) coordinate (u3);
\path (11,2) coordinate (u4);
\path (12,1) coordinate (u5);
\draw (w1)--(v1)--(v2)--(w2);
\draw (u1)--(u2)--(u5)--(u4)--(u1);
\draw (u1)--(u3)--(u5);
\draw (w1)--(u2);
\draw (w2)--(u4);
\draw (w1)--(w2);
\draw [style=dashed,rounded corners] (8.5,-0.5) rectangle (12.5,2.5);
\draw (v1) [fill=white] circle (\vr);
\draw (v2) [fill=black] circle (\vr);
\draw (w1) [fill=black] circle (\vr);
\draw (w2) [fill=white] circle (\vr);
\draw (u1) [fill=white] circle (\vr);
\draw (u2) [fill=white] circle (\vr);
\draw (u3) [fill=black] circle (\vr);
\draw (u4) [fill=black] circle (\vr);
\draw (u5) [fill=white] circle (\vr);
\draw (10,-1.25) node {{\small (b)  $F_{2}$}};
\path (15,0) coordinate (v01);
\path (16,0) coordinate (v1);
\path (17,0) coordinate (v2);
\path (17,2) coordinate (v3);
\path (16,2) coordinate (v4);
\path (15,2) coordinate (v04);
\path (17,0) coordinate (w1);
\path (17,2) coordinate (w2);
\path (18,1) coordinate (u1);
\path (19,0) coordinate (u2);
\path (19,1) coordinate (u3);
\path (19,2) coordinate (u4);
\path (20,1) coordinate (u5);
\draw (w1)--(v1)--(v2);
\draw (v3)--(v4)--(w2);
\draw (v1)--(v4);
\draw (u1)--(u2)--(u5)--(u4)--(u1);
\draw (u1)--(u3)--(u5);
\draw (w1)--(u2);
\draw (w2)--(u4);
\draw (w1)--(w2);
\draw (v1)--(v01);
\draw (v4)--(v04);
\draw [style=dashed,rounded corners] (16.5,-0.5) rectangle (20.5,2.5);
\draw (v01) [fill=black] circle (\vr);
\draw (v04) [fill=white] circle (\vr);
\draw (v1) [fill=white] circle (\vr);
\draw (v2) [fill=black] circle (\vr);
\draw (v3) [fill=white] circle (\vr);
\draw (v4) [fill=black] circle (\vr);
\draw (w1) [fill=black] circle (\vr);
\draw (w2) [fill=white] circle (\vr);
\draw (u1) [fill=white] circle (\vr);
\draw (u2) [fill=white] circle (\vr);
\draw (u3) [fill=black] circle (\vr);
\draw (u4) [fill=black] circle (\vr);
\draw (u5) [fill=white] circle (\vr);
\draw (17,-1.25) node {{\small (c) $F_{3}$}};
\end{tikzpicture}
\caption{Examples of graphs achieving equality in the upper bound of Theorem~\ref{thm:main2}}
\label{fig:example-2}
\end{center}
\end{figure}

\noindent
\textbf{Example~2.} If $G$ is the graph shown in Figure~\ref{fig:example-2}(a), then $i(G) = 5$ and the shaded vertices indicate an $i$-set in~$G$. In this example, $\w(G) = 10 \times 3 + 2 \times 4 = 38$ and $\Theta(G) = 2\tc(G) = 2 \times 1 = 2$ (where a $T_2$-configuration is indicated by the subgraph in the dashed box), and so $8i(G) = 40 = 38 + 2 = \w(G) + \Theta(G) = \Omega(G)$.

If $G$ is the graph shown in Figure~\ref{fig:example-2}(b), then $i(G) = 4$ and the shaded vertices indicate an $i$-set in~$G$. In this example, $\w(G) = 6 \times 3 + 3 \times 4 = 30$ and $\Theta(G) = 2\tc(G) = 2 \times 1 = 2$ (where a $T_2$-configuration is indicated by the subgraph in the dashed box), and so $8i(G) = 32 = 30 + 2 = \w(G) + \Theta(G) = \Omega(G)$.

If $G$ is the graph shown in Figure~\ref{fig:example-2}(c), then $i(G) = 5$ and the shaded vertices indicate an $i$-set in~$G$. In this example, $\w(G) = 8 \times 3 + 4 + 2 \times 5 = 38$ and $\Theta(G) = 2\tc(G) = 2 \times 1 = 2$ (where a $T_2$-configuration is indicated by the subgraph in the dashed box), and so $8i(G) = 40 = 38 + 2 = \w(G) + \Theta(G) = \Omega(G)$.

\begin{figure}[htb]
\begin{center}
\begin{tikzpicture}[scale=.75,style=thick,x=0.75cm,y=0.75cm]
\def\vr{2.25pt}
\def\vrn{1.5pt}
\path (0,1) coordinate (a1);
\path (0,2.75) coordinate (a2);
\path (0,4.25) coordinate (a3);
\path (0,6) coordinate (a4);
\path (1,0) coordinate (b1);
\path (1,1) coordinate (b2);
\path (1,2) coordinate (b3);
\path (1,5) coordinate (b4);
\path (1,6) coordinate (b5);
\path (1,7) coordinate (b6);
\path (2,1) coordinate (c1);
\path (2,2.75) coordinate (c2);
\path (2,4.25) coordinate (c3);
\path (2,6) coordinate (c4);
\path (5,-1) coordinate (d1);
\path (5,0) coordinate (d2);
\path (5,1.75) coordinate (d3);
\path (5,5.25) coordinate (d4);
\path (5,7) coordinate (d5);
\path (5,8) coordinate (d6);
\path (6,0.75) coordinate (e1);
\path (6,1.75) coordinate (e2);
\path (6,2.75) coordinate (e3);
\path (6,4.25) coordinate (e4);
\path (6,5.25) coordinate (e5);
\path (6,6.25) coordinate (e6);
\path (7,-1) coordinate (f1);
\path (7,0) coordinate (f2);
\path (7,1.75) coordinate (f3);
\path (7,5.25) coordinate (f4);
\path (7,7) coordinate (f5);
\path (7,8) coordinate (f6);
\path (10,1) coordinate (g1);
\path (10,2.75) coordinate (g2);
\path (10,4.25) coordinate (g3);
\path (10,6) coordinate (g4);
\path (11,0) coordinate (h1);
\path (11,1) coordinate (h2);
\path (11,2) coordinate (h3);
\path (11,5) coordinate (h4);
\path (11,6) coordinate (h5);
\path (11,7) coordinate (h6);
\path (12,1) coordinate (i1);
\path (12,2.75) coordinate (i2);
\path (12,4.25) coordinate (i3);
\path (12,6) coordinate (i4);
\draw (a1)--(a2)--(a3)--(a4);
\draw (b1)--(b2)--(b3);
\draw (b4)--(b5)--(b6);
\draw (c1)--(c2)--(c3)--(c4);
\draw (a1)--(b1)--(c1)--(b3)--(a1);
\draw (a2)--(c2);
\draw (a3)--(c3);
\draw (a4)--(b4)--(c4)--(b6)--(a4);
\draw (d1)--(d2)--(d3);
\draw (d4)--(d5)--(d6);
\draw (e1)--(e2)--(e3);
\draw (e4)--(e5)--(e6);
\draw (f1)--(f2)--(f3);
\draw (f4)--(f5)--(f6);
\draw (d1)--(f1);
\draw (d2)--(f2);
\draw (d5)--(f5);
\draw (d6)--(f6);
\draw (d3)--(e1)--(f3)--(e3)--(d3);
\draw (d4)--(e4)--(f4)--(e6)--(d4);
\draw (g1)--(g2)--(g3)--(g4);
\draw (h1)--(h2)--(h3);
\draw (h4)--(h5)--(h6);
\draw (i1)--(i2)--(i3)--(i4);
\draw (g1)--(h1)--(i1)--(h3)--(g1);
\draw (g2)--(i2);
\draw (g3)--(i3);
\draw (g4)--(h4)--(i4)--(h6)--(g4);
\draw (d1) to[out=180,in=-55, distance=1cm] (b2);
\draw (f1) to[out=0,in=235, distance=1cm] (h2);
\draw (d6) to[out=180,in=75, distance=1cm] (b5);
\draw (f6) to[out=0,in=110, distance=1cm] (h5);
\draw [style=dashed,rounded corners] (3.75,-0.5) rectangle (8,3.25);
\draw [style=dashed,rounded corners] (3.75,3.75) rectangle (8,7.5);
\draw (a1) [fill=white] circle (\vr);
\draw (a2) [fill=black] circle (\vr);
\draw (a3) [fill=white] circle (\vr);
\draw (a4) [fill=black] circle (\vr);
\draw (b1) [fill=black] circle (\vr);
\draw (b2) [fill=white] circle (\vr);
\draw (b3) [fill=black] circle (\vr);
\draw (b4) [fill=white] circle (\vr);
\draw (b5) [fill=white] circle (\vr);
\draw (b6) [fill=white] circle (\vr);
\draw (c1) [fill=white] circle (\vr);
\draw (c2) [fill=white] circle (\vr);
\draw (c3) [fill=white] circle (\vr);
\draw (c4) [fill=black] circle (\vr);
\draw (d1) [fill=white] circle (\vr);
\draw (d2) [fill=black] circle (\vr);
\draw (d3) [fill=white] circle (\vr);
\draw (d4) [fill=black] circle (\vr);
\draw (d5) [fill=white] circle (\vr);
\draw (d6) [fill=black] circle (\vr);
\draw (e1) [fill=white] circle (\vr);
\draw (e2) [fill=black] circle (\vr);
\draw (e3) [fill=white] circle (\vr);
\draw (e4) [fill=white] circle (\vr);
\draw (e5) [fill=black] circle (\vr);
\draw (e6) [fill=white] circle (\vr);
\draw (f1) [fill=white] circle (\vr);
\draw (f2) [fill=white] circle (\vr);
\draw (f3) [fill=black] circle (\vr);
\draw (f4) [fill=white] circle (\vr);
\draw (f5) [fill=black] circle (\vr);
\draw (f6) [fill=white] circle (\vr);
\draw (g1) [fill=white] circle (\vr);
\draw (g2) [fill=black] circle (\vr);
\draw (g3) [fill=white] circle (\vr);
\draw (g4) [fill=white] circle (\vr);
\draw (h1) [fill=white] circle (\vr);
\draw (h2) [fill=black] circle (\vr);
\draw (h3) [fill=white] circle (\vr);
\draw (h4) [fill=black] circle (\vr);
\draw (h5) [fill=white] circle (\vr);
\draw (h6) [fill=black] circle (\vr);
\draw (i1) [fill=black] circle (\vr);
\draw (i2) [fill=white] circle (\vr);
\draw (i3) [fill=black] circle (\vr);
\draw (i4) [fill=white] circle (\vr);
\end{tikzpicture}
\caption{An example of a graph achieving equality in the upper bound of Theorem~\ref{thm:main2}}
\label{fig:example-3}
\end{center}
\end{figure}
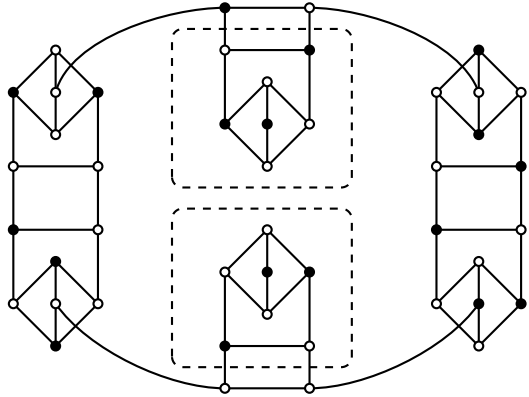

\noindent
\textbf{Example~3.} If $G$ is the graph shown in Figure~\ref{fig:example-3}, then $i(G) = 18$ and the shaded vertices indicate an $i$-set in~$G$. In this example, $\w(G) = 44 \times 3 + 2 \times 4 = 140$ and $\Theta(G) = 2\tc(G) = 2 \times 2 = 4$ (where the $T_2$-configurations are indicated by the subgraphs in the dashed boxes), and so $8i(G) = 144 = 140 + 4 = \w(G) + \Theta(G) = \Omega(G)$.

\section{Proof of key result}
\label{S:main-proof}

Recall that the vertex weight of a subcubic graph $G$ is $\w(G) = 8n_0(G) + 5n_1(G) + 4n_2(G) + 3n_3(G)$, and the structural weight of $G$ is $\Theta(G) = 2\tc(G) + 2b(G)$. Moreover the total weight of $G$, called the $\Omega$-weight of $G$, is the sum of its vertex weight and structural weight, that is, $\Omega(G) = \w(G) + \Theta(G)$. In this section, we present a proof of our key result, namely Theorem~\ref{thm:main2}. Its statement in terms of the $\Omega$-weight of $G$ is as follows:

\medskip
\noindent \textbf{Theorem~\ref{thm:main2}}. \emph{If $G$ is a subcubic graph of order~$n$ that contains no $K_{3,3}$-component and no $(C_5 \, \Box \, K_2)$-component, then $8i(G) \le \Omega(G)$.
}

\noindent
\proof Let $G$ be a subcubic graph of order~$n$ that contains no $K_{3,3}$-component and no $(C_5 \, \Box \, K_2)$-component, and let $V = V(G)$. We wish to prove that $8i(G) \le \Omega(G)$. Suppose that $G$ is a counterexample to the inequality with minimum order, and let $G$ have order~$n$. Thus, $8i(G) > \Omega(G)$ but every subcubic graph $G'$ of order less than~$n$ that contains no $K_{3,3}$-component and no $(C_5 \, \Box \, K_2)$-component satisfies $8i(G') \le \Omega(G')$. Since $\Theta(G) \ge 0$, we note that $\w(G) \le \Omega(G)$. Hence, if $8i(G) \le \w(G)$, then $8i(G) \le \Omega(G)$, a contradiction. Therefore, $8i(G) > \w(G)$.

Since our proof of Theorem~\ref{thm:main2} is long, we break the proof into four parts to enhance clarity and exposition of the proof. The first part of the proof establishes important structural properties of the counterexample $G$. The second part of the proof considers the case when the minimum degree of $G$ equals~$1$. The third part of the proof considers the case when the minimum degree of $G$ equals~$2$. The fourth and final part of the proof considers the case when the graph $G$ is cubic, that is, when $G$ is a $3$-regular graph.

\subsection{Part 1: Structural properties of $G$}
\label{S:part1}

In this section, namely Part~$1$ of our proof, we establish structural properties of the counterexample~$G$. A component of a graph that belongs to the family $\cB$ is a \emph{bad component} of the graph.

\begin{claim}
\label{claim:K23}
The subcubic graph $G$ contains at least one subgraph isomorphic to $K_{2,3}$.
\end{claim}
\proof If $G$ contains no subgraph isomorphic to $K_{2,3}$, then by Theorem~\ref{thm:Dorbec1}, we have $8i(G) \le \w(G)$, a contradiction.~\smallqed

\medskip
By Claim~\ref{claim:K23}, we note that $n \ge 5$. If $n = 5$, then either $G = K_{2,3}$, in which case $8i(G) = 16 < 18 = \w(G)$, or $G$ is obtained from $K_{2,3}$ by adding an edge, in which case $8i(G) = 16 = \w(G)$, contradicting the fact that $G$ is a counterexample. Hence, $n \ge 6$.

\begin{claim}
\label{claim:Gconn}
The graph $G$ is connected.
\end{claim}
\proof If $G$ is not connected, then, by the minimality of $n(G)$, the theorem holds for every component of $G$, and therefore also for $G$, a contradiction.~\smallqed

\begin{claim}
\label{claim:no-bad}
$G \notin \cB$.
\end{claim}
\proof If $G \in \cB$, then $b(G) = 1$ and $\Theta(G) = 2b(G) = 2$, and so $\Omega(G) = \w(G) + 2$. By Proposition~\ref{prop:cB}(a), $8i(G) = \w(G) + 2 = \w(G) + \Theta(G) = \Omega(G)$, a contradiction.~\smallqed

\medskip
By Claim~\ref{claim:no-bad}, we have $b(G) = 0$, and so $\Theta(G) = 2\tc(G)$.

\begin{claim}
\label{claim:no-T-subgraph}
The graph $G$ contains no troublesome configuration.
\end{claim}
\proof Suppose, to the contrary, that $G$ contains a troublesome configuration, say $T$. Thus, $T = T_{v_1,v_2}$ is obtained from a bad graph $B \in \cB_1$ with root vertex~$v_1$ (of degree~$1$ in $B$) by adding a new vertex~$v_2$ to~$B$, adding the edges $v_1v_2$ and $v_2w_2$, where $w_2$ is a vertex of degree~$2$ in $B$ different from~$v_1$. Furthermore, there are at most two exit edges that join $T$ to vertices in $V(G) \setminus V(T)$ via one edge incident with~$v_2$ and possibly one additional exit edge incident with~$v_1$.  Let $u$ be the vertex in $B$ adjacent to~$v_1$ in $G$.

Let $B_u$ be the subgraph $B - v_1$, and let $B_u$ contain $k \ge 1$ vertex disjoint copies of $K_{2,3}$. The weight of $B_u$ is $\w(B_u) = 16k + 2$. Moreover, $i(B_u) = 2k$, which follows by induction and using Lemma~\ref{lem:lem1}. Indeed, the non-canonical independent set of $B$, which we denote by $I_u$, is an $i$-set of $B_u$. Let $G' = G - V(B_u)$. We note that $v_1v_2$ is an edge of $G'$. Further, $v_1$ has degree~$1$ or~$2$ in $G'$ and $v_2$ has degree~$2$ in $G'$. Hence every $i$-set of $G'$ can be extended to an ID-set of $G$ by adding to it the $2k$ vertices from the set $I_u$, and so $i(G) \le i(G') + |I_u| = i(G') + 2k$. Moreover, $\w(G) = \w(G') + \w(B_u) - 4$ noting that the two edges $uv_1$ and $v_2w_2$ joining $V(B_u)$ to $V(G')$ decrease the weights of each of the vertices $u$, $v_1$, $v_2$ and $w_2$ by~$1$. Thus, $\w(G') = \w(G) - \w(B_u) + 4 = \w(G) - (16k+2) + 4 = \w(G) - 16k + 2$.

We do not create a bad component in $G'$ since such a component would contain the adjacent vertices $v_1$ and $v_2$ where we note that both $v_1$ and $v_2$ have degree at most~$2$ in $G'$. This contradicts the property that at least one of every two adjacent vertices in a bad component is a vertex of degree~$3$. We therefore infer that $b(G') = 0$.

We note that the troublesome configuration $T$ of $G$ is no longer a troublesome configuration of $G'$. Further we note that $v_1$ and $v_2$ are adjacent vertices in $G'$ where $v_1$ has degree~$1$ or~$2$ and $v_2$ has degree~$2$ in $G'$. Any new troublesome configuration $T'$ of $G'$ that is not a troublesome configuration of $G$ necessarily contains both $v_1$ and $v_2$ as its link vertices. However in this case, neither $v_1$ nor $v_2$ is adjacent to a vertex not in the $T'$-configuration noting that both $v_1$ and $v_2$ have the same degrees in $T'$ and $G'$. Hence we do not create a new troublesome configuration. Thus, $\tc(G') = \tc(G) - 1$. As observed earlier, $b(G') = 0$. Thus, $\Theta(G') = 2\tc(G') = 2\tc(G) - 2 = \Theta(G) - 2$, implying that $\Omega(G') = \w(G') + \Theta(G') \le (\w(G) - 16k + 2) + (\Theta(G) - 2) = \Omega(G) - 16k$. Hence, $8i(G) \le 8(i(G') + 2k) \le \Omega(G') + 16k \le \Omega(G)$, a contradiction.~\smallqed

\begin{claim}
\label{claim:theta-zero}
$\Theta(G) = 0$.
\end{claim}
\proof By Claim~\ref{claim:no-bad}, $b(G) = 0$. By Claim~\ref{claim:no-T-subgraph}, there are no troublesome configurations in $G$, implying that $\tc(G) = 0$. Hence, $\Theta(G) = b(G) + \tc(G) = 0$.~\smallqed

\begin{claim}
\label{claim:no-K23-two-deg2}
The graph $G$ does not contain a $K_{2,3}$-subgraph that contains two vertices of degree~$2$ in $G$.
\end{claim}
\proof Suppose, to the contrary, that $G$ contains a $K_{2,3}$-subgraph, $F$, that contains two vertices of degree~$2$ in $G$. Let $u_1$ and $u_2$ be the two vertices of degree~$2$ in $G$ that belong to $F$, and let $w_1$ and $w_2$ be the two common neighbors of $u_1$ and $u_2$. Further, let $u_3$ be the third common neighbor of $w_1$ and $w_2$. We note that the vertex $u_3$ has degree~$2$ in $F$ and has degree~$3$ in $G$. Let $v_1$ be the third neighbor of $u_3$. Let $G' = G - V(F)$. Every $i$-set of $G'$ can be extended to an ID-set of $G$ by adding to it the vertices $w_1$ and $w_2$, and so $i(G) \le i(G') + 2$.

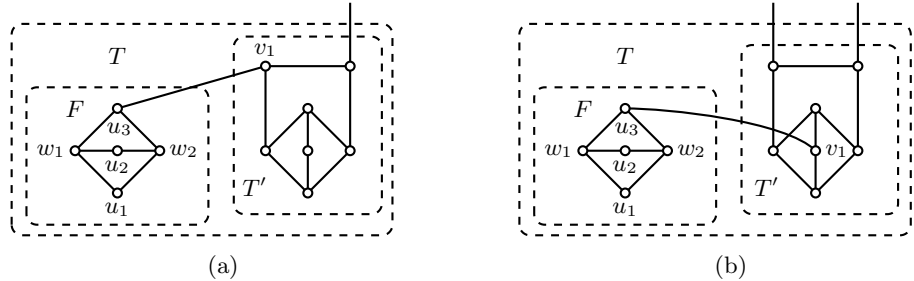
\begin{figure}[htb]
\begin{center}
\begin{tikzpicture}[scale=.75,style=thick,x=0.75cm,y=0.75cm]
\def\vr{2.25pt}
\def\vrn{1.25pt}
\path (-1.5,4.5) coordinate (b1);
\path (-0.5,3.5) coordinate (b2);
\path (-0.5,4.5) coordinate (b3);
\path (-0.5,5.5) coordinate (b4);
\path (0.5,4.5) coordinate (b5);
\path (-1.5,6.5) coordinate (e1);
\path (0.5,6.5) coordinate (e5);
\path (0.5,8) coordinate (e51);
\path (-6,4.5) coordinate (a1);
\path (-5,3.5) coordinate (a2);
\path (-5,4.5) coordinate (a3);
\path (-5,5.5) coordinate (a4);
\path (-4.95,5.4) coordinate (a4p);
\path (-4,4.5) coordinate (a5);
%
%
\draw (b1)--(e1);
\draw (e1)--(a4);
\draw (b5)--(e5)--(e51);
\draw (e1)--(e5);
\draw (b1)--(b2)--(b3)--(b4)--(b5)--(b2);
\draw (b1)--(b4);
\draw (a1)--(a2)--(a5)--(a4)--(a1);
\draw (a1)--(a3)--(a5);
\draw [style=dashed,rounded corners] (-7.5,2.5) rectangle (1.5,7.5);
\draw [style=dashed,rounded corners] (-7.15,2.75) rectangle (-2.85,6);
\draw [style=dashed,rounded corners] (-2.25,3) rectangle (1.25,7.2);
\draw (e1) [fill=white] circle (\vr);
\draw (e5) [fill=white] circle (\vr);
\draw (b1) [fill=white] circle (\vr);
\draw (b2) [fill=white] circle (\vr);
\draw (b3) [fill=white] circle (\vr);
\draw (b4) [fill=white] circle (\vr);
\draw (b5) [fill=white] circle (\vr);
\draw (a1) [fill=white] circle (\vr);
\draw (a2) [fill=white] circle (\vr);
\draw (a3) [fill=white] circle (\vr);
\draw (a4) [fill=white] circle (\vr);
\draw (a5) [fill=white] circle (\vr);
\draw[anchor = south] (e1) node {{\small $v_1$}};
\draw[anchor = east] (a1) node {{\small $w_1$}};
\draw[anchor = west] (a5) node {{\small $w_2$}};
\draw[anchor = north] (a2) node {{\small $u_1$}};
\draw[anchor = north] (a3) node {{\small $u_2$}};
\draw[anchor = north] (a4p) node {{\small $u_3$}};
\draw (-6,5.5) node {{\small $F$}};
\draw (-5,6.75) node {{\small $T$}};
\draw (-1.75,3.65) node {{\small $T'$}};
\path (10.5,4.5) coordinate (b1);
\path (11.5,3.5) coordinate (b2);
\path (11.5,4.5) coordinate (b3);
\path (11.5,5.5) coordinate (b4);
\path (12.5,4.5) coordinate (b5);
\path (10.5,6.5) coordinate (e1);
\path (10.5,8) coordinate (e11);
\path (12.5,6.5) coordinate (e5);
\path (12.5,8) coordinate (e51);
\path (6,4.5) coordinate (a1);
\path (7,3.5) coordinate (a2);
\path (7,4.5) coordinate (a3);
\path (7,5.5) coordinate (a4);
\path (7.05,5.4) coordinate (a4p);
\path (8,4.5) coordinate (a5);
%
%
\draw (b1)--(e1);
\draw (e1)--(e11);
\draw (b5)--(e5)--(e51);
\draw (e1)--(e5);
\draw (b1)--(b2)--(b3)--(b4)--(b5)--(b2);
\draw (b1)--(b4);
\draw (a1)--(a2)--(a5)--(a4)--(a1);
\draw (a1)--(a3)--(a5);
\draw (a4) to[out=0,in=135, distance=0.75cm] (b3);
%
\draw [style=dashed,rounded corners] (4.5,2.5) rectangle (13.75,7.5);
\draw [style=dashed,rounded corners] (4.85,2.75) rectangle (9.15,6);
\draw [style=dashed,rounded corners] (9.75,3) rectangle (13.45,7);
\draw (e1) [fill=white] circle (\vr);
\draw (e5) [fill=white] circle (\vr);
\draw (b1) [fill=white] circle (\vr);
\draw (b2) [fill=white] circle (\vr);
\draw (b3) [fill=white] circle (\vr);
\draw (b4) [fill=white] circle (\vr);
\draw (b5) [fill=white] circle (\vr);
\draw (a1) [fill=white] circle (\vr);
\draw (a2) [fill=white] circle (\vr);
\draw (a3) [fill=white] circle (\vr);
\draw (a4) [fill=white] circle (\vr);
\draw (a5) [fill=white] circle (\vr);
\draw[anchor = west] (b3) node {{\small $v_1$}};
\draw[anchor = east] (a1) node {{\small $w_1$}};
\draw[anchor = west] (a5) node {{\small $w_2$}};
\draw[anchor = north] (a2) node {{\small $u_1$}};
\draw[anchor = north] (a3) node {{\small $u_2$}};
\draw[anchor = north] (a4p) node {{\small $u_3$}};
\draw (6,5.5) node {{\small $F$}};
\draw (7,6.75) node {{\small $T$}};
\draw (10.35,3.65) node {{\small $T'$}};
%
\draw (-2.5,1.75) node {{\small (a)}};
\draw (9.5,1.75) node {{\small (b)}};
\end{tikzpicture}
\vskip -0.25 cm
\caption{Possible subgraphs in the proof of Claim~\ref{claim:no-K23-two-deg2}}
\label{fig:no-K23-two-deg2}
\end{center}
\end{figure}

We note that $\w(G) = \w(G') + 16$. By Claim~\ref{claim:theta-zero}, $\Theta(G) = 0$ and therefore $\Omega(G) = \w(G) = \w(G') + 16$.  We show next that $\Theta(G') = 0$. If $G'$ belongs to $\cB$, then so too does the original graph $G$. Thus, $G' \notin \cB$, and so $b(G') = 0$. By Claim~\ref{claim:theta-zero}, the graph $G$ contains no troublesome configuration. Suppose that $G'$ contains a troublesome configuration $T'$. Necessarily such a subgraph $T'$ contains the vertex~$v_1$. Note that $v_1$ is not a link vertex of $T'$, since this would imply that $T'$ is troublesome configuration already in $G$ (see Figure~\ref{fig:no-K23-two-deg2}(a)). Hence $v_1$ is not a link vertex, see Figure~\ref{fig:no-K23-two-deg2}(b)
as an illustration in this case. But then adding the subgraph $F$ to $T'$, together with the edge $u_3v_1$, produces a troublesome configuration $T$ of $G$, and so $\tc(G) \ge 1$, contradicting Claim~\ref{claim:theta-zero}. Hence, $\tc(G') = 0$. Thus, $\Theta(G') = 0$, and so $\Omega(G') = \w(G') = \Omega(G) - 16$. As observed earlier, $i(G) \le i(G') + 2$. Hence, $8i(G) \le 8(i(G') + 2) \le \Omega(G') + 16 = \Omega(G)$, a contradiction.~\smallqed

We shall frequently use the following property that removing an exit edge increases the total weight by at most~$3$.

\begin{claim}
\label{claim:cost-of-exit-edge}
If $X \subset V(G)$ and $G' = G - X$, then removing an $X$-exit edge when constructing $G'$ increases the total weight by at most~$3$.
\end{claim}
\proof Let $X \subset V(G)$ and $G' = G - X$, and consider an exit edge $vv' \in E(G)$ where $v \in X$ and $v' \in V(G')$. If $v'$ is isolated in $G'$, then the vertex~$v'$ has degree~$1$ in $G$ and degree~$0$ in $G'$, and so $\w_{G'}(v') - \w_G(v') = 8 - 5 = 3$. If $\deg_{G'}(v') \ge 1$, then $\w_{G'}(v') - \w_G(v') = 1$. However, in this case removing the edge $vv'$ may create a bad component or a troublesome configuration, resulting in an increase in the structural weight by~$2$. Thus the total weight increases by~$3$ if either we create an isolated vertex - in this case, the vertex weight increases by~$3$ and the structural weight by~$0$ - or if we create a bad component or a new troublesome configuration that contains the vertex incident with the exit edge - in this case, the vertex weight increases by~$1$ and the structural weight increases by~$2$, resulting in a total increase in the weight by~$3$. In both cases the removal of the $X$-exit edge $vv'$ increases the total weight by at most~$3$.~\smallqed

\medskip
By Claim~\ref{claim:no-K23-two-deg2}, every $K_{2,3}$-subgraph in $G$ contains at most one vertex of degree~$2$ in $G$.

\begin{claim}
\label{claim:no-K23-with-d2-vertex}
The graph $G$ does not contain a $K_{2,3}$-subgraph that contains a vertex of degree~$2$ in $G$.
\end{claim}
\proof Suppose, to the contrary, that $G$ contains a $K_{2,3}$-subgraph, $F$, that contains a vertex, say $u$, of degree~$2$ in $G$. Let $w_1$ and $w_2$ be the two neighbors of $u$ in $F$, and so $w_1$ and $w_2$ have degree~$3$ in $F$. Let $u_1$ and $u_2$ be the common neighbors of $w_1$ and $w_2$ in $F$ different from~$u$. Thus, $N_G(w_1) = N_G(w_2) = \{u,u_1,u_2\}$. By Claim~\ref{claim:no-K23-two-deg2}, both $u_1$ and $u_2$ have degree~$3$ in $G$. If $u_1u_2 \in E(G)$, then the graph $G$ is determined, and so $n = 5$, a contradiction. Hence, $u_1u_2 \notin E(G)$. Let $v_1$ and $v_2$ be the neighbors of $u_1$ and $u_2$, respectively, that do not belong to $F$.

We show firstly that $v_1 \ne v_2$. Suppose, to the contrary, that $v_1 = v_2$. If $\deg_G(v_1) = 2$, then the graph $G$ is determined. In this case, we have $i(G) = 2$, $\w(G) = 20$ and $\Theta(G) = 0$. Thus, $\Omega(G) = 20$, and so $8i(G) < \Omega(G)$, a contradiction.  Hence, $\deg_G(v_1) = 3$. Let $v$ be the neighbor of $v_1$ different from $u_1$ and $u_2$. If $\deg_G(v) = 1$, then the graph $G$ is determined, and as before $8i(G) < \Omega(G)$ (noting that $\{u,v_1\}$ is an $i$-set of $G$), a contradiction. Hence, $\deg_G(v) \ge 2$. We now let $X = V(F) \cup \{v,v_1\}$ and consider the graph $G' = G - X$.
We note that there are at most two $X$-exit edges in $G$ (and such edges are incident with~$v$). By Claim~\ref{claim:cost-of-exit-edge}, the removal of an $X$-exit edge when constructing $G'$ increases the total weight by at most~$3$. If there are two $X$-exit edges in $G$, then $\w_G(X) = 22$, and we infer that $\Omega(G') \le \Omega(G) - 22 + 6 = \Omega(G) - 16$. On the other hand, if there is only one $X$-exit edge in $G$, then $\w_G(X) = 23$, implying that $\Omega(G') \le \Omega(G) - 23 + 3 = \Omega(G) - 20$. Every $i$-set of $G'$ can be extended to an ID-set of $G$ by adding to it the vertices~$u$ and $v_1$, and so $i(G) \le i(G') + 2$. Hence, $8i(G) \le 8i(G') + 16 \le \Omega(G') + 16 \le \Omega(G)$, a contradiction. Hence, $v_1 \ne v_2$.

Suppose that $v_1v_2 \in E(G)$. If at least one of $v_1$ and $v_2$ has degree~$3$ in $G$, then the subgraph of $G$ induced by $V(F) \cup \{v_1,v_2\}$ is a troublesome configuration in $G$, contradicting Claim~\ref{claim:no-T-subgraph}. Hence, both $v_1$ and $v_2$ have degree~$2$ in $G$, and so the graph $G$ is determined. In this case, $i(G) = 3$ and $\Omega(G) = \w(G) = 24$, and so $8i(G) = \Omega(G)$, a contradiction. Hence, $v_1v_2 \notin E(G)$.

We now consider the graph $G'$ obtained from $G - V(F)$ by adding the edge $e = v_1v_2$. Let $I'$ be an $i$-set of $G'$. If neither $v_1$ nor $v_2$ belongs to $I'$, then let $I = I' \cup \{w_1,w_2\}$. If $v_1 \in I'$, then let $I = I' \cup \{u,u_2\}$. If $v_2 \in I'$, then let $I = I' \cup \{u,u_1\}$. In all cases, the set $I$ is an ID-set of $G$ and $|I| = |I'| + 2$. Thus, $i(G) \le |I| \le |I'| + 2 = i(G') + 2$. Moreover, $\w(G') = \w(G) - 16 = \Omega(G) - 16$. If $\Theta(G') = 0$, then $\Omega(G') = \w(G') =  \Omega(G) - 16$, and so $8i(G) \le 8(i(G') + 2) \le \Omega(G') + 16 = \Omega(G)$, a contradiction. Hence, $\Theta(G') > 0$, implying that either $G' \in \cB$ (noting that $G'$ is connected) or there is a troublesome configuration in $G'$ that contains the (adjacent) vertices $v_1$ and $v_2$. In both cases, we have $\Theta(G') = 2$, and so $\Omega(G') = \w(G') + \Theta(G') =  (\Omega(G) - 16) + 2 = \Omega(G) - 14$.

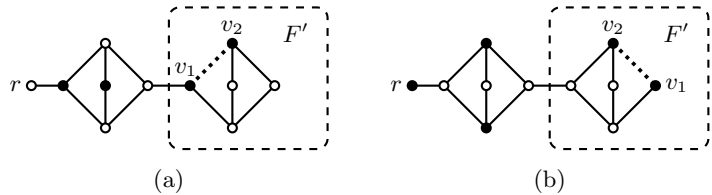
\begin{figure}[htb]
\begin{center}
\begin{tikzpicture}[scale=.75,style=thick,x=0.75cm,y=0.75cm]
\def\vr{2.25pt}
\path (3.25,1) coordinate (v6);
\path (4,1) coordinate (v1);
\path (5,0) coordinate (v2);
\path (5,1) coordinate (v3);
\path (5,2) coordinate (v4);
\path (6,1) coordinate (v5);
\path (6.9,1) coordinate (u1p);
\path (7,1) coordinate (u1);
\path (8,0) coordinate (u2);
\path (8,1) coordinate (u3);
\path (8,2) coordinate (u4);
\path (9,1) coordinate (u5);
\draw (v1)--(v2)--(v5)--(v4)--(v1);
\draw (v2)--(v3)--(v4);
\draw (u1)--(u2)--(u5)--(u4);
\draw (u2)--(u3)--(u4);
\draw (v5)--(u1);
\draw (v1)--(v6);
\draw[dotted, line width=0.05cm] (u1)--(u4);
%
\draw (v1) [fill=black] circle (\vr);
\draw (v2) [fill=white] circle (\vr);
\draw (v3) [fill=black] circle (\vr);
\draw (v4) [fill=white] circle (\vr);
\draw (v5) [fill=white] circle (\vr);
\draw (v6) [fill=white] circle (\vr);
\draw (u1) [fill=black] circle (\vr);
\draw (u2) [fill=white] circle (\vr);
\draw (u3) [fill=white] circle (\vr);
\draw (u4) [fill=black] circle (\vr);
\draw (u5) [fill=white] circle (\vr);
\draw (6.5,-1.25) node {{\small (a)}};
\draw[anchor = east] (v6) node {{\small $r$}};
\draw[anchor = south] (u1p) node {{\small $v_1$}};
\draw[anchor = south] (u4) node {{\small $v_2$}};
\draw [style=dashed,rounded corners] (6.5,-0.5) rectangle (10.25,2.85);
\draw (9.5,2.25) node {{\small $F'$}};
\path (12.25,1) coordinate (v6);
\path (13,1) coordinate (v1);
\path (14,0) coordinate (v2);
\path (14,1) coordinate (v3);
\path (14,2) coordinate (v4);
\path (15,1) coordinate (v5);
\path (16,1) coordinate (u1);
\path (17,0) coordinate (u2);
\path (17,1) coordinate (u3);
\path (17,2) coordinate (u4);
\path (18,1) coordinate (u5);
\draw (v1)--(v2)--(v5)--(v4)--(v1);
\draw (v2)--(v3)--(v4);
\draw (u1)--(u2)--(u5);
\draw (u4)--(u1);
\draw (u2)--(u3)--(u4);
\draw (v5)--(u1);
\draw (v1)--(v6);
\draw[dotted, line width=0.05cm] (u4)--(u5);
%
\draw (v1) [fill=white] circle (\vr);
\draw (v2) [fill=black] circle (\vr);
\draw (v3) [fill=white] circle (\vr);
\draw (v4) [fill=black] circle (\vr);
\draw (v5) [fill=white] circle (\vr);
\draw (v6) [fill=black] circle (\vr);
\draw (u1) [fill=white] circle (\vr);
\draw (u2) [fill=white] circle (\vr);
\draw (u3) [fill=white] circle (\vr);
\draw (u4) [fill=black] circle (\vr);
\draw (u5) [fill=black] circle (\vr);
\draw (15.5,-1.25) node {{\small (b)}};
\draw[anchor = east] (v6) node {{\small $r$}};
\draw[anchor = west] (u5) node {{\small $v_1$}};
\draw[anchor = south] (u4) node {{\small $v_2$}};
\draw [style=dashed,rounded corners] (15.5,-0.5) rectangle (19.25,2.85);
\draw (18.5,2.25) node {{\small $F'$}};
\end{tikzpicture}
\vskip -0.35 cm
\caption{Possible graphs $G'$ in the proof of Claim~\ref{claim:no-K23-with-d2-vertex}}
\label{fig:no-K23-with-d2-vertex}
\end{center}
\end{figure}

Suppose that $G' \in \cB$. We note that $G'$ contains at least one copy of $K_{2,3}$ that contains two vertices of degree~$2$ in $G'$. By Claim~\ref{claim:no-K23-two-deg2}, $G' - v_1v_2$ does not contain a $K_{2,3}$-subgraph that contains two vertices of degree~$2$ in $G$. We therefore infer that $G' \in \cB_1$ and that $G'$ contains exactly one copy of $K_{2,3}$, say $F'$, that contains two vertices of degree~$2$ in $G'$. Moreover, the added edge $v_1v_2$ belongs to $F'$. Thus, there are two possibilities that may occur. The first case is that both $v_1$ and $v_2$ have degree~$3$ in $G'$, as illustrated in Figure~\ref{fig:no-K23-with-d2-vertex}(a), and the second case is that one of $v_1$ and $v_2$ has degree~$2$ in $G'$, as illustrated in Figure~\ref{fig:no-K23-with-d2-vertex}(b). In both illustrations, we indicate the added edge $v_1v_2$ by the dotted edge, which recall exists in $G'$ but not in $G$.

By properties of the bad family $\cB$ (see Section~\ref{S:familyB}), we infer that there is an $i$-set, say $I'$, of $G' - v_1v_2$ that contains both $v_1$ and $v_2$ and is such that $|I'| = i(G') - 1$ in the first case (as indicated by the shaded vertices in Figure~\ref{fig:no-K23-with-d2-vertex}(a)) and $|I'| = i(G')$ in the second case (as indicated by the shaded vertices in Figure~\ref{fig:no-K23-with-d2-vertex}(b)). In both cases, we have that $I'$ is an ID-set of $G' - v_1v_2$. Moreover, $\{v_1,v_2\} \subset I'$ and $|I'| \le i(G')$. The set $I'$ can be extended to an ID-set of $G$ by adding to it the vertex~$u$. Thus, $i(G) \le i(G') + 1$. We note that $\w(G) = \w(G') + 16$ and $\Theta(G) = \Theta(G') - 2$, and so $\Omega(G) = \w(G) + 16 +\Theta(G) = \w(G') + \Theta(G') - 2 = \Omega(G') + 14$. Therefore, $8i(G) \le 8(i(G') + 1) \le \Omega(G') + 8 < \Omega(G)$, a contradiction. Hence, $G' \notin \cB$, and so $b(G') = 0$.

Thus, there is a troublesome configuration, say $T'$, in $G'$ that contains the (adjacent) vertices $v_1$ and~$v_2$. Let $z_1$ and $z_2$ be the two link vertices in $T'$. Recall that $e = v_1v_2$ is the edge that was added to construct $G'$. If $e$ does not belong to a copy of $K_{2,3}$ in $T'$ and if $e \ne z_1z_2$, then removing the edge $e$ from $T'$ and adding back the deleted $K_{2,3}$-subgraph $F$ (and the edges $u_1v_1$ and $u_2v_2$) produces a troublesome configuration in $G$, a contradiction. Hence either the edge $e$ belongs to a copy of $K_{2,3}$ in $T'$ or $v_1$ and $v_2$ are the two link vertices in $T'$ (that is, $\{v_1,v_2\} = \{z_1,z_2\}$). Let $T'$ be obtained from a bad graph $B' \in \cB$, where we may assume (renaming $z_1$ and $z_2$ if necessary) that $z_1$ is the root vertex of~$B'$.

\begin{subclaim}
\label{claim:no-K23-with-d2-vertex.1}
The edge $e$ joins the two link vertices in $T'$, that is, $\{v_1,v_2\} = \{z_1,z_2\}$.
\end{subclaim}
\proof
Suppose that $e$ belongs to a copy of $K_{2,3}$, say $F_e$, in $T'$. Let $T'$ have $k' \ge 1$ vertex disjoint copies of $K_{2,3}$. Let $I'$ be the non-canonical independent set of $B'$, and so $|I'| = 2k'$ (and neither $z_1$ nor $z_2$ belong to~$I'$). Moreover, let $I_e'$ be obtained from $I'$ by deleting the two vertices in $F_e$ that belong to~$I'$ and replacing them with the vertices~$v_1$ and~$v_2$. We note that $|I_e'| = |I'| = 2k'$.

\begin{figure}[htb]
\begin{center}
\begin{tikzpicture}[scale=.75,style=thick,x=0.75cm,y=0.75cm]
\def\vr{2.25pt}
\def\vrn{1.25pt}
\path (0,4.5) coordinate (b1);
\path (1.5,3.5) coordinate (b2);
\path (1.5,4.5) coordinate (b3);
\path (1.5,5.5) coordinate (b4);
\path (3.5,4.5) coordinate (b5);
\path (0,7) coordinate (e1);
\path (0,8) coordinate (e11);
\path (3.5,7) coordinate (e5);
\path (3.5,8) coordinate (e51);
\path (1.5,0.5) coordinate (a1);
\path (2.5,-0.5) coordinate (a2);
\path (2.5,0.5) coordinate (a3);
\path (2.5,1.5) coordinate (a4);
\path (3.5,0.5) coordinate (a5);
%
\draw (b1)--(e1);
\draw (e1)--(e11);
\draw (b5)--(e5)--(e51);
\draw (e1)--(e5);
\draw (b1)--(b2);
\draw (b3)--(b4)--(b5);
\draw (b1)--(b4);
\draw (b2)--(b5);
\draw (a1)--(a2)--(a3)--(a4)--(a5)--(a2);
\draw (a1)--(a4);
\draw (a1)--(b2);
%
\draw[dotted, line width=0.05cm] (b2)--(b3);
\draw (b3) to[out=0,in=90, distance=1.25cm] (a5);
\draw [style=dashed,rounded corners] (-2,2.75) rectangle (5,7.5);
\draw [style=dashed,rounded corners] (-1,3) rectangle (4.5,6);
\draw [style=dashed,rounded corners] (-0.5,-1.25) rectangle (4.5,2.25);
\draw [style=dashed,rounded corners] (-1,6.5) rectangle (4.5,9);
\draw (e1) [fill=white] circle (\vr);
\draw (e5) [fill=white] circle (\vr);
\draw (b1) [fill=white] circle (\vr);
\draw (b2) [fill=black] circle (\vr);
\draw (b3) [fill=black] circle (\vr);
\draw (b4) [fill=white] circle (\vr);
\draw (b5) [fill=white] circle (\vr);
\draw (a1) [fill=white] circle (\vr);
\draw (a2) [fill=white] circle (\vr);
\draw (a3) [fill=black] circle (\vr);
\draw (a4) [fill=white] circle (\vr);
\draw (a5) [fill=white] circle (\vr);
\draw[anchor = east] (b2) node {{\small $v_1$}};
\draw[anchor = east] (b3) node {{\small $v_2$}};
\draw[anchor = east] (e1) node {{\small $z_1$}};
\draw[anchor = west] (e5) node {{\small $z_2$}};
\draw[anchor = east] (a1) node {{\small $u_1$}};
\draw[anchor = west] (a5) node {{\small $u_2$}};
\draw[anchor = west] (a3) node {{\small $u$}};
\draw[anchor = north] (a2) node {{\small $w_1$}};
\draw[anchor = south] (a4) node {{\small $w_2$}};
\draw (-1.5,6.5) node {{\small $T'$}};
\draw (-0.5,3.5) node {{\small $F_e$}};
\draw (0,-0.5) node {{\small $F$}};
\draw (1.75,8.25) node {{\small $G^*$}};
\draw (1.75,4.05) node {{\small $e$}};
\draw (1.5,-2.5) node {{\small (a)}};
\path (9,9) coordinate (z0);
\path (9,5) coordinate (z1);
\path (12.25,5) coordinate (z2);
\path (12.25,3.35) coordinate (z3);
\path (13,3.35) coordinate (z4);
\path (13,-0.75) coordinate (z5);
\path (16,-0.75) coordinate (z6);
\path (16,9) coordinate (z7);
\path (10,4.5) coordinate (b1);
\path (11.5,3.5) coordinate (b2);
\path (11.5,4.5) coordinate (b3);
\path (11.5,5.5) coordinate (b4);
\path (13.5,4.5) coordinate (b5);
\path (10,7) coordinate (e1);
\path (10,8) coordinate (e11);
\path (13.5,7) coordinate (e5);
\path (13.5,8) coordinate (e51);
\path (11.5,0.5) coordinate (a1);
\path (12.5,-0.5) coordinate (a2);
\path (12.5,0.5) coordinate (a3);
\path (12.5,1.5) coordinate (a4);
\path (13.5,0.5) coordinate (a5);
%
\draw (b1)--(e1);
\draw (e1)--(e11);
\draw (b5)--(e5)--(e51);
\draw (e1)--(e5);
\draw (b1)--(b2)--(b3)--(b4)--(b5);
\draw (b1)--(b4);
\draw (a1)--(a2)--(a3)--(a4)--(a5)--(a2);
\draw (a1)--(a4);
\draw (a1)--(b2);
\draw (a5)--(b5);
\draw[dotted, line width=0.05cm] (b2)--(b5);
\draw[dotted, line width=0.05cm] (z0)--(z1)--(z2)--(z3)--(z4)--(z5)--(z6)--(z7)--(z0);
%
\draw [style=dashed,rounded corners] (7.5,2.75) rectangle (15,7.5);
\draw [style=dashed,rounded corners] (8.5,3) rectangle (14.5,6);
\draw [style=dashed,rounded corners] (9.5,-1.5) rectangle (14.5,2.25);
%
\draw (e1) [fill=white] circle (\vr);
\draw (e5) [fill=white] circle (\vr);
\draw (b1) [fill=white] circle (\vr);
\draw (b2) [fill=black] circle (\vr);
\draw (b3) [fill=white] circle (\vr);
\draw (b4) [fill=white] circle (\vr);
\draw (b5) [fill=white] circle (\vr);
\draw (a1) [fill=white] circle (\vr);
\draw (a2) [fill=white] circle (\vr);
\draw (a3) [fill=black] circle (\vr);
\draw (a4) [fill=white] circle (\vr);
\draw (a5) [fill=white] circle (\vr);
\draw[anchor = east] (b2) node {{\small $v_1$}};
\draw[anchor = west] (b5) node {{\small $v_2$}};
\draw[anchor = east] (e1) node {{\small $z_1$}};
\draw[anchor = west] (e5) node {{\small $z_2$}};
\draw[anchor = east] (a1) node {{\small $u_1$}};
\draw[anchor = west] (a5) node {{\small $u_2$}};
\draw[anchor = east] (a3) node {{\small $u$}};
\draw[anchor = north] (a2) node {{\small $w_1$}};
\draw[anchor = south] (a4) node {{\small $w_2$}};
\draw (8,6.5) node {{\small $T'$}};
\draw (9,3.5) node {{\small $F_e$}};
\draw (10,-0.5) node {{\small $F$}};
\draw (11.75,8.25) node {{\small $G^*$}};
\draw (12.75,4.5) node {{\small $e$}};
\draw (11.5,-2.5) node {{\small (b)}};
\end{tikzpicture}
\vskip -0.25 cm
\caption{Possible subgraphs in the proof of Claim~\ref{claim:no-K23-with-d2-vertex.1}}
\label{fig:no-K23-with-d2-vertex-2}
\end{center}
\end{figure}
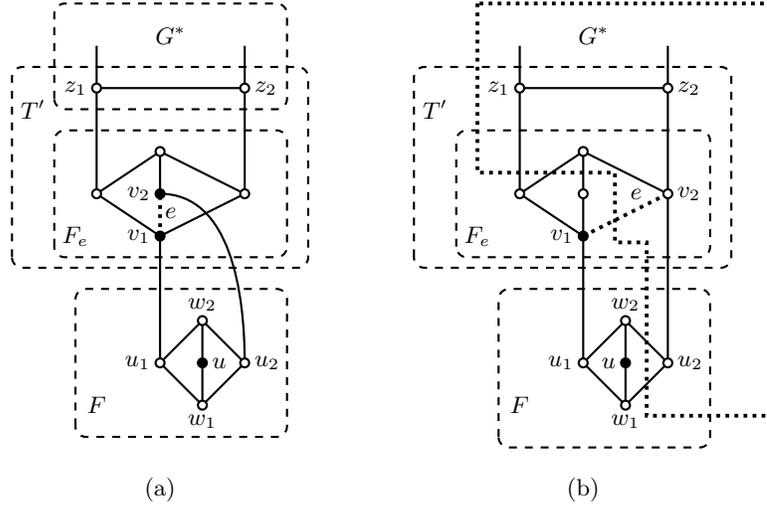

Suppose firstly that neither $v_1$ nor $v_2$ is adjacent to~$z_1$ or~$z_2$. In this case, we let $X^* = V(F) \cup V(B')\setminus\{z_1\}$ and we consider the graph $G^* = G - X^*$. In the case when $k' = 1$, the graph illustrated in Figure~\ref{fig:no-K23-with-d2-vertex-2}(a) is an example of a subgraph of $G$, where the subgraphs $F$, $F_e$, $T'$ and $G^*$ are indicated by the dashed boxes. Since $\tc(G) = 0$, any troublesome configuration in $G^*$ necessarily contains the adjacent vertices~$z_1$ and~$z_2$. However, both $z_1$ and~$z_2$ have degree at most~$2$ in $G^*$, and therefore cannot belong to a troublesome configuration in $G^*$. Moreover, since no bad graph contains two adjacent vertices, both of which has degree at most~$2$, and since $b(G) = 0$, the connected graph $G^*$ does not belong to the bad family~$\cB$, that is, $b(G^*) = 0$. Thus, $\Theta(G^*) = 0$ and $\Omega(G^*) = \w(G^*) = \w(G) - \w_G(X^*) + 2 = \w(G) - 16(k'+1) + 2 = \Omega(G) - 16k' - 14$. We now let $I^* = I_e' \cup \{u\}$. Thus, $|I^*| = 2k' + 1$. In the illustration in Figure~\ref{fig:no-K23-with-d2-vertex-2}, the set $I^*$ is indicated by the shaded vertices. Every $i$-set of $G^*$ can be extended to an ID-set of $G$ by adding to it the set $I^*$, and so $i(G) \le i(G^*) + |X^*| = i(G^*) + 2k' + 1$. Thus, $8i(G) \le 8(i(G^*) + 2k' + 1) \le \Omega(G^*) + 16k' + 8 = \Omega(G) - 6 < \Omega(G)$, a contradiction.

Hence renaming vertices if necessary, we may assume that $v_2z_2 \in E(G)$. In this case, we let $X^* = N_G[v_1] \cup N_G[u]$ and we consider the graph $G^* = G - X^*$. In the case when $k' = 1$, the graph illustrated in Figure~\ref{fig:no-K23-with-d2-vertex-2}(b) is an example of a subgraph of $G$, where the subgraphs $F$, $F_e$, $T'$ and $G^*$ are indicated by the dashed boxes. Since $b(G) = \tc(G) = 0$, we infer from the structure of the graph $G^*$ that $b(G^*) = \tc(G^*) = 0$, and so  $\Theta(G^*) = 0$ and $\Omega(G^*) = \w(G^*) = \w(G) - \w_G(X^*) + 5 = \w(G) - 23 + 5 = \Omega(G) - 18$. Every $i$-set of $G^*$ can be extended to an ID-set of $G$ by adding to it the vertices~$u$ and~$v_1$, and so $i(G) \le i(G^*) + 2$. Thus, $8i(G) \le 8(i(G^*) + 2) \le \Omega(G^*) + 16 = \Omega(G) - 2 < \Omega(G)$, a contradiction.~\smallqed

\medskip
By Claim~\ref{claim:no-K23-with-d2-vertex.1}, the edge $e$ joins the two link vertices in $T'$, that is, $\{v_1,v_2\} = \{z_1,z_2\}$. Recall that $F$ is an arbitrary $K_{2,3}$-subgraph that contains exactly one vertex of degree~$2$ in $G$. Moreover, $u_1$ and $u_2$ are the two vertices of degree~$2$ in $F$ that have degree~$3$ in $G$, where $v_1$ and $v_2$ are the neighbors of $u_1$ and $u_2$, respectively, that do not belong to~$F$. We have shown that $v_1 \ne v_2$ and $v_1$ and $v_2$ are not adjacent. Moreover, we have shown that if $G'$ is obtained from $G - V(F)$ by adding the edge~$e = v_1v_2$, then $G'$ contains a troublesome configuration, say $T'$, that necessarily contains $v_1$ and $v_2$ as the two link vertices of $F'$. This property, which we call property~($*$), holds for every $K_{2,3}$-subgraph that contains exactly one vertex of degree~$2$ in $G$.

We note that every $K_{2,3}$-subgraph in $T'$ contains only one vertex of degree~$2$ in $G$, for otherwise, $T'$ would contain a $K_{2,3}$-subgraph with two vertices of degree~$2$ in $G$, contradicting Claim~\ref{claim:no-K23-two-deg2}. Suppose that $T'$ contains two or more copies of $K_{2,3}$. In this case, we could have chosen the $K_{2,3}$-subgraph $F$ to be adjacent to another $K_{2,3}$-subgraph. With this choice of $F$, and with the vertices $u_1$, $u_2$, $v_1$, and $v_2$ as defined earlier (notably, $u_1$ and $u_2$ are vertices of degree $3$ whose neighbors outside $F$ are $v_1$ and $v_2$, respectively), at least one of $v_1$ or $v_2$ belongs to a copy of $K_{2,3}$ in $T'$, contradicting property~($*$) that $v_1$ and $v_2$ are necessarily the two link vertices of $T'$ (and noting that the link vertices of $T'$ do not belong to a copy of $K_{2,3}$ in $T'$).

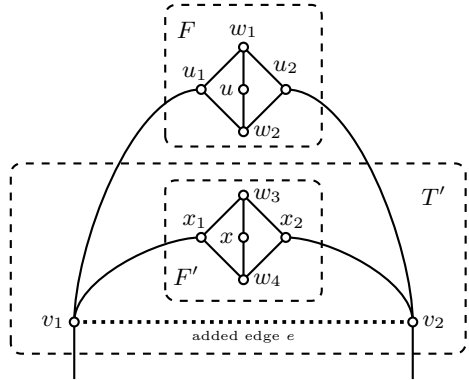
\begin{figure}[htb]
\begin{center}
\begin{tikzpicture}[scale=.75,style=thick,x=0.75cm,y=0.75cm]
\def\vr{2.25pt}
\def\vrn{1.25pt}
\path (0,1) coordinate (v1);
\path (0,-0.35) coordinate (y1);
\path (3,3) coordinate (a1);
\path (2.85,3) coordinate (a1p);
\path (4,2) coordinate (a2);
\path (4,3) coordinate (a3);
\path (4,4) coordinate (a4);
\path (5,3) coordinate (a5);
\path (5.15,3) coordinate (a5p);
\path (3,6.5) coordinate (b1);
\path (2.85,6.5) coordinate (b1p);
\path (4,5.5) coordinate (b2);
\path (4,6.5) coordinate (b3);
\path (4,7.5) coordinate (b4);
\path (5,6.5) coordinate (b5);
\path (5,6.65) coordinate (b5p);
\path (8,1) coordinate (v2);
\path (8,-0.35) coordinate (y2);
\draw (y1)--(v1);
%
\draw (y2)--(v2);
\draw (a1)--(a2)--(a3)--(a4)--(a5)--(a2);
\draw (a1)--(a4);
\draw (b1)--(b2)--(b3)--(b4)--(b5)--(b2);
\draw (b1)--(b4);
\draw (v1) to[out=90,in=180, distance=0.75cm] (a1);
\draw (v1) to[out=90,in=180, distance=1.25cm] (b1);
\draw (v2) to[out=90,in=0, distance=0.75cm] (a5);
\draw (v2) to[out=90,in=0, distance=1.25cm] (b5);
\draw[dotted, line width=0.05cm] (v1)--(v2);
\draw [style=dashed,rounded corners] (2.15,5.15) rectangle (5.85,8.5);
\draw [style=dashed,rounded corners] (2.15,1.5) rectangle (5.85,4.35);
\draw [style=dashed,rounded corners] (-1.5,0.25) rectangle (9.25,4.75);
%
\draw (v1) [fill=white] circle (\vr);
%
\draw (v2) [fill=white] circle (\vr);
%
\draw (a1) [fill=white] circle (\vr);
\draw (a2) [fill=white] circle (\vr);
\draw (a3) [fill=white] circle (\vr);
\draw (a4) [fill=white] circle (\vr);
\draw (a5) [fill=white] circle (\vr);
\draw (b1) [fill=white] circle (\vr);
\draw (b2) [fill=white] circle (\vr);
\draw (b3) [fill=white] circle (\vr);
\draw (b4) [fill=white] circle (\vr);
\draw (b5) [fill=white] circle (\vr);
\draw[anchor = east] (v1) node {{\small $v_1$}};
\draw[anchor = west] (v2) node {{\small $v_2$}};
%
\draw[anchor = south] (b1p) node {{\small $u_1$}};
\draw[anchor = south] (b5p) node {{\small $u_2$}};
\draw[anchor = west] (b2) node {{\small $w_2$}};
\draw[anchor = east] (b3) node {{\small $u$}};
\draw[anchor = south] (b4) node {{\small $w_1$}};
\draw[anchor = south] (a1p) node {{\small $x_1$}};
\draw[anchor = south] (a5p) node {{\small $x_2$}};
\draw[anchor = west] (a2) node {{\small $w_4$}};
\draw[anchor = east] (a3) node {{\small $x$}};
\draw[anchor = west] (a4) node {{\small $w_3$}};
\draw (2.65,7.9) node {{\small $F$}};
\draw (2.65,2.1) node {{\small $F'$}};
\draw (8.5,4) node {{\small $T'$}};
\draw (4,0.65) node {{\tiny added edge $e$}};
\end{tikzpicture}
\vskip -0.25 cm
\caption{A possible subgraph in the proof of Claim~\ref{claim:no-K23-with-d2-vertex}}
\label{fig:no-K23-two-deg2-b}
\end{center}
\end{figure}

Hence, $T'$ contains exactly one copy of $K_{2,3}$, say $F'$. Let $x$ be the vertex of $F'$ of degree~$2$ in $G$, and let $w_3$ and $w_4$ be the two neighbors of $x$ in $F'$. Let $x_1$ and $x_2$ be the common neighbors of $w_3$ and $w_4$ in $F'$ different from~$x$. Thus, $N_G(w_3) = N_G(w_4) = \{x,x_1,x_2\}$. Renaming $x_1$ and $x_2$ if necessary, we may assume that $v_1x_1$ and $v_2x_2$ are edges of $G$. Thus the graph illustrated in Figure~\ref{fig:no-K23-two-deg2-b} is a subgraph of $G$, where the added edge $e$ of $G'$ is indicated by the dotted line joining $v_1$ and $v_2$.

Since $T'$ is a troublesome configuration of $G'$, at least one of $v_1$ and $v_2$ has degree~$3$ in $G'$. By symmetry, we may assume that $\deg_G(v_1) = 3$ and we let $y_1$ be the neighbor of $v_1$ in $G$ different from $u_1$ and $x_1$. Suppose that $\deg_G(v_2) = 2$, and so $N_G(v_2) = \{u_2,x_2\}$. If $\deg_G(y_1) = 1$, then the graph $G$ is determined. In this case, $i(G) = |\{u,x,v_1,v_2\}| = 4$ and $\Omega(G) = \w(G) = 16 + 16 + 5 + 4 + 3 = 44$. Thus, $8i(G) < \Omega(G)$, a contradiction. Hence, $\deg_G(y_1) \ge 2$. In this case, we let $X^* = V(F) \cup V(F') \cup \{v_1,v_2,y_1\}$ and we consider the graph $G^* = G - X^*$. We note that every $X^*$-exit edge is incident with $y_1$. Moreover, since $\deg_G(y_1) \ge 2$, there is at least one $X^*$-exit edge. By Claim~\ref{claim:cost-of-exit-edge}, the removal of an $X^*$-exit edge when constructing $G^*$ increases the total weight by at most~$3$. Since there are at most two $X^*$-exit edges, we have that  $\Omega(G^*) \le \Omega(G) - \w_G(X^*) + 6 \le \Omega(G) - (16 + 16 + 4 + 3 + 3) + 6 = \Omega(G) - 36$. Every $i$-set of $G^*$ can be extended to an ID-set of $G$ by adding to it the set $\{u,x,v_1,v_2\}$, and so $i(G) \le i(G^*) + 4$. Thus, $8i(G) \le 8(i(G^*) + 4) \le \Omega(G^*) + 32 = \Omega(G) - 4 < \Omega(G)$, a contradiction. Hence, $\deg_G(v_2) = 3$. Let $y_2$ be the neighbor of $v_2$ in $G$ different from $u_2$ and $x_2$.

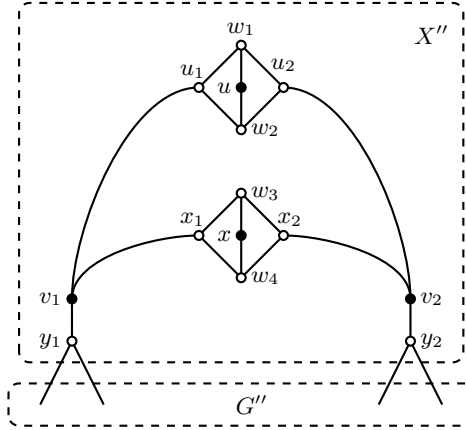
\begin{figure}[htb]
\begin{center}
\begin{tikzpicture}[scale=.75,style=thick,x=0.75cm,y=0.75cm]
\def\vr{2.25pt}
\def\vrn{1.25pt}
\path (0,1.5) coordinate (v1);
\path (0,0.5) coordinate (y1);
\path (-0.75,-1) coordinate (y11);
\path (0.75,-1) coordinate (y12);
\path (3,3) coordinate (a1);
\path (2.85,3) coordinate (a1p);
\path (4,2) coordinate (a2);
\path (4,3) coordinate (a3);
\path (4,4) coordinate (a4);
\path (5,3) coordinate (a5);
\path (5.15,3) coordinate (a5p);
\path (3,6.5) coordinate (b1);
\path (2.85,6.5) coordinate (b1p);
\path (4,5.5) coordinate (b2);
\path (4,6.5) coordinate (b3);
\path (4,7.5) coordinate (b4);
\path (5,6.5) coordinate (b5);
\path (5,6.65) coordinate (b5p);
\path (8,1.5) coordinate (v2);
\path (8,0.5) coordinate (y2);
\path (7.25,-1) coordinate (y21);
\path (8.75,-1) coordinate (y22);
\draw (y11)--(y1)--(y12);
\draw (y1)--(v1);
\draw (y21)--(y2)--(y22);
\draw (y2)--(v2);
\draw (a1)--(a2)--(a3)--(a4)--(a5)--(a2);
\draw (a1)--(a4);
\draw (b1)--(b2)--(b3)--(b4)--(b5)--(b2);
\draw (b1)--(b4);
\draw (v1) to[out=90,in=180, distance=0.75cm] (a1);
\draw (v1) to[out=90,in=180, distance=1.25cm] (b1);
\draw (v2) to[out=90,in=0, distance=0.75cm] (a5);
\draw (v2) to[out=90,in=0, distance=1.25cm] (b5);
%
%
\draw [style=dashed,rounded corners] (-1.25,0) rectangle (9.25,8.5);
\draw [style=dashed,rounded corners] (-1.5,-1.5) rectangle (9.5,-0.5);
%
\draw (v1) [fill=black] circle (\vr);
\draw (y1) [fill=white] circle (\vr);
%
\draw (v2) [fill=black] circle (\vr);
\draw (y2) [fill=white] circle (\vr);
%
\draw (a1) [fill=white] circle (\vr);
\draw (a2) [fill=white] circle (\vr);
\draw (a3) [fill=black] circle (\vr);
\draw (a4) [fill=white] circle (\vr);
\draw (a5) [fill=white] circle (\vr);
\draw (b1) [fill=white] circle (\vr);
\draw (b2) [fill=white] circle (\vr);
\draw (b3) [fill=black] circle (\vr);
\draw (b4) [fill=white] circle (\vr);
\draw (b5) [fill=white] circle (\vr);
\draw[anchor = east] (v1) node {{\small $v_1$}};
\draw[anchor = west] (v2) node {{\small $v_2$}};
\draw[anchor = east] (y1) node {{\small $y_1$}};
\draw[anchor = west] (y2) node {{\small $y_2$}};
\draw[anchor = south] (b1p) node {{\small $u_1$}};
\draw[anchor = south] (b5p) node {{\small $u_2$}};
\draw[anchor = west] (b2) node {{\small $w_2$}};
\draw[anchor = east] (b3) node {{\small $u$}};
\draw[anchor = south] (b4) node {{\small $w_1$}};
\draw[anchor = south] (a1p) node {{\small $x_1$}};
\draw[anchor = south] (a5p) node {{\small $x_2$}};
\draw[anchor = west] (a2) node {{\small $w_4$}};
\draw[anchor = east] (a3) node {{\small $x$}};
\draw[anchor = west] (a4) node {{\small $w_3$}};
\draw (8.5,7.75) node {{\small $X''$}};
\draw (4.25,-1) node {{\small $G''$}};
\end{tikzpicture}
\vskip -0.15 cm
\caption{A possible subgraph in the proof of Claim~\ref{claim:no-K23-with-d2-vertex}}
\label{fig:no-K23-two-deg2-c}
\end{center}
\end{figure}

  Suppose that $y_1 = y_2$. If $\deg_G(y_1) = 2$, then the graph $G$ is determined, and in this case $i(G) = |\{u,x,v_1,v_2\}| = 4$ and $\Omega(G) = \w(G) = 16 + 16 + 4 + 3 + 3 = 42$. Thus, $8i(G) < \Omega(G)$, a contradiction. Hence, $\deg_G(y_1) = 3$. As before, we let $X^* = V(F) \cup V(F') \cup \{v_1,v_2,y_1\}$ and we consider the graph $G^* = G - X^*$. In this case, there is only one $X^*$-exit edge, implying analogously as before that $\Omega(G^*) \le \Omega(G) - \w_G(X^*) + 3 = \Omega(G) - (16 + 16 + 3 + 3 + 3) + 3 = \Omega(G) - 38$. As before, $i(G) \le i(G^*) + 4$, yielding the contradiction $8i(G) < \Omega(G)$. Hence, $y_1 \ne y_2$. We now let $X'' = V(F) \cup V(F') \cup \{v_1,v_2,y_1,y_2\}$ and we consider the graph $G'' = G - X''$. Thus the graph illustrated in Figure~\ref{fig:no-K23-two-deg2-c} is a subgraph of $G$.

\vskip 0.25 cm

We note that there are at most four $X''$-exit edges, namely at most two edges incident with~$y_1$ and at most two edges incident with~$y_2$. By Claim~\ref{claim:cost-of-exit-edge}, removing each $X''$-exit edge when constructing $G''$ increases the total weight by at most~$3$. Thus, $\Omega(G'') \le \Omega(G) - \w_G(X'') + 4 \times 3 = \Omega(G) - (16 + 16 + 12) + 12 = \Omega(G) - 32$. Every $i$-set of $G''$ can be extended to an ID-set of $G$ by adding to it the set $\{u,x,v_1,v_2\}$ (as indicated by the shaded vertices in Figure~\ref{fig:no-K23-two-deg2-c}), and so $i(G) \le i(G'') + 4$. Thus, $8i(G) \le 8(i(G'') + 4) \le \Omega(G'') + 32 = \Omega(G)$, a contradiction.~\smallqed

\medskip
By Claim~\ref{claim:no-K23-with-d2-vertex}, every vertex that belongs to a $K_{2,3}$-subgraph has degree~$3$ in $G$.

\begin{claim}
\label{claim:induced-K23}
Every $K_{2,3}$-subgraph in $G$ is an induced subgraph in $G$.
\end{claim}
\proof Suppose, to the contrary, that $G$ contains a $K_{2,3}$-subgraph, $F$,  that is not an induced subgraph of $G$. Let the subgraph $F$ have partite sets $\{w_1,w_2\}$ and $\{u_1,u_2,u_3\}$. Renaming vertices if necessary, we may assume that $u_2u_3 \in E(G)$. If $\deg_G(u_1) = 2$, then the graph $G$ is determined, and so $n = 5$, a contradiction. Hence, $\deg_G(u_1) = 3$. We now consider the graph $G' = G - \{u_2,u_3,w_1,w_2\}$. We note that the vertex $u_1$ has degree~$1$ in $G'$, and so $u_1$ does not belong to a troublesome configuration in $G'$. Moreover, by Claim~\ref{claim:no-K23-two-deg2}, the (connected) graph $G' \notin \cB$. Thus, $\Omega(G') = \Omega(G) - 12 + 2 = \Omega(G) - 10$. Every $i$-set of $G'$ can be extended to an ID-set of $G$ by adding to it the vertex~$u_2$, and so $i(G) \le i(G') + 1$. Thus, $8i(G) \le 8(i(G') + 1) \le \Omega(G') + 8 < \Omega(G)$, a contradiction.~\smallqed

\begin{claim}
\label{claim:no-bad-in-Ge}
$b(G - e) = 0$ for all $e \in E(G)$.
\end{claim}
\proof Suppose, to the contrary, that the claim is false. Thus there exists an edge $e$ in $G$ whose removal creates a bad component. Let $B \in \cB$ be such a bad component in $G-e$. The component $B$ contains at least one copy of $K_{2,3}$, say $F$, with two (nonadjacent) vertices of degree~$2$ in $B$. By Claim~\ref{claim:induced-K23}, $F$ is an induced subgraph of $G$. Thus the edge $e$ is incident with at most one of the two vertices in $F$ that have degree~$2$ in $B$. Therefore, $F$ contains at least one vertex of degree~$2$ in $G$, contradicting Claim~\ref{claim:no-K23-with-d2-vertex}.~\smallqed

\begin{claim}
\label{claim:no-bad-in-Gee}
$b(G - e_1 - e_2) = 0$ for all $e_1,e_2 \in E(G)$.
\end{claim}
\proof Suppose, to the contrary, that the claim is false. Thus there exists edges $e_1$ and $e_2$ in $G$ whose removal creates a bad component. Let $B \in \cB$ be such a bad component in $G - e_1 - e_2$, and let $v$ be the root vertex of~$B$. Let $B$ contain $k \ge 1$ copies of $K_{2,3}$. If $B \in \cB_3$, then $k \ge 3$ and there are $k+3 \ge 6$ vertices of degree~$2$ in $B$. Thus, $B$ contains at least one copy of $K_{2,3}$ that contains a vertex of degree~$2$ in $G$, contradicting Claim~\ref{claim:no-K23-with-d2-vertex}. Hence, $B \in \cB_1 \cup \cB_2$.

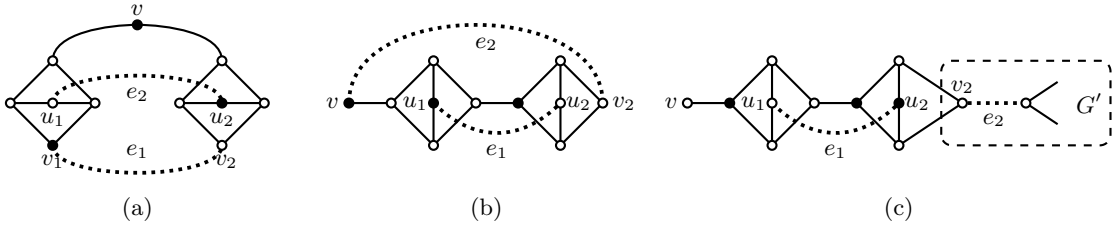
\begin{figure}[htb]
\begin{center}
\begin{tikzpicture}[scale=.75,style=thick,x=0.75cm,y=0.75cm]
\def\vr{2.25pt}
\def\vrn{1.25pt}
\path (0,1) coordinate (a1);
\path (1,0) coordinate (a2);
\path (1,1) coordinate (a3);
\path (1,2) coordinate (a4);
\path (2,1) coordinate (a5);
\path (4,1) coordinate (b1);
\path (5,0) coordinate (b2);
\path (5.1,0) coordinate (b2p);
\path (5,1) coordinate (b3);
\path (5,2) coordinate (b4);
\path (6,1) coordinate (b5);
\path (3,2.85) coordinate (v);
\draw (a1)--(a3)--(a5);
\draw (a1)--(a2)--(a5)--(a4)--(a1);
\draw (b1)--(b3)--(b5);
\draw (b1)--(b2)--(b5)--(b4)--(b1);
\draw (a4) to[out=90,in=180, distance=0.5cm] (v);
\draw (b4) to[out=90,in=0, distance=0.5cm] (v);
\draw[dotted, line width=0.05cm] (a2) to[out=-90,in=-90, distance=0.65cm] (b2);
\draw[dotted, line width=0.05cm] (a3) to[out=90,in=90, distance=0.65cm] (b3);
\draw (a1) [fill=white] circle (\vr);
\draw (a2) [fill=black] circle (\vr);
\draw (a3) [fill=white] circle (\vr);
\draw (a4) [fill=white] circle (\vr);
\draw (a5) [fill=white] circle (\vr);
\draw (b1) [fill=white] circle (\vr);
\draw (b2) [fill=white] circle (\vr);
\draw (b3) [fill=black] circle (\vr);
\draw (b4) [fill=white] circle (\vr);
\draw (b5) [fill=white] circle (\vr);
\draw (v) [fill=black] circle (\vr);
\draw[anchor = south] (v) node {{\small $v$}};
\draw[anchor = north] (a3) node {{\small $u_1$}};
\draw[anchor = north] (b3) node {{\small $u_2$}};
\draw[anchor = north] (a2) node {{\small $v_1$}};
\draw[anchor = north] (b2p) node {{\small $v_2$}};
\draw (3,-0.15) node {{\small $e_1$}};
\draw (3,1.2) node {{\small $e_2$}};
\draw (3,-1.5) node {{\small (a)}};
\path (8,1) coordinate (v6);
\path (9,1) coordinate (v1);
\path (10,0) coordinate (v2);
\path (10,1) coordinate (v3);
\path (10.1,1) coordinate (v3p);
\path (10,2) coordinate (v4);
\path (11,1) coordinate (v5);
\path (12,1) coordinate (u1);
\path (13,0) coordinate (u2);
\path (13,1) coordinate (u3);
\path (12.9,1) coordinate (u3p);
\path (13,2) coordinate (u4);
\path (14,1) coordinate (u5);
\draw (v1)--(v2)--(v5)--(v4)--(v1);
\draw (v2)--(v3)--(v4);
\draw (u1)--(u2)--(u5)--(u4)--(u1);
\draw (u2)--(u3)--(u4);
\draw (v5)--(u1);
\draw (v1)--(v6);
\draw (v6)[dotted, line width=0.05cm]  to[out=90,in=90, distance=1.85cm] (u5);
\draw (v3)[dotted, line width=0.05cm]  to[out=-45,in=-135, distance=1cm] (u3);
\draw (v1) [fill=white] circle (\vr);
\draw (v2) [fill=white] circle (\vr);
\draw (v3) [fill=black] circle (\vr);
\draw (v4) [fill=white] circle (\vr);
\draw (v5) [fill=white] circle (\vr);
\draw (v6) [fill=black] circle (\vr);
\draw (u1) [fill=black] circle (\vr);
\draw (u2) [fill=white] circle (\vr);
\draw (u3) [fill=white] circle (\vr);
\draw (u4) [fill=white] circle (\vr);
\draw (u5) [fill=white] circle (\vr);
\draw (11.25,-1.5) node {{\small (b)}};
\draw[anchor = east] (v6) node {{\small $v$}};
\draw[anchor = west] (u5) node {{\small $v_2$}};
\draw[anchor = east] (v3p) node {{\small $u_1$}};
\draw[anchor = west] (u3p) node {{\small $u_2$}};
\draw (11.5,-0.15) node {{\small $e_1$}};
\draw (11.25,2.4) node {{\small $e_2$}};
\path (16,1) coordinate (v6);
\path (17,1) coordinate (v1);
\path (18,0) coordinate (v2);
\path (18,1) coordinate (v3);
\path (18.1,1) coordinate (v3p);
\path (18,2) coordinate (v4);
\path (19,1) coordinate (v5);
\path (20,1) coordinate (u1);
\path (21,0) coordinate (u2);
\path (21,1) coordinate (u3);
\path (20.9,1) coordinate (u3p);
\path (21,2) coordinate (u4);
\path (22.5,1) coordinate (u5);
\path (24,1) coordinate (x);
\path (24.75,0.5) coordinate (x1);
\path (24.75,1.5) coordinate (x2);
\draw (v1)--(v2)--(v5)--(v4)--(v1);
\draw (v2)--(v3)--(v4);
\draw (u1)--(u2)--(u5)--(u4)--(u1);
\draw (u2)--(u3)--(u4);
\draw (v5)--(u1);
\draw (v1)--(v6);
\draw (x1)--(x)--(x2);
\draw[dotted, line width=0.05cm]   (x)--(u5);
%
\draw (v3)[dotted, line width=0.05cm]  to[out=-45,in=-135, distance=1cm] (u3);
\draw (v1) [fill=black] circle (\vr);
\draw (v2) [fill=white] circle (\vr);
\draw (v3) [fill=white] circle (\vr);
\draw (v4) [fill=white] circle (\vr);
\draw (v5) [fill=white] circle (\vr);
\draw (v6) [fill=white] circle (\vr);
\draw (u1) [fill=black] circle (\vr);
\draw (u2) [fill=white] circle (\vr);
\draw (u3) [fill=black] circle (\vr);
\draw (u4) [fill=white] circle (\vr);
\draw (u5) [fill=white] circle (\vr);
\draw (x) [fill=white] circle (\vr);
\draw (21,-1.5) node {{\small (c)}};
\draw[anchor = east] (v6) node {{\small $v$}};
\draw[anchor = south] (u5) node {{\small $v_2$}};
\draw[anchor = east] (v3p) node {{\small $u_1$}};
\draw[anchor = west] (u3p) node {{\small $u_2$}};
\draw (19.5,-0.15) node {{\small $e_1$}};
\draw (23.25,0.6) node {{\small $e_2$}};
\draw (25.5,1) node {{\small $G'$}};
\draw [style=dashed,rounded corners] (22,0) rectangle (26,2);
\end{tikzpicture}
\vskip -0.15 cm
\caption{Possible subgraphs in the proof of Claim~\ref{claim:no-bad-in-Gee}}
\label{fig:no-bad-in-Gee1}
\end{center}
\end{figure}

Suppose that $B \in \cB_2$. In this case, $k \ge 2$ and there are $k + 2 \ge 4$ vertices of degree~$2$ in $B$ that belong to $K_{2,3}$-subgraphs. If $k \ge 3$, then there are at least five vertices of degree~$2$ in $B$ that belong to $K_{2,3}$-subgraphs, implying that $B$ contains at least one copy of $K_{2,3}$ that contains a vertex of degree~$2$ in $G$, a contradiction. Hence, $k = 2$. Let $F_1$ and $F_2$ be the two copies of $K_{2,3}$ in $B$, and let $u_i$ and $v_i$ be the two vertices in $F_i$ that have degree~$2$ in $B$ for $i \in [2]$. By Claims~\ref{claim:no-K23-with-d2-vertex} and~\ref{claim:induced-K23}, we infer renaming vertices if necessary that $e_1 = u_1u_2$ and $e_2 = v_1v_2$. However, then the graph $G$ is determined and is illustrated in Figure~\ref{fig:no-bad-in-Gee1}(a), where the edges $e_1$ and $e_2$ are dotted. In this case, $i(G) = 3$ (an $i$-set of $G$ is indicated by the three shaded vertices in Figure~\ref{fig:no-bad-in-Gee1}(a)) and $\Omega(G) = \w(G) = 34$. Thus, $8i(G) < \Omega(G)$, a contradiction.

Hence, $B \in \cB_1$, and so $k \ge 1$ and there are $k + 1 \ge 2$ vertices of degree~$2$ in $B$ that belong to $K_{2,3}$-subgraphs. If $k \ge 4$, then there are at least five vertices of degree~$2$ in $B$ that belong to $K_{2,3}$-subgraphs, implying that $B$ contains at least one copy of $K_{2,3}$ that contains a vertex of degree~$2$ in $G$, a contradiction. Therefore, $k \le 3$.

Suppose that $k = 3$. In this case, there are four vertices of degree~$2$ in $B$ that belong to $K_{2,3}$-subgraphs. Since no copy of a $K_{2,3}$-subgraph in $B$ contains a vertex of degree~$2$ in $G$, we infer that these four vertices of degree~$2$ in $B$ are incident with the edges $e_1$ and $e_2$ in $G$. The graph $G$ is now determined and there are two non-isomorphic graphs to which $G$ can be isomorphic. By Claim~\ref{claim:induced-K23}, the edge $e_i$ joins vertices that belong to different copies of $K_{2,3}$ in $B$ for $i \in [2]$. The set consisting of the three core vertices in $B$, together with one end of $e_1$ and one end of $e_2$, is an ID-set of $G$, and so $i(G) \le 5$. Moreover, $\Omega(G) = \w(G) = 15 \times 3 + 5 = 50$, and so $8i(G) < \Omega(G)$, a contradiction. Therefore, $k \le 2$.

Suppose that $k = 2$. In this case, there are three vertices of degree~$2$ in $B$ that belong to $K_{2,3}$-subgraphs. Let $F_1$ and $F_2$ be the two copies of $K_{2,3}$ in $B$, where $F_2$ is the copy that contains two vertices of degree~$2$ in $B$. Let $u_1$ be the vertex in $F_1$ of degree~$2$ in $B$ and let $u_2$ and $v_2$ be the two vertices in $F_2$ of degree~$2$ in $B$.  Since no copy of a $K_{2,3}$-subgraph in $B$ contains a vertex of degree~$2$ in $G$, we infer that the vertices $u_1$, $u_2$ and $v_2$ are incident with the edges $e_1$ and $e_2$ in $G$. Moreover by Claim~\ref{claim:induced-K23}, vertices $u_2$ and $v_2$ are not adjacent in $G$. Renaming $u_2$ and $v_2$ if necessary, we may assume that $e_1 = u_1u_2$ and that $e_2$ is incident with~$v_2$.

Suppose that vertex~$v$ is incident with $e_2$, and so $e_2 = vv_2$. Thus, the graph $G$ is determined as is illustrated in Figure~\ref{fig:no-bad-in-Gee1}(b). In this case, $i(G) = 3$ (an $i$-set of $G$ is indicated by the three shaded vertices in Figure~\ref{fig:no-bad-in-Gee1}(b)) and $\Omega(G) = \w(G) = 34$. Thus,  $8i(G) < \Omega(G)$, a contradiction. Hence, the vertex~$v$ is not incident with $e_2$. We now consider the graph $G' = G - (V(B) \setminus \{v_2\})$ (see Figure~\ref{fig:no-bad-in-Gee1}(c), where the graph $G'$ is indicated by the dashed box). We note that $\Omega(G') = \w(G') = \w(G) - 32 + 2 = \Omega(G) - 30$. Every $i$-set of $G'$ can be extended to an ID-set of $G$ by adding to it the two core vertices of $B$ together with the vertex~$u_2$ (as indicated by the shaded vertices in Figure~\ref{fig:no-bad-in-Gee1}(c). Thus, $i(G) \le i(G') + 3$, and so $8i(G) < \Omega(G)$, a contradiction.

Hence, $k=1$. Let $F$ be the copy of $K_{2,3}$ in $B$, and let $u_1$ and $v_1$ be the two vertices of degree~$2$ in $B$ that belong to $F$. Since every vertex in $F$ has degree~$3$ in $G$, the vertices $u_1$ and $v_1$ are incident with the edges $e_1$ and $e_2$ in $G$. Moreover by Claim~\ref{claim:induced-K23}, vertices $u_1$ and $v_1$ are not adjacent in $G$. Renaming edges if necessary, we may assume that the edge $e_1$ is incident with~$u_1$ and the edge $e_2$ is incident with~$v_1$. Since $G \ne K_{3,3}$, the vertex $v$ is adjacent to at most one of~$u_1$ or $v_1$. If $v$ is adjacent to either~$u_1$ or $v_1$, then $v$ would be a vertex of degree~$2$ in $G$ that belongs to a copy of $K_{2,3}$, contradicting Claim~\ref{claim:no-K23-with-d2-vertex}. Hence, $v$ has degree~$1$ in $G$. Let $u$ be the neighbor of $v$. 
We now consider the graph $G' = G - N[u]$. We note that both $u_1$ and $v_1$ have degree~$1$ in $G'$. Moreover, $\Omega(G') = \w(G') = \w(G) - 14 + 4 = \Omega(G) - 10$. Every $i$-set of $G'$ can be extended to an ID-set of $G$ by adding to it the vertex~$u$, and so $i(G) \le i(G') + 1$. Thus, $8i(G) < \Omega(G)$, a contradiction.~\smallqed

\begin{claim}
\label{claim:no-TT-in-Ge}
$\tc(G - e) \le 1$ for all $e \in E(G)$.
\end{claim}
\proof Suppose, to the contrary, that the claim is false. Thus there exists an edge $e = uv$ in $G$ whose removal creates two troublesome configurations. Let $T_w = T_{w_1,w_2}$ and $T_v = T_{v_1,v_2}$ be the two new troublesome configuration in $G-e$, where $w_1$ and $w_2$ are the link vertices of $T_w$ and $v_1$ and $v_2$ are the link vertices of $T_v$. Let $T_w$ be obtained from the bad graph $B_w \in \cB_1$ with root vertex~$w_1$ (of degree~$1$ in $B_w$), and let $x_1$ and $x_2$ be the vertices of $B_w$ adjacent to $w_1$ and $w_2$, respectively, in $G$. Further, let $T_v$ be obtained from the bad graph $B_v \in \cB_1$ with root vertex~$v_1$ (of degree~$1$ in $B_v$), and let $y_1$ and $y_2$ be the vertices of $B_v$ adjacent to $v_1$ and $v_2$, respectively, in $G$. We may assume that $u \in V(B_w)$ and $v \in V(B_v)$. Let $u_1$ and $u_2$ be the two neighbors of $u$ in $B_w$, and let $a_1$ and $a_2$ be the two neighbors of $v$ in $B_v$. By Claim~\ref{claim:no-K23-with-d2-vertex}, we infer that each of $B_w$ and $B_v$ contain exactly one $K_{2,3}$-subgraph, which we name $F_w$ and $F_v$, respectively. Thus, $V(F_w) = \{u,u_1,u_2,x_1,x_2\}$ and $V(F_v) = \{v,a_1,a_2,y_1,y_2\}$. If $\deg_G(v_1) = 3$, then let $z_1$ be the neighbor of $v_1$ different from~$y_1$ and $v_2$. Moreover, let $z_2$ be the neighbor of $v_2$ different from~$y_2$ and $v_1$.

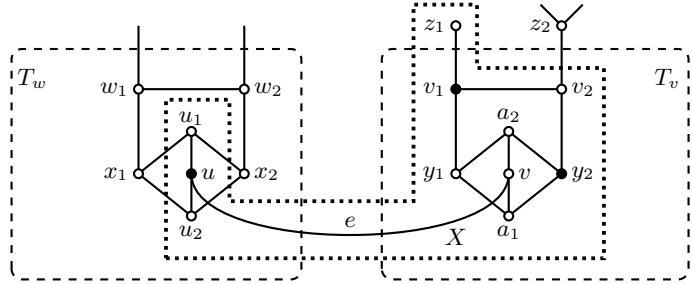
\begin{figure}[htb]
\begin{center}
\begin{tikzpicture}[scale=.75,style=thick,x=0.75cm,y=0.75cm]
\def\vr{2.25pt}
\def\vrn{1.25pt}
\path (1.65,2.5) coordinate (z0);
\path (12,2.5) coordinate (z1);
\path (12,7) coordinate (z2);
\path (9,7) coordinate (z3);
\path (9,8.5) coordinate (z4);
\path (7.5,8.5) coordinate (z5);
\path (7.5,3.85) coordinate (z6);
\path (3.15,3.85) coordinate (z7);
\path (3.15,6.25) coordinate (z8);
\path (1.65,6.25) coordinate (z9);
\path (1,6.5) coordinate (e1);
\path (1,8) coordinate (e11);
\path (1,4.5) coordinate (f1);
\path (2.25,3.5) coordinate (f2);
\path (2.25,4.5) coordinate (f3);
\path (2.25,5.5) coordinate (f4);
\path (2.25,5.425) coordinate (f4p);
\path (3.5,4.5) coordinate (f5);
%
\path (3.5,6.5) coordinate (g5);
\path (3.5,8) coordinate (g51);
\path (8.5,4.5) coordinate (c1);
\path (9.75,3.5) coordinate (c2);
\path (9.75,4.5) coordinate (c3);
\path (9.75,5.5) coordinate (c4);
\path (11,4.5) coordinate (c5);
\path (8.5,6.5) coordinate (h1);
\path (8.5,8) coordinate (h11);
%
\path (11,6.5) coordinate (q5);
\path (11,8) coordinate (q51);
\path (10.5,8.5) coordinate (q511);
\path (11.5,8.5) coordinate (q512);
\draw (f1)--(e1);
\draw (e1)--(g5);
\draw (e1)--(e11);
%
%
%
\draw (f5)--(g5)--(g51);
%
\draw (f1)--(f2)--(f3)--(f4)--(f5)--(f2);
\draw (f1)--(f4);
%
%
\draw [style=dashed,rounded corners] (-2,2) rectangle (4.85,7.45);
\draw[dotted, line width=0.05cm] (z0)--(z1)--(z2)--(z3)--(z4)--(z5)--(z6)--(z7)--(z8)--(z9)--(z0);
\draw (c1)--(h1)--(h11);
\draw (h1)--(q5);
\draw (c1)--(c2)--(c3)--(c4)--(c5)--(c2);
\draw (c1)--(c4);
\draw (q511)--(q51)--(q512);
\draw (c5)--(q5)--(q51);
%
%
\draw (f3) to[out=270,in=270, distance=1.45cm] (c3);
%
\draw [style=dashed,rounded corners] (6.75,2) rectangle (14,7.45);
\draw (e1) [fill=white] circle (\vr);
%
%
%
\draw (g5) [fill=white] circle (\vr);
\draw (f1) [fill=white] circle (\vr);
\draw (f2) [fill=white] circle (\vr);
\draw (f3) [fill=black] circle (\vr);
\draw (f4) [fill=white] circle (\vr);
\draw (f5) [fill=white] circle (\vr);
\draw (h1) [fill=black] circle (\vr);
\draw (h11) [fill=white] circle (\vr);
\draw (c1) [fill=white] circle (\vr);
\draw (c2) [fill=white] circle (\vr);
\draw (c3) [fill=white] circle (\vr);
\draw (c4) [fill=white] circle (\vr);
\draw (c5) [fill=black] circle (\vr);
%
%
\draw (q5) [fill=white] circle (\vr);
\draw (q51) [fill=white] circle (\vr);
%
%
\draw[anchor = east] (e1) node {{\small $w_1$}};
\draw[anchor = west] (g5) node {{\small $w_2$}};
\draw[anchor = west] (f5) node {{\small $x_2$}};
\draw[anchor = west] (f3) node {{\small $u$}};
\draw[anchor = east] (f1) node {{\small $x_1$}};
\draw[anchor = east] (h11) node {{\small $z_1$}};
\draw[anchor = east] (q51) node {{\small $z_2$}};
\draw[anchor = east] (h1) node {{\small $v_1$}};
\draw[anchor = west] (q5) node {{\small $v_2$}};
\draw[anchor = west] (c5) node {{\small $y_2$}};
\draw[anchor = east] (c1) node {{\small $y_1$}};
\draw[anchor = west] (c3) node {{\small $v$}};
\draw[anchor = north] (c2) node {{\small $a_1$}};
\draw[anchor = south] (c4) node {{\small $a_2$}};
\draw[anchor = north] (f2) node {{\small $u_2$}};
\draw[anchor = south] (f4p) node {{\small $u_1$}};
\draw (-1.5,6.75) node {{\small $T_w$}};
\draw (13.5,6.75) node {{\small $T_v$}};
\draw (8.5,3) node {{\small $X$}};
\draw (6,3.35) node {{\small $e$}};
\end{tikzpicture}
\caption{An illustration of a subgraph of $G$ in the proof of Claim~\ref{claim:no-TT-in-Ge}}
\label{fig:tsG-1a}
\end{center}
\end{figure}

Suppose that $\deg_G(v_1) = 3$ and $\deg_G(z_1) = 1$. In this case, we let $X = N_G[u] \cup (V(T_v) \cup \{z_1\})$ and we let $G' = G - X$. Thus, the graph shown in Figure~\ref{fig:tsG-1a} is a subgraph of $G$, where $T_v$ and $T_w$ are the troublesome configurations indicated in the dashed boxes and where $X$ is the set indicated by the vertices in the dotted region. The removal of the four $X$-exit edges incident with~$u_1$ and $u_2$ when constructing $G'$ increases the total weight by~$4$, noting that the component in $G'$ that contains the vertices $w_1$ and $w_2$ contains the path $x_1w_1w_2x_2$, where $x_1$ and $x_2$ have degree~$1$ in $G'$ (and therefore $w_1$ and $w_2$ do not belong to a bad component nor a troublesome configurations). Moreover by Claim~\ref{claim:cost-of-exit-edge}, the removal of the $X$-exit edge $v_2z_2$ when constructing $G'$ increases the total weight by at most~$3$. Therefore, $\Omega(G') \le \w(G) - \w_G(X) + 7 = \Omega(G) - 35 + 7 = \Omega(G) - 28$. Every $i$-set of $G'$ can be extended to an ID-set of $G$ by adding to it the vertices $u$, $v_1$ and $y_2$ (indicated by the shaded vertices in Figure~\ref{fig:tsG-1a}), and so $i(G) \le i(G') + 3$. Hence, $8i(G) \le 8(i(G') + 3) \le \Omega(G') + 24  < \Omega(G)$, a contradiction.

Hence either $\deg_G(v_1) = 2$ or $\deg_G(v_1) = 3$ and $\deg_G(z_1) \ge 2$. We now let $X = N_G[u] \cup (V(T_v) \setminus \{v_2,y_2\})$ and let $G' = G - X$. Thus, the graph shown in Figure~\ref{fig:tsG-1b} is a subgraph of $G$ (where in this illustration $\deg_G(v_1) = 3$), where $T_v$ and $T_w$ are the troublesome configurations indicated in the dashed boxes and where $X$ is the set indicated by the vertices in the dotted region. Every $i$-set of $G'$ can be extended to an ID-set of $G$ by adding to it the vertices $u$ and $y_1$ (indicated by the shaded vertices in Figure~\ref{fig:tsG-1b}), and so $i(G) \le i(G') + 2$. If $\deg_G(v_1) = 2$, then there are seven $X$-exit edges, and so $\w(G') = \w(G) - \w_G(X) + 7 = \w(G) - 25 + 7 = \w(G) - 18$. If $\deg_G(v_1) = 3$, then there are eight $X$-exit edges, and so $\w(G') = \w(G) - \w_G(X) + 8 = \w(G) - 24 + 8 = \w(G) - 16$.
The component in $G'$ that contains the vertices $w_1$ and $w_2$ contains at least two vertices of degree~$1$, namely the vertices $x_1$ and $x_2$. Further, we note that the vertex $y_2$ has degree~$1$ in $G'$ and its neighbor $v_2$ has degree~$2$ in $G'$. Moreover, if $v_1$ has degree~$3$ in $G$ and its neighbor $z_1$ belongs to a component of $G'$ that contains neither $w_1$ nor $v_2$, then such a component is not a bad component (that belongs to $\cB$), by Claim~\ref{claim:no-bad-in-Ge}. We therefore infer that $b(G') = 0$. We also note that the component in $G'$ that contains the vertex $w_1$ (and $w_2$) does not contain a troublesome configuration and neither does the component in $G'$ that contains the vertex $v_2$ (and $y_2$). Hence either $G'$ has no troublesome configuration or $G'$ has exactly one troublesome configuration and such a troublesome configuration necessarily contains the vertex $z_1$.

\begin{figure}[htb]
\begin{center}
\begin{tikzpicture}[scale=.75,style=thick,x=0.75cm,y=0.75cm]
\def\vr{2.25pt}
\def\vrn{1.25pt}
\path (1.65,2.5) coordinate (z0);
\path (10.5,2.5) coordinate (z2);
\path (10.5,7) coordinate (z3);
\path (7.5,7) coordinate (z6);
\path (7.5,3.85) coordinate (z7);
\path (3.15,3.85) coordinate (z8);
\path (3.15,6.25) coordinate (z9);
\path (1.65,6.25) coordinate (z10);
\path (1,6.5) coordinate (e1);
\path (1,8) coordinate (e11);
\path (1,4.5) coordinate (f1);
\path (2.25,3.5) coordinate (f2);
\path (2.25,4.5) coordinate (f3);
\path (2.25,5.5) coordinate (f4);
\path (2.25,5.425) coordinate (f4p);
\path (3.5,4.5) coordinate (f5);
%
\path (3.5,6.5) coordinate (g5);
\path (3.5,8) coordinate (g51);
\path (8.5,4.5) coordinate (c1);
\path (9.75,3.5) coordinate (c2);
\path (9.75,4.5) coordinate (c3);
\path (9.75,5.5) coordinate (c4);
\path (11,4.5) coordinate (c5);
\path (8.5,6.5) coordinate (h1);
\path (8.5,8) coordinate (h11);
%
\path (11,6.5) coordinate (q5);
\path (11,8) coordinate (q51);
\draw (f1)--(e1);
\draw (e1)--(g5);
\draw (e1)--(e11);
%
%
%
\draw (f5)--(g5)--(g51);
%
\draw (f1)--(f2)--(f3)--(f4)--(f5)--(f2);
\draw (f1)--(f4);
%
%
\draw [style=dashed,rounded corners] (-2,2) rectangle (4.85,7.45);
\draw[dotted, line width=0.05cm] (z0)--(z2)--(z3)--(z6)--(z7)--(z8)--(z9)--(z10)--(z0);
\draw (c1)--(h1)--(h11);
\draw (h1)--(q5);
\draw (c1)--(c2)--(c3)--(c4)--(c5)--(c2);
\draw (c1)--(c4);
%
%
\draw (c5)--(q5)--(q51);
%
%
\draw (f3) to[out=270,in=270, distance=1.45cm] (c3);
%
\draw [style=dashed,rounded corners] (6.75,2) rectangle (14,7.45);
\draw (e1) [fill=white] circle (\vr);
%
%
%
\draw (g5) [fill=white] circle (\vr);
\draw (f1) [fill=white] circle (\vr);
\draw (f2) [fill=white] circle (\vr);
\draw (f3) [fill=black] circle (\vr);
\draw (f4) [fill=white] circle (\vr);
\draw (f5) [fill=white] circle (\vr);
\draw (h1) [fill=white] circle (\vr);
\draw (c1) [fill=black] circle (\vr);
\draw (c2) [fill=white] circle (\vr);
\draw (c3) [fill=white] circle (\vr);
\draw (c4) [fill=white] circle (\vr);
\draw (c5) [fill=white] circle (\vr);
%
%
\draw (q5) [fill=white] circle (\vr);
%
%
\draw[anchor = east] (e1) node {{\small $w_1$}};
\draw[anchor = west] (g5) node {{\small $w_2$}};
\draw[anchor = west] (f5) node {{\small $x_2$}};
\draw[anchor = west] (f3) node {{\small $u$}};
\draw[anchor = east] (f1) node {{\small $x_1$}};
\draw[anchor = east] (h1) node {{\small $v_1$}};
\draw[anchor = west] (q5) node {{\small $v_2$}};
\draw[anchor = west] (c5) node {{\small $y_2$}};
\draw[anchor = east] (c1) node {{\small $y_1$}};
\draw[anchor = west] (c3) node {{\small $v$}};
\draw[anchor = north] (c2) node {{\small $a_1$}};
\draw[anchor = south] (c4) node {{\small $a_2$}};
\draw[anchor = north] (f2) node {{\small $u_2$}};
\draw[anchor = south] (f4p) node {{\small $u_1$}};
\draw (-1.5,6.75) node {{\small $T_w$}};
\draw (13.5,6.75) node {{\small $T_v$}};
\draw (8.5,3) node {{\small $X$}};
\draw (6,3.35) node {{\small $e$}};
\end{tikzpicture}
\caption{An illustration of a subgraph of $G$ in the proof of Claim~\ref{claim:no-TT-in-Ge}}
\label{fig:tsG-1b}
\end{center}
\end{figure}

Suppose that $G'$ has no troublesome configuration; that is, $\tc(G') = 0$. As observed earlier, $b(G') = 0$, and so $\Theta(G') = 0$. Thus, $\Omega(G') = \w(G') \le \w(G) - 16$. Hence, $8i(G) \le 8(i(G') + 2) \le \Omega(G') + 16  \le \w(G) = \Omega(G)$, a contradiction. Therefore, $G'$ has exactly one troublesome configuration, namely a troublesome configuration, $T_{z_1}$ say, that contains the vertex $z_1$. In this case, we note that $\deg_G(v_1) = 3$. Interchanging the roles of $v_1$ and $v_2$, by analogous arguments (letting $X = N_G[u] \cup (V(T_v) \setminus \{v_1,y_1\})$) we infer that $G - v_2z_2$ contains a troublesome configuration, $T_{z_2}$ say, that contains the vertex $z_2$. Thus the graph shown in Figure~\ref{fig:tsG-2} is a subgraph of $G$. Let $Z$ be the set of vertices in $T_{z_1}$ different from the two link vertices of $T_{z_1}$ in $G - v_1z_1$. We note that $G[Z]$ is a copy of $K_{2,3}$. Let $I_z$ be the non-canonical independent set of $G[Z]$, and so $I_z$ consists of the two neighbors of $z_1$ in $T_{z_1}$. Let $Y = N_G[u] \cup V(T_v) \cup Z \cup \{z_2\}$. The set $Y$ is indicated by the vertices in the dotted region in Figure~\ref{fig:tsG-2}. We now consider the graph $G'' = G - Y$.

\begin{figure}[htb]
\begin{center}
\begin{tikzpicture}[scale=.75,style=thick,x=0.75cm,y=0.75cm]
\def\vr{2.25pt}
\def\vrn{1.25pt}
\path (1.65,2.5) coordinate (z0);
\path (15.75,2.5) coordinate (z2);
\path (15.75,10.5) coordinate (z3);
\path (14.25,10.5) coordinate (z31);
\path (14.25,9.5) coordinate (z32);
\path (15.25,9.5) coordinate (z33);
\path (15.25,8.5) coordinate (z34);
\path (11.75,8.5) coordinate (z4);
\path (9.5,8.5) coordinate (z6);
\path (9.5,11.35) coordinate (z61);
\path (6,11.35) coordinate (z62);
\path (6,3.85) coordinate (z7);
\path (3.15,3.85) coordinate (z8);
\path (3.15,6.25) coordinate (z9);
\path (1.65,6.25) coordinate (z10);
\path (1,6.5) coordinate (e1);
\path (1,8) coordinate (e11);
\path (1,4.5) coordinate (f1);
\path (2.25,3.5) coordinate (f2);
\path (2.25,4.5) coordinate (f3);
\path (2.25,5.5) coordinate (f4);
\path (2.25,5.425) coordinate (f4p);
\path (3.5,4.5) coordinate (f5);
\path (3.5,6.5) coordinate (g5);
\path (3.5,8) coordinate (g51);
\path (10,4.5) coordinate (c1);
\path (11.25,3.5) coordinate (c2);
\path (11.25,4.5) coordinate (c3);
\path (11.25,5.5) coordinate (c4);
\path (12.5,4.5) coordinate (c5);
\path (10,6.5) coordinate (h1);
\path (10,8) coordinate (h11);
%
\path (12.5,6.5) coordinate (q5);
\path (12.5,8) coordinate (q51);
\path (6.5,10) coordinate (d1);
\path (7.75,9) coordinate (d2);
\path (7.75,10) coordinate (d3);
\path (7.8,10) coordinate (d3p);
\path (7.75,11) coordinate (d4);
\path (9,10) coordinate (d5);
\path (6.5,12) coordinate (j1);
\path (6.5,13.5) coordinate (j11);
\path (9,12) coordinate (j2);
\path (9,13.5) coordinate (j21);
\path (13.25,10) coordinate (k1);
\path (14.75,9) coordinate (k2);
\path (14.75,10) coordinate (k3);
\path (14.7,10) coordinate (k3p);
\path (14.75,11) coordinate (k4);
\path (16.25,10) coordinate (k5);
\path (13.25,12) coordinate (l1);
\path (13.25,13.5) coordinate (l11);
\path (16.25,12) coordinate (l2);
\path (16.25,13.5) coordinate (l21);

\draw (f1)--(e1);
\draw (e1)--(g5);
\draw (e1)--(e11);
\draw (f5)--(g5)--(g51);
\draw (f1)--(f2)--(f3)--(f4)--(f5)--(f2);
\draw (f1)--(f4);
\draw [style=dashed,rounded corners] (-2,2) rectangle (4.85,7.45);
\draw[dotted, line width=0.05cm] (z0)--(z2)--(z3)--(z31)--(z32)--(z33)--(z34)--(z4)--(z6)--(z61)--(z62)--(z7)--(z8)--(z9)--(z10)--(z0);
\draw (h1)--(q5);
\draw (c1)--(c2)--(c3)--(c4)--(c5)--(c2);
\draw (c1)--(c4);
\draw (c1)--(h1);
\draw (d1)--(d2)--(d3)--(d4)--(d5)--(d2);
\draw (d1)--(d4);
\draw (d1)--(j1)--(j11);
\draw (d5)--(j2)--(j1);
\draw (j2)--(j21);
\draw (h1) to[out=90,in=0, distance=0.75cm] (d3);
\draw (k1)--(k2)--(k3)--(k4)--(k5)--(k2);
\draw (k1)--(k4);
\draw (k1)--(l1)--(l11);
\draw (k5)--(l2)--(l1);
\draw (l2)--(l21);
\draw (q5) to[out=90,in=180, distance=0.75cm] (k3);
\draw (c5)--(q5);
%
%
\draw (f3) to[out=270,in=270, distance=1.45cm] (c3);
\draw [style=dashed,rounded corners] (9,2) rectangle (15,7.45);
\draw [style=dashed,rounded corners] (5,8) rectangle (10,12.5);
\draw [style=dashed,rounded corners] (12.5,8) rectangle (17.75,12.5);
\draw (5.75,11.85) node {{\small $T_{z_1}$}};
\draw (17,11.5) node {{\small $T_{z_2}$}};
\draw (e1) [fill=white] circle (\vr);
%
%
%
\draw (g5) [fill=white] circle (\vr);
\draw (f1) [fill=white] circle (\vr);
\draw (f2) [fill=white] circle (\vr);
\draw (f3) [fill=black] circle (\vr);
\draw (f4) [fill=white] circle (\vr);
\draw (f5) [fill=white] circle (\vr);
\draw (h1) [fill=white] circle (\vr);
\draw (c1) [fill=black] circle (\vr);
\draw (c2) [fill=white] circle (\vr);
\draw (c3) [fill=white] circle (\vr);
\draw (c4) [fill=white] circle (\vr);
\draw (c5) [fill=white] circle (\vr);
\draw (d1) [fill=white] circle (\vr);
\draw (d2) [fill=black] circle (\vr);
\draw (d3) [fill=white] circle (\vr);
\draw (d4) [fill=black] circle (\vr);
\draw (d5) [fill=white] circle (\vr);
\draw (j1) [fill=white] circle (\vr);
\draw (j2) [fill=white] circle (\vr);
\draw (k1) [fill=white] circle (\vr);
\draw (k2) [fill=white] circle (\vr);
\draw (k3) [fill=white] circle (\vr);
\draw (k4) [fill=white] circle (\vr);
\draw (k5) [fill=white] circle (\vr);
\draw (l1) [fill=white] circle (\vr);
\draw (l2) [fill=white] circle (\vr);
\draw (q5) [fill=black] circle (\vr);
%
%
\draw[anchor = east] (e1) node {{\small $w_1$}};
\draw[anchor = west] (g5) node {{\small $w_2$}};
\draw[anchor = west] (f5) node {{\small $x_2$}};
\draw[anchor = west] (f3) node {{\small $u$}};
\draw[anchor = east] (f1) node {{\small $x_1$}};
\draw[anchor = east] (h1) node {{\small $v_1$}};
\draw[anchor = west] (q5) node {{\small $v_2$}};
\draw[anchor = west] (c5) node {{\small $y_2$}};
\draw[anchor = east] (c1) node {{\small $y_1$}};
\draw[anchor = west] (c3) node {{\small $v$}};
\draw[anchor = north] (f2) node {{\small $u_2$}};
\draw[anchor = south] (f4p) node {{\small $u_1$}};
\draw[anchor = east] (d3p) node {{\small $z_1$}};
\draw[anchor = west] (k3p) node {{\small $z_2$}};
%
\draw (-1.5,6.75) node {{\small $T_w$}};
\draw (14.5,6.75) node {{\small $T_v$}};
\draw (7,6) node {{\small $Y$}};
\draw (5.75,3.35) node {{\small $e$}};
\end{tikzpicture}
\caption{An illustration of a subgraph of $G$ in the proof of Claim~\ref{claim:no-TT-in-Ge}}
\label{fig:tsG-2}
\end{center}
\end{figure}

Every $i$-set of $G''$ can be extended to an ID-set of $G$ by adding to it the set $I_z \cup \{u,v_2,y_1\}$ (as indicated by the five shaded vertices in Figure~\ref{fig:tsG-2}), and so $i(G) \le i(G'') + 5$. We note that $\w_G(Y) = 48$. Since there are eight $Y$-exit edges, we have $\w(G'') = \w(G) - \w_G(Y) + 8 = \w(G) - 48 + 8 = \w(G) - 40$. Because $\tc(G) = b(G) = 0$ and by the structure of the graph $G''$ (recalling also the properties of bad components and troublesome configurations in Propositions~\ref{prop:cB}(f) and~\ref{prop:troublesome-subgraph}(f)), we infer $\tc(G'') = b(G'') = 0$. Thus, $\Theta(G'') = 0$, and so $\Omega(G'') = \w(G'') = \w(G'') - 40 = \Omega(G) - 40$. Hence, $8i(G) \le 8(i(G'') + 5) \le \Omega(G'') + 40 = \Omega(G)$, a contradiction.~\smallqed

\begin{claim}
\label{claim:no-T-in-Ge}
$\tc(G - e) = 0$ for all $e \in E(G)$.
\end{claim}
\proof Suppose, to the contrary, that $\tc(G - e) \ge 1$ for some edge $e \in E(G)$. By Claim~\ref{claim:no-TT-in-Ge}, $\tc(G - e) \le 1$. Thus, $\tc(G - e) = 1$, that is, there exists an edge $e$ in $G$ whose removal creates a troublesome configuration. Let $T = T_{v_1,v_2}$ be a troublesome configuration in $G-e$ with link vertices $v_1$ and $v_2$, where $\deg_G(v_2) = 3$. Thus, $T$ is obtained from a bad graph $B \in \cB_1$ with root vertex~$v_1$ (of degree~$1$ in $B$). Further, let $u_1$ and $u_2$ be the vertices of $B$ adjacent to $v_1$ and $v_2$, respectively, in $G$. If $v_1$ has degree~$3$ in $G$, then let $z_1$ be the neighbor of $v_1$ different from~$u_1$ and $v_2$. Moreover, let $z_2$ be the neighbor of $v_2$ different from~$u_2$ and $v_1$.

\begin{figure}[htb]
\begin{center}
\begin{tikzpicture}[scale=.75,style=thick,x=0.75cm,y=0.75cm]
\def\vr{2.25pt}
\def\vrn{1.25pt}
%
\path (12,6.5) coordinate (h1);
A
\path (12,4.5) coordinate (p1);
\path (13.5,3.3) coordinate (p2);
\path (13.5,4.5) coordinate (p3);
\path (13.5,5.7) coordinate (p4);
\path (15,4.5) coordinate (p5);
\path (15,6.5) coordinate (q5);
\path (15,8) coordinate (q51);
\path (14.25,8.75) coordinate (q511);
\path (15.75,8.75) coordinate (q512);
\draw (h1)--(q5);
\draw (p5)--(q5)--(q51);
\draw (p1)--(p2)--(p3)--(p4)--(p5)--(p2);
\draw (p1)--(p4);
\draw (h1)--(p1);
\draw (q511)--(q51)--(q512);
\draw (p3) to[out=180,in=270, distance=0.5cm] (h1);
\draw [style=dashed,rounded corners] (9.5,2.75) rectangle (16.25,7.45);
\draw (h1) [fill=black] circle (\vr);
\draw (q5) [fill=white] circle (\vr);
\draw (q51) [fill=white] circle (\vr);
\draw (p1) [fill=white] circle (\vr);
\draw (p2) [fill=white] circle (\vr);
\draw (p3) [fill=white] circle (\vr);
\draw (p4) [fill=white] circle (\vr);
\draw (p5) [fill=black] circle (\vr);
\draw[anchor = east] (h1) node {{\small $v_1$}};
\draw[anchor = west] (q5) node {{\small $v_2$}};
\draw[anchor = west] (p5) node {{\small $u_2$}};
\draw[anchor = east] (p1) node {{\small $u_1$}};
\draw[anchor = west] (p3) node {{\small $u$}};
\draw[anchor = west] (q51) node {{\small $z_2$}};
%
\draw (10.15,6.85) node {{\small $T$}};
\draw (12.5,5.75) node {{\small $e$}};
\end{tikzpicture}
\caption{An illustration of a subgraph of $G$ in the proof of Claim~\ref{claim:no-T-in-Ge.0}}
\label{fig:tsG-0}
\end{center}
\end{figure}
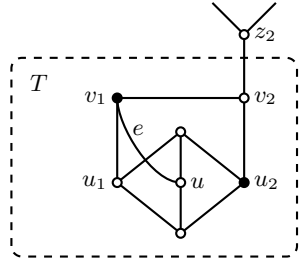

\begin{subclaim}
\label{claim:no-T-in-Ge.0}
The vertex~$v_1$ is not incident with the edge~$e$.
\end{subclaim}
\proof Suppose that $v_1$ is incident with the edge~$e$. In this case, by Claims~\ref{claim:no-K23-with-d2-vertex} and~\ref{claim:induced-K23}, the bad graph $B$ contains exactly one copy of $K_{2,3}$, say $F$, and so $B \in \cB_1$. Moreover, the end of the edge~$e$, say~$u$, different from~$v_1$ belongs to~$F$, for otherwise once again $G$ would contain a copy of $K_{2,3}$ that contains a vertex of degree~$2$ in $G$, a contradiction. Thus, the graph illustrated in Figure~\ref{fig:tsG-0} is a subgraph of $G$, where $T$ is the troublesome configuration indicated in the dashed box and where $e = uv_1$. We now let $X = V(T)$ and we let $G' = G - X$. The removal of the (unique) $X$-exit edge, namely $v_2z_2$ when constructing $G'$ increases the total weight by at most~$3$ by Claim~\ref{claim:cost-of-exit-edge}.  By Claim~\ref{claim:no-bad-in-Ge}, $b(G') = 0$. Thus, $\Omega(G') \le \w(G) - \w_G(X) + 3 = \Omega(G) - 21 + 3 = \Omega(G) - 18$. Every $i$-set of $G'$ can be extended to an $i$-set of $G$ by adding to it the vertices~$v_1$ and $u_2$, and so $i(G) \le i(G') + 2$. Hence, $8i(G) \le 8(i(G') + 2) \le \Omega(G') + 16 = \Omega(G) - 2 < \Omega(G)$, a contradiction.~\smallqed

\medskip
By Claim~\ref{claim:no-T-in-Ge.0}, the vertex~$v_1$ is not incident with the edge~$e$. Let $e = uv$, where $u \in V(B)$.

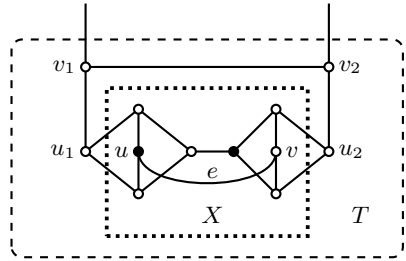
\begin{figure}[htb]
\begin{center}
\begin{tikzpicture}[scale=.75,style=thick,x=0.75cm,y=0.75cm]
\def\vr{2.25pt}
\def\vrn{1.25pt}
\path (9,2.5) coordinate (z0);
\path (13.75,2.5) coordinate (z1);
\path (13.75,6) coordinate (z2);
\path (9,6) coordinate (z3);
\path (8.5,4.5) coordinate (c1);
\path (9.75,3.5) coordinate (c2);
\path (9.75,4.5) coordinate (c3);
\path (9.75,5.5) coordinate (c4);
\path (11,4.5) coordinate (c5);
\path (8.5,6.5) coordinate (h1);
\path (8.5,8) coordinate (h11);
\path (12,4.5) coordinate (p1);
\path (13,3.5) coordinate (p2);
\path (13,4.5) coordinate (p3);
\path (13,5.5) coordinate (p4);
\path (14.25,4.5) coordinate (p5);
\path (14.25,6.5) coordinate (q5);
\path (14.25,8) coordinate (q51);
\draw[dotted, line width=0.05cm] (z0)--(z1)--(z2)--(z3)--(z0);
\draw (c1)--(h1)--(h11);
\draw (h1)--(q5);
\draw (c1)--(c2)--(c3)--(c4)--(c5)--(c2);
\draw (c1)--(c4);
%
%
\draw (p5)--(q5)--(q51);
%
\draw (p1)--(p2)--(p3)--(p4)--(p5)--(p2);
\draw (p1)--(p4);
\draw (c5)--(p1);
\draw (c3) to[out=270,in=270, distance=0.75cm] (p3);
\draw [style=dashed,rounded corners] (6.75,2) rectangle (16,7.15);
\draw (h1) [fill=white] circle (\vr);
%
\draw (c1) [fill=white] circle (\vr);
\draw (c2) [fill=white] circle (\vr);
\draw (c3) [fill=black] circle (\vr);
\draw (c4) [fill=white] circle (\vr);
\draw (c5) [fill=white] circle (\vr);
%
%
\draw (q5) [fill=white] circle (\vr);
%
\draw (p1) [fill=black] circle (\vr);
\draw (p2) [fill=white] circle (\vr);
\draw (p3) [fill=white] circle (\vr);
\draw (p4) [fill=white] circle (\vr);
\draw (p5) [fill=white] circle (\vr);
%
\draw[anchor = east] (h1) node {{\small $v_1$}};
\draw[anchor = west] (q5) node {{\small $v_2$}};
\draw[anchor = west] (p5) node {{\small $u_2$}};
\draw[anchor = east] (c1) node {{\small $u_1$}};
\draw[anchor = east] (c3) node {{\small $u$}};
\draw[anchor = west] (p3) node {{\small $v$}};
%
\draw (15,3) node {{\small $T$}};
\draw (11.5,3) node {{\small $X$}};
\draw (11.5,4) node {{\small $e$}};
\end{tikzpicture}
\caption{An illustration of a subgraph of $G$ in the proof of Claim~\ref{claim:no-T-in-Ge.1}}
\label{fig:tsG-3}
\end{center}
\end{figure}

\begin{subclaim}
\label{claim:no-T-in-Ge.1}
$v \notin V(B)$.
\end{subclaim}
\proof Suppose that $v \in V(B)$, and so $\{u,v\} \subset V(B)$. By Claims~\ref{claim:no-K23-with-d2-vertex} and~\ref{claim:induced-K23}, the bad graph $B$ contains two copies of $K_{2,3}$. Let $F_u$ and $F_v$ be the copies of $K_{2,3}$ that contain~$u$ and~$v$, respectively, where we may assume that $u_1 \in V(F_u)$ (and so, $u_2 \in V(F_v)$). Let $X = V(T) \setminus \{u_1,u_2,v_1,v_2\}$ and let $G' = G - X$. Thus, the graph illustrated in Figure~\ref{fig:tsG-3} is a subgraph of $G$, where $T$ is the troublesome configuration indicated in the dashed box and the set $X$ is indicated by the vertices in the dotted region. Every $i$-set of $G'$ can be extended to an ID-set of $G$ by adding to it two vertices (indicated by the two shaded vertices in Figure~\ref{fig:tsG-3}), and so $i(G) \le i(G') + 2$. From the structure of the graph $G$ and since $b(G) = \tc(G) = 0$, we note that $\Theta(G') = 0$, and so $\Omega(G') = \w(G') = \w(G) - 24 + 4 = \w(G) - 20 = \Omega(G) - 20$. Thus, $8i(G) \le 8(i(G') + 2) \le \Omega(G') + 16 < \Omega(G)$, a contradiction.~\smallqed

\medskip
By Claim~\ref{claim:no-T-in-Ge.1}, we have $v \notin V(B)$. By Claim~\ref{claim:no-K23-with-d2-vertex} and Claim~\ref{claim:induced-K23}, the bad graph $B$ contains only one copy of $K_{2,3}$, say $F_u$. Let $w_1$ and $w_2$ be the two neighbors of $u$ in $B$, and so $N_B(w_1) = N_B(w_2) = \{u,u_1,u_2\}$.

\begin{figure}[htb]
\begin{center}
\begin{tikzpicture}[scale=.75,style=thick,x=0.75cm,y=0.75cm]
\def\vr{2.25pt}
\def\vrn{1.25pt}
\path (10.5,2.5) coordinate (z0);
\path (13.8,2.5) coordinate (z1);
\path (13.8,4.15) coordinate (z2);
\path (12.8,4.15) coordinate (z3);
\path (12.8,4.85) coordinate (z4);
\path (13.95,4.85) coordinate (z5);
\path (13.95,7) coordinate (z6);
\path (10.5,7) coordinate (z7);
%
\path (12,6.5) coordinate (h1);
A
\path (12,4.5) coordinate (p1);
\path (13.5,3.3) coordinate (p2);
\path (13.5,4.5) coordinate (p3);
\path (13.5,5.7) coordinate (p4);
\path (15,4.5) coordinate (p5);
\path (15,6.5) coordinate (q5);
\path (15,8) coordinate (q51);
%
\path (17.5,6.5) coordinate (r);
\path (17,7.25) coordinate (r1);
\path (18,7.25) coordinate (r2);
\draw[dotted, line width=0.05cm] (z0)--(z1)--(z2)--(z3)--(z4)--(z5)--(z6)--(z7)--(z0);
\draw (h1)--(q5);
\draw (p5)--(q5)--(q51);
%
\draw (p1)--(p2)--(p3)--(p4)--(p5)--(p2);
\draw (p1)--(p4);
\draw (h1)--(p1);
\draw (r1)--(r)--(r2);
\draw (p3) to[out=-45,in=270, distance=1.65cm] (r);
\draw [style=dashed,rounded corners] (9.5,2) rectangle (16.25,7.45);
\draw (h1) [fill=white] circle (\vr);
\draw (q5) [fill=black] circle (\vr);
%
\draw (p1) [fill=black] circle (\vr);
\draw (p2) [fill=white] circle (\vr);
\draw (p3) [fill=white] circle (\vr);
\draw (p4) [fill=white] circle (\vr);
\draw (p5) [fill=white] circle (\vr);
\draw (r) [fill=white] circle (\vr);
\draw[anchor = east] (h1) node {{\small $v_1$}};
\draw[anchor = west] (q5) node {{\small $v_2$}};
\draw[anchor = west] (p5) node {{\small $u_2$}};
\draw[anchor = east] (p1) node {{\small $u_1$}};
\draw[anchor = west] (r) node {{\small $v$}};
\draw[anchor = east] (p3) node {{\small $u$}};
\draw[anchor = east] (p2) node {{\small $w_1$}};
\draw[anchor = east] (p4) node {{\small $w_2$}};
%
\draw (17.5,3) node {{\small $T$}};
\draw (11.25,3) node {{\small $X$}};
\draw (17.15,4.75) node {{\small $e$}};
\end{tikzpicture}
\caption{An illustration of a subgraph of $G$ in the proof of Claim~\ref{claim:no-T-in-Ge.2}}
\label{fig:tsG-5}
\end{center}
\end{figure}
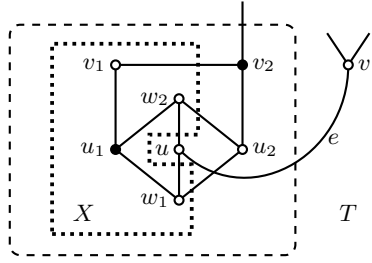

\begin{subclaim}
\label{claim:no-T-in-Ge.2}
$\deg_G(v_1) = 3$.
\end{subclaim}
\proof
Suppose that $\deg_G(v_1) = 2$. Let $X = N_G[u_1] = \{u_1,v_1,w_1,w_2\}$ and let $G' = G - X$. Thus, the graph illustrated in Figure~\ref{fig:tsG-5} is a subgraph of $G$, where $T$ is the troublesome configuration indicated in the dashed box and the set $X$ is indicated by the vertices in the dotted region. Every $i$-set of $G'$ can be extended to an ID-set of $G$ by adding to it the vertex $u_1$ (indicated by the shaded vertes in Figure~\ref{fig:tsG-5}), and so $i(G) \le i(G') + 1$. From the structure of the graph $G$ and since $b(G) = \tc(G) = 0$, we note that $\Theta(G') = 0$, and so $\Omega(G') = \w(G') = \w(G) - 13 + 5 = \w(G) - 8 = \Omega(G) - 8$. Thus, $8i(G) \le 8(i(G') + 1) \le \Omega(G') + 8 = \Omega(G)$, a contradiction.~\smallqed

By Claim~\ref{claim:no-T-in-Ge.2}, we have $\deg_G(v_1) = 3$.

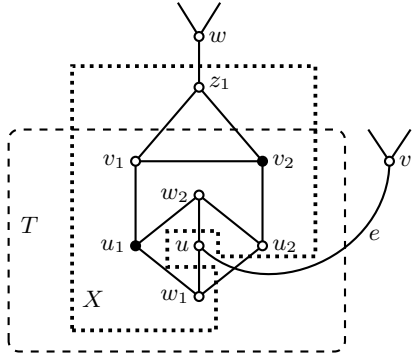
\begin{figure}[htb]
\begin{center}
\begin{tikzpicture}[scale=.75,style=thick,x=0.75cm,y=0.75cm]
\def\vr{2.25pt}
\def\vrn{1.25pt}
\path (10.5,2.5) coordinate (z0);
\path (13.9,2.5) coordinate (z1);
\path (13.9,4) coordinate (z2);
\path (12.75,4) coordinate (z21);
\path (12.75,4.85) coordinate (z22);
\path (13.95,4.85) coordinate (z23);
\path (13.95,4.25) coordinate (z3);
\path (16.25,4.25) coordinate (z4);
\path (16.25,8.75) coordinate (z5);
\path (10.5,8.75) coordinate (z6);
%
\path (12,6.5) coordinate (h1);
A
\path (12,4.5) coordinate (p1);
\path (13.5,3.3) coordinate (p2);
\path (13.5,4.5) coordinate (p3);
\path (13.5,5.7) coordinate (p4);
\path (15,4.5) coordinate (p5);
\path (15,6.5) coordinate (q5);
\path (13.5,8.25) coordinate (q51);
\path (13.5,8.35) coordinate (q51n);
\path (13.5,9.45) coordinate (q52);
\path (13,10.25) coordinate (q53);
\path (14,10.25) coordinate (q54);
%
\path (18,6.5) coordinate (r);
\path (17.5,7.25) coordinate (r1);
\path (18.5,7.25) coordinate (r2);
\draw[dotted, line width=0.05cm] (z0)--(z1)--(z2)--(z21)--(z22)--(z23)--(z3)--(z4)--(z5)--(z6)--(z0);
\draw (h1)--(q5);
\draw (p5)--(q5)--(q51)--(h1);
\draw (q52)--(q51);
\draw (q53)--(q52)--(q54);
\draw (p1)--(p2)--(p3)--(p4)--(p5)--(p2);
\draw (p1)--(p4);
\draw (h1)--(p1);
\draw (r1)--(r)--(r2);
\draw (p3) to[out=-45,in=270, distance=1.65cm] (r);
\draw [style=dashed,rounded corners] (9,2) rectangle (16.95,7.25);
\draw (h1) [fill=white] circle (\vr);
\draw (q5) [fill=black] circle (\vr);
\draw (q51) [fill=white] circle (\vr);
\draw (q52) [fill=white] circle (\vr);
\draw (p1) [fill=black] circle (\vr);
\draw (p2) [fill=white] circle (\vr);
\draw (p3) [fill=white] circle (\vr);
\draw (p4) [fill=white] circle (\vr);
\draw (p5) [fill=white] circle (\vr);
\draw (r) [fill=white] circle (\vr);
\draw[anchor = east] (h1) node {{\small $v_1$}};
\draw[anchor = west] (q5) node {{\small $v_2$}};
\draw[anchor = west] (p5) node {{\small $u_2$}};
\draw[anchor = east] (p1) node {{\small $u_1$}};
\draw[anchor = west] (r) node {{\small $v$}};
\draw[anchor = east] (p3) node {{\small $u$}};
\draw[anchor = west] (q51n) node {{\small $z_1$}};
\draw[anchor = west] (q52) node {{\small $w$}};
\draw[anchor = east] (p2) node {{\small $w_1$}};
\draw[anchor = east] (p4) node {{\small $w_2$}};

\draw (9.5,5) node {{\small $T$}};
\draw (11,3.25) node {{\small $X$}};
\draw (17.65,4.75) node {{\small $e$}};
\end{tikzpicture}
\caption{An illustration of a subgraph of $G$ in the proof of Claim~\ref{claim:no-T-in-Ge.3}}
\label{fig:tsG-7}
\end{center}
\end{figure}

\begin{subclaim}
\label{claim:no-T-in-Ge.3}
$z_1 \ne z_2$.
\end{subclaim}
\proof
Suppose that $z_1 = z_2$, and so $\{v_1,v_2\} \subseteq N_G(z_1)$.  Let $X = (V(T) \setminus \{u\}) \cup \{z_1\}$ and let $G' = G - X$. Every $i$-set of $G'$ can be extended to an ID-set of $G$ by adding to it the vertices $u_1$ and $v_2$, and so $i(G) \le i(G') + 2$. If $N_G(z_1) = \{v_1,v_2\}$, then $\deg_G(z_1) = 2$, $\w_G(X) = 22$ and $\Omega(G') = \w(G') = \w(G) - \w_G(X) + 2 = \w(G) - 22 + 2 = \w(G) - 20$. Thus, $8i(G) \le 8(i(G') + 2) \le \Omega(G') + 16 = \w(G) - 4 < \Omega(G)$, a contradiction. Hence, $\deg_G(z_1) = 3$. Let $w$ be the neighbor of $z_1$ different from $v_1$ and $v_2$. Thus the graph shown in Figure~\ref{fig:tsG-7} is a subgraph of $G$, where $T$ is the troublesome configuration indicated in the dashed box and the set $X$ is indicated by the vertices in the dotted region. As noted earlier, $i(G) \le i(G') + 2$. The removal of the $X$-exit edge $z_1w$ when constructing $G'$ increases the total weight by at most~$3$ by Claim~\ref{claim:cost-of-exit-edge}. Moreover, the removal of the two $X$-exit edges incident with~$u$ when constructing $G'$ increases the total weight by~$2$ noting that the vertex~$u$ does not belong to a bad component or a troublesome configuration in $G'$ but its vertex weight increased by~$2$. Thus, $\Omega(G') \le \Omega(G) - \w_G(X) + 3 + 2 = \Omega(G) - 21 + 3 + 2 = \Omega(G) - 16$. Hence, $8i(G) \le 8(i(G') + 2) \le \Omega(G') + 16 = \Omega(G)$, a contradiction.~\smallqed

\medskip
By Claim~\ref{claim:no-T-in-Ge.3}, we have $z_1 \ne z_2$.

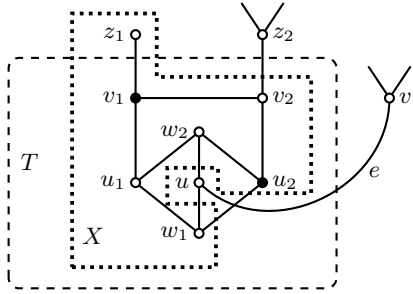
\begin{figure}[htb]
\begin{center}
\begin{tikzpicture}[scale=.75,style=thick,x=0.75cm,y=0.75cm]
\def\vr{2.25pt}
\def\vrn{1.25pt}
\path (10.5,2.5) coordinate (z0);
\path (13.9,2.5) coordinate (z1);
\path (13.9,4) coordinate (z2);
\path (12.75,4) coordinate (z21);
\path (12.75,4.85) coordinate (z22);
\path (13.95,4.85) coordinate (z23);
\path (13.95,4.25) coordinate (z3);
\path (16.15,4.25) coordinate (z4);
\path (16.15,7) coordinate (z5);
\path (12.5,7) coordinate (z51);
\path (12.5,8.5) coordinate (z52);
\path (10.5,8.5) coordinate (z6);

\path (12,6.5) coordinate (h1);
\path (12,8) coordinate (h11);
\path (12,4.5) coordinate (p1);
\path (13.5,3.3) coordinate (p2);
\path (13.5,4.5) coordinate (p3);
\path (13.5,5.7) coordinate (p4);
\path (15,4.5) coordinate (p5);
\path (15,6.5) coordinate (q5);
\path (15,8) coordinate (q51);
\path (14.5,8.75) coordinate (q52);
\path (15.5,8.75) coordinate (q53);
\path (18,6.5) coordinate (r);
\path (17.5,7.25) coordinate (r1);
\path (18.5,7.25) coordinate (r2);
\draw[dotted, line width=0.05cm] (z0)--(z1)--(z2)--(z21)--(z22)--(z23)--(z3)--(z4)--(z5)--(z51)--(z52)--(z6)--(z0);
\draw (p1)--(h1)--(h11);
\draw (h1)--(q5);
\draw (p5)--(q5)--(q51);
\draw (q52)--(q51)--(q53);
\draw (p1)--(p2)--(p3)--(p4)--(p5)--(p2);
\draw (p1)--(p4);
\draw (r1)--(r)--(r2);
\draw (p3) to[out=-45,in=270, distance=1.65cm] (r);
\draw [style=dashed,rounded corners] (9,2) rectangle (16.75,7.45);
\draw (h1) [fill=black] circle (\vr);
\draw (h11) [fill=white] circle (\vr);
\draw (q5) [fill=white] circle (\vr);
\draw (q51) [fill=white] circle (\vr);
\draw (p1) [fill=white] circle (\vr);
\draw (p2) [fill=white] circle (\vr);
\draw (p3) [fill=white] circle (\vr);
\draw (p4) [fill=white] circle (\vr);
\draw (p5) [fill=black] circle (\vr);
\draw (r) [fill=white] circle (\vr);
\draw[anchor = east] (h1) node {{\small $v_1$}};
\draw[anchor = west] (q5) node {{\small $v_2$}};
\draw[anchor = west] (p5) node {{\small $u_2$}};
\draw[anchor = east] (p1) node {{\small $u_1$}};
\draw[anchor = west] (r) node {{\small $v$}};
\draw[anchor = east] (p3) node {{\small $u$}};
\draw[anchor = west] (q51) node {{\small $z_2$}};
\draw[anchor = east] (h11) node {{\small $z_1$}};
\draw[anchor = east] (p2) node {{\small $w_1$}};
\draw[anchor = east] (p4) node {{\small $w_2$}};
\draw (9.5,5) node {{\small $T$}};
\draw (11,3.25) node {{\small $X$}};
\draw (17.65,4.75) node {{\small $e$}};
\end{tikzpicture}
\caption{An illustration of a subgraph of $G$ in the proof of Claim~\ref{claim:no-T-in-Ge.4}}
\label{fig:tsG-8}
\end{center}
\end{figure}

\begin{subclaim}
\label{claim:no-T-in-Ge.4}
$\deg_G(z_i) \ge 2$ for $i \in [2]$.
\end{subclaim}
\proof
Suppose that $\deg_G(z_1) = 1$ or $\deg_G(z_2) = 1$. By symmetry, we may assume that $\deg_G(z_1) = 1$. Let $X = (V(T) \setminus \{u\}) \cup \{z_1\}$ and let $G' = G - X$. Every $i$-set of $G'$ can be extended to an ID-set of $G$ by adding to it the vertices $v_1$ and $u_2$, and so $i(G) \le i(G') + 2$. The graph shown in Figure~\ref{fig:tsG-8} is a subgraph of $G$, where $T$ is the troublesome configuration indicated in the dashed box and the set $X$ is indicated by the vertices in the dotted region. The removal of the $X$-exit edge $z_2v_2$ when constructing $G'$ increases the total weight by at most~$3$. Moreover, the removal of the two $X$-exit edges incident with~$u$ when constructing $G'$ increases the total weight by~$2$. Thus, $\Omega(G') \le \Omega(G) - \w_G(X) + 3 + 2 = \Omega(G) - 23 + 3 + 2 = \Omega(G) - 18$. Hence, $8i(G) \le 8(i(G') + 2) \le \Omega(G') + 16 < \Omega(G)$, a contradiction.~\smallqed

\medskip
By Claim~\ref{claim:no-T-in-Ge.4}, we have $\deg_G(z_i) \ge 2$ for $i \in [2]$.

\begin{subclaim}
\label{claim:no-T-in-Ge.5}
$z_1z_2 \in E(G)$.
\end{subclaim}
\proof
Suppose that $z_1z_2 \notin E(G)$. Let $X = V(T) \setminus \{u\}$ and let $G'$ be obtained from $G - X$ by adding the edge $f = z_1z_2$. Let $I'$ be an $i$-set of $G'$. If $I' \cap \{z_1,z_2\} = \emptyset$, then let $I_u = \{u_1,u_2\}$. If $z_1 \in I'$, then let $I_u = \{u_1,v_2\}$. If $z_2 \in I'$, then let $I_u = \{v_1,u_2\}$. In all cases, we have $|I_u| = 2$ and the set $I'$ can be extended to an ID-set of $G$ by adding to it the set $I_u$. Thus, $i(G) \le |I'| + |I_u| = i(G') + 2$. We note that $\w_G(X) = 18$ and $\w(G') = \w(G) - \w_G(X) + 2 = \w(G) - 18 + 2 = \w(G) - 16$. If $\Theta(G') = 0$, then $\Omega(G') = \w(G')$, and so $8i(G) \le 8(i(G') + 2) \le \Omega(G') + 16 = \w(G') + 16 = \w(G) = \Omega(G)$, a contradiction. Hence, $\Theta(G') > 0$. We therefore infer that the component, $G_f'$ say, of $G'$ that contains the edge~$f$ is either a bad component or contains a troublesome configuration.

\begin{subsubclaim}
\label{claim:no-T-in-Ge.5.1}
$G_f' \notin \cB$.
\end{subsubclaim}
\proof
Suppose that $G_f' \in \cB$. Let $G_f'$ have root vertex~$v_f$ (of degree~$1$ in $B_f$). We note that the degrees of the vertices in $G_f'$ are the same as their degrees in $G$, except possibly for the vertex~$u$ in the case when it belongs to $G_f'$. By Claims~\ref{claim:no-K23-with-d2-vertex} and~\ref{claim:induced-K23}, we infer that $G_f' \in \cB_1$, and so the root $v_f$ of $G_f'$ has degree~$1$ in $G$. Further, we infer from these claims that $G_f'$ contains exactly one copy of $K_{2,3}$, say $H_f$.
By Claim~\ref{claim:no-T-in-Ge.4}, $\deg_G(z_i) \ge 2$ for $i \in [2]$.

Suppose that $u$ belongs to $G_f'$, implying that $u$ is the root vertex~$v_f$. Since the vertex~$u$ is adjacent to neither~$v_1$ nor~$v_2$ in $G$ and since $z_i$ is adjacent to $v_i$ in $G$ for $i \in [2]$, we note that $u \notin \{z_1,z_2\}$. Thus, neither $z_1$ nor $z_2$ is the root vertex~$v_f$ of $G_f'$, and so the added edge $f = z_1z_2$ is an edge of $H_f$. The graph $G$ is now determined and $V(G) = V(T) \cup V(H_f)$. In particular, $G$ has order~$n = 12$. There are two possibilities up to isomorphism, depending on whether $v \in \{z_1,z_2\}$ or $v \notin \{z_1,z_2\}$. In both cases, we have $\Omega(G) = \w(G) = 38$ and the set $\{w_1,w_2,z_1,z_2\}$ is an ID-set of $G$, and so $i(G) \le 4$. Thus, $8i(G) < \Omega(G)$, a contradiction. Hence, $u \notin V(G_f')$, implying that $G'$ has two components, namely the component $G_f'$ and the component containing the vertices~$u$ and~$v$.

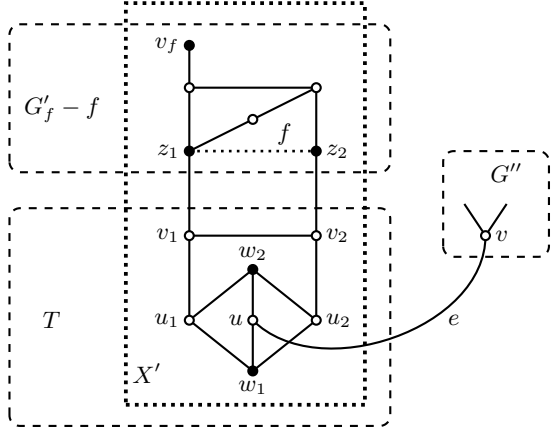
\begin{figure}[htb]
\begin{center}
\begin{tikzpicture}[scale=.75,style=thick,x=0.75cm,y=0.75cm]
\def\vr{2.25pt}
\def\vrn{1.25pt}
\path (10.5,2.5) coordinate (z0);
\path (16.15,2.5) coordinate (z1);
\path (16.15,4.25) coordinate (z4);
\path (16.15,12) coordinate (z5);
\path (10.5,12) coordinate (z6);
%
\path (12,6.5) coordinate (h1);
\path (12,8.5) coordinate (h11);
\path (12,10) coordinate (h12);
\path (12,11) coordinate (h13);
\path (13.5,9.25) coordinate (m);
\path (12,4.5) coordinate (p1);
\path (13.5,3.3) coordinate (p2);
\path (13.5,4.5) coordinate (p3);
\path (13.5,5.7) coordinate (p4);
\path (15,4.5) coordinate (p5);
\path (15,6.5) coordinate (q5);
\path (15,8.5) coordinate (q51);
\path (15,10) coordinate (q52);
\path (19,6.5) coordinate (r);
\path (18.5,7.25) coordinate (r1);
\path (19.5,7.25) coordinate (r2);
\draw[dotted, line width=0.05cm] (z0)--(z1)--(z4)--(z5)--(z6)--(z0);
\draw[dotted, line width=0.035cm] (h11)--(q51);
\draw (p1)--(h1)--(h11)--(h12);
\draw (h1)--(q5);
\draw (h13)--(h12)--(q52);
\draw (h11)--(m)--(q52);
\draw (p5)--(q5)--(q51);
\draw (q52)--(q51);
\draw (p1)--(p2)--(p3)--(p4)--(p5)--(p2);
\draw (p1)--(p4);
\draw (r1)--(r)--(r2);
\draw (p3) to[out=-45,in=270, distance=1.65cm] (r);
\draw [style=dashed,rounded corners] (7.75,2) rectangle (16.75,7.15);
\draw [style=dashed,rounded corners] (7.75,8) rectangle (16.75,11.5);
\draw [style=dashed,rounded corners] (18,6) rectangle (20.5,8.5);
\draw (h1) [fill=white] circle (\vr);
\draw (h11) [fill=black] circle (\vr);
\draw (h12) [fill=white] circle (\vr);
\draw (h13) [fill=black] circle (\vr);
\draw (q5) [fill=white] circle (\vr);
\draw (q51) [fill=black] circle (\vr);
\draw (q52) [fill=white] circle (\vr);
\draw (m) [fill=white] circle (\vr);
\draw (p1) [fill=white] circle (\vr);
\draw (p2) [fill=black] circle (\vr);
\draw (p3) [fill=white] circle (\vr);
\draw (p4) [fill=black] circle (\vr);
\draw (p5) [fill=white] circle (\vr);
\draw (r) [fill=white] circle (\vr);
\draw[anchor = east] (h1) node {{\small $v_1$}};
\draw[anchor = west] (q5) node {{\small $v_2$}};
\draw[anchor = west] (p5) node {{\small $u_2$}};
\draw[anchor = east] (p1) node {{\small $u_1$}};
\draw[anchor = west] (r) node {{\small $v$}};
\draw[anchor = east] (p3) node {{\small $u$}};
\draw[anchor = west] (q51) node {{\small $z_2$}};
\draw[anchor = east] (h11) node {{\small $z_1$}};
\draw[anchor = east] (h13) node {{\small $v_f$}};
\draw[anchor = north] (p2) node {{\small $w_1$}};
\draw[anchor = south] (p4) node {{\small $w_2$}};
\draw (8.75,4.5) node {{\small $T$}};
\draw (9,9.5) node {{\small $G_f' - f$}};
\draw (11,3.25) node {{\small $X'$}};
\draw (19.5,8) node {{\small $G''$}};
\draw (18.25,4.5) node {{\small $e$}};
\draw (14.25,8.9) node {{\small $f$}};
\end{tikzpicture}
\caption{An illustration of a subgraph of $G$ in the proof of Claim~\ref{claim:no-T-in-Ge.5.1}}
\label{fig:tsG-9}
\end{center}
\end{figure}

As observed earlier, $\deg_G(z_i) \ge 2$ for $i \in [2]$, and so the added edge $f = z_1z_2$ is an edge of $H_f$. Let $X' = V(G_f') \cup V(T)$ and let $G'' = G - X'$. Thus the graph shown in Figure~\ref{fig:tsG-9} is an example of a subgraph of $G$, where $T$ is the troublesome configuration indicated in the dashed box, $G_f' - f$ is the subgraph indicated in the dashed box, and the set $X'$ is indicated by the vertices in the dotted region. We note that there exists an ID-set, $I_f$ say, of $G_f' - f$ that contains both vertices $z_1$ and $z_2$ and is such that $|I_f| \le 3$. Every $i$-set of $G''$ can be extended to an ID-set of $G$ by adding to it the set $I_f \cup \{w_1,w_2\}$ (as indicated by the shaded vertices in Figure~\ref{fig:tsG-9}), and so $i(G) \le i(G'') + 5$. By Claims~\ref{claim:no-K23-with-d2-vertex} and~\ref{claim:induced-K23}, we infer that $b(G'') = 0$ (and note that $\tc(G'') = 0$). Thus, $\Omega(G'') = \w(G'') = \w(G) - \w_G(X) + 1 = \Omega(G) - 43 + 1 = \Omega(G) - 42$. Hence, $8i(G) \le 8(i(G') + 5) \le \Omega(G') + 40 < \Omega(G)$, a contradiction.~\smallqed

\medskip
We now return to the proof of Claim~\ref{claim:no-T-in-Ge.5}. By Claim~\ref{claim:no-T-in-Ge.5.1}, $G_f' \notin \cB$, implying that $G_f'$ contains a troublesome configuration, $T_f$ say. Necessarily, $T_f$ contains the added edge $f$.  By Claims~\ref{claim:no-K23-with-d2-vertex} and~\ref{claim:induced-K23}, we infer that $T_f$ contains exactly one copy of $K_{2,3}$, and this copy of $K_{2,3}$ contains the added edge $f = z_1z_2$.  We may assume, by symmetry, that $z_1$ has degree~$3$ in the copy of $K_{2,3}$ in $T_f$. We now let $X^* = (V(T) \setminus \{u\}) \cup \{z_1\}$ and let $G^* = G - X^*$. Every $i$-set of $G^*$ can be extended to an ID-set of $G$ by adding to it the vertices $v_1$ and $u_2$, and so $i(G) \le i(G^*) + 2$.  We note that $\deg_{G^*}(u)=1$ and the degree of a neighbor of $z_1$ in $G^*$ is also $1$. Thus by the structure of the graph $G^*$ so far determined, we infer that $b(G^*) = 0$ and $\tc(G^*) = 0$. Thus, $\Omega(G^*) = \w(G^*) = \w(G) - \w_G(X^*) + 5 = \Omega(G) - 21 + 5 = \Omega(G) - 16$. Hence, $8i(G) \le 8(i(G^*) + 2) \le \Omega(G^*) + 16 = \Omega(G)$, a contradiction.~\smallqed

\medskip
By Claim~\ref{claim:no-T-in-Ge.5}, $z_1z_2 \in E(G)$.

\begin{subclaim}
\label{claim:no-T-in-Ge.6}
$\deg_G(z_i) = 3$ for $i \in [2]$.
\end{subclaim}
\proof
Suppose that $\deg_G(z_i) = 2$ for some $i \in [2]$. Suppose firstly that $\deg_G(z_1) = \deg_G(z_2) = 2$. Thus, $N_G(z_i) = \{v_i,z_{3-i}\}$ for $i \in [2]$. We now let $X = (V(T) \setminus \{u\}) \cup \{z_1,z_2\}$ and let $G' = G - X$. Every $i$-set of $G'$ can be extended to an ID-set of $G$ by adding to it the vertices in the set $\{u_1,u_2,z_1\}$, and so $i(G) \le i(G') + 3$. We note that $\Theta(G') = 0$ and $\w_G(X) = 26$, and so $\Omega(G') = \w(G') = \w(G) - \w_G(X) + 2 = \Omega(G) - 26 + 2 = \Omega(G) - 24$, and so $8i(G) \le 8(i(G') + 3) \le \Omega(G') + 24 = \Omega(G)$, a contradiction. Hence, at least one of $z_1$ and $z_2$ has degree~$3$ in $G$. By symmetry, we may assume that $\deg_G(z_2) = 3$, and so, by supposition, $\deg_G(z_1) = 2$. In this case, we let $X = (V(T) \setminus \{u\}) \cup \{z_1\}$ and let $G' = G - X$. Every $i$-set of $G'$ can be extended to an ID-set of $G$ by adding to it the vertices $v_1$ and $u_2$, and so $i(G) \le i(G') + 2$. We note that $\Theta(G') = 0$ and $\w_G(X) = 22$, and so $\Omega(G') = \w(G') = \w(G) - \w_G(X) + 4 = \Omega(G) - 22 + 4 = \Omega(G) - 18$, and so $8i(G) \le 8(i(G') + 2) \le \Omega(G') + 16 < \Omega(G)$, a contradiction.~\smallqed

\medskip
By Claim~\ref{claim:no-T-in-Ge.6}, $\deg_G(z_i) = 3$ for $i \in [2]$. Let $y_i$ be the neighbor of $z_i$ different from~$v_i$ and $z_{3-i}$ for $i \in [2]$.

\begin{subclaim}
\label{claim:no-T-in-Ge.7}
The following properties hold. \\ [-20pt]
\begin{enumerate}
\item[{\rm (a)}] $\deg_G(y_i) \ge 2$ for $i \in [2]$.
\item[{\rm (b)}] $y_1 \ne y_2$.
\end{enumerate}
\end{subclaim}
\proof (a) Suppose that at least one of $y_1$ and $y_2$ has degree~$1$ in $G$. By symmetry, we may assume that $\deg_G(y_1) = 1$. In this case, we let $X = N_G[z_1] = \{v_1,y_1,z_1,z_2\}$ and let $G' = G - X$. Every $i$-set of $G'$ can be extended to an ID-set of $G$ by adding to it the vertex $z_1$, and so $i(G) \le i(G') + 1$. The cost of removing the $X$-exit edge $y_2z_2$ when constructing $G'$ increases the total weight by at most~$3$ by Claim~\ref{claim:cost-of-exit-edge}. Moreover, the removal of the three $X$-exit edges different from~$y_2z_2$ increases the total weight by~$3$ since the removal of these three edges does not create a bad component or a troublesome configuration. Thus, $\Omega(G') \le \Omega(G) - \w_G(X) + 3 + 3 = \Omega(G) - 14 + 3 + 3  = \Omega(G) - 8$. Hence, $8i(G) \le 8(i(G') + 1) \le \Omega(G') + 8 = \Omega(G)$, a contradiction. This proves part~(a).

(b) Suppose that $y_1 = y_2$. In this case, we let $X = (V(T) \setminus \{u\}) \cup \{y_1,z_1,z_2\}$ and let $G' = G - X$. Every $i$-set of $G'$ can be extended to an ID-set of $G$ by adding to it the vertices in the set $\{u_1,u_2,z_1\}$, and so $i(G) \le i(G') + 3$. We note that $\Theta(G') = 0$ and $\w_G(X) = 27$, and so $\Omega(G') = \w(G') \le \w(G) - \w_G(X) + 3 = \Omega(G) - 27 + 3 = \Omega(G) - 24$, and so $8i(G) \le 8(i(G') + 3) \le \Omega(G') + 24 = \Omega(G)$, a contradiction. This proves part~(b).~\smallqed

\medskip
By Claim~\ref{claim:no-T-in-Ge.7}, $\deg_G(y_i) \ge 2$ for $i \in [2]$ and $y_1 \ne y_2$. Suppose that $y_1$ does not belong to a troublesome configuration in $G - y_1z_1$. In this case, we let $X = (V(T) \setminus \{u\}) \cup \{z_1\}$ and let $G' = G - X$. Every $i$-set of $G'$ can be extended to an ID-set of $G$ by adding to it the vertices $v_1$ and $u_2$, and so $i(G) \le i(G') + 2$. We note that both $u$ and $z_2$ have degree~$1$ in $G'$ and do not belong to a bad component or a troublesome configuration in $G'$. By supposition, $y_1$ does not belong to a troublesome configuration in $G'$. Moreover, removing the edge $y_1z_1$ does not create a bad component containing~$y_1$. Hence, $\Theta(G') = 0$ and $\w_G(X) = 21$, and so $\Omega(G') = \w(G') \le \w(G) - \w_G(X) + 5 = \Omega(G) - 21 + 5 = \Omega(G) - 16$, and so $8i(G) \le 8(i(G') + 2) \le \Omega(G') + 16 = \Omega(G)$, a contradiction. Hence, $y_1$ belongs to a troublesome configuration in $G - y_1z_1$. By symmetry, $y_2$ belongs to a troublesome configuration in $G - y_2z_2$. Let $T_{y_i}$ be the troublesome configuration in $G - y_iz_i$ that contains the vertex~$y_i$ for $i \in [2]$.

\begin{figure}[htb]
\begin{center}
\begin{tikzpicture}[scale=.75,style=thick,x=0.75cm,y=0.75cm]
\def\vr{2.25pt}
\def\vrn{1.25pt}
\path (6.75,2) coordinate (z0);
\path (11.75,2) coordinate (z1);
\path (11.75,3) coordinate (z11);
\path (10.6,3) coordinate (z2);
\path (10.6,3.85) coordinate (z3);
\path (11.75,3.85) coordinate (z4);
\path (11.75,3.25) coordinate (z5);
\path (13.75,3.25) coordinate (z6);
\path (13.75,7) coordinate (z7);
\path (8.5,7) coordinate (z8);
\path (8.5,11.35) coordinate (z9);
\path (6.75,11.35) coordinate (z10);
\path (10,3.5) coordinate (c1);
\path (11.25,2.5) coordinate (c2);
\path (11.25,3.5) coordinate (c3);
\path (11.25,4.5) coordinate (c4);
\path (12.5,3.5) coordinate (c5);
\path (10,5.5) coordinate (h1);
\path (10,6.5) coordinate (h11);
\path (12.5,5.5) coordinate (q5);
\path (12.5,6.5) coordinate (q51);
\path (6.25,10) coordinate (d1);
\path (7.75,9) coordinate (d2);
\path (7.75,10) coordinate (d3);
\path (7.9,10) coordinate (d3p);
\path (7.75,11) coordinate (d4);
\path (9,10) coordinate (d5);
\path (6.25,12) coordinate (j1);
\path (6.25,13.5) coordinate (j11);
\path (9,12) coordinate (j2);
\path (9,13.5) coordinate (j21);
\path (13.25,10) coordinate (k1);
\path (14.75,9) coordinate (k2);
\path (14.75,10) coordinate (k3);
\path (14.7,10) coordinate (k3p);
\path (14.75,11) coordinate (k4);
\path (16.25,10) coordinate (k5);
\path (13.25,12) coordinate (l1);
\path (13.25,13.5) coordinate (l11);
\path (16.25,12) coordinate (l2);
\path (16.25,13.5) coordinate (l21);
\path (19,6.5) coordinate (r);
\path (18.5,7.25) coordinate (r1);
\path (19.5,7.25) coordinate (r2);
%
\draw[dotted, line width=0.05cm] (z0)--(z1)--(z11)--(z2)--(z3)--(z4)--(z5)--(z6)--(z7)--(z8)--(z9)--(z10)--(z0);
\draw (h1)--(q5);
\draw (c1)--(c2)--(c3)--(c4)--(c5)--(c2);
\draw (c1)--(c4);
\draw (c1)--(h1)--(h11)--(q51)--(q5);
\draw (d1)--(d2)--(d3)--(d4)--(d5)--(d2);
\draw (d1)--(d4);
\draw (d1)--(j1)--(j11);
\draw (d5)--(j2)--(j1);
\draw (j2)--(j21);
\draw (h11) to[out=90,in=0, distance=0.75cm] (d3);
\draw (k1)--(k2)--(k3)--(k4)--(k5)--(k2);
\draw (k1)--(k4);
\draw (k1)--(l1)--(l11);
\draw (k5)--(l2)--(l1);
\draw (l2)--(l21);
\draw (q51) to[out=90,in=180, distance=0.75cm] (k3);
\draw (c5)--(q5);
\draw (r1)--(r)--(r2);
\draw (c3) to[out=-45,in=270, distance=1.65cm] (r);
\draw [style=dashed,rounded corners] (8.75,1.5) rectangle (15,7.45);
\draw [style=dashed,rounded corners] (4.75,8) rectangle (10,12.5);
\draw [style=dashed,rounded corners] (12.5,8) rectangle (18,12.5);
\draw (5.5,11.85) node {{\small $T_{y_1}$}};
\draw (17.25,11.85) node {{\small $T_{y_2}$}};
\draw (h1) [fill=white] circle (\vr);
\draw (h11) [fill=white] circle (\vr);
\draw (c1) [fill=black] circle (\vr);
\draw (c2) [fill=white] circle (\vr);
\draw (c3) [fill=white] circle (\vr);
\draw (c4) [fill=white] circle (\vr);
\draw (c5) [fill=white] circle (\vr);
\draw (d1) [fill=white] circle (\vr);
\draw (d2) [fill=white] circle (\vr);
\draw (d3) [fill=black] circle (\vr);
\draw (d4) [fill=white] circle (\vr);
\draw (d5) [fill=white] circle (\vr);
\draw (j1) [fill=white] circle (\vr);
\draw (j2) [fill=white] circle (\vr);
\draw (k1) [fill=white] circle (\vr);
\draw (k2) [fill=white] circle (\vr);
\draw (k3) [fill=white] circle (\vr);
\draw (k4) [fill=white] circle (\vr);
\draw (k5) [fill=white] circle (\vr);
\draw (l1) [fill=white] circle (\vr);
\draw (l2) [fill=white] circle (\vr);
\draw (q5) [fill=black] circle (\vr);
\draw (q51) [fill=white] circle (\vr);
\draw (r) [fill=white] circle (\vr);
\draw[anchor = east] (h1) node {{\small $v_1$}};
\draw[anchor = west] (q5) node {{\small $v_2$}};
\draw[anchor = east] (h11) node {{\small $z_1$}};
\draw[anchor = west] (q51) node {{\small $z_2$}};
\draw[anchor = west] (c5) node {{\small $u_2$}};
\draw[anchor = east] (c1) node {{\small $u_1$}};
\draw[anchor = east] (c3) node {{\small $u$}};
\draw[anchor = east] (d3p) node {{\small $y_1$}};
\draw[anchor = west] (k3p) node {{\small $y_2$}};
\draw[anchor = west] (r) node {{\small $v$}};

\draw (14.5,6.75) node {{\small $T$}};
\draw (7.5,5) node {{\small $X$}};
\end{tikzpicture}
\caption{An illustration of a subgraph of $G$ in the proof of Claim~\ref{claim:no-T-in-Ge}}
\label{fig:tsG-end}
\end{center}
\end{figure}
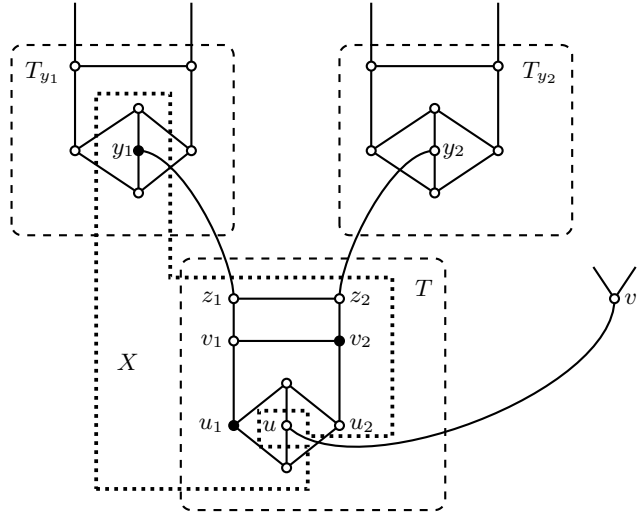

We now let $X = (V(T) \setminus \{u\}) \cup \{z_2\} \cup N_G[y_1]$ and let $G' = G - X$. Thus the graph shown in Figure~\ref{fig:tsG-end} is a subgraph of $G$, where $T$, $T_{y_1}$ and $T_{y_2}$ are indicated in the dashed box, and the set $X$ is indicated by the vertices in the dotted region. Every $i$-set of $G'$ can be extended to an ID-set of $G$ by adding to it the vertices in the set $\{u_1,v_2,y_1\}$ (indicated by the shaded vertices in Figure~\ref{fig:tsG-end}), and so $i(G) \le i(G') + 3$. The cost of removing the $X$-exit edge $y_2z_2$ when constructing $G'$ increases the total weight by~$3$. Moreover, the removal of the six $X$-exit edges different from~$y_2z_2$ increases the total weight by~$6$ since the removal of these six edges does not create a bad component or a troublesome configuration. Thus, $\Omega(G') \le \Omega(G) - \w_G(X) + 3 + 6 = \Omega(G) - 33 + 3 + 6  = \Omega(G) - 24$. Hence, $8i(G) \le 8(i(G') + 3) \le \Omega(G') + 24 = \Omega(G)$, a contradiction. This completes the proof of Claim~\ref{claim:no-T-in-Ge}.~\smallqed

\smallskip
We present next two consequences of Claim~\ref{claim:no-T-in-Ge}.

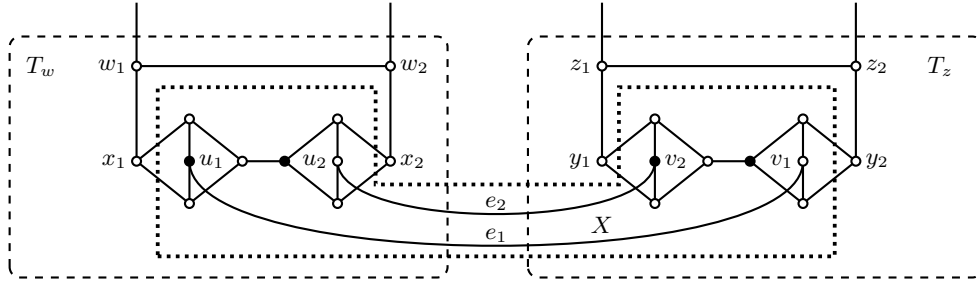
\begin{figure}[htb]
\begin{center}
\begin{tikzpicture}[scale=.75,style=thick,x=0.75cm,y=0.75cm]
\def\vr{2.25pt}
\def\vrn{1.25pt}
\path (-2,2.25) coordinate (z0);
\path (14,2.25) coordinate (z2);
\path (14,6.25) coordinate (z3);
\path (8.9,6.25) coordinate (z6);
\path (8.9,3.95) coordinate (z7);
\path (3.15,3.95) coordinate (z8);
\path (3.15,6.25) coordinate (z9);
\path (-2,6.25) coordinate (z10);
\path (-2.5,6.75) coordinate (e1);
\path (-2.5,8.25) coordinate (e11);
\path (-2.5,4.5) coordinate (r1);
\path (-1.25,3.5) coordinate (r2);
\path (-1.25,4.5) coordinate (r3);
\path (-1.25,5.5) coordinate (r4);
\path (-1.25,5.425) coordinate (r4p);
\path (0,4.5) coordinate (r5);
\path (1,4.5) coordinate (f1);
\path (2.25,3.5) coordinate (f2);
\path (2.25,4.5) coordinate (f3);
\path (2.25,5.5) coordinate (f4);
\path (2.25,5.425) coordinate (f4p);
\path (3.5,4.5) coordinate (f5);
%
\path (3.5,6.75) coordinate (g5);
\path (3.5,8.25) coordinate (g51);
\path (8.5,4.5) coordinate (s1);
\path (9.75,3.5) coordinate (s2);
\path (9.75,4.5) coordinate (s3);
\path (9.75,5.5) coordinate (s4);
\path (11,4.5) coordinate (s5);
\path (8.5,6.75) coordinate (h1);
\path (8.5,8.25) coordinate (h11);
\path (12,4.5) coordinate (c1);
\path (13.25,3.5) coordinate (c2);
\path (13.25,4.5) coordinate (c3);
\path (13.25,5.5) coordinate (c4);
\path (14.5,4.5) coordinate (c5);
%
\path (14.5,6.75) coordinate (q5);
\path (14.5,8.25) coordinate (q51);
\draw (r1)--(e1);
\draw (e1)--(g5);
\draw (e1)--(e11);
\draw (f1)--(r5);
\draw (f5)--(g5)--(g51);
\draw (f1)--(f2)--(f3)--(f4)--(f5)--(f2);
\draw (f1)--(f4);
\draw (r1)--(r2)--(r3)--(r4)--(r5)--(r2);
\draw (r1)--(r4);
\draw [style=dashed,rounded corners] (-5.5,1.75) rectangle (4.85,7.45);
\draw[dotted, line width=0.05cm] (z0)--(z2)--(z3)--(z6)--(z7)--(z8)--(z9)--(z10)--(z0);
\draw (s1)--(h1)--(h11);
\draw (h1)--(q5);
\draw (c1)--(c2)--(c3)--(c4)--(c5)--(c2);
\draw (c1)--(c4);
\draw (c1)--(s5);
\draw (s1)--(s2)--(s3)--(s4)--(s5)--(s2);
\draw (s1)--(s4);
\draw (c5)--(q5)--(q51);
%
%
\draw (f3) to[out=270,in=270, distance=1.25cm] (s3);
\draw (r3) to[out=270,in=270, distance=2cm] (c3);
\draw [style=dashed,rounded corners] (6.75,1.75) rectangle (17.5,7.45);
\draw (e1) [fill=white] circle (\vr);
%
%
%
\draw (g5) [fill=white] circle (\vr);
\draw (f1) [fill=black] circle (\vr);
\draw (f2) [fill=white] circle (\vr);
\draw (f3) [fill=white] circle (\vr);
\draw (f4) [fill=white] circle (\vr);
\draw (f5) [fill=white] circle (\vr);
\draw (r1) [fill=white] circle (\vr);
\draw (r2) [fill=white] circle (\vr);
\draw (r3) [fill=black] circle (\vr);
\draw (r4) [fill=white] circle (\vr);
\draw (r5) [fill=white] circle (\vr);
\draw (s1) [fill=white] circle (\vr);
\draw (s2) [fill=white] circle (\vr);
\draw (s3) [fill=black] circle (\vr);
\draw (s4) [fill=white] circle (\vr);
\draw (s5) [fill=white] circle (\vr);
\draw (h1) [fill=white] circle (\vr);
\draw (c1) [fill=black] circle (\vr);
\draw (c2) [fill=white] circle (\vr);
\draw (c3) [fill=white] circle (\vr);
\draw (c4) [fill=white] circle (\vr);
\draw (c5) [fill=white] circle (\vr);
%
%
\draw (q5) [fill=white] circle (\vr);
%
%
\draw[anchor = east] (e1) node {{\small $w_1$}};
\draw[anchor = west] (g5) node {{\small $w_2$}};
\draw[anchor = west] (f5) node {{\small $x_2$}};
\draw[anchor = east] (r1) node {{\small $x_1$}};
\draw[anchor = east] (h1) node {{\small $z_1$}};
\draw[anchor = west] (q5) node {{\small $z_2$}};
\draw[anchor = west] (c5) node {{\small $y_2$}};
\draw[anchor = east] (s1) node {{\small $y_1$}};
\draw[anchor = east] (c3) node {{\small $v_1$}};
%
\draw[anchor = west] (s3) node {{\small $v_2$}};
%
\draw[anchor = west] (r3) node {{\small $u_1$}};
%
\draw[anchor = east] (f3) node {{\small $u_2$}};
\draw (-4.75,6.75) node {{\small $T_w$}};
\draw (16.5,6.75) node {{\small $T_z$}};
\draw (8.5,3) node {{\small $X$}};
\draw (6,3.5) node {{\small $e_2$}};
\draw (6,2.75) node {{\small $e_1$}};
\end{tikzpicture}
\caption{An illustration of a subgraph of $G$ in the proof of Claim~\ref{claim:no-T-in-Gee}}
\label{fig:tsG-e1-e2}
\end{center}
\end{figure}

\begin{claim}
\label{claim:no-T-in-Gee}
$\tc(G - e_1 - e_2) \le 1$ for all $e_1,e_2 \in E(G)$.
\end{claim}
\proof Suppose, to the contrary, that the claim is false. Thus there exists edges $e_1 = u_1v_1$ and $e_2 = u_2v_2$ in $G$ whose removal creates two troublesome configurations. Let $T_w = T_{w_1,w_2}$ and $T_z = T_{z_1,z_2}$ be the two troublesome configuration in $G - e_1 - e_2$, where $w_1$ and $w_2$ are the link vertices of $T_w$ and $z_1$ and $z_2$ are the link vertices of $T_z$. By Claims~\ref{claim:no-K23-with-d2-vertex},~\ref{claim:induced-K23}, and~\ref{claim:no-T-in-Ge}, we infer that the edges $e_1$ and $e_2$ both have one end in $T_w$ and one end in $T_z$. Furthermore, each of $T_w$ and $T_z$ contains exactly two copies of $K_{2,3}$ and each copy of $K_{2,3}$ contains an end of one of the edges $e_1$ and $e_2$.
We now let $X = (V(T_w) \setminus \{w_1,w_2,x_1,x_2\}) \cup (V(T_z) \setminus \{y_1,y_2,z_1,z_2\})$ and let $G' = G - X$. Renaming vertices if necessary, we may assume that the graph illustrated in Figure~\ref{fig:tsG-e1-e2} is a subgraph of $G$, where $T_w$ and $T_z$ are the troublesome configurations indicated in the dashed boxes and where $X$ is the set indicated by the vertices in the dotted region. We note that $b(G') = \tc(G') = 0$, and so $\Theta(G') = 0$. Since there are eight $X$-exit edges, we have $\Omega(G') = \w(G') = \w(G) - \w_G(X) + 8 = \w(G) - 48 + 8 = \Omega(G) - 40$. Every $i$-set of $G'$ can be extended to an ID-set of $G$ by adding to it four vertices from the set $X$ (indicated, for example, by the shaded vertices in Figure~\ref{fig:tsG-e1-e2}), and so $i(G) \le i(G') + 4$. Thus, $8i(G) \le 8(i(G') + 4) \le \Omega(G') + 32 = \Omega(G) - 8 < \Omega(G)$, a contradiction.~\smallqed

\begin{figure}[htb]
\begin{center}
\begin{tikzpicture}[scale=.75,style=thick,x=0.8cm,y=0.8cm]
\def\vr{2.5pt} 
\path (0.5,-0.1) coordinate (v1);
\path (3.5,-0.1) coordinate (v2);
\path (2,1) coordinate (u3);
\path (2,1.75) coordinate (u);
\path (2,2.5) coordinate (u4);
\path (0.5,1.75) coordinate (u1);
\path (3.5,1.75) coordinate (u2);
%
\draw (u1)--(u3)--(u)--(u4)--(u2)--(u3);
\draw (v1)--(u1);
\draw (v2)--(u2);
\draw (u1)--(u4);
\draw (v1)--(v2);
\draw (v1) [fill=white] circle (\vr);
\draw (v2) [fill=white] circle (\vr);
\draw (u) [fill=white] circle (\vr);
\draw (u1) [fill=white] circle (\vr);
\draw (u2) [fill=white] circle (\vr);
\draw (u3) [fill=white] circle (\vr);
\draw (u4) [fill=white] circle (\vr);
%
%
\draw[anchor = east] (v1) node {{\small $u_1$}};
\draw[anchor = west] (v2) node {{\small $u_2$}};
\draw[anchor = east] (u) node {{\small $v_3$}};
\draw[anchor = east] (u1) node {{\small $v_1$}};
\draw[anchor = west] (u2) node {{\small $v_2$}};
\draw[anchor = north] (u3) node {{\small $x_1$}};
\draw[anchor = south] (u4) node {{\small $x_2$}};
\end{tikzpicture}
\caption{The graph $G_7$}
\label{fig:G7G8}
\end{center}
\end{figure}
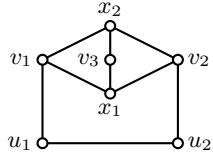

\medskip
Let $G_7$ be the graph shown in Figure~\ref{fig:G7G8}. As a special case of Claim~\ref{claim:no-T-in-Ge}, we have the following result.

\begin{claim}
\label{no-G7}
The graph $G$ does not contain $G_7$ as a subgraph.
\end{claim}
\proof Suppose, to the contrary, that $G$ contains $G_7$ as a subgraph. Let $F$ be such a $G_7$-configuration and let the vertices of $F$ be named as in Figure~\ref{fig:G7G8}. By Claim~\ref{claim:no-K23-with-d2-vertex}, there is no $K_{2,3}$-subgraph in $G$ that contains a vertex of degree~$2$ in $G$. Hence, $\deg_G(v_3) = 3$. Let $v$ be the neighbor of~$v_3$ different from~$x_1$ and~$x_2$. By Claim~\ref{claim:induced-K23}, either $v \in \{u_1,u_2\}$ or $v \notin V(F)$. We show that at least one of $u_1$ and $u_2$ has a neighbor in $G$ that does not belong to~$V(F)$. If $v = u_1$ and $\deg_G(u_2) = 2$, then $G = G_7 + u_1v_3$. In this case, $i(G) = 2$ and $\w(G) = 22$, and so $8i(G) < \Omega(G)$, a contradiction. Hence if $v = u_1$, then $\deg_G(u_2) = 3$, and so $u_2$ has a neighbor in $G$ that does not belong to~$V(F)$. By symmetry, if $v = u_2$, then $u_1$ has a neighbor in $G$ that does not belong to~$V(F)$. If $v \notin V(F)$ and $\deg_G(u_1) = \deg_G(u_2) = 2$, then let $X' = V(F) \setminus \{v_3\}$ and let $G' = G - X'$. Every $i$-set of $G'$ can be extended to an ID-set of $G$ by adding to it the vertices~$v_1$ and~$v_2$, and so $i(G) \le i(G') + 2$. We note that $\Theta(G') = 0$ and $\w_G(X') = 20$, and so $\Omega(G') = \w(G') = \w(G) - \w_G(X) + 2 = \Omega(G) - 20 + 2 = \Omega(G) - 18$, and so $8i(G) \le 8(i(G') + 2) \le \Omega(G') + 16 < \Omega(G)$, a contradiction. Hence if $v \notin V(F)$, then at least one of $u_1$ and $u_2$ has a neighbor in $G$ that does not belong to~$V(F)$. In all cases, at least one of $u_1$ and $u_2$ has a neighbor in $G$ that does not belong to~$V(F)$. Thus, removing the edge $v_3v$ creates a troublesome configuration, contradicting Claim~\ref{claim:no-T-in-Ge}.~\smallqed

\subsection{Part 2: The minimum degree of $G$ equals~$1$}
\label{S:part2}

In this section, we consider the case when the minimum degree of $G$ equals~$1$, that is, when $\delta(G) = 1$. We prove a number of claims that will culminate in a contradiction, thereby showing that this case cannot occur. A \emph{support vertex} in $G$ is a vertex that has at least one neighbor that has degree~$1$, and we define a \emph{strong support vertex} as a vertex that has at least two neighbors of degree~$1$ in $G$.

\begin{claim}
\label{claim:min-deg2}
$\delta(G) \ge 2$.
\end{claim}
\proof Suppose, to the contrary, that $G$ contains a vertex of degree~$1$.

\begin{subclaim}
\label{claim:no-strong-support}
There is no strong support vertex in $G$.
\end{subclaim}
\proof Suppose, to the contrary, that $G$ contains a strong support vertex~$v$. Let $N_G(v) = \{v_1,v_2,v_3\}$ where $\deg_G(v_1) = \deg_G(v_2) = 1$. Since $n \ge 6$, we note that $\deg_G(v_3) \ge 2$. Suppose that $\deg_G(v_3) = 3$ and the two neighbors of $v_3$ different from $v$ both have degree~$1$ in $G$. In this case, the graph $G$ is determined (and is a double star $S(2,2)$) and $\w(G) = 26$. Moreover, $i(G) = 3$, and so $8i(G) < \w(G)$, a contradiction. Hence, at most one neighbor of $v_3$ has degree~$1$.

Suppose that $v_3$ has a neighbor, say $v_4$, of degree~$1$. Since $n \ge 6$, we note that $\deg_G(v_3) = 3$. Let $v_5$ be the third neighbor of $v_3$ different from $v$ and $v_4$. By our earlier observations, $\deg_G(v_5) \ge 2$. We now let $X = N_G[v] \cup \{v_4\}$ and let $G' = G - X$. By Claim~\ref{claim:no-bad-in-Ge}, $b(G') = 0$ and by Claim~\ref{claim:no-T-in-Ge}, $\tc(G') = 0$. Hence, $\Theta(G') = 0$, and so since there is exactly one $X$-exit edge, namely $v_3v_5$, we have $\Omega(G') = \w(G') = \w(G) - \w_G(X) + 1 = \Omega(G) - 21 + 1 = \Omega(G) - 20$. Every $i$-set of $G'$ can be extended to an ID-set of $G$ by adding to it the vertices $v$ and $v_4$, and so $i(G) \le i(G') + 2$. Thus, $8i(G) \le 8(i(G') + 2) \le \Omega(G') + 16 = \Omega(G) - 4 < \Omega(G)$, a contradiction. Hence, no neighbor of $v_3$ has degree~$1$ in $G$.

We now let $X = N_G[v]$ and let $G' = G - X$. By Claim~\ref{claim:no-bad-in-Gee}, $b(G') = 0$ and by Claim~\ref{claim:no-T-in-Gee}, $\tc(G') \le 1$. Hence, $\Theta(G') = 2\tc(G') \le 2$. Since there are two $X$-exit edges, neither of which is incident to a vertex of degree~$1$ in $G$, we have $\Omega(G') = \w(G') + \Theta(G') \le (\w(G) - \w_G(X) + 2) + 2 \le \Omega(G) - 16 + 2 + 2 = \Omega(G) - 12$. Every $i$-set of $G'$ can be extended to an ID-set of $G$ by adding to it the vertex $v$, and so $i(G) \le i(G') + 1$. Thus, $8i(G) \le 8(i(G') + 1) \le \Omega(G') + 8 = \Omega(G) - 4 < \Omega(G)$, a contradiction.~\smallqed

\medskip
By Claim~\ref{claim:no-strong-support}, every support vertex in $G$ has exactly one leaf neighbor.

\begin{subclaim}
\label{claim:no-deg2-support}
No vertex of degree~$2$ is a support vertex.
\end{subclaim}
\proof Suppose, to the contrary, that $G$ contains a support vertex $v$ of degree~$2$. Let $v_1$ be the leaf neighbor of $v$, and let $v_2$ be the second neighbor of $v$. Since $n \ge 6$, we note that $\deg_G(v_2) \ge 2$. We now consider the connected subcubic graph $G' = G - \{v,v_1\}$.  By Claims~\ref{claim:no-bad-in-Ge} and~\ref{claim:no-T-in-Ge}, we have $\Theta(G') = 0$. Thus, $\Omega(G') = \w(G') = \w(G) - 9 + 1 = \Omega(G) - 8$. Thus, $8i(G) \le 8(i(G') + 1) \le \Omega(G') + 8 = \Omega(G)$, a contradiction.~\smallqed

\medskip
Let $P \colon v_1 \ldots v_k$ be the longest path in $G$ such that vertex $v_i$ is adjacent to a vertex, say $u_i$, of degree~$1$ for all $i \in [k]$. Possibly, $k = 1$. By Claim~\ref{claim:no-deg2-support}, $\deg_G(v_i) = 3$ for all $i \in [k]$.

\begin{subclaim}
\label{claim:min-deg2.a}
$k \le 2$.
\end{subclaim}
\proof  Suppose, to the contrary, that $k \ge 3$. If $v_1v_k \in E(G)$, then the graph $G$ is determined and $G$ is the corona of a cycle of length~$k$, and so $i(G) = k$ and $\w(G) = 8k$, whence $8i(G) = 8k = \w(G)$, a contradiction. Hence, $v_1v_k \notin E(G)$. Let $G'$ be obtained from $G - \{u_2,v_2\}$ by adding the edge $v_1v_3$. We note that $\w(G) = \w(G') + 8$. If $G'$ contains a troublesome configuration, then such a configuration has a copy of $K_{2,3}$ that contains a vertex of degree~$2$ in $G$, contradicting Claim~\ref{claim:no-K23-with-d2-vertex}. Thus, $\tc(G') = 0$. Moreover since the connected graph $G'$ contains at least two vertices of degree~$1$, we note that $G' \notin \cB$, and so $b(G') = 0$. Thus, $\Theta(G') = 0$, and so $\Omega(G') = \w(G') = \w(G) - 8 = \Omega(G) - 8$. Every $i$-set of $G'$ can be extended to an ID-set of $G$ by adding to it the vertex~$u_2$, and so $i(G) \le i(G') + 1$. Hence, $8i(G) \le 8(i(G') + 1) \le \Omega(G') + 8 = \Omega(G)$, a contradiction.

\medskip
By Claim~\ref{claim:min-deg2.a}, we have $k \le 2$. If $k = 1$, then we let $N_G(v_1) = \{a,b,u_1\}$, and if $k = 2$, then we let $N_G(v_1) = \{a,u_1,v_2\}$ and $N_G(v_2) = \{b,u_2,v_1\}$, where possibly $a = b$.

\begin{figure}[htb]
\begin{center}
\begin{tikzpicture}[scale=.75,style=thick,x=0.75cm,y=0.75cm]
\def\vr{2.25pt}
\def\vrn{1.25pt}
\path (1.5,2.5) coordinate (a);
\path (1.5,3.25) coordinate (a1);
\path (1.5,3.4) coordinate (a1p);
\path (0.8,4) coordinate (b1);
\path (2.2,4) coordinate (b2);
\path (0,0) coordinate (u1);
\path (3,0) coordinate (u2);
%
\path (0,1) coordinate (v1);
\path (3,1) coordinate (v2);
%
\path (1.75,1) coordinate (n1);
\path (2,1) coordinate (n2);
\path (2.25,1) coordinate (n3);
\draw (u1)--(v1);
\draw (u2)--(v2);
\draw (v1)--(v2);
\draw (a)--(a1);
\draw (b1)--(a1);
\draw[densely dashed] (a1)--(b2);
\draw (v1) to[out=90,in=180, distance=0.75cm] (a);
\draw (v2) to[out=90,in=0, distance=0.75cm] (a);
\draw (v1) [fill=white] circle (\vr);
\draw (v2) [fill=white] circle (\vr);
%
\draw (u1) [fill=white] circle (\vr);
\draw (u2) [fill=white] circle (\vr);
%
\draw (a) [fill=white] circle (\vr);
\draw (a1) [fill=white] circle (\vr);
\draw[anchor = east] (u1) node {{\small $u_1$}};
\draw[anchor = west] (u2) node {{\small $u_2$}};
%
\draw[anchor = east] (v1) node {{\small $v_1$}};
\draw[anchor = west] (v2) node {{\small $v_2$}};
\draw[anchor = north] (a) node {{\small $a$}};
\draw[anchor = south] (a1p) node {{\small $a_1$}};
\draw (1.5,-1) node {{\small (a)}};
\path (7,2.25) coordinate (a);
\path (7,3.05) coordinate (a1);
\path (7.65,2.25) coordinate (a2);
\path (9,2.25) coordinate (b);
\path (9,3.05) coordinate (b1);
\path (8.35,2.25) coordinate (b2);
\path (7,0) coordinate (u1);
\path (9,0) coordinate (u2);
%
\path (7,1) coordinate (v1);
\path (9,1) coordinate (v2);
%
%
\draw (u1)--(v1);
\draw (u2)--(v2);
\draw (v1)--(v2);
%
\draw (a)--(a1);
\draw (v1)--(a);
\draw[densely dashed] (a)--(a2);
\draw (b)--(b1);
\draw (v2)--(b);
\draw[densely dashed] (b)--(b2);
%
%
\draw (v1) [fill=white] circle (\vr);
\draw (v2) [fill=white] circle (\vr);
%
\draw (u1) [fill=white] circle (\vr);
\draw (u2) [fill=white] circle (\vr);
%
\draw (a) [fill=white] circle (\vr);
\draw (b) [fill=white] circle (\vr);
%
\draw[anchor = east] (u1) node {{\small $u_1$}};
\draw[anchor = west] (u2) node {{\small $u_2$}};
%
\draw[anchor = east] (v1) node {{\small $v_1$}};
\draw[anchor = west] (v2) node {{\small $v_2$}};
\draw[anchor = east] (a) node {{\small $a$}};
\draw[anchor = west] (b) node {{\small $b$}};
%
\draw (8,-1) node {{\small (b)}};
\path (13,2.25) coordinate (a);
\path (13,3.05) coordinate (a1);
\path (13.65,2.25) coordinate (a2);
\path (15,2.25) coordinate (b);
\path (15,3.05) coordinate (b1);
\path (14.35,2.25) coordinate (b2);
\path (14,0) coordinate (u2);
%
\path (14,1) coordinate (v2);
%
%
\draw (u2)--(v2);
%
\draw (a)--(a1);
\draw (v2)--(a);
\draw[densely dashed] (a)--(a2);
\draw (b)--(b1);
\draw (v2)--(b);
\draw[densely dashed] (b)--(b2);
%
%
\draw (v2) [fill=white] circle (\vr);
%
\draw (u2) [fill=white] circle (\vr);
%
\draw (a) [fill=white] circle (\vr);
\draw (b) [fill=white] circle (\vr);
%
\draw[anchor = west] (u2) node {{\small $u_1$}};
%
\draw[anchor = west] (v2) node {{\small $v_1$}};
\draw[anchor = east] (a) node {{\small $a$}};
\draw[anchor = west] (b) node {{\small $b$}};
%
\draw (14,-1) node {{\small (c)}};
\end{tikzpicture}
\caption{Subgraphs in the proof of Claim~\ref{claim:min-deg2}}
\label{fig:deg2}
\end{center}
\end{figure}
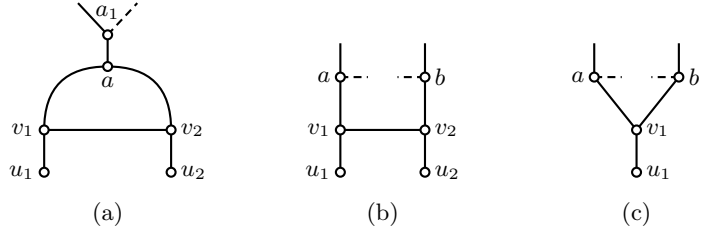

\begin{subclaim}
\label{claim:min-deg2.1}
$a \ne b$.
\end{subclaim}
\proof  Suppose, to the contrary, that $a = b$, implying that $k = 2$. Since $n \ge 6$ (and $G$ contains a copy of $K_{2,3}$), we note that $\deg_G(a) = 3$. Let $N_G(a) = \{a_1,v_1,v_2\}$. [See Figure~\ref{fig:deg2}(a).] By the maximality of the path $P$, the vertex $a_1$ has degree at least~$2$ in $G$. Let $X = \{a,u_1,u_2,v_1,v_2\}$ and let $G' = G - X$. By Claims~\ref{claim:no-bad-in-Ge} and~\ref{claim:no-T-in-Ge}, we have $\Theta(G') = 0$. Thus, $\Omega(G') = \w(G') = \w(G) - \w_G(X) + 1 = \Omega(G) - 19 + 1 = \Omega(G) - 18$. Every $i$-set of $G'$ can be extended to an ID-set of $G$ by adding to it the vertices $v_1$ and $u_2$, and so $i(G) \le i(G') + 2$. Thus, $8i(G) \le 8(i(G') + 2) \le \Omega(G') + 16 = \Omega(G) - 2 < \Omega(G)$, a contradiction.~\smallqed

\medskip
By Claim~\ref{claim:min-deg2.1}, $a \ne b$. [See Figure~\ref{fig:deg2}(b) and~\ref{fig:deg2}(c).] Relabeling vertices if necessary, we may assume that $\deg_G(a) \le \deg_G(b)$. By Claim~\ref{claim:no-strong-support}, we have $2 \le \deg_G(a) \le \deg_G(b)$. If $k = 1$, let $X = \{u_1,v_1\}$ and if $k = 2$, let $X = \{u_1,u_2,v_1,v_2\}$. In both cases, let $X' = X \cup \{a\}$.

\begin{subclaim}
\label{claim:min-deg2.2}
The graph $G - X'$ is isolate-free.
\end{subclaim}
\proof Suppose, to the contrary, that $G - X'$ has an isolated vertex, $z$ say. By Claim~\ref{claim:no-deg2-support} and the maximality of the path $P$, the vertex $z$ has degree at least~$2$ in $G$, implying that $\deg_G(z) = 2$ and $N_G(z) = \{a,v_k\}$. Hence, $z = b$. Thus, $2 \le \deg_G(a) \le \deg_G(b) = 2$, and so $\deg_G(a) = \deg_G(b)$ and $ab \in E(G)$. The graph $G$ is therefore determined. However, $G$ does not contain a subgraph isomorphic to $K_{2,3}$, contradicting Claim~\ref{claim:K23}.~\smallqed

\begin{subclaim}
\label{claim:min-deg2.3}
If $k = 2$, then $ab \in E(G)$.
\end{subclaim}
\proof Suppose, to the contrary, that $k = 2$ and the vertices $a$ and $b$ are not adjacent. We now consider the graph $G'$ obtained from $G - X$ by adding the edge $e = ab$; that is, $G' = (G - X) + e$. The graph $G'$ is a connected subcubic graph. Recall that $2 \le \deg_G(a) \le \deg_G(b)$. By construction, the degrees of vertices~$a$ and $b$ in $G'$ are the same as their degrees in~$G$. Let $I'$ be an $i$-set of $G'$. If $a \in I'$, then let $I = I' \cup \{u_1,v_2\}$. If $b \in I'$, then let $I = I' \cup \{v_1,u_2\}$. If neither $a$ nor $b$ belongs to $I'$, then let $I = I' \cup \{u_1,u_2\}$. In all cases, the set $I$ is an ID-set of $G$, and so $i(G) \le i(G') + 2$.

\begin{figure}[htb]
\begin{center}
\begin{tikzpicture}[scale=.75,style=thick,x=0.75cm,y=0.75cm]
\def\vr{2.25pt}
\def\vrn{1.25pt}
\path (0,1) coordinate (a1);
\path (1,1) coordinate (a2);
\path (2,0) coordinate (a3);
\path (2,1) coordinate (a4);
\path (2,2) coordinate (a5);
\path (3,1) coordinate (a6);
\path (4.75,0) coordinate (v2);
\path (4.75,1) coordinate (v1);
\path (6,0) coordinate (u2);
\path (6,1) coordinate (u1);
\draw (a1)--(a2)--(a3)--(a4)--(a5)--(a2);
\draw (a3)--(v2)--(v1)--(a6)--(a5);
\draw (v1)--(u1);
\draw (v2)--(u2);
\draw [style=dashed,rounded corners] (-1,-1) rectangle (3.5,2.75);
\draw [style=dashed,rounded corners] (4,-1) rectangle (7.25,2.75);
\draw[dotted, line width=0.05cm] (a3)--(a6);
\draw (a1) [fill=white] circle (\vr);
\draw (a2) [fill=black] circle (\vr);
\draw (a3) [fill=white] circle (\vr);
\draw (a4) [fill=black] circle (\vr);
\draw (a5) [fill=white] circle (\vr);
\draw (a6) [fill=white] circle (\vr);
\draw (v1) [fill=black] circle (\vr);
\draw (v2) [fill=white] circle (\vr);
\draw (u1) [fill=white] circle (\vr);
\draw (u2) [fill=black] circle (\vr);
\draw[anchor = south] (a6) node {{\small $a$}};
\draw[anchor = south] (v1) node {{\small $v_1$}};
\draw[anchor = south] (u1) node {{\small $u_1$}};
\draw[anchor = north] (a3) node {{\small $b$}};
\draw[anchor = north] (v2) node {{\small $v_2$}};
\draw[anchor = north] (u2) node {{\small $u_2$}};
\draw (2.85,0.5) node {{\small $e$}};
\draw (-0.5,2.25) node {{\small $G'$}};
\draw (6.75,2.25) node {{\small $X$}};
\end{tikzpicture}
\caption{The graph $G$ in the proof of Claim~\ref{claim:min-deg2.3}}
\label{fig:deg2.3.1a}
\end{center}
\end{figure}
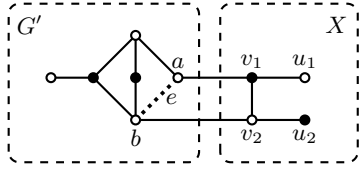

We show firstly that $b(G') = 0$. Suppose, to the contrary, that $b(G') \ge 1$. Since $G'$ is a connected graph, we infer that $G' \in \cB$ and $b(G') = 1$. If $G'$ contains two or more copies of $K_{2,3}$, then at least one of these copies contains a vertex of degree~$2$ in $G$, contradicting Claim~\ref{claim:no-K23-with-d2-vertex}. Hence, $G' \in \cB_1$ and $G'$ contains exactly one copy of $K_{2,3}$, say $F$. If $F$ contains at most one of~$a$ and~$b$, then again we contradict Claim~\ref{claim:no-K23-with-d2-vertex}. Hence, $F$ contains both vertices~$a$ and~$b$, and the root vertex of $G'$ has degree~$1$ in $G$ (and is distinct from~$a$ and~$b$). If $a$ or~$b$ is adjacent to the root vertex, then we would contradict the maximality of the path~$P$. Hence the graph $G$ is as illustrated in Figure~\ref{fig:deg2.3.1a}, where the added edge $e$ is indicated by the dotted line. Thus, $i(G) = 4$ (the shaded vertices in Figure~\ref{fig:deg2.3.1a} are an example of an $i$-set in $G$) and $\Omega(G) = \w(G) = 38$, and so $8i(G) < \Omega(G)$, a contradiction. Hence, $b(G') = 0$.

We show next that $\tc(G') = 0$. Suppose, to the contrary, that $\tc(G') \ge 1$. Thus, adding the edge $e = ab$ to the graph $G - X$ creates a new troublesome configuration, which we call $T_e$, that necessarily contains the edge~$e$. Let $r_1$ and $r_2$ be the two link vertices of $T_e$. We may assume that $\deg_G(r_1) \le \deg_G(r_2)$, and so $2 \le \deg_G(r_1) \le \deg_G(r_2) \le 3$. Thus, $T_e$ is obtained from a bad graph $B_e \in \cB$ with root vertex~$r_1$. By Claims~\ref{claim:no-K23-with-d2-vertex} and~\ref{claim:induced-K23}, the bad graph $B_e$ contains exactly one copy of $K_{2,3}$, say $F_e$, and so $B_e \in \cB_1$. Moreover, $F_e$ contains both vertices~$a$ and~$b$, for otherwise $F_e$ would contain a vertex of degree~$2$ in $G$, a contradiction.

\begin{figure}[htb]
\begin{center}
\begin{tikzpicture}[scale=.75,style=thick,x=0.75cm,y=0.75cm]
\def\vr{2.25pt}
\def\vrn{1.25pt}
\path (12,6.75) coordinate (h1);
\path (12,8) coordinate (h2);
\path (12,4.75) coordinate (p1);
\path (13.5,3.3) coordinate (p2);
\path (13.5,4.75) coordinate (p3);
\path (13.5,5.75) coordinate (p4);
\path (15,4.75) coordinate (p5);
\path (15,6.75) coordinate (q5);
\path (15,8) coordinate (q51);
%
\path (9.5,3.3) coordinate (v2);
\path (9.5,4.75) coordinate (v1);
\path (8,3.3) coordinate (u2);
\path (8,4.75) coordinate (u1);
\draw (h1)--(q5);
\draw (h1)--(h2);
\draw (p5)--(q5)--(q51);
\draw (p1)--(p2);
\draw (p3)--(p4)--(p5)--(p2);
\draw (p1)--(p4);
\draw (h1)--(p1);
%
%
\draw (v2)--(v1);
\draw (v1)--(u1);
\draw (u2)--(v2)--(p2);
\draw [dotted, line width=0.05cm] (p2)--(p3);
\draw  (v1) to[out=-45,in=225, distance=0.75cm] (p3);
\draw [style=dashed,rounded corners] (11,2.15) rectangle (17,7.45);
\draw [style=dashed,rounded corners] (6.75,2.5) rectangle (15.5,6.2);
\draw (h1) [fill=white] circle (\vr);
\draw (q5) [fill=white] circle (\vr);
%
\draw (p1) [fill=white] circle (\vr);
\draw (p2) [fill=white] circle (\vr);
\draw (p3) [fill=white] circle (\vr);
\draw (p4) [fill=black] circle (\vr);
\draw (p5) [fill=white] circle (\vr);
\draw (v1) [fill=white] circle (\vr);
\draw (v2) [fill=black] circle (\vr);
\draw (u1) [fill=black] circle (\vr);
\draw (u2) [fill=white] circle (\vr);
%
\draw[anchor = east] (h1) node {{\small $r_1$}};
\draw[anchor = west] (q5) node {{\small $r_2$}};
\draw[anchor = west] (p3) node {{\small $a$}};
\draw[anchor = west] (p2) node {{\small $b$}};
%
\draw[anchor = south] (v1) node {{\small $v_1$}};
\draw[anchor = south] (u1) node {{\small $u_1$}};
\draw[anchor = north] (v2) node {{\small $v_2$}};
\draw[anchor = north] (u2) node {{\small $u_2$}};
\draw (16.5,6.85) node {{\small $T_e$}};
\draw (7.5,4) node {{\small $X^*$}};
\draw (13.75,4.2) node {{\small $e$}};
\draw (12,1) node {{\small (a)}};
\path (24,6.75) coordinate (h1);
\path (24,8) coordinate (h2);
\path (24,4.75) coordinate (p1);
\path (25.5,3.3) coordinate (p2);
\path (25.5,4.75) coordinate (p3);
\path (25.5,5.75) coordinate (p4);
\path (27,4.75) coordinate (p5);
\path (27,6.75) coordinate (q5);
\path (27,8) coordinate (q51);
\path (21.5,3.3) coordinate (v2);
\path (21.5,4.75) coordinate (v1);
\path (20,3.3) coordinate (u2);
\path (20,4.75) coordinate (u1);
\draw (h1)--(q5);
\draw (h1)--(h2);
\draw (p5)--(q5)--(q51);
\draw (p3)--(p2);
\draw (p3)--(p4)--(p5)--(p2);
\draw (p1)--(p4);
\draw (h1)--(p1);
%
%
\draw (v2)--(v1);
\draw (v1)--(u1);
\draw (u2)--(v2)--(p2);
\draw [dotted, line width=0.05cm] (p1)--(p2);
%
\draw (v1)--(p1);
\draw [style=dashed,rounded corners] (23,2.15) rectangle (29,7.45);
\draw [style=dashed,rounded corners] (18.75,2.5) rectangle (27.5,6.2);
\draw (h1) [fill=white] circle (\vr);
\draw (q5) [fill=white] circle (\vr);
%
\draw (p1) [fill=white] circle (\vr);
\draw (p2) [fill=white] circle (\vr);
\draw (p3) [fill=white] circle (\vr);
\draw (p4) [fill=black] circle (\vr);
\draw (p5) [fill=white] circle (\vr);
\draw (v1) [fill=white] circle (\vr);
\draw (v2) [fill=black] circle (\vr);
\draw (u1) [fill=black] circle (\vr);
\draw (u2) [fill=white] circle (\vr);
%
\draw[anchor = east] (h1) node {{\small $r_1$}};
\draw[anchor = west] (q5) node {{\small $r_2$}};
\draw[anchor = north] (p1) node {{\small $a$}};
\draw[anchor = west] (p2) node {{\small $b$}};
%
\draw[anchor = south] (v1) node {{\small $v_1$}};
\draw[anchor = south] (u1) node {{\small $u_1$}};
\draw[anchor = north] (v2) node {{\small $v_2$}};
\draw[anchor = north] (u2) node {{\small $u_2$}};
\draw (28.5,6.85) node {{\small $T_e$}};
\draw (19.5,4) node {{\small $X^*$}};
\draw (25,4.25) node {{\small $e$}};
\draw (24,1) node {{\small (b)}};
\end{tikzpicture}
\caption{A subgraph of $G$ in the proof of Claim~\ref{claim:min-deg2.3}}
\label{fig:min-deg2.3.2-0}
\end{center}
\end{figure}
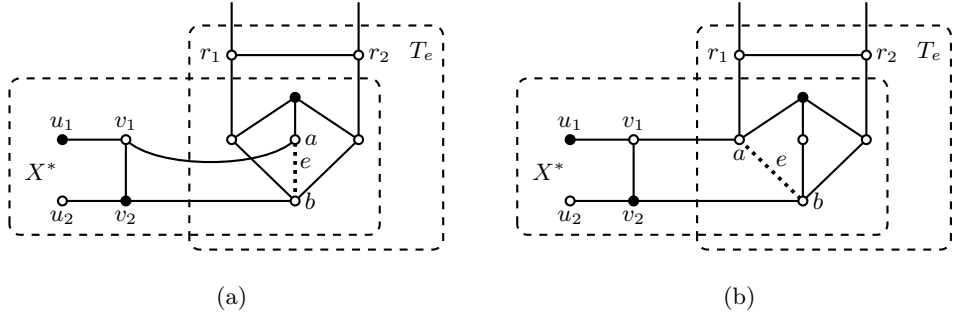

We now let $X^* = X \cup (V(T_e) \setminus \{r_1,r_2\})$ and we let $G^* = G - X^*$. If $\deg_G(a) = 2$, then the graph illustrated in Figure~\ref{fig:min-deg2.3.2-0}(a) is a subgraph of $G$, while if $\deg_G(a) = 3$, then the graph illustrated in Figure~\ref{fig:min-deg2.3.2-0}(b) is a subgraph of $G$ where in both cases $T_e$ is the troublesome configuration indicated in the dashed box, the set $X^*$ is indicated in the dashed box, and where $e = ab$ is indicated by the dotted edge. Every $i$-set of $G^*$ can be extended to an ID-set of $G$ by adding to it the vertices $u_1$ and $v_2$ and the neighbor of $a$ in $F_e$ different from~$b$ (as indicated by the shaded vertices in Figures~\ref{fig:min-deg2.3.2-0}(a) and~\ref{fig:min-deg2.3.2-0}(b)). Thus, $i(G) \le i(G^*) + 3$. We note that $\Theta(G^*) = 0$, and so $\Omega(G^*) = \w(G^*) = \w(G) - 32 + 2 = \Omega(G) - 30$. Hence, $8i(G) \le 8(i(G^*) + 3) \le \Omega(G') + 24 < \Omega(G)$, a contradiction. Hence, $\tc(G') = 0$.

We now return to the proof of Claim~\ref{claim:min-deg2.3}. By our earlier observations, $b(G') = 0$ and $\tc(G') = 0$. Hence, $\Theta(G') = 0$, and so $\Omega(G') = \w(G') = \w(G) - \w_G(X) = \Omega(G) - 16$. Recall that $i(G) \le i(G') + 2$. Thus, $8i(G) \le 8(i(G') + 2) \le \Omega(G') + 16 = \Omega(G)$, a contradiction.~\smallqed

\begin{subclaim}
\label{claim:min-deg2.4}
$k = 1$.
\end{subclaim}
\proof Suppose, to the contrary, that $k = 2$. By Claim~\ref{claim:min-deg2.3}, $ab \in E(G)$. Recall that $2 \le \deg_G(a) \le \deg_G(b)$. If $\deg_G(a) = \deg_G(b) = 2$, then the graph $G$ is determined. In this case, $G$ does not contain a copy of $K_{2,3}$, a contradiction. Hence, $\deg_G(b) = 3$.

We show firstly that $\deg_G(a) = 3$. Suppose, to the contrary, that $\deg_G(a) = 2$. In this case, we let $G' = G - X'$ where recall that $X = \{u_1,u_2,v_1,v_2\}$ and $X' = X \cup \{a\}$. Every $i$-set of $G'$ can be extended to an ID-set of $G$ by adding to it the vertices $v_1$ and $u_2$, and so $i(G) \le i(G') + 2$. We note that $G'$ is a connected subcubic graph and $\deg_{G'}(b) = 1$. Moreover, $\Theta(G') = 0$, and so $\Omega(G') = \w(G') = \w(G) - \w_G(X) + 2 = \Omega(G) - 20 + 2 = \Omega(G) - 18$. Thus, $8i(G) \le 8(i(G') + 2) \le \Omega(G') + 16 < \Omega(G)$, a contradiction. Hence, $\deg_G(a) = 3$. Let $N_G(a) = \{a_1,b,v_1\}$ and $N_G(b) = \{a,b_1,v_2\}$.

We show next that $a_1 \ne b_1$. Suppose, to the contrary, that $a_1 = b_1$. If $\deg_G(a_1) = 2$, then the graph $G$ is determined. In this case, $G$ does not contain a copy of $K_{2,3}$, a contradiction. Hence, $\deg_G(a_1) = 3$. Let $a_2$ be the neighbor of $a_1$ different from $a$ and $b$. If $\deg_G(a_2) = 1$, then the graph $G$ is determined, and again we contradict the supposition that $G$ contain a copy of $K_{2,3}$. Hence, $\deg_G(a_2) \ge 2$. We now let $X^* = X \cup \{a,b,a_1\}$ and consider the connected subcubic graph $G^* = G - X^*$. Every $i$-set of $G^*$ can be extended to an ID-set of $G$ by adding to it the set $\{b,v_1,u_2\}$, and so $i(G) \le i(G^*) + 3$. By Claim~\ref{claim:no-bad-in-Ge}, $b(G^*) = 0$. By Claim~\ref{claim:no-T-in-Ge}, $\tc(G^*) = 0$. Thus, $\Theta(G^*) = 0$, and so noting that $a_1a_2$ is the unique $X^*$-exit edge and $\deg_G(a_2) \ge 2$, we have $\Omega(G^*) = \w(G^*) = \w(G) - \w_G(X) + 1 = \Omega(G) - 25 + 1 = \Omega(G) - 24$. Thus, $8i(G) \le 8(i(G') + 3) \le \Omega(G') + 24 = \Omega(G)$, a contradiction. Hence, $a_1 \ne b_1$.

Recall by maximality of the path~$P$ that $\deg_G(a_1) \ge 2$ and  $\deg_G(b_1) \ge 2$. We now consider the graph $G' = G - X'$, where recall that $X = \{u_1,u_2,v_1,v_2\}$ and $X' = X \cup \{a\}$. Every $i$-set of $G'$ can be extended to an ID-set of $G$ by adding to it the vertices $v_1$ and $u_2$, and so $i(G) \le i(G') + 2$. We note that $\deg_{G'}(b) = 1$, and so there is no troublesome configuration containing~$b$ in~$G'$. If the vertex~$b$ belongs to a bad component, then such a component has a copy of $K_{2,3}$ that contains a vertex of degree~$2$ in $G$, contradicting Claim~\ref{claim:no-K23-with-d2-vertex}. Hence, $b$ belongs to no bad component in $G'$. If the removal of the edge $aa_1$ when constructing $G'$ creates a new troublesome configuration, then by Claim~\ref{claim:no-T-in-Ge} we infer that such a troublesome configuration must contain the three vertices~$a_1$,~$b$ and~$b_1$, which is not possible noting that $\deg_{G'}(b) = 1$. Hence, $\tc(G') = 0$. Thus, $\Theta(G') = 0$, and so noting that there are three $X'$-exit edges and $\deg_G(a_1) \ge 2$, we have $\Omega(G') = \w(G') = \w(G) - \w_G(X') + 3 = \Omega(G) - 19 + 3 = \Omega(G) - 16$. Thus, $8i(G) \le 8(i(G') + 2) \le \Omega(G') + 16 = \Omega(G)$, a contradiction.~\smallqed

\medskip
By Claim~\ref{claim:min-deg2.4}, we have $k = 1$. Recall that $2 \le \deg_G(a) \le \deg_G(b)$. Recall that $n \ge 6$. In particular, we note that $\deg_G(b) = 3$.

\begin{subclaim}
\label{claim:min-deg2.5a}
$ab \notin E(G)$.
\end{subclaim}
\proof Suppose, to the contrary, that $ab \in E(G)$. In this case, we let $X^* = \{a,b,u_1,v_1\}$ and consider the graph $G^* = G -  X^*$. Every $i$-set of $G^*$ can be extended to an ID-set of $G$ by adding to it the vertex~$v_1$, and so $i(G) \le i(G^*) + 1$. Suppose that $\deg_G(a) = 2$, and so there is exactly one $X^*$-exit edge. By Claim~\ref{claim:no-bad-in-Ge}, $b(G^*) = 0$. By Claim~\ref{claim:no-T-in-Ge}, $\tc(G^*) = 0$. Thus, $\Theta(G^*) = 0$, and so $\Omega(G^*) = \w(G^*) = \w(G) - \w_G(X) + 1 = \Omega(G) - 15 + 1 = \Omega(G) - 14$. Thus, $8i(G) \le 8(i(G^*) + 1) \le \Omega(G^*) + 8 = \Omega(G)$, a contradiction. Hence, $\deg_G(a) = 3$. Let $N_G(a) = \{a_1,b,v_1\}$ and $N_G(b) = \{a,b_1,v_1\}$. By the maximality of the path~$P$, we note that $\deg_G(a_1) \ge 2$ and $\deg_G(b_1) \ge 2$. If $a_1$ is isolated in $G^*$, then the graph $G$ is determined. In this case, $a_1 = b_1$ and $N_G(a_1) = \{a,b\}$. However, then, $n = 5$, a contradiction. Hence, $a_1$ is not isolated in~$G^*$. By symmetry, $b_1$ is not isolated in~$G^*$. Thus, $\w(G^*) = \w(G) - \w_G(X^*) + 2 = \w(G) - 14 + 2 = \w(G) - 12$. By Claim~\ref{claim:no-bad-in-Ge}, we have $b(G^*) = 0$, and by Claim~\ref{claim:no-T-in-Gee}, we have $\tc(G^*) \le 1$. Thus, $\Theta(G^*) \le 2$, and so $\Omega(G^*) = \w(G^*) + \Theta(G^*) \le (\w(G) - 12) + 2 = \Omega(G) - 12 + 2 = \Omega(G) - 10$. Thus, $8i(G) \le 8(i(G^*) + 1) \le \Omega(G^*) + 8 < \Omega(G)$, a contradiction.~\smallqed

\medskip
By Claim~\ref{claim:min-deg2.5a}, $ab \notin E(G)$.

\begin{subclaim}
\label{claim:min-deg2.2}
The graph $G - \{a,b,u_1,v_1\}$ is isolate-free.
\end{subclaim}
\proof Suppose, to the contrary, that $G - \{a,b,u_1,v_1\}$ has an isolated vertex. By Claim~\ref{claim:no-deg2-support} and the maximality of the path $P$, such a vertex has degree~$2$ in $G$ with $a$ and $b$ as its two neighbors. If $G - \{a,b,u_1,v_1\}$ contains two isolated vertices, then the graph $G$ is determined and is the graph $B_1$ shown in Figure~\ref{fig:B1}. However, then $G \in \cB$, a contradiction. Hence, $G - \{a,b,u_1,v_1\}$ contains exactly one isolated vertex, say $c$. As observed earlier, $N_G(c) = \{a,b\}$. Since $n \ge 6$, we note that $\deg_G(b) = 3$. Let $N_G(b) = \{b_1,c,v_1\}$. By the maximality of the path $P$, we note that $\deg_G(b_1) \ge 2$. We now let $X'' = \{a,b,c,u_1,v_1\}$.

Suppose that $\deg_G(a) = 2$. In this case, we consider the graph $G'' = G - X''$. Every $i$-set of $G''$ can be extended to an ID-set of $G$ by adding to it the vertices $v_1$ and $c$, and so $i(G) \le i(G'') + 2$. By Claim~\ref{claim:no-bad-in-Ge}, $b(G'') = 0$. By Claim~\ref{claim:no-T-in-Ge}, $\tc(G'') = 0$. Thus, $\Theta(G'') = 0$, and so $\Omega(G'') = \w(G'') = \w(G) - \w_G(X'') + 1 = \Omega(G) - 19 + 1 = \Omega(G) - 18$. Thus, $8i(G) \le 8(i(G^*) + 2) \le \Omega(G^*) + 16 < \Omega(G)$, a contradiction. Hence, $\deg_G(a) = 3$. Let $N_G(a) = \{a_1,c,v_1\}$. By the maximality of the path $P$, we note that $\deg_G(a_1) \ge 2$. If $a_1 = b_1$, then $G[\{a,a_1,b,c,v_1\}]$ is a copy of $K_{2,3}$ that contains a vertex of degree~$2$ in $G$, namely vertex~$c$, contradicting Claim~\ref{claim:no-K23-with-d2-vertex}. Hence, $a_1 \ne b_1$. 
Renaming vertices if necessary, we may assume by symmetry that $2\le\deg_G(a_1) \le \deg_G(b_1)$.

We show next that $a_1b_1 \in E(G)$. Suppose, to the contrary, that $a_1b_1 \notin E(G)$. Recall that $X'' = \{a,b,c,u_1,v_1\}$. In this case, we consider the graph $G'$ obtained from $G - X''$ by adding the edge $e = a_1b_1$; that is, $G' = (G - X'') + e$. The graph $G'$ is a connected subcubic graph. Let $I''$ be an $i$-set of $G'$. If $a_1 \in I''$, then let $I = I'' \cup \{b,u_1\}$. If $b_1 \in I''$, then let $I = I'' \cup \{a,u_1\}$. If neither $a_1$ nor $b_1$ belongs to $I''$, then let $I = I'' \cup \{c,v_1\}$. In all cases, the set $I$ is an ID-set of $G$, and so $i(G) \le i(G') + 2$. 
By construction, the degrees of vertices~$a_1$ and $b_1$ in $G'$ are the same as their degrees in~$G$. However adding the edge~$e$ may have created a bad graph or a troublesome configuration. Thus, $b(G') + \tc(G') \le 1$, implying that $\Theta(G') \le 2$, and so $\Omega(G')  =\w(G') + \Theta(G') \le (\w(G) - \w_G(X'')) + 2 = \w(G) - 18 + 2 = \Omega(G) - 16$. Thus, $8i(G) \le 8(i(G') + 2) \le \Omega(G') + 16 = \Omega(G)$, a contradiction. Hence, $a_1b_1 \in E(G)$.

In this case, we consider the graph $G'' = G - X''$. Every $i$-set of $G''$ can be extended to an ID-set of $G$ by adding to it the vertices $v_1$ and $c$, and so $i(G) \le i(G'') + 2$. By Claim~\ref{claim:no-bad-in-Gee}, $b(G'') = 0$. We note that $a_1$ and $b_1$ are adjacent vertices in $G''$ and both $a_1$ and $b_1$ have degree at most~$2$ in $G''$. We therefore infer that $a_1$ and $b_1$ do not belong to a troublesome configuration in $G''$, and so $\tc(G'') = 0$. Thus, $\Theta(G'') = 0$, and so $\Omega(G'') = \w(G'') = \w(G) - \w_G(X'') + 2 = \Omega(G) - 18 + 2 = \Omega(G) - 16$. Thus, $8i(G) \le 8(i(G'') + 2) \le \Omega(G'') + 16 = \Omega(G)$, a contradiction.~\smallqed

\medskip
By Claim~\ref{claim:min-deg2.2}, the graph $G^* = G - X^*$ is isolate-free, where recall that $X^* = \{a,b,u_1,v_1\}$. Further recall that $2 \le \deg_G(a) \le \deg_G(b)$ and that $ab \notin E(G)$. Every $i$-set of $G^*$ can be extended to an ID-set of $G$ by adding to it the vertex $v_1$, and so $i(G) \le i(G^*) + 1$. By Claim~\ref{claim:no-bad-in-Ge}, if there is a bad component in $G^*$, then such a component is incident with at least three $X^*$-edges. By Claims~\ref{claim:no-T-in-Ge} and~\ref{claim:no-T-in-Gee}, if there is a troublesome configuration in $G^*$, then such a component is incident with at least two $X^*$-edges. Thus since there are at most four $X^*$-edges, namely two edges incident with vertex~$a$ and two incident with vertex~$b$, we infer that either $b(G^*) = \tc(G^*) = 0$ or $b(G^*) = 1$ and $\tc(G^*) = 0$ or $b(G^*) = 0$ and $\tc(G^*) \le 2$. In all cases, $b(G^*) + \tc(G^*) \le 2$.

Suppose that $b(G^*) + \tc(G^*) \le 1$, and so $\Theta(G^*) \le 2$ and $\Omega(G^*) \le \w(G^*) + 2$. If $\deg_G(a) = 2$, then $\w(G^*) = \w(G) - \w_G(X^*) + 3 = \w(G) - 15 + 3 = \Omega(G) - 12$, and so $\Omega(G^*) \le \w(G^*)  + 2 \le \Omega(G) - 10$. If $\deg_G(a) = 3$, then $\w(G^*) = \w(G) - \w_G(X^*) + 4 = \w(G) - 14 + 4 = \Omega(G) - 10$, and so $\Omega(G^*) \le \w(G^*)  + 2 \le \Omega(G) - 8$. In both cases, $\Omega(G^*) + 8 \le \Omega(G)$. Thus, $8i(G) \le 8(i(G^*) + 1) \le \Omega(G^*) + 8 \le \Omega(G)$, a contradiction.

Hence, $b(G^*) + \tc(G^*) = 2$, implying by our earlier observations that $b(G^*) = 0$ and $\tc(G^*) = 2$. Thus, $G^*$ contains two troublesome configurations, say $T_w$ and~$T_z$, each of which contain two neighbors of~$a$ or~$b$. In particular, we note that $\deg_G(a) = \deg_G(b) = 3$. Let $N_G(a) = \{a_1,a_2,v_1\}$ and let $N_G(b) = \{b_1,b_2,v_1\}$. Let $w_1$ and $w_2$ be the link vertices of $T_w$ and let $z_1$ and $z_2$ be the link vertices of $T_z$. Let $T_w$ be obtained from the bad graph $B_w \in \cB_1$ with root vertex~$w_1$, and let $T_z$ be obtained from the bad graph $B_z \in \cB_1$ with root vertex~$z_1$. Moreover, let $x_1$ and $x_2$ be the vertices of $B_w$ adjacent to $w_1$ and $w_2$, respectively, and let $y_1$ and $y_2$ be the vertices of $B_z$ adjacent to $z_1$ and $z_2$, respectively. By Claims~\ref{claim:no-K23-with-d2-vertex} and~\ref{claim:induced-K23}, each of $B_w$ and $B_z$ contains exactly two copies of $K_{2,3}$. Further, each copy of $K_{2,3}$ in $B_w$ and $B_z$ contains a vertex incident with an $X^*$-exit edge.

We now let $X'' = X^* \cup (V(T_w) \setminus \{x_2,w_2\}) \cup (V(T_z) \setminus \{y_2,z_2\})$ and let $G'' = G - X''$. In the case when $\{a_1,b_1\} \subset V(T_w)$ and $\{a_2,b_2\} \subset V(T_z)$, the graph illustrated in Figure~\ref{fig:claim:min-deg2.2E} is a subgraph of $G$, where $T_w$ and $T_z$ are the troublesome configurations indicated in the dashed boxes and where $X^*$ and $X''$ are the sets indicated by the vertices in the dotted region.

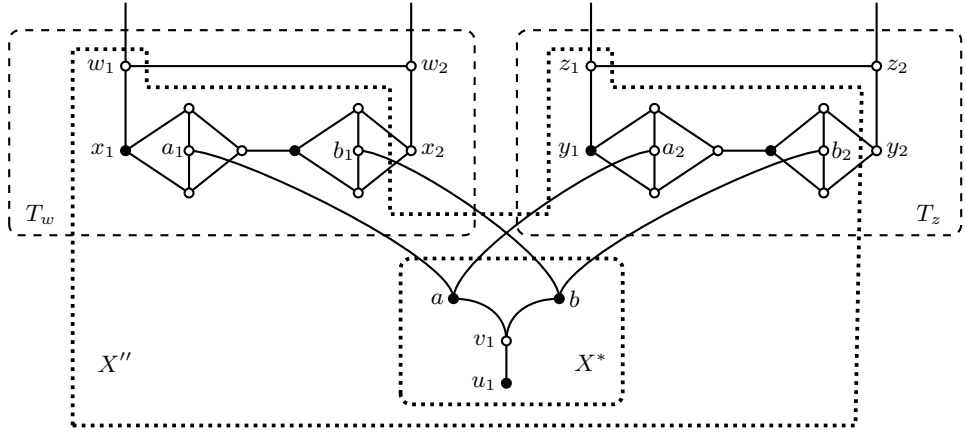
\begin{figure}[htb]
\begin{center}
\begin{tikzpicture}[scale=.75,style=thick,x=0.75cm,y=0.75cm]
\def\vr{2.25pt}
\def\vrn{1.25pt}
\path (1,3.5) coordinate (z0);
\path (19.5,3.5) coordinate (z1);
\path (19.65,11.5) coordinate (z2);
\path (13.75,11.5) coordinate (z3);
\path (13.75,12.4) coordinate (z4);
\path (12.25,12.4) coordinate (z5);
\path (12.25,8.5) coordinate (z6);
\path (8.5,8.5) coordinate (z7);
\path (8.5,11.5) coordinate (z8);
\path (2.75,11.5) coordinate (z9);
\path (2.75,12.4) coordinate (z10);
\path (1,12.4) coordinate (z11);
\path (10,6.5) coordinate (h11);
\path (12.5,6.5) coordinate (q51);

\path (11.25,5.5) coordinate (v1);
\path (11.25,4.5) coordinate (u1);
\path (2.25,10) coordinate (a1);
\path (3.75,9) coordinate (a2);
\path (3.75,10) coordinate (a3);
\path (3.9,10) coordinate (a3p);
\path (3.75,11) coordinate (a4);
\path (5,10) coordinate (a5);
\path (6.25,10) coordinate (d1);
\path (7.75,9) coordinate (d2);
\path (7.75,10) coordinate (d3);
\path (7.9,10) coordinate (d3p);
\path (7.75,11) coordinate (d4);
\path (9,10) coordinate (d5);
\path (2.25,12) coordinate (j1);
\path (2.25,13.5) coordinate (j11);
\path (9,12) coordinate (j2);
\path (9,13.5) coordinate (j21);
\path (13.25,10) coordinate (k1);
\path (14.75,9) coordinate (k2);
\path (14.75,10) coordinate (k3);
\path (14.7,10) coordinate (k3p);
\path (14.75,11) coordinate (k4);
\path (16.25,10) coordinate (k5);
\path (17.5,10) coordinate (b1);
\path (18.75,9) coordinate (b2);
\path (18.75,10) coordinate (b3);
\path (18.7,10) coordinate (b3p);
\path (18.75,11) coordinate (b4);
\path (20,10) coordinate (b5);
\path (13.25,12) coordinate (l1);
\path (13.25,13.5) coordinate (l11);
\path (20,12) coordinate (l2);
\path (20,13.5) coordinate (l21);
%
%
\draw[dotted, line width=0.05cm] (z0)--(z1)--(z2)--(z3)--(z4)--(z5)--(z6)--(z7)--(z8)--(z9)--(z10)--(z11)--(z0);
\draw (d1)--(d2)--(d3)--(d4)--(d5)--(d2);
\draw (d1)--(d4);
\draw (a1)--(a2)--(a3)--(a4)--(a5)--(a2);
\draw (a1)--(a4);
\draw (a5)--(d1);
\draw (a1)--(j1)--(j11);
\draw (d5)--(j2)--(j1);
\draw (j2)--(j21);
\draw (u1)--(v1);
\draw (h11) to[out=90,in=0, distance=0.75cm] (a3);
\draw (h11) to[out=90,in=180, distance=0.75cm] (k3);
\draw (v1) to[out=90,in=0, distance=0.5cm] (h11);
\draw (v1) to[out=90,in=180, distance=0.5cm] (q51);
\draw (k1)--(k2)--(k3)--(k4)--(k5)--(k2);
\draw (k1)--(k4);
\draw (b1)--(b2)--(b3)--(b4)--(b5)--(b2);
\draw (b1)--(b4);
\draw (b1)--(k5);
\draw (k1)--(l1)--(l11);
\draw (b5)--(l2)--(l1);
\draw (l2)--(l21);
\draw (q51) to[out=90,in=0, distance=0.75cm] (d3);
\draw (q51) to[out=90,in=180, distance=0.75cm] (b3);
\draw[dotted, line width=0.05cm] [rounded corners] (8.75,4) rectangle (14,7.45);
\draw [style=dashed,rounded corners] (-0.5,8) rectangle (10.5,12.85);
\draw [style=dashed,rounded corners] (11.5,8) rectangle (22,12.85);
\draw (0.25,8.5) node {{\small $T_w$}};
\draw (21.25,8.5) node {{\small $T_z$}};
\draw (a1) [fill=black] circle (\vr);
\draw (a2) [fill=white] circle (\vr);
\draw (a3) [fill=white] circle (\vr);
\draw (a4) [fill=white] circle (\vr);
\draw (a5) [fill=white] circle (\vr);
\draw (b1) [fill=black] circle (\vr);
\draw (b2) [fill=white] circle (\vr);
\draw (b3) [fill=white] circle (\vr);
\draw (b4) [fill=white] circle (\vr);
\draw (b5) [fill=white] circle (\vr);
\draw (d1) [fill=black] circle (\vr);
\draw (d2) [fill=white] circle (\vr);
\draw (d3) [fill=white] circle (\vr);
\draw (d4) [fill=white] circle (\vr);
\draw (d5) [fill=white] circle (\vr);
\draw (j1) [fill=white] circle (\vr);
\draw (j2) [fill=white] circle (\vr);
\draw (k1) [fill=black] circle (\vr);
\draw (k2) [fill=white] circle (\vr);
\draw (k3) [fill=white] circle (\vr);
\draw (k4) [fill=white] circle (\vr);
\draw (k5) [fill=white] circle (\vr);
\draw (l1) [fill=white] circle (\vr);
\draw (l2) [fill=white] circle (\vr);
\draw (q51) [fill=black] circle (\vr);
\draw (h11) [fill=black] circle (\vr);
\draw (v1) [fill=white] circle (\vr);
\draw (u1) [fill=black] circle (\vr);
\draw[anchor = east] (h11) node {{\small $a$}};
\draw[anchor = west] (q51) node {{\small $b$}};
\draw[anchor = east] (d3p) node {{\small $b_1$}};
\draw[anchor = west] (k3p) node {{\small $a_2$}};
\draw[anchor = east] (a3p) node {{\small $a_1$}};
\draw[anchor = west] (b3p) node {{\small $b_2$}};
\draw[anchor = east] (v1) node {{\small $v_1$}};
\draw[anchor = east] (u1) node {{\small $u_1$}};
\draw[anchor = east] (a1) node {{\small $x_1$}};
\draw[anchor = west] (d5) node {{\small $x_2$}};
\draw[anchor = east] (k1) node {{\small $y_1$}};
\draw[anchor = west] (b5) node {{\small $y_2$}};
\draw[anchor = east] (j1) node {{\small $w_1$}};
\draw[anchor = west] (j2) node {{\small $w_2$}};
\draw[anchor = east] (l1) node {{\small $z_1$}};
\draw[anchor = west] (l2) node {{\small $z_2$}};
\draw (13.25,5) node {{\small $X^*$}};
\draw (2,5) node {{\small $X''$}};
\end{tikzpicture}
\caption{An illustration of a subgraph of $G$ in the proof of Claim~\ref{claim:min-deg2}}
\label{fig:claim:min-deg2.2E}
\end{center}
\end{figure}

If $\deg_G(w_1) = 3$, then let $e_w = ww_1$ be the edge incident with $w_1$ not in $T_w$, and if $\deg_G(z_1) = 3$, then let $e_z = zz_1$ be the edge incident with $z_1$ not in $T_z$. The vertices~$x_2$ and $y_2$ have degree~$1$ in $G''$, and the neighbors of $x_2$ and $y_2$, namely $w_2$ and $z_2$, in $G''$ have degree~$2$. Hence, there is no bad component or troublesome configuration in $G''$ containing these four vertices. By Claim~\ref{claim:no-bad-in-Ge},~\ref{claim:no-bad-in-Gee},~\ref{claim:no-T-in-Ge} and~\ref{claim:no-T-in-Gee}, we have $b(G'') = 0$ and $\tc(G'') \le 1$. Moreover if $\tc(G'') = 1$, then the edges $e_w$ and $e_z$ exist and their removal increases the total weight by~$2 + 2 = 4$ (namely, an increase in~$2$ from the vertex weight and an increase in~$2$ in the structural weight arising from a troublesome configuration containing~$w$ and~$z$). If $\tc(G'') = 0$, then the removal of $e_w$ and $e_z$, if these edges exist, increases the total weight by at most~$2 \times 3 = 6$ in the worst case (which can only occur if $w$ and $z$ are isolated vertices in $G''$). As observed earlier, the removal of the six $X''$-exit edges different from~$e_w$ and~$e_z$ do not increase the structural weight and therefore only increase the vertex weight (by~$6$). We therefore infer that $\Omega(G'') \le \Omega(G) - \w_G(X'') + 6 + 6  \le \Omega(G) - 74 + 6 + 6 = \Omega(G) - 62$.

Every $i$-set of $G''$ can be extended to an ID-set of $G''$ by adding to it the vertices~$a$,~$b$ and~$u_1$, together with two vertices from each of~$T_w$ and~$T_z$ (as illustrated by the shaded vertices in Figure~\ref{fig:claim:min-deg2.2E}). Thus, $i(G) \le i(G'') + 7$. As observed earlier, $\Omega(G'') + 62 \le \Omega(G)$. Therefore,  $8i(G) \le 8(i(G'') + 7) \le \Omega(G'') + 56 < \Omega(G)$, a contradiction. We therefore deduce that our supposition that $G$ contains a vertex of degree~$1$ is false. Hence, $\delta(G) \ge 2$, completing the proof of Claim~\ref{claim:min-deg2}.~\smallqed

\subsection{Part 3: The minimum degree of $G$ is at least two}
\label{S:part3}

By Claim~\ref{claim:min-deg2}, $\delta(G) \ge 2$. Recall that $n \ge 6$ and that $G$ contains $K_{2,3}$ as a subgraph. By supposition, $G \ne K_{3,3}$. As in the proof of Part~2, in what follows we will frequently use the structural property stated in Claim~\ref{claim:no-K23-with-d2-vertex} that there is no $K_{2,3}$-subgraph in $G$ that contains a vertex of degree~$2$ in $G$.

We prove next a key structural property of the graph $G$. For this purpose, we define a \emph{troublesome subgraph} of $G$ as a subgraph $T$ such that there exists a set of edges incident with vertices of $T$ whose removal creates a troublesome configuration in $G$ with vertex set $V(T)$. For example, the subgraph $T$ illustrated in Figure~\ref{fig:tsG-0} is a troublesome subgraph of $G$ since the removal of the edge~$e$ produces a troublesome configuration in $G - e$. As a further example, the subgraph $T$ illustrated in Figure~\ref{fig:tsG-3} or in Figure~\ref{fig:tsG-5} is a troublesome subgraph of $G$ since the removal of the edge~$e$ (in both figures) produces a troublesome configuration in $G - e$.

\begin{claim}
\label{claim:no-troublesome-subgraph}
There is no troublesome subgraph in $G$.
\end{claim}
\proof Suppose, to the contrary, that $G$ contains a troublesome subgraph, say $T$. Hence there exists a set, say $E_T \subset E(G)$, of edges incident with vertices of $T$ such that $G' = G - E_T$ contains a troublesome configuration, say $T'$, with vertex set $V(T') = V(T)$. By Claim~\ref{claim:no-T-subgraph}, the graph $G$ contains no troublesome configuration, that is, $\tc(G) = 0$. By Claim~\ref{claim:no-T-in-Ge}, $\tc(G - e) = 0$ for all $e \in E(G)$. Hence, $|E_T| \ge 2$, that is, at least two edges incident with vertices in $T$ were removed when constructing~$T'$.

Let the troublesome configuration $T'$ in $G'$ have link vertices $v_1$ and $v_2$, where $\deg_G(v_2) = 3$. Thus, $T'$ is obtained from the bad graph $B \in \cB$ with root vertex~$v_1$ (of degree~$1$ in $B$). Further, let $u_1$ and $u_2$ be the vertices of $B$ adjacent to $v_1$ and $v_2$, respectively, in $G$. Let $z_2$ be the neighbor of $v_2$ different from~$u_2$ and $v_1$. Let the bad graph $B$ have $k$ copies of $K_{2,3}$ in $G'$. Since $|E_T| \ge 2$, we note that $k \ge 2$. Let $A = \{a_1,a_2,\ldots,a_k\}$ be the set of vertices in $B$ incident with edges in $E_T$.

\begin{subclaim}
\label{claim:no-troublesome-subgraph.1}
If the link vertex~$v_1$ is not incident with an edge of $E_T$, then $8i(G) \le \Omega(G)$.
\end{subclaim}
\proof Suppose that the link vertex~$v_1$ is not incident with an edge of $E_T$. If $v_1$ has degree~$3$ in $G$, then let $z_1$ be the neighbor of $v_1$ not in $T$. We now let $X = V(T) \setminus (A \cup \{u_2,v_2\})$ and we consider the graph $G' = G - X$. We note that each vertex in $A$ has degree~$1$ in $G'$, and belongs to neither a bad component nor a troublesome configuration in $G'$. Moreover, the vertex~$u_2$ has degree~$1$ in $G'$ and its neighbor, $v_2$, in $G'$ has degree~$2$ in $G'$. Hence, $u_2$ and $v_2$ do not belong to a bad component or a troublesome configuration in $G'$. If $\deg_G(v_1) = 3$, then we infer by Claims~\ref{claim:no-bad-in-Ge} and~\ref{claim:no-T-in-Ge} that the vertex $z_1$ does not belong to a bad component or a troublesome configuration in $G'$. Hence, $b(G') = \tc(G') = 0$, and so $\Theta(G') = 0$. We note that there are at most~$2k+4$ $X$-exit edges. Moreover, $|X| = 4k$ and every vertex of $X$, except for possibly the vertex~$v_1$, has degree~$3$ in $G$. Therefore, $\Omega(G') = \w(G') = \w(G) - \w_G(X) + 2k+4 \le \w(G) - 3 \times 4k + 2k+4 = \w(G) - 10k + 4 = \Omega(G) - 10k + 4$.

Every $i$-set of $G'$ can be extended to an ID-set of $G$ by adding to it the core independent set (see Section~\ref{S:familyB}) of the bad graph $B$, and so $i(G) \le i(G') + k$. (For example, if the graph shown in Figure~\ref{fig:claim:min-deg2.2F}(a) is a subgraph of $G$, where the dotted edges are the edges that belong to the set $E_T$, then $k = 6$ and the subgraph $G[X]$ is illustrated in Figure~\ref{fig:claim:min-deg2.2F}(b). Furthermore, the core independent set of $B$ is indicated by the $k = 6$ shaded vertices in Figure~\ref{fig:claim:min-deg2.2F}(b).) Thus noting that $k \ge 2$, we have $8i(G) \le 8(i(G') + k ) \le \Omega(G') + 8k \le (\Omega(G) - 10k + 4) + 8k = \Omega(G) - 2k + 4 \le \Omega(G)$, a contradiction.~\smallqed

\begin{figure}[htb]
\begin{center}
\begin{tikzpicture}[scale=.75,style=thick,x=0.75cm,y=0.75cm]
\def\vr{2.25pt}
\def\vrn{1.25pt}

\path (2.25,10) coordinate (a1);
\path (3.75,9) coordinate (a2);
\path (3.75,10) coordinate (a3);
\path (3.9,10) coordinate (a3p);
\path (3.75,11) coordinate (a4);
\path (5,10) coordinate (a5);
\path (6.5,10) coordinate (d1);
\path (7.75,9) coordinate (d2);
\path (7.75,10) coordinate (d3);
\path (7.9,10) coordinate (d3p);
\path (7.75,11) coordinate (d4);
\path (9,10) coordinate (d5);
\path (9,9) coordinate (d33);
\path (6.5,7) coordinate (f1);
\path (7.75,6) coordinate (f2);
\path (7.75,7) coordinate (f3);
\path (7.85,7) coordinate (f3p);
\path (7.75,8) coordinate (f4);
\path (9,7) coordinate (f5);
\path (9,5) coordinate (f33);
\path (10.75,7) coordinate (g1);
\path (12.25,6) coordinate (g2);
\path (12.25,7) coordinate (g3);
\path (12.35,7) coordinate (g3p);
\path (12.25,8) coordinate (g4);
\path (13.75,7) coordinate (g5);
\path (13.7,7) coordinate (g5p);
\path (13.75,5) coordinate (g33);
\path (15,5) coordinate (g55);
\path (2.25,12.5) coordinate (j1);
\path (2.25,14) coordinate (j11);
\path (1.75,14.75) coordinate (j111);
\path (2.75,14.75) coordinate (j112);
\path (9,12) coordinate (j2);
\path (9,13.5) coordinate (j21);
\path (10.75,10) coordinate (k1);
\path (12.25,9) coordinate (k2);
\path (12.25,10) coordinate (k3);
\path (12.2,10) coordinate (k3p);
\path (12.25,11) coordinate (k4);
\path (13.75,10) coordinate (k5);
\path (13.75,8) coordinate (k33);
\path (15,10) coordinate (b1);
\path (16.25,9) coordinate (b2);
\path (16.25,10) coordinate (b3);
\path (16.35,10) coordinate (b3p);
\path (16.25,11) coordinate (b4);
\path (17.5,10) coordinate (b5);
\path (17.25,8) coordinate (b33);
%
\path (17.5,12.5) coordinate (l2);
\path (17.5,14) coordinate (l21);
\path (17,14.75) coordinate (l211);
\path (18,14.75) coordinate (l212);
%
%
\draw (d1)--(d2)--(d3)--(d4)--(d5)--(d2);
\draw (d1)--(d4);
\draw (f1)--(f2)--(f3)--(f4)--(f5)--(f2);
\draw (f1)--(f4);
\draw (a1)--(a2)--(a3)--(a4)--(a5)--(a2);
\draw (a1)--(a4);
\draw (a5)--(d1);
\draw (a1)--(j1);
%
\draw (k1)--(k2)--(k3)--(k4)--(k5)--(k2);
\draw (k1)--(k4);
\draw (d5)--(k1);
\draw (j1)--(j11);
\draw (j111)--(j11)--(j112);
\draw (g1)--(g2)--(g3)--(g4)--(g5)--(g2);
\draw (g1)--(g4);
\draw (f5)--(g1);
\draw (b1)--(b2)--(b3)--(b4)--(b5)--(b2);
\draw (b1)--(b4);
\draw (b1)--(k5);
%
\draw (b5)--(l2)--(j1);
\draw (l2)--(l21);
\draw (l211)--(l21)--(l212);
\draw (a3) to[out=0,in=90, distance=0.75cm] (f1);
\draw[dotted, line width=0.05cm]  (d3) to[out=90,in=90, distance=0.75cm] (k3);
\draw[dotted, line width=0.05cm] (f3) to[out=0,in=90, distance=0.75cm] (f33);
\draw[dotted, line width=0.05cm] (g3) to[out=0,in=90, distance=0.75cm] (g33);
\draw[dotted, line width=0.05cm] (g5) to[out=0,in=90, distance=0.75cm] (g55);
\draw[dotted, line width=0.05cm] (b3) to[out=0,in=90, distance=0.75cm] (b33);

%
\draw [style=dashed,rounded corners] (1,5.5) rectangle (19,13.25);
\draw (1.5,6.5) node {{\small $T$}};
%
\draw (a1) [fill=white] circle (\vr);
\draw (a2) [fill=white] circle (\vr);
\draw (a3) [fill=white] circle (\vr);
\draw (a4) [fill=white] circle (\vr);
\draw (a5) [fill=white] circle (\vr);
\draw (b1) [fill=white] circle (\vr);
\draw (b2) [fill=white] circle (\vr);
\draw (b3) [fill=white] circle (\vr);
\draw (b4) [fill=white] circle (\vr);
\draw (b5) [fill=white] circle (\vr);
\draw (d1) [fill=white] circle (\vr);
\draw (d2) [fill=white] circle (\vr);
\draw (d3) [fill=white] circle (\vr);
\draw (d4) [fill=white] circle (\vr);
\draw (d5) [fill=white] circle (\vr);
\draw (f1) [fill=white] circle (\vr);
\draw (f2) [fill=white] circle (\vr);
\draw (f3) [fill=white] circle (\vr);
\draw (f4) [fill=white] circle (\vr);
\draw (f5) [fill=white] circle (\vr);
\draw (g1) [fill=white] circle (\vr);
\draw (g2) [fill=white] circle (\vr);
\draw (g3) [fill=white] circle (\vr);
\draw (g4) [fill=white] circle (\vr);
\draw (g5) [fill=white] circle (\vr);
\draw (j1) [fill=white] circle (\vr);
\draw (j11) [fill=white] circle (\vr);
\draw (k1) [fill=white] circle (\vr);
\draw (k2) [fill=white] circle (\vr);
\draw (k3) [fill=white] circle (\vr);
\draw (k4) [fill=white] circle (\vr);
\draw (k5) [fill=white] circle (\vr);
%
\draw (l2) [fill=white] circle (\vr);
\draw (l21) [fill=white] circle (\vr);
%
%
\draw[anchor = east] (d3p) node {{\small $a_1$}};
\draw[anchor = west] (k3) node {{\small $a_2$}};
\draw[anchor = east] (f3p) node {{\small $a_4$}};
\draw[anchor = east] (b3) node {{\small $a_3$}};
\draw[anchor = east] (g3p) node {{\small $a_5$}};
\draw[anchor = south] (g5) node {{\small $a_6$}};
%
\draw[anchor = east] (a1) node {{\small $u_1$}};
\draw[anchor = west] (b5) node {{\small $u_2$}};
\draw[anchor = east] (j1) node {{\small $v_1$}};
\draw[anchor = west] (l2) node {{\small $v_2$}};
\draw[anchor = west] (l21) node {{\small $z_2$}};
\draw[anchor = east] (j11) node {{\small $z_1$}};
\draw (10,4) node {{\small (a) A subgraph $T$ of $G$ where the dotted edges belong to $E_T$}};
\end{tikzpicture}
\vskip 0.5 cm
\begin{tikzpicture}[scale=.75,style=thick,x=0.75cm,y=0.75cm]
\def\vr{2.25pt}
\def\vrn{1.25pt}
\path (2.25,10) coordinate (a1);
\path (3.75,9) coordinate (a2);
\path (3.75,10) coordinate (a3);
\path (3.9,10) coordinate (a3p);
\path (3.75,11) coordinate (a4);
\path (5,10) coordinate (a5);
\path (6.5,10) coordinate (d1);
\path (7.75,9) coordinate (d2);
\path (7.75,11) coordinate (d4);
\path (9,10) coordinate (d5);
\path (6.5,7) coordinate (f1);
\path (7.75,6) coordinate (f2);
\path (7.75,8) coordinate (f4);
\path (9,7) coordinate (f5);
\path (10.75,7) coordinate (g1);
\path (12.25,6) coordinate (g2);
\path (12.25,8) coordinate (g4);
\path (13.75,7) coordinate (g5);
\path (13.7,7) coordinate (g5p);
\path (2.25,12) coordinate (j1);
\path (10.75,10) coordinate (k1);
\path (12.25,9) coordinate (k2);
\path (12.25,10) coordinate (k3);
\path (12.2,10) coordinate (k3p);
\path (12.25,11) coordinate (k4);
\path (13.75,10) coordinate (k5);
\path (13.75,8) coordinate (k33);
\path (15,10) coordinate (b1);
\path (16.25,9) coordinate (b2);
\path (16.25,11) coordinate (b4);
\path (17.5,10) coordinate (b5);
%
%
\draw (d1)--(d2);
\draw (d4)--(d5)--(d2);
\draw (d1)--(d4);
\draw (f1)--(f2);
\draw (f4)--(f5)--(f2);
\draw (f1)--(f4);
\draw (a1)--(a2)--(a3)--(a4)--(a5)--(a2);
\draw (a1)--(a4);
\draw (a5)--(d1);
\draw (a1)--(j1);
\draw (g1)--(f5);
%
\draw (k1)--(k2);
\draw (k4)--(k5)--(k2);
\draw (k1)--(k4);
\draw (d5)--(k1);
\draw (g1)--(g2);
\draw (g1)--(g4);
%
\draw (b1)--(b2);
\draw (b1)--(b4);
\draw (b1)--(k5);
%
\draw (a3) to[out=0,in=90, distance=0.75cm] (f1);
%
%
\draw (a1) [fill=black] circle (\vr);
\draw (a2) [fill=white] circle (\vr);
\draw (a3) [fill=white] circle (\vr);
\draw (a4) [fill=white] circle (\vr);
\draw (a5) [fill=white] circle (\vr);
\draw (b1) [fill=black] circle (\vr);
\draw (b2) [fill=white] circle (\vr);
\draw (b4) [fill=white] circle (\vr);
\draw (d1) [fill=black] circle (\vr);
\draw (d2) [fill=white] circle (\vr);
\draw (d4) [fill=white] circle (\vr);
\draw (d5) [fill=white] circle (\vr);
\draw (f1) [fill=black] circle (\vr);
\draw (f2) [fill=white] circle (\vr);
\draw (f4) [fill=white] circle (\vr);
\draw (f5) [fill=white] circle (\vr);
\draw (g1) [fill=black] circle (\vr);
\draw (g2) [fill=white] circle (\vr);
\draw (g4) [fill=white] circle (\vr);
%
\draw (j1) [fill=white] circle (\vr);
%
\draw (k1) [fill=black] circle (\vr);
\draw (k2) [fill=white] circle (\vr);
\draw (k4) [fill=white] circle (\vr);
\draw (k5) [fill=white] circle (\vr);
%
%
%
\draw[anchor = east] (j1) node {{\small $v_1$}};
\draw (10,5) node {{\small (b) The subgraph $G[X]$ where $X = V(T) \setminus \{a_1,a_2,a_3,a_4,a_5,a_6,u_2,v_2\}$}};
\end{tikzpicture}
\caption{An illustration of a subgraph of $G$ in the proof of Claim~\ref{claim:no-troublesome-subgraph.1}}
\label{fig:claim:min-deg2.2F}
\end{center}
\end{figure}

\begin{subclaim}
\label{claim:no-troublesome-subgraph.2}
If the link vertex~$v_1$ is incident with an edge of $E_T$, then $8i(G) < \Omega(G)$.
\end{subclaim}
\proof Suppose that the link vertex~$v_1$ is incident with an edge of $E_T$. Renaming vertices if necessary, we may assume that $e = a_1v_1$. We now let $X = V(T) \setminus (A \setminus \{a_1\})$ and we consider the graph $G' = G - X$. We note that each vertex in $A \setminus \{a_1\}$ has degree~$1$ in $G'$, and belongs to neither a bad component nor a troublesome configuration in $G'$. Moreover, the vertex~$z_2$ has degree~$1$ or~$2$  in $G'$, and by Claims~\ref{claim:no-bad-in-Ge} and~\ref{claim:no-T-in-Ge} we infer that $z_2$ does not belong to a bad component or a troublesome configuration in $G'$. Hence, $b(G') = \tc(G') = 0$, and so $\Theta(G') = 0$. Thus, $\Omega(G') = \w(G') = \w(G) - \w_G(X) + 2k-1 = \w(G) - 3 \times (4k+3) + 2k-1 = \w(G) - 10k - 10 = \Omega(G) - 10k - 10$.

Let $I_B$ be the core independent set (see Section~\ref{S:familyB}) of the bad graph $B$, and so $|I_B| = k$. We define the modified core independent set of $B$ as follows. We say that two copies of $K_{2,3}$ are adjacent if they are joined by an edge. By the structure of the bad graph $G$, there exists a sequence $F_1, \ldots, F_\ell$ of copies of $K_{2,3}$ where $u_1 \in V(F_1)$, $u_2 \in V(F_\ell)$ and if $\ell \ge 2$, then $F_{i+1}$ is adjacent to $F_i$ for all $i \in [\ell - 1]$. Let $c_i$ be the vertex in $F_i$ that belongs to $I_B$ for $i \in [\ell]$. In particular, we note that $c_1 = u_1$. Let $C = \{c_1,\ldots,c_{\ell}\}$.  We define $d_1 = v_1$, and if $\ell \ge 2$, then we define $d_{i+1}$ as the vertex of $F_{i}$ that is adjacent to $c_{i+1}$ in $B$ for $i \in [\ell - 1]$. Moreover, we define $d_{\ell + 1} = u_2$, and we let $D = \{d_1,d_2,\ldots,d_{\ell + 1}\}$. We now define the modified core independent set of $B$ by
\[
I_B^* = (I_B \setminus C) \cup D.
\]

As an illustration, suppose that $T$ is the troublesome configuration in Figure~\ref{fig:claim:min-deg2.2Fb}, where $B$ is the associated bad graph used to construct $T$. In this example, the core independent set $I_B$ is indicated by the shaded vertices in Figure~\ref{fig:claim:min-deg2.2Fb}(a), and the modified core independent set $I_B^*$ is indicated by the shaded vertices in Figure~\ref{fig:claim:min-deg2.2Fb}(b). We note that $|I_B| = k$ and $|I_B^*| = k+1$, where here $k = 6$ is the number of copies of $K_{2,3}$ in the bad graph $B$.

\begin{figure}[htb]
\begin{center}
\begin{tikzpicture}[scale=.75,style=thick,x=0.75cm,y=0.75cm]
\def\vr{2.25pt}
\def\vrn{1.25pt}
\path (-0.5,5.5) coordinate (z0);
\path (20,5.5) coordinate (z1);
\path (20,11.5) coordinate (z2);
\path (3,11.5) coordinate (z3);
\path (3,13.15) coordinate (z4);
\path (-0.5,13.15) coordinate (z5);
\path (2.25,10) coordinate (a1);
\path (3.75,9) coordinate (a2);
\path (3.75,10) coordinate (a3);
\path (3.9,10) coordinate (a3p);
\path (3.75,11) coordinate (a4);
\path (5,10) coordinate (a5);
\path (6.5,10) coordinate (d1);
\path (7.75,9) coordinate (d2);
\path (7.75,10) coordinate (d3);
\path (7.9,10) coordinate (d3p);
\path (7.75,11) coordinate (d4);
\path (9,10) coordinate (d5);
\path (9,9) coordinate (d33);
\path (6.5,7) coordinate (f1);
\path (7.75,6) coordinate (f2);
\path (7.75,7) coordinate (f3);
\path (7.85,7) coordinate (f3p);
\path (7.75,8) coordinate (f4);
\path (9,7) coordinate (f5);
\path (9,5) coordinate (f33);
\path (10.75,7) coordinate (g1);
\path (12.25,6) coordinate (g2);
\path (12.25,7) coordinate (g3);
\path (12.35,7) coordinate (g3p);
\path (12.25,8) coordinate (g4);
\path (13.75,7) coordinate (g5);
\path (13.7,7) coordinate (g5p);
\path (13.75,5) coordinate (g33);
\path (15,5) coordinate (g55);
\path (2.25,12.5) coordinate (j1);
\path (2.25,14) coordinate (j11);
\path (1.75,14.75) coordinate (j111);
\path (2.75,14.75) coordinate (j112);
\path (9,12) coordinate (j2);
\path (9,13.5) coordinate (j21);
\path (10.75,10) coordinate (k1);
\path (12.25,9) coordinate (k2);
\path (12.25,10) coordinate (k3);
\path (12.2,10) coordinate (k3p);
\path (12.25,11) coordinate (k4);
\path (13.75,10) coordinate (k5);
\path (13.75,8) coordinate (k33);
\path (15,10) coordinate (b1);
\path (16.25,9) coordinate (b2);
\path (16.25,10) coordinate (b3);
\path (16.35,10) coordinate (b3p);
\path (16.25,11) coordinate (b4);
\path (17.5,10) coordinate (b5);
\path (17.25,8) coordinate (b33);
%
\path (17.5,12.5) coordinate (l2);
\path (17.5,14) coordinate (l21);
\path (17,14.75) coordinate (l211);
\path (18,14.75) coordinate (l212);
%
%
\draw[dotted, line width=0.05cm] (z0)--(z1)--(z2)--(z3)--(z4)--(z5)--(z0);
\draw (d1)--(d2)--(d3)--(d4)--(d5)--(d2);
\draw (d1)--(d4);
\draw (f1)--(f2)--(f3)--(f4)--(f5)--(f2);
\draw (f1)--(f4);
\draw (a1)--(a2)--(a3)--(a4)--(a5)--(a2);
\draw (a1)--(a4);
\draw (a5)--(d1);
\draw (a1)--(j1);
\draw (k1)--(k2)--(k3)--(k4)--(k5)--(k2);
\draw (k1)--(k4);
\draw (d5)--(k1);
\draw (j1)--(j11);
\draw (j111)--(j11)--(j112);
\draw (g1)--(g2)--(g3)--(g4)--(g5)--(g2);
\draw (g1)--(g4);
\draw (f5)--(g1);
\draw (b1)--(b2)--(b3)--(b4)--(b5)--(b2);
\draw (b1)--(b4);
\draw (b1)--(k5);
%
\draw (b5)--(l2)--(j1);
\draw (l2)--(l21);
\draw (l211)--(l21)--(l212);
\draw (a3) to[out=0,in=90, distance=0.75cm] (f1);
\draw (0,6.5) node {{\small $B$}};
\draw (19.5,12.5) node {{\small $T$}};
\draw [style=dashed,rounded corners] (-1,5) rectangle (20.5,13.5);
\draw (a1) [fill=black] circle (\vr);
\draw (a2) [fill=white] circle (\vr);
\draw (a3) [fill=white] circle (\vr);
\draw (a4) [fill=white] circle (\vr);
\draw (a5) [fill=white] circle (\vr);
\draw (b1) [fill=black] circle (\vr);
\draw (b2) [fill=white] circle (\vr);
\draw (b3) [fill=white] circle (\vr);
\draw (b4) [fill=white] circle (\vr);
\draw (b5) [fill=white] circle (\vr);
\draw (d1) [fill=black] circle (\vr);
\draw (d2) [fill=white] circle (\vr);
\draw (d3) [fill=white] circle (\vr);
\draw (d4) [fill=white] circle (\vr);
\draw (d5) [fill=white] circle (\vr);
\draw (f1) [fill=black] circle (\vr);
\draw (f2) [fill=white] circle (\vr);
\draw (f3) [fill=white] circle (\vr);
\draw (f4) [fill=white] circle (\vr);
\draw (f5) [fill=white] circle (\vr);
\draw (g1) [fill=black] circle (\vr);
\draw (g2) [fill=white] circle (\vr);
\draw (g3) [fill=white] circle (\vr);
\draw (g4) [fill=white] circle (\vr);
\draw (g5) [fill=white] circle (\vr);
\draw (j1) [fill=white] circle (\vr);
\draw (j11) [fill=white] circle (\vr);
\draw (k1) [fill=black] circle (\vr);
\draw (k2) [fill=white] circle (\vr);
\draw (k3) [fill=white] circle (\vr);
\draw (k4) [fill=white] circle (\vr);
\draw (k5) [fill=white] circle (\vr);
%
\draw (l2) [fill=white] circle (\vr);
\draw (l21) [fill=white] circle (\vr);
%
%
\draw[anchor = south] (a5) node {{\small $d_2$}};
\draw[anchor = south] (d5) node {{\small $d_3$}};
\draw[anchor = south] (k5) node {{\small $d_4$}};
\draw[anchor = east] (a1) node {{\small $u_1 = c_1$}};
\draw[anchor = south] (d1) node {{\small $c_2$}};
\draw[anchor = south] (k1) node {{\small $c_3$}};
\draw[anchor = south] (b1) node {{\small $c_4$}};
\draw[anchor = west] (b5) node {{\small $u_2 = d_5$}};
\draw[anchor = east] (j1) node {{\small $v_1 = d_1$}};
\draw[anchor = west] (l2) node {{\small $v_2$}};
%
\draw (10,4) node {{\small (a) The core independent set, $I_B$, of $B$}};
\end{tikzpicture}
\vskip 0.5 cm
\begin{tikzpicture}[scale=.75,style=thick,x=0.75cm,y=0.75cm]
\def\vr{2.25pt}
\def\vrn{1.25pt}
\path (-0.5,5.5) coordinate (z0);
\path (20,5.5) coordinate (z1);
\path (20,11.5) coordinate (z2);
\path (3,11.5) coordinate (z3);
\path (3,13.15) coordinate (z4);
\path (-0.5,13.15) coordinate (z5);
\path (2.25,10) coordinate (a1);
\path (3.75,9) coordinate (a2);
\path (3.75,10) coordinate (a3);
\path (3.9,10) coordinate (a3p);
\path (3.75,11) coordinate (a4);
\path (5,10) coordinate (a5);
\path (6.5,10) coordinate (d1);
\path (7.75,9) coordinate (d2);
\path (7.75,10) coordinate (d3);
\path (7.9,10) coordinate (d3p);
\path (7.75,11) coordinate (d4);
\path (9,10) coordinate (d5);
\path (9,9) coordinate (d33);
\path (6.5,7) coordinate (f1);
\path (7.75,6) coordinate (f2);
\path (7.75,7) coordinate (f3);
\path (7.85,7) coordinate (f3p);
\path (7.75,8) coordinate (f4);
\path (9,7) coordinate (f5);
\path (9,5) coordinate (f33);
\path (10.75,7) coordinate (g1);
\path (12.25,6) coordinate (g2);
\path (12.25,7) coordinate (g3);
\path (12.35,7) coordinate (g3p);
\path (12.25,8) coordinate (g4);
\path (13.75,7) coordinate (g5);
\path (13.7,7) coordinate (g5p);
\path (13.75,5) coordinate (g33);
\path (15,5) coordinate (g55);
\path (2.25,12.5) coordinate (j1);
\path (2.25,14) coordinate (j11);
\path (1.75,14.75) coordinate (j111);
\path (2.75,14.75) coordinate (j112);
\path (9,12) coordinate (j2);
\path (9,13.5) coordinate (j21);
\path (10.75,10) coordinate (k1);
\path (12.25,9) coordinate (k2);
\path (12.25,10) coordinate (k3);
\path (12.2,10) coordinate (k3p);
\path (12.25,11) coordinate (k4);
\path (13.75,10) coordinate (k5);
\path (13.75,8) coordinate (k33);
\path (15,10) coordinate (b1);
\path (16.25,9) coordinate (b2);
\path (16.25,10) coordinate (b3);
\path (16.35,10) coordinate (b3p);
\path (16.25,11) coordinate (b4);
\path (17.5,10) coordinate (b5);
\path (17.25,8) coordinate (b33);
%
\path (17.5,12.5) coordinate (l2);
\path (17.5,14) coordinate (l21);
\path (17,14.75) coordinate (l211);
\path (18,14.75) coordinate (l212);
%
%
%
\draw[dotted, line width=0.05cm] (z0)--(z1)--(z2)--(z3)--(z4)--(z5)--(z0);
\draw [style=dashed,rounded corners] (-1,5) rectangle (20.5,13.5);
\draw (0,6.5) node {{\small $B$}};
\draw (19.5,12.5) node {{\small $T$}};
\draw (d1)--(d2)--(d3)--(d4)--(d5)--(d2);
\draw (d1)--(d4);
\draw (f1)--(f2)--(f3)--(f4)--(f5)--(f2);
\draw (f1)--(f4);
\draw (a1)--(a2)--(a3)--(a4)--(a5)--(a2);
\draw (a1)--(a4);
\draw (a5)--(d1);
\draw (a1)--(j1);
%
\draw (k1)--(k2)--(k3)--(k4)--(k5)--(k2);
\draw (k1)--(k4);
\draw (d5)--(k1);
\draw (j1)--(j11);
\draw (j111)--(j11)--(j112);
\draw (g1)--(g2)--(g3)--(g4)--(g5)--(g2);
\draw (g1)--(g4);
\draw (f5)--(g1);
\draw (b1)--(b2)--(b3)--(b4)--(b5)--(b2);
\draw (b1)--(b4);
\draw (b1)--(k5);
%
\draw (b5)--(l2)--(j1);
\draw (l2)--(l21);
\draw (l211)--(l21)--(l212);
\draw (a3) to[out=0,in=90, distance=0.75cm] (f1);
%
\draw (a1) [fill=white] circle (\vr);
\draw (a2) [fill=white] circle (\vr);
\draw (a3) [fill=white] circle (\vr);
\draw (a4) [fill=white] circle (\vr);
\draw (a5) [fill=black] circle (\vr);
\draw (b1) [fill=white] circle (\vr);
\draw (b2) [fill=white] circle (\vr);
\draw (b3) [fill=white] circle (\vr);
\draw (b4) [fill=white] circle (\vr);
\draw (b5) [fill=black] circle (\vr);
\draw (d1) [fill=white] circle (\vr);
\draw (d2) [fill=white] circle (\vr);
\draw (d3) [fill=white] circle (\vr);
\draw (d4) [fill=white] circle (\vr);
\draw (d5) [fill=black] circle (\vr);
\draw (f1) [fill=black] circle (\vr);
\draw (f2) [fill=white] circle (\vr);
\draw (f3) [fill=white] circle (\vr);
\draw (f4) [fill=white] circle (\vr);
\draw (f5) [fill=white] circle (\vr);
\draw (g1) [fill=black] circle (\vr);
\draw (g2) [fill=white] circle (\vr);
\draw (g3) [fill=white] circle (\vr);
\draw (g4) [fill=white] circle (\vr);
\draw (g5) [fill=white] circle (\vr);
\draw (j1) [fill=black] circle (\vr);
\draw (j11) [fill=white] circle (\vr);
\draw (k1) [fill=white] circle (\vr);
\draw (k2) [fill=white] circle (\vr);
\draw (k3) [fill=white] circle (\vr);
\draw (k4) [fill=white] circle (\vr);
\draw (k5) [fill=black] circle (\vr);
%
\draw (l2) [fill=white] circle (\vr);
\draw (l21) [fill=white] circle (\vr);
%
%
\draw[anchor = south] (a5) node {{\small $d_2$}};
\draw[anchor = south] (d5) node {{\small $d_3$}};
\draw[anchor = south] (k5) node {{\small $d_4$}};
\draw[anchor = east] (a1) node {{\small $u_1 = c_1$}};
\draw[anchor = south] (d1) node {{\small $c_2$}};
\draw[anchor = south] (k1) node {{\small $c_3$}};
\draw[anchor = south] (b1) node {{\small $c_4$}};
\draw[anchor = west] (b5) node {{\small $u_2 = d_5$}};
\draw[anchor = east] (j1) node {{\small $v_1 = d_1$}};
\draw[anchor = west] (l2) node {{\small $v_2$}};
\draw (10,4) node {{\small (b) The modified core independent set, $I_B^*$, of $B$}};
\end{tikzpicture}
\caption{An illustration of a modified core independent set in a bad graph}
\label{fig:claim:min-deg2.2Fb}
\end{center}
\end{figure}

We now return to the proof of Claim~\ref{claim:no-troublesome-subgraph.2}. Every $i$-set of $G'$ can be extended to an ID-set of $G$ by adding to it the modified core independent set $I_B^*$ of $B$, and so $i(G) \le i(G') + |I_B^*| = i(G') + k + 1$. (For example, if the graph shown in Figure~\ref{fig:claim:min-deg2.2E}(a) is a subgraph of $G$, where the dotted edges are the edges that belong to the set $E_T$, then $k = 6$ and the subgraph $G[X]$ is illustrated in Figure~\ref{fig:claim:min-deg2.2E}(b) where the shaded vertices indicate the vertices of the modified core independent set $I_B^*$ of $B$.) Recall that $\Omega(G') + 10k + 10 = \Omega(G)$. Thus, $8i(G) \le 8(i(G') + k + 1) \le \Omega(G') + 8k + 8 = \Omega(G) - 2k - 2 < \Omega(G)$, a contradiction.~\smallqed

\begin{figure}[htb]
\begin{center}
\begin{tikzpicture}[scale=.75,style=thick,x=0.75cm,y=0.75cm]
\def\vr{2.25pt}
\def\vrn{1.25pt}

\path (2.25,10) coordinate (a1);
\path (3.75,9) coordinate (a2);
\path (3.75,10) coordinate (a3);
\path (3.9,10) coordinate (a3p);
\path (3.75,11) coordinate (a4);
\path (5,10) coordinate (a5);
\path (6.5,10) coordinate (d1);
\path (7.75,9) coordinate (d2);
\path (7.75,10) coordinate (d3);
\path (7.9,10) coordinate (d3p);
\path (7.75,11) coordinate (d4);
\path (9,10) coordinate (d5);
\path (9,9) coordinate (d33);
\path (6.5,7) coordinate (f1);
\path (7.75,6) coordinate (f2);
\path (7.75,7) coordinate (f3);
\path (7.9,7) coordinate (f3p);
\path (7.75,8) coordinate (f4);
\path (9,7) coordinate (f5);
\path (9,5) coordinate (f33);
\path (10.75,7) coordinate (g1);
\path (12.25,6) coordinate (g2);
\path (12.25,7) coordinate (g3);
\path (12.2,7) coordinate (g3p);
\path (12.25,8) coordinate (g4);
\path (13.75,7) coordinate (g5);
\path (13.7,7) coordinate (g5p);
\path (13.75,5) coordinate (g33);
\path (15,5) coordinate (g55);
\path (2.25,12.5) coordinate (j1);
\path (10.75,10) coordinate (k1);
\path (12.25,9) coordinate (k2);
\path (12.25,10) coordinate (k3);
\path (12.2,10) coordinate (k3p);
\path (12.25,11) coordinate (k4);
\path (13.75,10) coordinate (k5);
\path (13.75,8) coordinate (k33);
\path (15,10) coordinate (b1);
\path (16.25,9) coordinate (b2);
\path (16.25,10) coordinate (b3);
\path (16.2,10) coordinate (b3p);
\path (16.25,11) coordinate (b4);
\path (17.5,10) coordinate (b5);
\path (17.25,8) coordinate (b33);
%
\path (17.5,12.5) coordinate (l2);
\path (17.5,14) coordinate (l21);
\path (17,14.75) coordinate (l211);
\path (18,14.75) coordinate (l212);
%
%
\draw (d1)--(d2)--(d3)--(d4)--(d5)--(d2);
\draw (d1)--(d4);
\draw (f1)--(f2)--(f3)--(f4)--(f5)--(f2);
\draw (f1)--(f4);
\draw (a1)--(a2)--(a3)--(a4)--(a5)--(a2);
\draw (a1)--(a4);
\draw (a5)--(d1);
\draw (a1)--(j1);
%
\draw (k1)--(k2)--(k3)--(k4)--(k5)--(k2);
\draw (k1)--(k4);
\draw (d5)--(k1);
\draw (g1)--(g2)--(g3)--(g4)--(g5)--(g2);
\draw (g1)--(g4);
\draw (f5)--(g1);
\draw (b1)--(b2)--(b3)--(b4)--(b5)--(b2);
\draw (b1)--(b4);
\draw (b1)--(k5);
%
\draw (b5)--(l2)--(j1);
\draw (l2)--(l21);
\draw (l211)--(l21)--(l212);
\draw (a3) to[out=0,in=90, distance=0.75cm] (f1);
\draw[dotted, line width=0.05cm]  (j1) to[out=0,in=90, distance=0.75cm] (d3);
\draw[dotted, line width=0.05cm] (f3) to[out=0,in=90, distance=0.75cm] (f33);
\draw[dotted, line width=0.05cm] (g3) to[out=0,in=90, distance=0.75cm] (g33);
\draw[dotted, line width=0.05cm] (g5) to[out=0,in=90, distance=0.75cm] (g55);
\draw[dotted, line width=0.05cm] (k3) to[out=0,in=90, distance=0.75cm] (k33);
\draw[dotted, line width=0.05cm] (b3) to[out=0,in=90, distance=0.75cm] (b33);

%
\draw [style=dashed,rounded corners] (1,5.5) rectangle (19,13.25);
\draw (1.5,6.5) node {{\small $T$}};
%
\draw (a1) [fill=white] circle (\vr);
\draw (a2) [fill=white] circle (\vr);
\draw (a3) [fill=white] circle (\vr);
\draw (a4) [fill=white] circle (\vr);
\draw (a5) [fill=white] circle (\vr);
\draw (b1) [fill=white] circle (\vr);
\draw (b2) [fill=white] circle (\vr);
\draw (b3) [fill=white] circle (\vr);
\draw (b4) [fill=white] circle (\vr);
\draw (b5) [fill=white] circle (\vr);
\draw (d1) [fill=white] circle (\vr);
\draw (d2) [fill=white] circle (\vr);
\draw (d3) [fill=white] circle (\vr);
\draw (d4) [fill=white] circle (\vr);
\draw (d5) [fill=white] circle (\vr);
\draw (f1) [fill=white] circle (\vr);
\draw (f2) [fill=white] circle (\vr);
\draw (f3) [fill=white] circle (\vr);
\draw (f4) [fill=white] circle (\vr);
\draw (f5) [fill=white] circle (\vr);
\draw (g1) [fill=white] circle (\vr);
\draw (g2) [fill=white] circle (\vr);
\draw (g3) [fill=white] circle (\vr);
\draw (g4) [fill=white] circle (\vr);
\draw (g5) [fill=white] circle (\vr);
\draw (j1) [fill=white] circle (\vr);
%
\draw (k1) [fill=white] circle (\vr);
\draw (k2) [fill=white] circle (\vr);
\draw (k3) [fill=white] circle (\vr);
\draw (k4) [fill=white] circle (\vr);
\draw (k5) [fill=white] circle (\vr);
%
\draw (l2) [fill=white] circle (\vr);
\draw (l21) [fill=white] circle (\vr);
%
%
\draw[anchor = east] (d3p) node {{\small $a_1$}};
\draw[anchor = east] (k3p) node {{\small $a_2$}};
\draw[anchor = east] (f3p) node {{\small $a_4$}};
\draw[anchor = east] (b3p) node {{\small $a_3$}};
\draw[anchor = east] (g3p) node {{\small $a_5$}};
\draw[anchor = east] (g5p) node {{\small $a_6$}};
%
\draw[anchor = east] (a1) node {{\small $u_1$}};
\draw[anchor = west] (b5) node {{\small $u_2$}};
\draw[anchor = east] (j1) node {{\small $v_1$}};
\draw[anchor = west] (l2) node {{\small $v_2$}};
\draw[anchor = west] (l21) node {{\small $z_2$}};
\draw (10,4) node {{\small (a) A subgraph $T$ of $G$ where the dotted edges belong to $E_T$}};
\end{tikzpicture}
\vskip 0.5 cm
\begin{tikzpicture}[scale=.75,style=thick,x=0.75cm,y=0.75cm]
\def\vr{2.25pt}
\def\vrn{1.25pt}
\path (2.25,10) coordinate (a1);
\path (3.75,9) coordinate (a2);
\path (3.75,10) coordinate (a3);
\path (3.9,10) coordinate (a3p);
\path (3.75,11) coordinate (a4);
\path (5,10) coordinate (a5);
\path (6.5,10) coordinate (d1);
\path (7.75,9) coordinate (d2);
\path (7.75,11) coordinate (d4);
\path (9,10) coordinate (d5);
\path (6.5,7) coordinate (f1);
\path (7.75,6) coordinate (f2);
\path (7.75,8) coordinate (f4);
\path (9,7) coordinate (f5);
\path (10.75,7) coordinate (g1);
\path (12.25,6) coordinate (g2);
\path (12.25,8) coordinate (g4);
\path (13.75,7) coordinate (g5);
\path (13.7,7) coordinate (g5p);
\path (2.25,12.5) coordinate (j1);
\path (10.75,10) coordinate (k1);
\path (12.25,9) coordinate (k2);
\path (12.25,10) coordinate (k3);
\path (12.2,10) coordinate (k3p);
\path (12.25,11) coordinate (k4);
\path (13.75,10) coordinate (k5);
\path (13.75,8) coordinate (k33);
\path (15,10) coordinate (b1);
\path (16.25,9) coordinate (b2);
\path (16.25,11) coordinate (b4);
\path (17.5,10) coordinate (b5);
%
\path (17.5,12.5) coordinate (l2);
%
\draw (d1)--(d2)--(d3)--(d4)--(d5)--(d2);
\draw (d1)--(d4);
\draw (f1)--(f2);
\draw (f4)--(f5)--(f2);
\draw (f1)--(f4);
\draw (a1)--(a2)--(a3)--(a4)--(a5)--(a2);
\draw (a1)--(a4);
\draw (a5)--(d1);
\draw (a1)--(j1);
\draw (g1)--(f5);
%
\draw (k1)--(k2);
\draw (k4)--(k5)--(k2);
\draw (k1)--(k4);
\draw (d5)--(k1);
\draw (g1)--(g2);
\draw (g1)--(g4);
%
\draw (b1)--(b2);
\draw (b4)--(b5)--(b2);
\draw (b1)--(b4);
\draw (b1)--(k5);
%
\draw (b5)--(l2)--(j1);
\draw (a3) to[out=0,in=90, distance=0.75cm] (f1);
\draw (j1) to[out=0,in=90, distance=0.75cm] (d3);
%
%
\draw (a1) [fill=white] circle (\vr);
\draw (a2) [fill=white] circle (\vr);
\draw (a3) [fill=white] circle (\vr);
\draw (a4) [fill=white] circle (\vr);
\draw (a5) [fill=black] circle (\vr);
\draw (b1) [fill=white] circle (\vr);
\draw (b2) [fill=white] circle (\vr);
\draw (b4) [fill=white] circle (\vr);
\draw (b5) [fill=black] circle (\vr);
\draw (d1) [fill=white] circle (\vr);
\draw (d2) [fill=white] circle (\vr);
\draw (d3) [fill=white] circle (\vr);
\draw (d4) [fill=white] circle (\vr);
\draw (d5) [fill=black] circle (\vr);
\draw (f1) [fill=black] circle (\vr);
\draw (f2) [fill=white] circle (\vr);
\draw (f4) [fill=white] circle (\vr);
\draw (f5) [fill=white] circle (\vr);
\draw (g1) [fill=black] circle (\vr);
\draw (g2) [fill=white] circle (\vr);
\draw (g4) [fill=white] circle (\vr);
%
\draw (j1) [fill=black] circle (\vr);
%
\draw (k1) [fill=white] circle (\vr);
\draw (k2) [fill=white] circle (\vr);
\draw (k4) [fill=white] circle (\vr);
\draw (k5) [fill=black] circle (\vr);
\draw (l2) [fill=white] circle (\vr);
\draw[anchor = east] (d3p) node {{\small $a_1$}};
%
\draw[anchor = east] (a1) node {{\small $u_1$}};
\draw[anchor = west] (b5) node {{\small $u_2$}};
\draw[anchor = east] (j1) node {{\small $v_1$}};
\draw[anchor = west] (l2) node {{\small $v_2$}};
\draw (10,5) node {{\small (b) The subgraph $G[X]$ where $X = V(T) \setminus \{a_2,a_3,a_4,a_5,a_6\}$}};
\end{tikzpicture}
\caption{An illustration of a subgraph of $G$ in the proof of Claim~\ref{claim:no-troublesome-subgraph.2}}
\label{fig:claim:min-deg2.2E}
\end{center}
\end{figure}

In both Claims~\ref{claim:no-troublesome-subgraph.1} and~\ref{claim:no-troublesome-subgraph.2}, we have $8i(G) \le \Omega(G)$, a contradiction. This completes the proof of Claim~\ref{claim:no-troublesome-subgraph}.~\smallqed

\medskip
By Claim~\ref{claim:no-troublesome-subgraph}, there is no troublesome subgraph in $G$. Hence, if $G'$ is an arbitrary induced subgraph of $G$, then $\tc(G') = 0$, and so $\Theta(G') = 2b(G')$. We state this formally as follows.

\begin{claim}
\label{structural-weight}
If $G'$ is an induced subgraph of $G$, then $\Theta(G') = 2b(G')$.
\end{claim}

By Claim~\ref{structural-weight}, the structural weight of an induced subgraph $G'$ of $G$ is determined only by the number of bad components of $G'$, if any, that belong to the bad family~$\cB$. Recall that $\cB = \cB_1 \cup \cB_2 \cup \cB_3$. Moreover, every graph in the family~$\cB_1$ has one vertex of degree~$1$ and at least two vertices of degree~$2$, while every graph in the family~$\cB_2$ has at least five vertices of degree~$2$. Furthermore, every graph in the family~$\cB_3$ has at least six vertices of degree~$2$. We prove next certain properties that must hold if the graph $G$ contains vertices of degree~$2$.

\begin{claim}
\label{claim:min-cubic.1}
The graph $G$ does not contain a path $v_1v_2v_3$ such that $\deg_G(v_i) = 2$ for all $i \in [3]$.
\end{claim}
\textbf{Proof.} Suppose, to the contrary, that $G$ contains a path $v_1v_2v_3$ such that $\deg_G(v_i) = 2$ for all $i \in [3]$. Let $X = \{v_1,v_2,v_3\}$ and let $G' = G - X$. We note that there are two $X$-exit edges. If $G'$ contains an isolated vertex, then $G = C_4$, a contradiction. Hence, $\delta(G') \ge 1$ and the graph $G'$ contains at most two vertices of degree~$1$. By Claim~\ref{claim:no-bad-in-Ge}, we infer that $b(G') = 0$. Thus, by Claim~\ref{structural-weight}, we have $\Theta(G') = 0$. Thus, $\Omega(G') = \w(G') = \w(G) - \w_G(X) + 2 = \Omega(G) - 12 + 2 = \Omega(G) - 10$. Every $i$-set of $G'$ can be extended to an ID-set of $G$ by adding to it the vertex~$v_2$, and so $i(G) \le i(G') + 1$. Hence, $8i(G) \le 8(i(G') + 1) \le \Omega(G') + 8 <\Omega(G)$, a contradiction.~\smallqed

\begin{claim}
\label{claim:min-cubic.2}
No two adjacent vertices in $G$ both have degree~$2$.
\end{claim}
\proof Suppose, to the contrary, that $G$ contains two adjacent vertices, say $v_1$ and $v_2$, of degree~$2$. Suppose firstly that $v_1$ and $v_2$ have a common neighbor~$v$, and so $vv_1v_2v$ is a $3$-cycle in $G$. Let $u$ be the third neighbor of $v$ different from $v_1$ and $v_2$. Let $X = \{v,v_1,v_2\}$ and let $G' = G - X$. The graph $G'$ is a connected subcubic graph and all vertices in $G'$ have degree at least~$2$, except possibly for the vertex~$u$ which has degree at least~$1$. By Claims~\ref{claim:no-bad-in-Ge} and~\ref{structural-weight}, we infer that $\Theta(G') = 0$. Thus, $\Omega(G') = \w(G') = \w(G) - \w_G(X) + 1 = \Omega(G) - 11 + 1 = \Omega(G) - 10$. Every $i$-set of $G'$ can be extended to an ID-set of $G$ by adding to it the vertex~$v_1$, and so $i(G) \le i(G') + 1$. Hence, $8i(G) \le 8(i(G') + 1) \le \Omega(G') + 8 = \Omega(G)$, a contradiction.

The vertices $v_1$ and $v_2$ therefore have no common neighbor. Let $u_i$ be the neighbor of $v_i$ different from $v_{3-i}$ for $i \in [2]$. Since $v_1$ and $v_2$ have no common neighbor, $u_1 \ne u_2$. By Claim~\ref{claim:min-cubic.1}, both vertices $u_1$ and $u_2$ have degree~$3$ in $G$. We show firstly that $u_1u_2$ is not an edge in $G$. Suppose, to the contrary, that $u_1u_2 \in E(G)$. In this case, we let $X = \{v_1,v_2,u_1\}$ and consider the graph $G' = G - X$. We note that $G'$ is isolate-free and the vertex $u_2$ has degree~$1$ in $G'$.  By Claim~\ref{structural-weight}, $\tc(G')=0$. By Claims~\ref{claim:no-bad-in-Ge} and~\ref{claim:no-bad-in-Gee}, we infer that if $G'$ has more than one component, none of them is bad. In addition, if $G'$ is connected, then by Claim~\ref{claim:no-K23-with-d2-vertex} we infer $G' \notin {\cal B}$. Thus, $\Theta(G') = 0$. Thus, $\Omega(G') = \w(G') = \w(G) - \w_G(X) + 3 = \Omega(G) - 11 + 3 = \Omega(G) - 8$. Every $i$-set of $G'$ can be extended to an ID-set of $G$ by adding to it the vertex~$v_1$, and so $i(G) \le i(G') + 1$. Hence, $8i(G) \le 8(i(G') + 1) \le \Omega(G') + 8 = \Omega(G)$, a contradiction. Therefore, $u_1u_2 \notin E(G)$.

We now let $X = \{v_1,v_2\}$ and consider the graph $G'$ obtained from $G - X$ by adding the edge $e = u_1u_2$. We note that $G'$ is a connected subcubic graph satisfying $\delta(G') \ge 2$. Furthermore, $u_1$ and $u_2$ are adjacent vertices of degree~$3$ in $G'$. Let $I'$ be an $i$-set of $G'$. If $u_1 \in I'$, then let $I = I' \cup \{v_2\}$. If $u_2 \in I'$ or if neither $u_1$ nor $u_2$ belong to $I'$, then let $I = I' \cup \{v_1\}$. In all cases, $|I| = |I'| + 1$ and the set $I$ is an ID-set of $G$, and so $i(G) \le i(G') + 1$. If $\Theta(G') = 0$, then $\Omega(G') = \w(G') = \w(G) - \w_G(X) = \Omega(G) - 8$, and so $8i(G) \le 8(i(G') + 1) \le \Omega(G') + 8 = \Omega(G)$, a contradiction. Hence, $\Theta(G') \ge 1$, implying that either $G' \in \cB$ or the edge $u_1u_2$ belongs to a troublesome configuration in $G'$. If $G' \in \cB$, then $G$ would contain a copy of $K_{2,3}$ that contains two vertices of degree~$2$ in $G$ (noting that in this case, the root vertex in $G'$ has degree at least~$2$), a contradiction. Hence, $G' \notin \cB$.

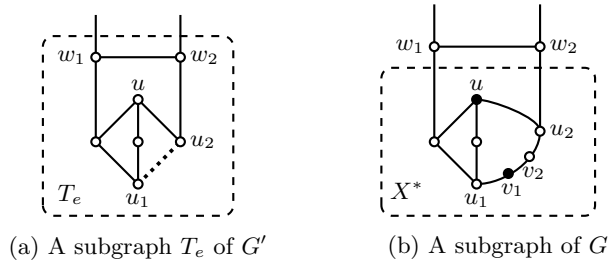
\begin{figure}[htb]
\begin{center}
\begin{tikzpicture}[scale=.75,style=thick,x=0.75cm,y=0.75cm]
\def\vr{2.25pt}
\def\vrn{1.25pt}
\path (-2,4.5) coordinate (b1);
\path (-1,3.5) coordinate (b2);
\path (-1,4.5) coordinate (b3);
\path (-1,5.5) coordinate (b4);
\path (0,4.5) coordinate (b5);
\path (-2,6.5) coordinate (e1);
\path (-2,7.5) coordinate (e11);
\path (0,6.5) coordinate (g5);
\path (0,7.5) coordinate (g51);
\draw (b1)--(e1)--(e11);
\draw (e1)--(g5);
\draw (b1)--(b2)--(b3)--(b4)--(b5);
\draw (b1)--(b4);
\draw (b5)--(g5)--(g51);
\draw [dotted, line width=0.05cm] (b2)--(b5);
\draw [style=dashed,rounded corners] (-3.25,2.75) rectangle (1.25,7);
\draw (-2.6,3.25) node {{\small $T_e$}};
\draw (e1) [fill=white] circle (\vr);
\draw (b1) [fill=white] circle (\vr);
\draw (b2) [fill=white] circle (\vr);
\draw (b3) [fill=white] circle (\vr);
\draw (b4) [fill=white] circle (\vr);
\draw (b5) [fill=white] circle (\vr);
\draw (g5) [fill=white] circle (\vr);
\draw[anchor = west] (b5) node {{\small $u_2$}};
\draw[anchor = north] (b2) node {{\small $u_1$}};
\draw[anchor = south] (b4) node {{\small $u$}};
\draw[anchor = west] (g5) node {{\small $w_2$}};
\draw[anchor = east] (e1) node {{\small $w_1$}};
\draw (-1,2) node {{\small (a) A subgraph $T_e$ of $G'$}};
\path (6,4.5) coordinate (b1);
\path (7,3.5) coordinate (b2);
\path (7,4.5) coordinate (b3);
\path (7,5.5) coordinate (b4);
\path (8.5,4.75) coordinate (b5);
\path (7.75,3.75) coordinate (v1);
\path (7.85,3.75) coordinate (v1p);
\path (8.25,4.15) coordinate (v2);
\path (8.35,4.15) coordinate (v2p);
\path (6,6.75) coordinate (e1);
\path (6,7.75) coordinate (e11);
\path (8.5,6.75) coordinate (g5);
\path (8.5,7.75) coordinate (g51);
\draw (b1)--(e1)--(e11);
\draw (e1)--(g5);
\draw (b1)--(b2)--(b3)--(b4);
\draw (b1)--(b4);
%
\draw (b5)--(g5)--(g51);
\draw [style=dashed,rounded corners] (4.75,2.75) rectangle (9.8,6.25);
\draw (5.35,3.35) node {{\small $X^*$}};
%
\draw (b4) to[out=0,in=90, distance=0.25cm] (b5);
\draw (b2) to[out=0,in=270, distance=0.5cm] (b5);
\draw (e1) [fill=white] circle (\vr);
\draw (b1) [fill=white] circle (\vr);
\draw (b2) [fill=white] circle (\vr);
\draw (b3) [fill=white] circle (\vr);
\draw (b4) [fill=black] circle (\vr);
\draw (b5) [fill=white] circle (\vr);
\draw (g5) [fill=white] circle (\vr);
\draw (v1) [fill=black] circle (\vr);
\draw (v2) [fill=white] circle (\vr);
\draw[anchor = west] (b5) node {{\small $u_2$}};
\draw[anchor = north] (b2) node {{\small $u_1$}};
\draw[anchor = south] (b4) node {{\small $u$}};
\draw[anchor = west] (g5) node {{\small $w_2$}};
\draw[anchor = east] (e1) node {{\small $w_1$}};
\draw[anchor = north] (v1p) node {{\small $v_1$}};
\draw[anchor = north] (v2p) node {{\small $v_2$}};
\draw (7.5,2) node {{\small (b) A subgraph of $G$}};
\end{tikzpicture}
\caption{Subgraphs of $G'$ and $G$ in the proof of Claim~\ref{claim:min-cubic.2}}
\label{fig:claim:min-cubic.2f}
\end{center}
\end{figure}

Thus, adding the edge $e = u_1u_2$ to the graph $G - X$ creates a new troublesome configuration, which we call $T_e$, that necessarily contains the edge~$e$. Let $w_1$ and $w_2$ be the two link vertices of $T_e$. By combining Claim~\ref{claim:no-K23-with-d2-vertex} and Claim~\ref{claim:no-troublesome-subgraph}, we infer that $T_e$ contains exactly one copy of $K_{2,3}$, say $F_e$, and both $u_1$ and $u_2$ belong to $F_e$ (for otherwise, $F_e$ would contain a vertex of degree~$2$ in $G$, a contradiction). By symmetry and renaming vertices if necessary, we may assume that $u_2$ is adjacent to the link vertex $w_2$ of $T_e$. Let $u$ be the third neighbor of $u_2$ in $T_e$ different from $u_1$ and $w_2$, as illustrated in Figure~\ref{fig:claim:min-cubic.2f}(a). We now let $X^* = V(F_e) \cup X$ and let $G^* = G - X^*$. Every $i$-set of $G^*$ can be extended to an ID-set of $G$ by adding to it the vertices $v_1$ and $u$, as indicated by the shaded vertices in Figure~\ref{fig:claim:min-cubic.2f}(b). Thus, $i(G) \le i(G^*) + 2$. We note that $\Theta(G^*) = \Theta(G) = 0$, and so $\Omega(G^*) = \w(G^*) = \w(G) - \w_G(X^*) + 2 = \Omega(G) - 24 + 2 = \Omega(G) - 22$. Hence, $8i(G) \le 8(i(G^*) + 2) \le \Omega(G^*) + 16 < \Omega(G)$, a contradiction.~\smallqed

\begin{claim}
\label{claim:min-cubic.3}
A vertex of degree~$3$ in $G$ has at most one neighbor of degree~$2$.
\end{claim}
\proof Suppose, to the contrary, that $G$ contains a vertex, say $v$, of degree~$3$ that is adjacent to at least two vertices of degree~$2$. Let $X = N_G[v]$ and let $G' = G - X$. Every $i$-set of $G'$ can be extended to an ID-set of $G$ by adding to it the vertex $v$, and so $i(G) \le i(G') + 1$.

Suppose that all three neighbors of $v$ have degree~$2$ in $G$. By Claim~\ref{claim:min-cubic.2}, the set $N_G(v)$ is therefore an independent set. Suppose that $G'$ contains an isolated vertex, $u$ say. Since $\delta(G) \ge 2$, we note that $\deg_G(u) \ge 2$. If $\deg_G(u) = 3$, then $G = K_{2,3}$, contradicting the fact that $n \ge 6$. If $\deg_G(u) = 2$, then since $u$ is adjacent to two vertices in $N_G(v)$, the graph $G$ contains a path $P_3$ all of whose vertices have degree~$2$ in $G$, contradicting Claim~\ref{claim:min-cubic.1}. Hence, the graph $G'$ is isolate-free. Since there are three $X$-exit edges, we note that $\w(G') = (\w(G) - 15) + 3 = \w(G) - 12$. By Claim~\ref{claim:no-bad-in-Gee}, we infer that $b(G') \le 1$. Thus by Claim~\ref{structural-weight}, we have $\Theta(G') \le 2$, and so $\Omega(G') = \w(G') + \Theta(G') = (\w(G) - 12) + 2 = \Omega(G) - 10$. Every $i$-set of $G'$ can be extended to an ID-set of $G$ by adding to it the vertex~$v$, and so $i(G) \le i(G') + 1$. Hence, $8i(G) \le 8(i(G') + 1) \le \Omega(G') + 8 < \Omega(G)$, a contradiction.

Therefore, the vertex $v$ has two neighbors of degree~$2$ and one of degree~$3$. Let $v_1$ and $v_2$ be the two neighbors of~$v$ of degree~$2$ and let $v_3$ be the third neighbor of $v$ (of degree~$3$). Suppose that $G'$ contains an isolated vertex, $u$ say. If $\deg_G(u) = 3$, then $G$ contains a copy of $K_{2,3}$ that contains two vertices of degree~$2$ in $G$, contradicting Claim~\ref{claim:no-K23-with-d2-vertex}. If $\deg_G(u) = 2$, then since $u$ is adjacent to at least one of $v_1$ and $v_2$, the graph $G$ contains two adjacent vertices of degree~$2$, contradicting Claim~\ref{claim:min-cubic.2}. Hence, the graph $G'$ is isolate-free. Every $i$-set of $G'$ can be extended to an ID-set of $G$ by adding to it the vertex $v$, and so $i(G) \le i(G') + 1$. By Claim~\ref{claim:no-bad-in-Gee}, we infer that $b(G') \le 1$. Thus by Claim~\ref{structural-weight}, we have $\Theta(G') \le 2$, and so $\Omega(G') = \w(G') + \Theta(G') \le (\w(G) - \w_G(X) + 4) + 2 = \w(G) - 14 + 4 + 2 = \Omega(G) - 8$. Thus, $8i(G) \le 8(i(G') + 1) \le \Omega(G') + 8 = \Omega(G)$, a contradiction.~\smallqed

\medskip
By Claim~\ref{claim:min-cubic.2}, both neighbors of a vertex of degree~$2$ have degree~$3$ in $G$. By Claim~\ref{claim:min-cubic.3}, at most one neighbor of a vertex of degree~$3$ has degree~$2$ in $G$.

\begin{claim}
\label{claim:min-cubic.4}
No vertex of degree~$2$ belongs to a triangle.
\end{claim}
\proof Suppose, to the contrary, that $G$ contains a vertex, say $u$, of degree~$2$ that belongs to a triangle, $T$ say. Let $v_1$ and $v_2$ be the two neighbors of $u$, and so $v_1v_2 \in E(G)$ and $V(T) = \{u,v_1,v_2\}$. By our earlier observations, both $v_1$ and $v_2$ have degree~$3$ in $G$. Let $w_i$ be the neighbor of $v_i$ not in $T$ for $i \in [2]$. Since at most one neighbor of a vertex of degree~$3$ has degree~$2$ in $G$, both $w_1$ and $w_2$ have degree~$3$.

Suppose firstly that $w_1 = w_2$. In this case, we let $X = N_G[v_1]$ and consider the graph $G' = G - X$. Every $i$-set of $G'$ can be extended to an ID-set of $G$ by adding the vertex~$v_1$, and so $i(G) \le i(G') + 1$.  By Claims~\ref{claim:no-bad-in-Ge} and~\ref{structural-weight}, we infer that $\Theta(G') = 0$, and so $\Omega(G') = \w(G') = \w(G) - \w_G(X) + 1 = \w(G) - 13 + 1 = \Omega(G) - 12$. Thus, $8i(G) \le 8(i(G') + 1) \le \Omega(G') + 8 < \Omega(G)$, a contradiction. Hence, $w_1 \ne w_2$.

Suppose that $w_1w_2 \notin E(G)$. In this case, we let $X = N_G[u] = \{u,v_1,v_2\}$ and let $G'$ be obtained from $G - X$ by adding the edge $w_1w_2$. Let $I'$ be an $i$-set of $G'$. If $w_1 \in I'$, then let $I = I \cup \{v_2\}$, while if $w_1 \notin I'$, then let $I = I \cup \{v_1\}$. In both cases, $|I| = |I'| + 1$ and the set $I$ is an ID-set of $G$, and so $i(G) \le i(G') + 1$. The added edge $w_1w_2$ increases the structural weight by at most~$2$ (which occurs if $G' \in\cB$ or if the edge $w_1w_2$ belongs to a troublesome configuration). Hence, $\Theta(G') \le 2$, and so $\Omega(G') = \w(G') + \Theta(G') \le (\w(G) - 10) + 2 = \Omega(G) - 8$. Thus, $8i(G) \le 8(i(G') + 1) \le \Omega(G') + 8 = \Omega(G)$, a contradiction.

Hence, $w_1w_2 \in E(G)$. We now let $X = N_G[v_1] = \{u,v_1,v_2,w_1\}$ and consider the graph $G' = G - X$. Every $i$-set of $G'$ can be extended to an ID-set of $G$ by adding the vertex~$v_1$, and so $i(G) \le i(G') + 1$. 
By combining Claims~\ref{claim:no-K23-with-d2-vertex},~\ref{claim:no-bad-in-Ge} and~\ref{structural-weight}, we infer that $\Theta(G') = 0$, and so $\Omega(G') = \w(G') = \w(G) - \w_G(X) + 3 = \w(G) - 13 + 3 = \Omega(G) - 10$. Therefore, $8i(G) \le 8(i(G') + 1) \le \Omega(G') + 8 < \Omega(G)$, a contradiction.~\smallqed

\medskip
Recall that by Claim~\ref{claim:K23}, the graph $G$ contains at least one subgraph isomorphic to $K_{2,3}$, and by Claim~\ref{claim:Gconn}, the graph $G$ is connected. By supposition, $G \ne K_{3,3}$ and $G \ne C_5 \, \Box \, K_2$. We proceed further by proving additional structural properties of the graph $G$. A \emph{diamond} in $G$ is a copy of $K_4 - e$ where $e$ is an arbitrary edge of the complete graph $K_4$. A graph is \emph{diamond}-\emph{free} if it does not contain a diamond as an induced subgraph.

\begin{claim}
\label{no-diamond}
The graph $G$ is diamond-free.
\end{claim}
\proof Suppose, to the contrary, that $G$ contains a diamond $D$. Let $V(D) = \{v_1,v_2,v_3,v_4\}$ where $v_1v_2$ is the missing edge in the diamond. By Claim~\ref{claim:min-cubic.4}, no vertex of degree~$2$ belongs to a triangle. Thus, both $v_1$ and $v_2$ have degree~$3$ in $G$. We now let $X = V(D)$ and we consider the graph $G' = G - X$. Since $n \ge 6$, we note that $G'$ is isolate-free. Every $i$-set of $G'$ can be extended to an ID-set of $G$ by adding the vertex~$v_3$, and so $i(G) \le i(G') + 1$. By Claims~\ref{claim:no-bad-in-Gee} and~\ref{structural-weight}, we infer that $\Theta(G') = 0$, and so $\Omega(G') = \w(G') + \Theta(G') \le \w(G) - \w_G(X) + 2  = \w(G) - 12 + 2 = \Omega(G) - 10$. Thus, $8i(G) \le 8(i(G') + 1) \le \Omega(G') + 8 <\Omega(G)$, a contradiction.~\smallqed

\begin{claim}
\label{claim:no-K33-minus-e}
The graph $G$ does not contain $K_{3,3} - e$ as a subgraph, where $e$ is an edge of the $K_{3,3}$.
\end{claim}
\proof Suppose, to the contrary, that $G$ contains $K_{3,3} - e$ as a subgraph, where $e$ is an edge of the $K_{3,3}$. Let $F$ be such a subgraph in $G$, and let $F$ have partite sets $A = \{x_1,x_2,x_3\}$ and $B = \{v_1,v_2,v_3\}$, where $e = x_1v_1$ is the missing edge from the copy of $K_{3,3}$ in $F$. Since $G$ is diamond-free by Claim~\ref{no-diamond}, the set $B$ is an independent set in $G$. If $\deg_G(x_1) = 2$ or if $\deg_G(v_1) = 2$, then there is a $K_{2,3}$-subgraph in $G$ that contains a vertex of degree~$2$ in $G$, a contradiction. Hence, $\deg_G(x_1) = \deg_G(v_1) = 3$. Let $u_1$ be the neighbor of $v_1$ different from $x_2$ and $x_3$, and let $w_1$ be the neighbor of $x_1$ different from $v_2$ and $v_3$. Since $G \ne K_{3,3}$, the vertices $u_1$ and $x_1$ are distinct. If $x_1u_1 \in E(G)$, then $G$ contains $G_7$ as a subgraph (see Figure~\ref{fig:G7G8}), contradicting Claim~\ref{no-G7}. Hence, $u_1x_1 \notin E(G)$, and so the vertices $u_1$ and $w_1$ are distinct. We now let $X = V(F) \cup \{u_1,w_1\}$ and consider the graph $G' = G - X$.

\begin{figure}[htb]
\begin{center}
\begin{tikzpicture}[scale=.75,style=thick,x=0.8cm,y=0.8cm]
\def\vr{2.5pt} 
%
\path (2.85,1.875) coordinate (z2);
\path (4,1.875) coordinate (w);
\path (4,0.875) coordinate (u);
%
\path (2.85,0) coordinate (z4);
\path (4,0) coordinate (x1);
\path (5,1.875) coordinate (v2n);
\path (6.25,0) coordinate (u1);
\path (6.25,1.25) coordinate (u2);
\path (6.25,2.5) coordinate (u3);
\path (7.5,0.625) coordinate (w1);
\path (7.5,1.875) coordinate (w2);
\draw (x1)--(u1);
\draw (v2n)--(w);
\draw (u2)--(v2n)--(u3);
\draw (w1)--(u1)--(w2);
\draw (w)--(u)--(x1);
\draw (w1)--(u2)--(w2);
\draw (w1)--(u3)--(w2);
\draw (w)--(z2);
\draw (x1)--(z4);
\draw (x1) [fill=white] circle (\vr);
\draw (w) [fill=white] circle (\vr);
\draw (u) [fill=black] circle (\vr);
\draw (v2n) [fill=black] circle (\vr);
\draw (u1) [fill=black] circle (\vr);
\draw (u2) [fill=white] circle (\vr);
\draw (u3) [fill=white] circle (\vr);
\draw (w1) [fill=white] circle (\vr);
\draw (w2) [fill=white] circle (\vr);
\draw (5.25,-2) node {{\small (a)}};
\draw [style=dashed,rounded corners] (1.5,-1) rectangle (3.5,3.1);
\draw (2.255,1.25) node {{\small $G''$}};
%
\draw[anchor = north] (x1) node {{\small $u_1$}};
\draw[anchor = south] (w) node {{\small $w_1$}};
\draw[anchor = west] (u) node {{\small $u$}};
\draw[anchor = south] (v2n) node {{\small $x_1$}};
%
\draw[anchor = north] (u1) node {{\small $v_1$}};
\draw[anchor = north] (u2) node {{\small $v_2$}};
\draw[anchor = south] (u3) node {{\small $v_3$}};
\draw[anchor = west] (w1) node {{\small $x_2$}};
\draw[anchor = west] (w2) node {{\small $x_3$}};
\path (12.15,1.375) coordinate (z1);
\path (12.15,2.375) coordinate (z2);
\path (13,1.875) coordinate (w);
\path (12.15,-0.55) coordinate (z3);
\path (12.15,0.5) coordinate (z4);
\path (13,0) coordinate (x1);
\path (14,1.875) coordinate (v2n);
\path (15.25,0) coordinate (u1);
\path (15.25,1.25) coordinate (u2);
\path (15.25,2.5) coordinate (u3);
\path (16.5,0.625) coordinate (w1);
\path (16.5,1.875) coordinate (w2);
\draw (x1)--(u1);
\draw (v2n)--(w);
\draw (u2)--(v2n)--(u3);
\draw (w1)--(u1)--(w2);
\draw (w1)--(u2)--(w2);
\draw (w1)--(u3)--(w2);
\draw (z1)--(w)--(z2);
\draw (z3)--(x1)--(z4);
\draw (x1) [fill=white] circle (\vr);
\draw (w) [fill=white] circle (\vr);
\draw (v2n) [fill=black] circle (\vr);
\draw (u1) [fill=black] circle (\vr);
\draw (u2) [fill=white] circle (\vr);
\draw (u3) [fill=white] circle (\vr);
\draw (w1) [fill=white] circle (\vr);
\draw (w2) [fill=white] circle (\vr);
\draw (14.25,-2) node {{\small (b)}};
\draw [style=dashed,rounded corners] (10.5,-1) rectangle (12.5,3.1);
\draw (11.255,1.25) node {{\small $G'$}};
%
\draw[anchor = north] (x1) node {{\small $u_1$}};
\draw[anchor = south] (w) node {{\small $w_1$}};
\draw[anchor = south] (v2n) node {{\small $x_1$}};
%
\draw[anchor = north] (u1) node {{\small $v_1$}};
\draw[anchor = north] (u2) node {{\small $v_2$}};
\draw[anchor = south] (u3) node {{\small $v_3$}};
\draw[anchor = west] (w1) node {{\small $x_2$}};
\draw[anchor = west] (w2) node {{\small $x_3$}};
\end{tikzpicture}
\end{center}
\begin{center}
\vskip - 0.65 cm
\caption{Subgraphs of $G$ in the proof of Claim~\ref{claim:no-K33-minus-e}}
\label{fig:no-K33-minus-e}
\end{center}
\end{figure}
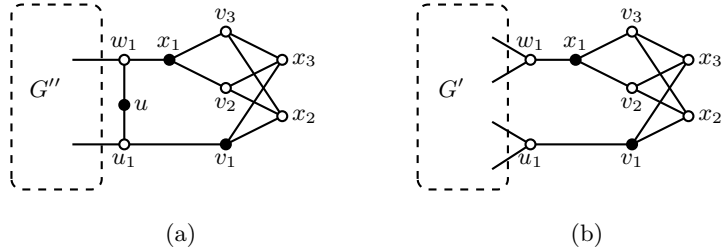

Suppose that $G'$ contains an isolated vertex, say~$u$. Thus, $N_G(u) = \{u_1,w_1\}$. By Claim~\ref{claim:min-cubic.2}, $\deg_G(u_1) = \deg_G(w_1) = 3$. By Claim~\ref{claim:min-cubic.4}, $u_1w_1 \notin E(G)$. In this case, we let $X' = X \cup \{u\}$ and consider the graph $G'' = G - X'$. By Claim~\ref{claim:min-cubic.3}, the neighbor of $u_1$ in $G''$ has degree~$3$ in $G$ and the neighbor of $w_1$ in $G''$ has degree~$3$ in $G$. Thus, $G''$ is isolate-free and there are two $X'$-exit edges. By Claims~\ref{claim:no-bad-in-Gee} and~\ref{structural-weight}, we infer that $\Theta(G'') = 0$, and so $\Omega(G'') = \w(G'') = \w(G) - \w_G(X') + 2 = \w(G) - 28 + 2 = \Omega(G) - 26$. Every $i$-set of $G''$ can be extended to an ID-set of $G$ by adding to it the vertices in the set $\{u,v_1,x_1\}$ (indicated by the shaded vertices in Figure~\ref{fig:no-K33-minus-e}(a)), and so $i(G) \le i(G'') + 3$. Thus, $8i(G) \le 8(i(G'') + 3) \le \Omega(G'') + 24 < \Omega(G)$, a contradiction.

Hence, $G'$ is isolate-free. Every $i$-set of $G'$ can be extended to an ID-set of $G$ by adding to it the vertices in the set $\{v_1,x_1\}$ (indicated by the shaded vertices in Figure~\ref{fig:no-K33-minus-e}(b)), and so $i(G) \le i(G') + 2$. We note that there are at most four $X$-exit edges. By Claims~\ref{claim:no-bad-in-Gee} and~\ref{structural-weight}, we infer that $b(G') \le 1$ and $\Theta(G') \le 2$, and so $\Omega(G') = \w(G') + \Theta(G') \le (\w(G) - \w_G(X) + 4) + 2 = \w(G) - 24 + 4 + 2 = \Omega(G) - 18$. Every $i$-set of $G'$ can be extended to an ID-set of $G$ by adding to it the vertices in the set $\{v_1,x_1\}$ (indicated by the shaded vertices in Figure~\ref{fig:no-K33-minus-e}(b)), and so $i(G) \le i(G') + 2$. Thus, $8i(G) \le 8(i(G') + 2) \le \Omega(G') + 16 < \Omega(G)$, a contradiction.~\smallqed

\begin{figure}[htb]
\begin{center}
\begin{tikzpicture}[scale=.75,style=thick,x=0.8cm,y=0.8cm]
\def\vr{2.5pt} 
\path (11.25,0.75) coordinate (q21);
\path (11.25,1.75) coordinate (q22);
\path (11.25,2.25) coordinate (p21);
\path (11.25,3.25) coordinate (p22);
\path (11.25,-0.75) coordinate (y21);
\path (11.25,0.25) coordinate (y22);
\path (12.5,-0.25) coordinate (x1);
\path (12.6,-0.25) coordinate (x1p);
%
\path (12.5,1.25) coordinate (v2);
\path (12.6,1.25) coordinate (v2p);
\path (12.5,2.75) coordinate (v3);
\path (14,-0.25) coordinate (u1);
\path (14,1.25) coordinate (u2);
\path (14,2.75) coordinate (u3);
\path (15.5,0.625) coordinate (w1);
\path (15.5,1.875) coordinate (w2);
\draw (x1)--(u1);
\draw (v2)--(u2);
\draw (v3)--(u3);
\draw (w1)--(u1)--(w2);
\draw (w1)--(u2)--(w2);
\draw (w1)--(u3)--(w2);
\draw (y21)--(x1)--(y22);
\draw (p21)--(v3)--(p22);
\draw (q21)--(v2)--(q22);
\draw (x1) [fill=white] circle (\vr);
\draw (v2) [fill=white] circle (\vr);
\draw (v3) [fill=white] circle (\vr);
\draw (u1) [fill=white] circle (\vr);
\draw (u2) [fill=white] circle (\vr);
\draw (u3) [fill=white] circle (\vr);
\draw (w1) [fill=white] circle (\vr);
\draw (w2) [fill=white] circle (\vr);
%
%
\draw[anchor = north] (x1p) node {{\small $u_3$}};
\draw[anchor = north] (v2p) node {{\small $u_2$}};
\draw[anchor = south] (v3) node {{\small $u_1$}};
\draw[anchor = north] (u1) node {{\small $v_3$}};
\draw[anchor = north] (u2) node {{\small $v_2$}};
\draw[anchor = south] (u3) node {{\small $v_1$}};
\draw[anchor = west] (w1) node {{\small $x_1$}};
\draw[anchor = west] (w2) node {{\small $x_2$}};
\draw (13.5,-2) node {{\small (a)}};
\draw [style=dashed,rounded corners] (11.85,-1) rectangle (13,4.25);
\draw (12.5,3.75) node {{\small $C$}};
\draw [style=dashed,rounded corners] (13.35,-1) rectangle (14.5,4.25);
\draw (14,3.75) node {{\small $B$}};
\draw [style=dashed,rounded corners] (15,-1) rectangle (16.5,4.25);
\draw (15.5,3.75) node {{\small $A$}};
%
\path (20.5,-0.25) coordinate (x1);
\path (20.6,-0.25) coordinate (x1p);
%
\path (20.5,1.25) coordinate (v2);
\path (20.6,1.25) coordinate (v2p);
\path (20.5,2.75) coordinate (v3);
\path (22,-0.25) coordinate (u1);
\path (22,1.25) coordinate (u2);
\path (22,2.75) coordinate (u3);
\path (23.5,0.625) coordinate (w1);
\path (23.5,1.875) coordinate (w2);
\draw (x1)--(u1);
\draw (v2)--(u2);
\draw (v3)--(u3);
\draw (w1)--(u1)--(w2);
\draw (w1)--(u2)--(w2);
\draw (w1)--(u3)--(w2);
\draw (x1) [fill=white] circle (\vr);
\draw (v2) [fill=white] circle (\vr);
\draw (v3) [fill=white] circle (\vr);
\draw (u1) [fill=white] circle (\vr);
\draw (u2) [fill=white] circle (\vr);
\draw (u3) [fill=white] circle (\vr);
\draw (w1) [fill=white] circle (\vr);
\draw (w2) [fill=white] circle (\vr);
%
%
\draw (21.5,-2) node {{\small (b) $G_{8.2}$}};
\end{tikzpicture}
\end{center}
\begin{center}
\vskip - 0.65 cm
\caption{Subgraphs of $G$ in the proof of Claim~\ref{claim:K23-config}}
\label{fig:K23-config}
\end{center}
\end{figure}
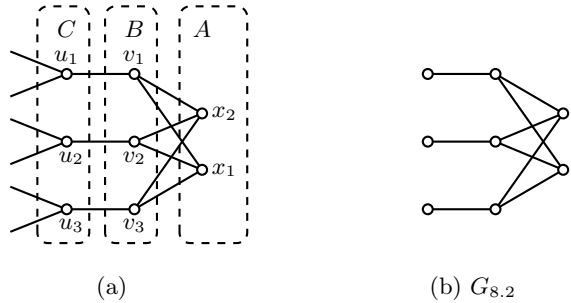

\begin{claim}
\label{claim:K23-config}
Every $K_{2,3}$-subgraph of $G$ belongs to a $G_{8.2}$-configuration that is an induced subgraph in $G$, where $G_{8.2}$ is the graph in Figure~\ref{fig:K23-config}(b).
\end{claim}
\proof Let $F$ be a $K_{2,3}$-subgraph in $G$ with partite sets $A = \{x_1,x_2\}$ and $B = \{v_1,v_2,v_3\}$. Since $G$ is diamond-free by Claim~\ref{no-diamond}, the set $B$ is an independent set in $G$. Let $u_i$ be the neighbor of $v_i$ different from $x_1$ and $x_2$ for $i \in [3]$. By Claim~\ref{claim:no-K33-minus-e}, the vertices $u_1$, $u_2$ and $u_3$ are distinct. Let $U = \{u_1,u_2,u_3\}$. Since $G_7$ is not a subgraph of $G$, the set $C$ is independent. Thus the subgraph $G[A \cup B \cup C]$ of $G$ induced by the sets $A \cup B \cup C$ (see Figure~\ref{fig:K23-config}(a)) is an induced subgraph of $G$.~\smallqed

\begin{claim}
\label{claim:G-cubic}
The graph $G$ is cubic.
\end{claim}
\proof Suppose, to the contrary, that $\delta(G) = 2$. Let $u$ be an arbitrary vertex of degree~$2$ in $G$. By Claim~\ref{claim:min-cubic.2}, both neighbors of~$u$ have degree~$3$ in $G$. Let $v$ be a neighbor of~$u$, and let $N_G(v) = \{u,x,y\}$. By Claim~\ref{claim:min-cubic.3}, $\deg_G(x) = \deg_G(y) = 3$. We now let $X = N_G[v]$ and we consider the graph $G' = G - X$.

\begin{subclaim}
\label{claim:G-cubic.1}
The graph $G'$ is isolate-free.
\end{subclaim}
\proof Suppose that $G'$ contains an isolated vertex, say~$w$. Thus, $N_G(w) \subseteq N_G(v)$. If $N_G(w) = N_G(v)$, then $G[\{u,v,w,x,y\}]$ is a copy of $K_{2,3}$ that contains a vertex of degree~$2$ in $G$, a contradiction. Hence, $w$ has degree~$2$ in $G$ and is adjacent to exactly two vertices in~$N_G(v)$. Since no two vertices of degree~$2$ in $G$ are adjacent by Claim~\ref{claim:min-cubic.2}, we have $N_G(w) = \{x,y\}$. Let $x_1$ and $y_1$ be the neighbors of $x$ and~$y$, respectively, different from~$v$ and~$w$. Since there is no $K_{2,3}$-subgraph that contains a vertex of degree~$2$ in $G$, we note that $x_1 \ne y_1$. By Claim~\ref{claim:min-cubic.3}, $\deg_G(x_1) = \deg_G(y_1) = 3$. Since $\deg_G(u) = 2$, the vertex~$u$ is adjacent to at most one of~$x_1$ and~$y_1$. By symmetry, we may assume that $u$ and $x_1$ are not adjacent.

We now let $X' = N_G[x] = \{v,w,x,x_1\}$ and consider the graph $G'' = G - X'$. We note that both vertices~$u$ and~$y$ have degree~$1$ in $G''$. If $u$ belongs to a bad component $B_u \in \cB$ in $G''$, then necessarily $B_u = B_1$ where $u$ is the vertex of degree~$1$ in~$B_u$ and where $B_u$ contains exactly one copy of $K_{2,3}$, which we denote by $F_u$. In this case, the vertex~$x_1$ is adjacent in $G$ to the two vertices of degree~$2$ in $F_u$. However, then $G[V(F_u) \cup \{x_1\}]$ is a copy of $K_{3,3}$ with an edge removed, contradicting Claim~\ref{claim:no-K33-minus-e}. Hence, $u$ does not belong to a bad component in $G''$. Analogously, $y$ does not belong to a bad component in $G''$. By Claim~\ref{claim:no-bad-in-Gee}, we therefore infer that removing the two $X'$-exit edges incident with~$x_1$ does not create a bad component. Hence, $b(G'') = 0$, and so by Claim~\ref{structural-weight} we have $\Theta(G'') = 0$, and so $\Omega(G'') = \w(G'') \le \w(G) - \w_G(X) + 5 = \w(G) - 13 + 5 = \Omega(G) - 8$. Every $i$-set of $G''$ can be extended to an ID-set of $G$ by adding to it the vertex~$x$, and so $i(G) \le i(G'') + 1$. Thus, $8i(G) \le 8(i(G'') + 1) \le \Omega(G'') + 8 = \Omega(G)$, a contradiction.~\smallqed

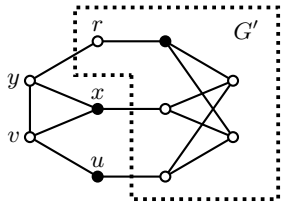
\begin{figure}[htb]
\begin{center}
\begin{tikzpicture}[scale=.75,style=thick,x=0.8cm,y=0.8cm]
\def\vr{2.5pt} 
\path (21.25,2) coordinate (r0);
\path (21.25,-0.75) coordinate (r1);
\path (24.5,-0.75) coordinate (r2);
\path (24.5,3.5) coordinate (r3);
\path (20,3.5) coordinate (r4);
\path (20,2) coordinate (r5);
\path (19,0.625) coordinate (z1);
\path (19,1.875) coordinate (z2);
\path (20.5,-0.25) coordinate (x1);
\path (20.6,-0.25) coordinate (x1p);
\path (20.5,1.25) coordinate (v2);
\path (20.6,1.25) coordinate (v2p);
\path (20.5,2.75) coordinate (v3);
\path (22,-0.25) coordinate (u1);
\path (22,1.25) coordinate (u2);
\path (22,2.75) coordinate (u3);
\path (23.5,0.625) coordinate (w1);
\path (23.5,1.875) coordinate (w2);
\draw (x1)--(u1);
\draw (v2)--(u2);
\draw (v3)--(u3);
\draw (w1)--(u1)--(w2);
\draw (w1)--(u2)--(w2);
\draw (w1)--(u3)--(w2);
\draw (z1)--(z2)--(v3);
\draw (x1)--(z1)--(v2);
\draw (z2)--(v2);
\draw[dotted, line width=0.05cm] (r0)--(r1)--(r2)--(r3)--(r4)--(r5)--(r0);
\draw (x1) [fill=black] circle (\vr);
\draw (v2) [fill=black] circle (\vr);
\draw (v3) [fill=white] circle (\vr);
\draw (u1) [fill=white] circle (\vr);
\draw (u2) [fill=white] circle (\vr);
\draw (u3) [fill=black] circle (\vr);
\draw (w1) [fill=white] circle (\vr);
\draw (w2) [fill=white] circle (\vr);
\draw (z1) [fill=white] circle (\vr);
\draw (z2) [fill=white] circle (\vr);
\draw[anchor = east] (z1) node {{\small $v$}};
\draw[anchor = east] (z2) node {{\small $y$}};
\draw[anchor = south] (v3) node {{\small $r$}};
\draw[anchor = south] (v2) node {{\small $x$}};
\draw[anchor = south] (x1) node {{\small $u$}};
\draw (23.8,3) node {{\small $G'$}};
\end{tikzpicture}
\end{center}
\begin{center}
\vskip - 0.35 cm
\caption{The graph $G$ in the proof of Claim~\ref{claim:G-cubic.1b}}
\label{fig:K23-config-h}
\end{center}
\end{figure}

\begin{subclaim}
\label{claim:G-cubic.1b}
$xy \notin E(G)$.
\end{subclaim}
\proof Suppose that $xy \in E(G)$. Thus there are exactly three $X$-exit edges. By Claim~\ref{claim:G-cubic.1}, the graph $G'$ is isolate-free. By Claim~\ref{claim:no-bad-in-Gee}, if there is a bad component in $G'$, then all three $X$-exit edges emanate from such a component and $G' =  B_1$. Let $r$ be the root vertex of $G'$. By Claim~\ref{claim:K23-config}, the $K_{2,3}$-subgraph in $G'$ belongs to a $G_{8.2}$-configuration (see Figure~\ref{fig:K23-config}(b)). The graph $G$ is therefore determined (up to isomorphism) and is illustrated in Figure~\ref{fig:K23-config-h}. In this case, $i(G) = 3$ (the shaded vertices in Figure~\ref{fig:K23-config-h} indicate an $i$-set of $G$) and $\w(G) = 32$, and so  $8i(G) < \Omega(G)$, a contradiction. Hence, $b(G') = 0$. By Claim~\ref{structural-weight}, we have $\Theta(G') = 0$, and so $\Omega(G') = \w(G') \le \w(G) - \w_G(X) + 3 = \w(G) - 13 + 3 = \Omega(G) - 10$. Every $i$-set of $G'$ can be extended to an ID-set of $G$ by adding to it the vertex~$v$, and so $i(G) \le i(G') + 1$. Thus, $8i(G) \le 8(i(G') + 1) \le \Omega(G') + 8 < \Omega(G)$, a contradiction.~\smallqed

By Claim~\ref{claim:G-cubic.1}, the graph $G'$ is isolate-free. By Claim~\ref{claim:G-cubic.1b}, $xy \notin E(G)$. Hence, there are at exactly five $X$-exit edges. By Claim~\ref{claim:no-bad-in-Gee}, if there is a bad component in $G'$, then at least three $X$-exit edges emanate from such a component, implying that $b(G') \le 1$. Every $i$-set of $G'$ can be extended to an ID-set of $G$ by adding to it the vertex~$v$, and so $i(G) \le i(G') + 1$. If $b(G') = 0$, then by Claim~\ref{structural-weight} we have $\Theta(G') = 0$, and so $\Omega(G') = \w(G') \le \w(G) - \w_G(X) + 5 = \w(G) - 13 + 5 = \Omega(G) - 8$. Thus, $8i(G) \le 8(i(G') + 1) \le \Omega(G') + 8 = \Omega(G)$, a contradiction. Hence, $b(G') = 1$. Let $B$ be the bad component in $G'$, and let $r$ be the root vertex of~$B$.

\begin{subclaim}
\label{claim:G-cubic.2}
$B \in \cB_1$.
\end{subclaim}
\proof Suppose that $B \in \cB_2 \cup \cB_3$. If $B \in \cB_3$, then at least six $X$-exit edges emanate from such a component, contradicting the fact that there are five $X$-exit edges. Hence, $B \in \cB_2$. Thus, the root vertex $r$ of $B$ has degree~$2$ in $B$ and there are either two or three $K_{2,3}$-subgraphs in $B$.
Suppose that there are three $K_{2,3}$-subgraphs in $B$. In this case the graph $G$ is determined and $V(G) = V(B) \cup X$. The core independent set of $B$ (see Section~\ref{S:familyB}) of cardinality~$3$ can be extended to an ID-set of $G$ by adding to it the set $N_G(v) = \{u,x,y\}$, and so $i(G) \le 6$. (For example, if $G$ is the graph illustrated in Figure~\ref{fig:G-cubic.2C}, then the ID-set of $G$ containing the  core independent set of $B$ together with the vertices in $\{u,x,y\}$ is indicated by the six shaded vertices.) Thus, $8i(G) \le 48 < 62 = \Omega(G)$, a contradiction.

\begin{figure}[htb]
\begin{center}
\begin{tikzpicture}[scale=.75,style=thick,x=0.75cm,y=0.75cm]
\def\vr{2.25pt}
\def\vrn{1.25pt}
\path (7.5,6.5) coordinate (h11);
\path (10,6.5) coordinate (p);
\path (12.5,6.5) coordinate (q51);
\path (10,5) coordinate (v);
\path (3.5,10) coordinate (a1);
\path (4.75,9) coordinate (a2);
\path (4.75,10) coordinate (a3);
\path (4.9,10) coordinate (a3p);
\path (4.75,11) coordinate (a4);
\path (6,10) coordinate (a5);
\path (7.5,10) coordinate (d1);
\path (8.75,9) coordinate (d2);
\path (8.75,10) coordinate (d3);
\path (8.9,10) coordinate (d3p);
\path (8.75,11) coordinate (d4);
\path (10,10) coordinate (d5);
\path (11.125,10) coordinate (c);
\path (12.5,10) coordinate (k1);
\path (13.75,9) coordinate (k2);
\path (13.75,10) coordinate (k3);
\path (13.7,10) coordinate (k3p);
\path (13.75,11) coordinate (k4);
\path (15.25,10) coordinate (k5);
\draw (d1)--(d2)--(d3)--(d4)--(d5)--(d2);
\draw (d1)--(d4);
\draw (a1)--(a2)--(a3)--(a4)--(a5)--(a2);
\draw (a1)--(a4);
\draw (a5)--(d1);
\draw (d5)--(c)--(k1);
\draw (v)--(p);
\draw (h11) to[out=180,in=-90, distance=0.75cm] (a1);
\draw (h11) to[out=90,in=180, distance=0.75cm] (d3);
\draw (p) to[out=90,in=180, distance=0.75cm] (k3);
\draw (p) to[out=90,in=0, distance=0.75cm] (a3);
\draw (k1)--(k2)--(k3)--(k4)--(k5)--(k2);
\draw (k1)--(k4);
\draw (q51) to[out=0,in=-90, distance=0.75cm] (k5);
\draw (v) to[out=180,in=-90, distance=0.5cm] (h11);
\draw (v) to[out=0,in=-90, distance=0.5cm] (q51);
\draw [style=dashed,rounded corners] (2.5,8) rectangle (16,12);
\draw (15,11.5) node {{\small $B$}};
\draw [style=dashed,rounded corners] (6.25,4.25) rectangle (13.25,7.25);
\draw (12.5,5) node {{\small $X$}};
\draw (a1) [fill=white] circle (\vr);
\draw (a2) [fill=white] circle (\vr);
\draw (a3) [fill=white] circle (\vr);
\draw (a4) [fill=white] circle (\vr);
\draw (a5) [fill=black] circle (\vr);
\draw (d1) [fill=white] circle (\vr);
\draw (d2) [fill=white] circle (\vr);
\draw (d3) [fill=white] circle (\vr);
\draw (d4) [fill=white] circle (\vr);
\draw (d5) [fill=black] circle (\vr);
\draw (k1) [fill=black] circle (\vr);
\draw (k2) [fill=white] circle (\vr);
\draw (k3) [fill=white] circle (\vr);
\draw (k4) [fill=white] circle (\vr);
\draw (k5) [fill=white] circle (\vr);
\draw (q51) [fill=black] circle (\vr);
\draw (h11) [fill=black] circle (\vr);
\draw (p) [fill=black] circle (\vr);
\draw (v) [fill=white] circle (\vr);
\draw (c) [fill=white] circle (\vr);
\draw[anchor = west] (h11) node {{\small $x$}};
\draw[anchor = east] (q51) node {{\small $u$}};
\draw[anchor = west] (p) node {{\small $y$}};
\draw[anchor = north] (v) node {{\small $v$}};
\draw[anchor = south] (c) node {{\small $r$}};
\end{tikzpicture}
\caption{An example of a graph $G$ in the proof of Claim~\ref{claim:G-cubic.2}}
\label{fig:G-cubic.2C}
\end{center}
\end{figure}

Hence, $B$ contains exactly two $K_{2,3}$-subgraphs. If $V(G) = V(B) \cup X$, then the graph $G$ is determined and as before the core independent set (of cardinality~$2$) of $B$ can be extended to an ID-set of $G$ by adding to it the set $N_G(v) = \{u,x,y\}$, and so $i(G) \le 5$. Moreover, $\Omega(G) = 46$, and so $8i(G) < \Omega(G)$, a contradiction. Hence, $V(G) \ne V(B) \cup X$, implying that there are exactly four $X$-exit edges that are incident with vertices in $B$ and the root vertex~$r$ of $B$ is not incident with any of these four $X$-exit edges.

Suppose the vertex~$u$ is not adjacent to a vertex of~$B$ in $G$. We now let $X^* = V(B) \cup \{v,x,y\}$ and consider the graph $G^* = G - X^*$. In this case, there is exactly one $X^*$-exit edge and the vertex~$u$ has degree~$1$ in $G^*$. By Claims~\ref{claim:no-bad-in-Ge} and~\ref{structural-weight}, we infer that $b(G^*) = 0$ and $\Theta(G^*) = 0$, and so $\Omega(G^*) = \w(G^*) = \w(G) - \w_G(X^*) + 1 = \w(G) - 43 + 1 = \Omega(G) - 42$. Every $i$-set of $G^*$ can be extended to an ID-set of $G$ by adding to it the core independent set of $B$ of cardinality~$2$ and the vertices $x$ and~$y$ (as indicated by the four shaded vertices in Figure~\ref{fig:G-cubic.2D}), and so $i(G^*) \le i(G) + 4$. Thus, $8i(G) \le 8(i(G^*) + 4) \le \Omega(G^*) + 32 < \Omega(G)$, a contradiction.

\begin{figure}[htb]
\begin{center}
\begin{tikzpicture}[scale=.75,style=thick,x=0.75cm,y=0.75cm]
\def\vr{2.25pt}
\def\vrn{1.25pt}
\path (8,6.5) coordinate (h11);
\path (10.5,6.5) coordinate (p);
\path (15,5) coordinate (q);
\path (16,5) coordinate (q1);
\path (16.75,4.5) coordinate (q11);
\path (16.75,5.5) coordinate (q12);
\path (9.25,5.5) coordinate (v);
\path (5.5,10) coordinate (d1);
\path (6.75,9) coordinate (d2);
\path (6.75,10) coordinate (d3);
\path (6.9,10) coordinate (d3p);
\path (6.75,11) coordinate (d4);
\path (8,10) coordinate (d5);
\path (9.125,10) coordinate (c);
\path (10.5,10) coordinate (k1);
\path (11.75,9) coordinate (k2);
\path (11.75,10) coordinate (k3);
\path (11.7,10) coordinate (k3p);
\path (11.75,11) coordinate (k4);
\path (13,10) coordinate (k5);
\draw (d1)--(d2)--(d3)--(d4)--(d5)--(d2);
\draw (d1)--(d4);
\draw (d5)--(c)--(k1);
\draw (v)--(p);
\draw (v)--(h11);
\draw (q)--(q1);
\draw (q11)--(q1)--(q12);
\draw (h11) to[out=180,in=-90, distance=0.75cm] (d1);
\draw (h11) to[out=90,in=180, distance=0.75cm] (k3);
\draw (p) to[out=180,in=-90, distance=0.75cm] (d3);
\draw (p) to[out=0,in=-90, distance=0.75cm] (k5);
\draw (v) to[out=-45,in=180, distance=0.75cm] (q);
\draw (k1)--(k2)--(k3)--(k4)--(k5)--(k2);
\draw (k1)--(k4);
%
%
\draw [style=dashed,rounded corners] (5,8) rectangle (13.5,11.5);
\draw (5.5,11) node {{\small $B$}};
\draw [style=dashed,rounded corners] (4.5,4.25) rectangle (14,11.75);
\draw (5.25,5) node {{\small $X^*$}};
\draw [style=dashed,rounded corners] (14.5,4) rectangle (18,6);
\draw (17.5,5) node {{\small $G^*$}};
\draw (d1) [fill=white] circle (\vr);
\draw (d2) [fill=white] circle (\vr);
\draw (d3) [fill=white] circle (\vr);
\draw (d4) [fill=white] circle (\vr);
\draw (d5) [fill=black] circle (\vr);
\draw (k1) [fill=black] circle (\vr);
\draw (k2) [fill=white] circle (\vr);
\draw (k3) [fill=white] circle (\vr);
\draw (k4) [fill=white] circle (\vr);
\draw (k5) [fill=white] circle (\vr);
\draw (q) [fill=white] circle (\vr);
\draw (q1) [fill=white] circle (\vr);
\draw (h11) [fill=black] circle (\vr);
\draw (p) [fill=black] circle (\vr);
\draw (v) [fill=white] circle (\vr);
\draw (c) [fill=white] circle (\vr);
\draw[anchor = north] (h11) node {{\small $x$}};
\draw[anchor = north] (q) node {{\small $u$}};
\draw[anchor = north] (p) node {{\small $y$}};
\draw[anchor = north] (v) node {{\small $v$}};
\draw[anchor = south] (c) node {{\small $r$}};
\end{tikzpicture}
\caption{An example of a graph $G$ in the proof of Claim~\ref{claim:G-cubic.2}}
\label{fig:G-cubic.2D}
\end{center}
\end{figure}

Hence, the vertex~$u$ is adjacent to a vertex of~$B$ in $G$. Renaming $x$ and $y$ if necessary, we may assume that~$x$ is adjacent to exactly one vertex, say~$x^*$, of~$B$ in $G$. In this case, we let $X^* = V(B) \cup X$ and consider the graph $G^* = G - X^*$. Since there is exactly one $X^*$-exit edge and the graph $G^*$ is isolate-free, by Claims~\ref{claim:no-bad-in-Ge} and~\ref{structural-weight} we infer that $b(G^*) = 0$ and $\Theta(G^*) = 0$, and so $\Omega(G^*) = \w(G^*) = \w(G) - \w_G(X^*) + 1 = \w(G) - 47 + 1 = \Omega(G) - 46$. Every $i$-set of $G'$ can be extended to an ID-set of $G$ by adding to it the core independent set of $B$ of cardinality~$2$, together with the vertices in the set $\{u,x^*,y\}$,
and so $i(G^*) \le i(G) + 5$. Thus, $8i(G) \le 8(i(G^*) + 5) \le \Omega(G^*) + 40 < \Omega(G)$, a contradiction. Since both cases produce a contradiction, this completes the proof of Claim~\ref{claim:G-cubic.2}.~\smallqed

By Claim~\ref{claim:G-cubic.2}, we have $B \in \cB_1$. Thus, the root vertex~$r$ in $B$ has degree~$1$ in $B$. Let $B$ have $k$ copies of $K_{2,3}$. We note that $k \le 3$.

\begin{subclaim}
\label{claim:G-cubic.3}
$k=1$.
\end{subclaim}
\proof Suppose that $k \in \{2,3\}$. If $k=3$, then the graph $G$ is determined and $V(G) = V(B) \cup X$. In this case, the core independent set of $B$ (see Section~\ref{S:familyB}) of cardinality~$3$ can be extended to an ID-set of $G$ by adding to it the set $N_G(v) = \{u,x,y\}$, and so $i(G) \le 6$. Thus, $8i(G) \le 48 < 62 = \Omega(G)$, a contradiction. Hence, $k=2$.

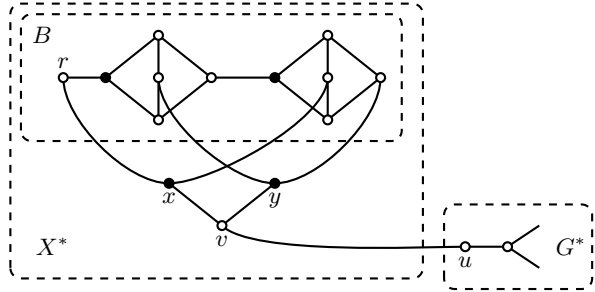
\begin{figure}[htb]
\begin{center}
\begin{tikzpicture}[scale=.75,style=thick,x=0.75cm,y=0.75cm]
\def\vr{2.25pt}
\def\vrn{1.25pt}
\path (8,7.5) coordinate (h11);
\path (10.5,7.5) coordinate (p);
\path (15,6) coordinate (q);
\path (16,6) coordinate (q1);
\path (16.75,5.5) coordinate (q11);
\path (16.75,6.5) coordinate (q12);
\path (9.25,6.5) coordinate (v);
\path (6.5,10) coordinate (d1);
\path (7.75,9) coordinate (d2);
\path (7.75,10) coordinate (d3);
\path (7.9,10) coordinate (d3p);
\path (7.75,11) coordinate (d4);
\path (9,10) coordinate (d5);
\path (5.5,10) coordinate (c);
\path (10.5,10) coordinate (k1);
\path (11.75,9) coordinate (k2);
\path (11.75,10) coordinate (k3);
\path (11.7,10) coordinate (k3p);
\path (11.75,11) coordinate (k4);
\path (13,10) coordinate (k5);
\draw (d1)--(d2)--(d3)--(d4)--(d5)--(d2);
\draw (d1)--(d4);
\draw (d1)--(c);
\draw (d5)--(k1);
\draw (v)--(p);
\draw (v)--(h11);
\draw (q)--(q1);
\draw (q11)--(q1)--(q12);
\draw (h11) to[out=180,in=-90, distance=0.75cm] (c);
\draw (h11) to[out=0,in=-90, distance=0.75cm] (k3);
\draw (p) to[out=180,in=-90, distance=0.75cm] (d3);
\draw (p) to[out=0,in=-90, distance=0.75cm] (k5);
\draw (v) to[out=-45,in=180, distance=0.75cm] (q);
\draw (k1)--(k2)--(k3)--(k4)--(k5)--(k2);
\draw (k1)--(k4);
%
%
\draw [style=dashed,rounded corners] (4.5,8.5) rectangle (13.5,11.5);
\draw (5,11) node {{\small $B$}};
\draw [style=dashed,rounded corners] (4.25,5.25) rectangle (14,11.75);
\draw (5.25,6) node {{\small $X^*$}};
\draw [style=dashed,rounded corners] (14.5,5) rectangle (18,7);
\draw (17.5,6) node {{\small $G^*$}};
\draw (d1) [fill=black] circle (\vr);
\draw (d2) [fill=white] circle (\vr);
\draw (d3) [fill=white] circle (\vr);
\draw (d4) [fill=white] circle (\vr);
\draw (d5) [fill=white] circle (\vr);
\draw (k1) [fill=black] circle (\vr);
\draw (k2) [fill=white] circle (\vr);
\draw (k3) [fill=white] circle (\vr);
\draw (k4) [fill=white] circle (\vr);
\draw (k5) [fill=white] circle (\vr);
\draw (q) [fill=white] circle (\vr);
\draw (q1) [fill=white] circle (\vr);
\draw (h11) [fill=black] circle (\vr);
\draw (p) [fill=black] circle (\vr);
\draw (v) [fill=white] circle (\vr);
\draw (c) [fill=white] circle (\vr);
\draw[anchor = north] (h11) node {{\small $x$}};
\draw[anchor = north] (q) node {{\small $u$}};
\draw[anchor = north] (p) node {{\small $y$}};
\draw[anchor = north] (v) node {{\small $v$}};
\draw[anchor = south] (c) node {{\small $r$}};
\end{tikzpicture}
\caption{An example of a graph $G$ in the proof of Claim~\ref{claim:G-cubic.3}}
\label{fig:G-cubic.2F}
\end{center}
\end{figure}

Suppose the vertex~$u$ is not adjacent to a vertex of~$B$ in $G$. We now let $X^* = V(B) \cup \{v,x,y\}$ and consider the graph $G^* = G - X^*$. In this case, there is exactly one $X^*$-exit edge and the vertex~$u$ has degree~$1$ in $G^*$. By Claims~\ref{claim:no-bad-in-Ge} and~\ref{structural-weight}, we infer that $b(G^*) = 0$ and $\Theta(G^*) = 0$, and so $\Omega(G^*) = \w(G^*) = \w(G) - \w_G(X^*) + 1 = \w(G) - 43 + 1 = \Omega(G) - 42$. Every $i$-set of $G'$ can be extended to an ID-set of $G$ by adding to it the core independent set of $B$ of cardinality~$2$ and the vertices $x$ and~$y$ (as indicated by the four shaded vertices in Figure~\ref{fig:G-cubic.2F}), and so $i(G^*) \le i(G) + 4$. Thus, $8i(G) \le 8(i(G^*) + 4) \le \Omega(G^*) + 32 < \Omega(G)$, a contradiction.

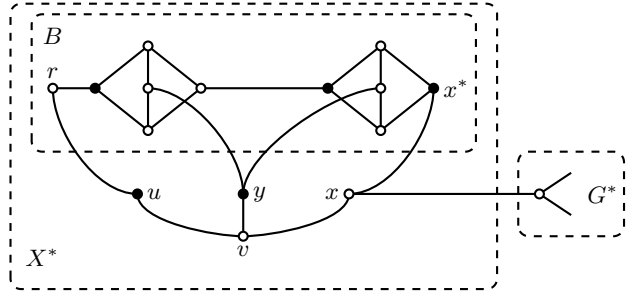
\begin{figure}[htb]
\begin{center}
\begin{tikzpicture}[scale=.75,style=thick,x=0.75cm,y=0.75cm]
\def\vr{2.25pt}
\def\vrn{1.25pt}
\path (7.5,7.5) coordinate (h11);
\path (10,7.5) coordinate (p);
\path (12.5,7.5) coordinate (q);
\path (17,7.5) coordinate (q1);
\path (17.75,8) coordinate (q11);
\path (17.75,7) coordinate (q12);
\path (10,6.5) coordinate (v);
\path (6.5,10) coordinate (a1);
\path (7.75,9) coordinate (a2);
\path (7.75,10) coordinate (a3);
\path (7.9,10) coordinate (a3p);
\path (7.75,11) coordinate (a4);
\path (9,10) coordinate (a5);
\path (5.5,10) coordinate (c);
\path (12,10) coordinate (k1);
\path (13.25,9) coordinate (k2);
\path (13.25,10) coordinate (k3);
\path (13.2,10) coordinate (k3p);
\path (13.25,11) coordinate (k4);
\path (14.5,10) coordinate (k5);
\draw (a1)--(a2)--(a3)--(a4)--(a5)--(a2);
\draw (a1)--(a4);
\draw (d1)--(c);
\draw (d5)--(k1);
\draw (v)--(p);
\draw (h11) to[out=180,in=-90, distance=0.75cm] (c);
\draw (p) to[out=90,in=180, distance=0.75cm] (k3);
\draw (p) to[out=90,in=0, distance=0.75cm] (a3);
\draw (q) to[out=0,in=-90, distance=0.75cm] (k5);
\draw (k1)--(k2)--(k3)--(k4)--(k5)--(k2);
\draw (k1)--(k4);
\draw (q)--(q1);
\draw (q11)--(q1)--(q12);
\draw (v) to[out=180,in=-90, distance=0.5cm] (h11);
\draw (v) to[out=0,in=-90, distance=0.5cm] (q);
\draw [style=dashed,rounded corners] (5,8.5) rectangle (15.5,11.75);
\draw (5.5,11.25) node {{\small $B$}};
\draw [style=dashed,rounded corners] (4.5,5.25) rectangle (16,12);
\draw (5.25,6) node {{\small $X^*$}};
\draw [style=dashed,rounded corners] (16.5,6.5) rectangle (19,8.5);
\draw (18.5,7.5) node {{\small $G^*$}};
\draw (a1) [fill=black] circle (\vr);
\draw (a2) [fill=white] circle (\vr);
\draw (a3) [fill=white] circle (\vr);
\draw (a4) [fill=white] circle (\vr);
\draw (a5) [fill=white] circle (\vr);
\draw (k1) [fill=black] circle (\vr);
\draw (k2) [fill=white] circle (\vr);
\draw (k3) [fill=white] circle (\vr);
\draw (k4) [fill=white] circle (\vr);
\draw (k5) [fill=black] circle (\vr);
\draw (q) [fill=white] circle (\vr);
\draw (q1) [fill=white] circle (\vr);
\draw (h11) [fill=black] circle (\vr);
\draw (p) [fill=black] circle (\vr);
\draw (v) [fill=white] circle (\vr);
\draw (c) [fill=white] circle (\vr);
\draw[anchor = west] (h11) node {{\small $u$}};
\draw[anchor = east] (q) node {{\small $x$}};
\draw[anchor = west] (p) node {{\small $y$}};
\draw[anchor = west] (k5) node {{\small $x^*$}};
\draw[anchor = north] (v) node {{\small $v$}};
\draw[anchor = south] (c) node {{\small $r$}};
\end{tikzpicture}
\caption{An example of a graph $G$ in the proof of Claim~\ref{claim:G-cubic.3}}
\label{fig:G-cubic.2G}
\end{center}
\end{figure}

Hence, the vertex~$u$ is adjacent to a vertex of~$B$ in $G$. Renaming $x$ and $y$ if necessary, we may assume that~$x$ is adjacent to exactly one vertex, say~$x^*$, of~$B$ in $G$. In this case, we let $X^* = V(B) \cup X$ and consider the graph $G^* = G - X^*$. Since there is exactly one $X^*$-exit edge and the graph $G^*$ is isolate-free, by Claims~\ref{claim:no-bad-in-Ge} and~\ref{structural-weight} we infer that $b(G^*) = 0$ and $\Theta(G^*) = 0$, and so $\Omega(G^*) = \w(G^*) = \w(G) - \w_G(X^*) + 1 = \w(G) - 47 + 1 = \Omega(G) - 46$. Every $i$-set of $G'$ can be extended to an ID-set of $G$ by adding to it the core independent set of $B$ of cardinality~$2$, together with the vertices in the set $\{u,x^*,y\}$ (as indicated by the five shaded vertices in Figure~\ref{fig:G-cubic.2G}), and so $i(G^*) \le i(G) + 5$. Thus, $8i(G) \le 8(i(G^*) + 5) \le \Omega(G^*) + 40 < \Omega(G)$, a contradiction. Since both cases produce a contradiction, this completes the proof of Claim~\ref{claim:G-cubic.2}.~\smallqed

\medskip
By Claim~\ref{claim:G-cubic.3}, we have $k=1$, and so $B$ contains exactly one copy of $K_{2,3}$. Let the copy of $K_{2,3}$ in $B$ have partite sets $\{x_1,x_1\}$ and $\{v_1,v_2,v_3\}$. Renaming vertices if necessary, we may assume that the root vertex of $B$ is adjacent to~$v_3$. By Claim~\ref{claim:K23-config}, every $K_{2,3}$-subgraph of $G$ belongs to a $G_{8.2}$-configuration that is an induced subgraph in $G$, where $G_{8.2}$ is the graph in Figure~\ref{fig:K23-config}(b). Let $w_1$ and $w_2$ be the neighbors of~$v_1$ and~$v_2$, respectively, not in $B$. We note by  Claim~\ref{claim:K23-config} that $\{r,w_1,w_2\}$ is an independent set.  Let $w_3$ be the third neighbor of~$v$, and so the root vertex~$r$ is adjacent to~$w_3$. Thus, $\deg_G(r) = 2$ and $N_G(r) = \{v_3,w_3\}$. We note that $N_G(v) = \{w_1,w_2,w_3\} = \{u,x,y\}$, where recall that $\deg_G(u) = 2$ and $\deg_G(x) = \deg_G(y) = 3$. Since there are no two adjacent vertices of degree~$2$ in $G$ by Claim~\ref{claim:min-cubic.2}, we infer that $\deg_G(w_3) = 3$, and so $u \in \{w_1,w_2\}$. Renaming vertices if necessary, we may assume that $u = w_2$, and so $\deg_G(w_2) = 2$ and $N_G(w_2) = \{v,v_2\}$.

\begin{figure}[htb]
\begin{center}
\begin{tikzpicture}[scale=.75,style=thick,x=0.8cm,y=0.8cm]
\def\vr{2.5pt} 
\path (11.25,0.75) coordinate (q21);
\path (11.25,1.75) coordinate (q22);
\path (11.25,2.25) coordinate (p21);
\path (11.25,3.25) coordinate (p22);
\path (11.25,-0.75) coordinate (y21);
\path (11.25,0.25) coordinate (y22);
\path (12.75,-0.25) coordinate (x1);
\path (9.5,1.25) coordinate (v);
\path (11,-0.25) coordinate (v1);
\path (11,1.25) coordinate (v2);
\path (11.1,1.25) coordinate (v2p);
\path (11,2.75) coordinate (v3);
\path (8,2.75) coordinate (z3);
\path (7,2.25) coordinate (z31);
\path (7,3.25) coordinate (z32);
\path (8,-0.25) coordinate (z1);
\path (7,-0.75) coordinate (z11);
\path (7,0.25) coordinate (z12);
\path (14,-0.25) coordinate (u1);
\path (14,1.25) coordinate (u2);
\path (14,2.75) coordinate (u3);
\path (15.75,0.625) coordinate (w1);
\path (15.75,1.875) coordinate (w2);
\draw (v1)--(v)--(v2);
\draw (v)--(v3);
\draw (v1)--(x1)--(u1);
\draw (v2)--(u2);
\draw (v3)--(u3);
\draw (w1)--(u1)--(w2);
\draw (w1)--(u2)--(w2);
\draw (w1)--(u3)--(w2);
\draw (z1)--(v1);
\draw (z3)--(v3);
\draw (z11)--(z1)--(z12);
\draw (z31)--(z3)--(z32);
\draw (x1) [fill=black] circle (\vr);
\draw (v1) [fill=white] circle (\vr);
\draw (v2) [fill=black] circle (\vr);
\draw (v3) [fill=white] circle (\vr);
\draw (u1) [fill=white] circle (\vr);
\draw (u2) [fill=white] circle (\vr);
\draw (u3) [fill=black] circle (\vr);
\draw (w1) [fill=white] circle (\vr);
\draw (w2) [fill=white] circle (\vr);
\draw (z1) [fill=white] circle (\vr);
\draw (z3) [fill=white] circle (\vr);
\draw (v) [fill=white] circle (\vr);
\draw [style=dashed,rounded corners] (12.25,-1.25) rectangle (17,3.5);
\draw (16.5,3) node {{\small $B$}};
\draw [style=dashed,rounded corners] (9,-1.5) rectangle (17.5,3.75);
\draw (9.5,-1) node {{\small $X''$}};
\draw [style=dashed,rounded corners] (5.5,-1.5) rectangle (8.5,3.75);
\draw (6.5,1.25) node {{\small $G''$}};
%
\draw[anchor = north] (x1) node {{\small $r$}};
\draw[anchor = north] (v) node {{\small $v$}};
\draw[anchor = north] (v1) node {{\small $w_3$}};
\draw[anchor = north] (v2) node {{\small $w_2$}};
\draw[anchor = south] (v3) node {{\small $w_1$}};
\draw[anchor = north] (u1) node {{\small $v_3$}};
\draw[anchor = north] (u2) node {{\small $v_2$}};
\draw[anchor = south] (u3) node {{\small $v_1$}};
\draw[anchor = west] (w1) node {{\small $x_1$}};
\draw[anchor = west] (w2) node {{\small $x_2$}};
%
%
\end{tikzpicture}
\end{center}
\begin{center}
\vskip - 0.35 cm
\caption{A subgraph of $G$ in the proof of Claim~\ref{claim:G-cubic}}
\label{fig:K23-configN}
\end{center}
\end{figure}
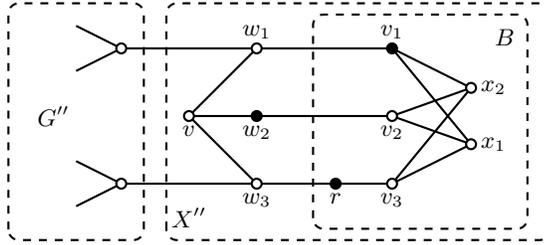

We now let $X'' = N_G[v] \cup V(B)$ and we consider the graph $G'' = G - X''$. Suppose that $G''$ contains an isolated vertex. In this case, such a vertex, say~$z$, has degree~$2$ in $G$ and is adjacent to~$w_1$ and $w_3$, implying that the vertex~$w_3$ has two neighbors of degree~$2$, namely~$r$ and~$z$, contradicting Claim~\ref{claim:min-cubic.3}. Hence, $G''$ is isolate-free. Thus, the graph shown in Figure~\ref{fig:K23-configN} is a subgraph of $G$, where in this illustration the neighbors of $w_1$ and $w_2$ not in $G''$ are distinct. By Claims~\ref{claim:no-bad-in-Gee} and~\ref{structural-weight} we infer that $b(G'') = 0$ and $\Theta(G'') = 0$, and so $\Omega(G'') = \w(G'') \le \w(G) - \w_G(X'') + 2 = \w(G) - 32 + 2 = \Omega(G) - 30$. Every $i$-set of $G'$ can be extended to an ID-set of $G$ by adding to it the vertices in the set $\{r,v_1,w_2\}$ (as indicated by the three shaded vertices in Figure~\ref{fig:K23-configN}), and so $i(G'') \le i(G) + 3$. Thus, $8i(G) \le 8(i(G'') + 3) \le \Omega(G') + 24 < \Omega(G)$, a contradiction. This completes the proof of Claim~\ref{claim:G-cubic}.~\smallqed

\subsection{Part 4: The graph $G$ is cubic}
\label{S:part3}

By Claim~\ref{claim:G-cubic}, the graph $G$ is a cubic graph of order~$n$. Recall that $n \ge 6$. By Claim~\ref{claim:Gconn}, the graph $G$ is connected. By supposition, $G \ne K_{3,3}$ and $G \ne C_5 \, \Box \, K_2$. Since $G$ is a cubic graph, we note that $\Theta(G) = 0$, and so $\Omega(G) = \w(G)$. We prove a number of claims that will culminate in a contradiction, thereby proving our main theorem.

\begin{claim}
\label{no-triangle}
The graph $G$ contains no triangle.
\end{claim}
\proof Suppose, to the contrary, that $G$ contains a triangle $T$, where $V(T) = \{v_1,v_2,v_3\}$. Let $u_i$ be the neighbor of $v_i$ not in $T$ for $i \in [3]$. By Claim~\ref{no-diamond}, the graph $G$ is diamond-free, implying that the vertices $u_1$, $u_2$ and $u_3$ are distinct. Let $U = \{u_1,u_2,u_3\}$. If $G[U] = K_3$, then $G$ is the $3$-prism $C_3 \, \Box \, K_2$. In this case, $i(G) = 2$ and $\w(G) = 18$, and so $8i(G) < \Omega(G)$, a contradiction. Hence renaming vertices if necessary, we may assume that $u_1u_2 \notin E(G)$. Let $G'$ be the graph obtained from $G - V(T)$ by adding the edge~$u_1u_2$. We note that every vertex of $G'$ has degree~$3$, except for the vertex~$u_3$ which has degree~$2$ in $G'$. We therefore infer that $b(G') = 0$.

Suppose that $\tc(G') = 0$. In this case, $\Theta(G') = 0$ and $\Omega(G') = \w(G') = \w(G) - 9 + 1 = \Omega(G) - 8$. Let $I'$ be an $i$-set of $G'$. If neither $u_1$ nor $u_2$ belong to~$I'$, then let $I = I' \cup \{v_1\}$. If $u_i \in I'$ for some $i \in [2]$, then let $I = I' \cup \{v_{3-i}\}$. In both cases, $|I| = |I'| + 1 = i(G') + 1$ and the set $I$ is an ID-set of $G$, and so $i(G) \le i(G') + 1$. Hence, $8i(G) \le 8(i(G') + 1) \le \Omega(G') + 8 = \Omega(G)$, a contradiction.

Hence, $\tc(G') = 1$, implying that there is a troublesome configuration, say $T'$, in $G'$ where $T'$ has exactly one copy of $K_{2,3}$, and this copy contains the vertex~$u_3$ as the vertex of degree~$2$ in $T'$, and where $T'$ contains the added edge $u_1u_2$
We now let $X^* = V(T) \cup \{u_3\}$ and $G^* = G - X^*$. In this case, $\tc(G^*) = 0$ and $b(G^*) = 0$, and so $\Theta(G^*) = 0$. Thus, $\Omega(G^*) = \w(G^*) = \w(G) - \w_G(X^*) + 4 = \Omega(G) - 12 + 4 = \Omega(G) - 8$. Every $i$-set of $G^*$ can be extended to an ID-set of $G$ by adding to it the vertex~$v_3$, and so $i(G) \le i(G^*) + 1$. Hence, $8i(G) \le 8(i(G^*) + 1) \le \Omega(G^*) + 8 = \Omega(G)$, a contradiction.

\begin{claim}
\label{claim:no-adj-K23}
The graph $G$ does not contain two vertex disjoint copies of $K_{2,3}$ joined by at least one edge.
\end{claim}
\proof Let $F_1$ and $F_2$ be two vertex disjoint copies of $K_{2,3}$, where $F_1$ has partite sets $X = \{x_1,x_2\}$ and $V_1 = \{v_1,v_2,v_3\}$, and $F_2$ partite sets $Y = \{y_1,y_2\}$ and $U_1 = \{u_1,u_2,u_3\}$. Suppose, to the contrary, that $F_1$ and $F_2$ are joined by at least one edge.

\begin{subclaim}
\label{claim:no-adj-K23.1}
$F_1$ and $F_2$ are joined by exactly one edge.
\end{subclaim}
\proof
If $F_1$ and $F_2$ are joined by three edges, then $G = G_{10.2}$, where $G_{10.2}$ is the graph illustrated in Figure~\ref{fig:G7-fig9b}(a). In this case, $i(G) = 3$ (an $i$-set is indicated by the shaded vertices in Figure~\ref{fig:G7-fig9b}(a)) and $\w(G) = 30$, and so $8i(G) < \w(G)$, a contradiction.

\begin{figure}[htb]
\begin{center}
\begin{tikzpicture}[scale=.75,style=thick,x=0.8cm,y=0.8cm]
\def\vr{2.5pt} 
%
\path (11.5,-0.25) coordinate (x1);
\path (11.6,-0.25) coordinate (x1p);
%
\path (11.5,1.25) coordinate (v2);
\path (11.6,1.25) coordinate (v2p);
\path (11.5,2.75) coordinate (v3);
\path (12.75,-0.25) coordinate (u1);
\path (12.75,1.25) coordinate (u2);
\path (12.75,2.75) coordinate (u3);
\path (14,0.625) coordinate (w1);
\path (14,1.875) coordinate (w2);
\path (10.25,0.625) coordinate (z1);
\path (10.25,1.875) coordinate (z2);
\draw (x1)--(u1);
\draw (v2)--(u2);
\draw (v3)--(u3);
\draw (w1)--(u1)--(w2);
\draw (w1)--(u2)--(w2);
\draw (w1)--(u3)--(w2);
\draw (z1)--(x1)--(z2);
\draw (z1)--(v2)--(z2);
\draw (z1)--(v3)--(z2);
\draw (x1) [fill=white] circle (\vr);
\draw (v2) [fill=white] circle (\vr);
\draw (v3) [fill=black] circle (\vr);
\draw (u1) [fill=black] circle (\vr);
\draw (u2) [fill=black] circle (\vr);
\draw (u3) [fill=white] circle (\vr);
\draw (w1) [fill=white] circle (\vr);
\draw (w2) [fill=white] circle (\vr);
\draw (z1) [fill=white] circle (\vr);
\draw (z2) [fill=white] circle (\vr);
%
%
\draw (12.25,-1.5) node {{\small (a) $G_{10.2}$}};
%
%
\path (22,-0.25) coordinate (x1);
\path (20.25,-0.25) coordinate (x11);
\path (22.1,-0.25) coordinate (x1p);
%
\path (22,1.25) coordinate (v2);
\path (22.2,1.25) coordinate (v2p);
\path (20.25,1.25) coordinate (v21);
\path (22,2.75) coordinate (v3);
\path (20.25,2.75) coordinate (v31);
\path (23.25,-0.25) coordinate (u1);
\path (23.25,1.25) coordinate (u2);
\path (23.25,2.75) coordinate (u3);
\path (24.5,0.625) coordinate (w1);
\path (24.5,1.875) coordinate (w2);
%
\path (21.25,2) coordinate (z2);
\path (21,2) coordinate (z2p);
\draw (x1)--(u1);
\draw (v2)--(u2);
\draw (v3)--(u3);
\draw (w1)--(u1)--(w2);
\draw (w1)--(u2)--(w2);
\draw (w1)--(u3)--(w2);
\draw (x1)--(z2);
\draw (v2)--(z2);
\draw (v3)--(z2);
\draw (x1)--(x11);
\draw (v2)--(v21);
\draw (v3)--(v31);
\draw (x1) [fill=white] circle (\vr);
\draw (v2) [fill=white] circle (\vr);
\draw (v3) [fill=white] circle (\vr);
\draw (u1) [fill=white] circle (\vr);
\draw (u2) [fill=white] circle (\vr);
\draw (u3) [fill=white] circle (\vr);
\draw (w1) [fill=black] circle (\vr);
\draw (w2) [fill=black] circle (\vr);
%
\draw (z2) [fill=black] circle (\vr);
%
%
\draw[anchor = north] (x1p) node {{\small $u_3$}};
\draw[anchor = north] (v2p) node {{\small $u_2$}};
\draw[anchor = south] (v3) node {{\small $u_1$}};
\draw[anchor = north] (u1) node {{\small $v_3$}};
\draw[anchor = north] (u2) node {{\small $v_2$}};
\draw[anchor = south] (u3) node {{\small $v_1$}};
\draw[anchor = west] (w1) node {{\small $x_1$}};
\draw[anchor = west] (w2) node {{\small $x_2$}};
%
\draw[anchor = west] (z2) node {{\small $w$}};
\draw (22.75,-1.5) node {{\small (b)}};
\draw [style=dashed,rounded corners] (19.25,-0.8) rectangle (20.65,3.4);
\draw (20,2) node {{\small $G'$}};
\end{tikzpicture}
\end{center}
\begin{center}
\vskip - 0.65 cm
\caption{Subgraphs of $G$ in the proof of Claim~\ref{claim:no-adj-K23}}
\label{fig:G7-fig9b}
\end{center}
\end{figure}

Suppose next that $F_1$ and $F_2$ are joined by exactly two edges. Renaming edges if necessary, we may assume that $u_2v_2$ and $u_3v_3$ are edges of $G$. Let $x$ be the neighbor of $v_1$ not in $X$, and let $y$ be the neighbor of $u_1$ not in $Y$. Possibly, $x = y$. Let $X' = (V(F_1) \cup V(F_2)) \setminus \{u_1,v_1\}$ and let $G' = G - X'$. Thus the graph in Figure~\ref{fig:adj-K23-two-edges}(a) is a subgraph of $G$, where in this illustration $x \ne y$. Every $i$-set of $G'$ can be extended to an ID-set of $G$ by adding to it the vertices in the set $\{u_2,v_3\}$ (indicated by the shaded vertices in Figure~\ref{fig:adj-K23-two-edges}(a)), and so $i(G) \le i(G') + 2$. Since every vertex of $G'$ has degree~$3$, except for the vertices $u_1$ and $v_1$ both of which have degree~$1$, we infer that $b(G') = 0$, and so, by Claim~\ref{structural-weight}, we have $\Theta(G') = 0$. Hence, $\Omega(G') = \w(G') \le \w(G) - \w_G(X') + 4 = \Omega(G) - 24 + 4 = \Omega(G) - 20$. Thus, $8i(G) \le 8(i(G') + 2) \le \Omega(G') + 16 < \Omega(G)$, a contradiction.~\smallqed

\begin{figure}[htb]
\begin{center}
\begin{tikzpicture}[scale=.75,style=thick,x=0.8cm,y=0.8cm]
\def\vr{2.5pt} 
\path (-1.5,-2.25) coordinate (z0);
\path (-1.5,3) coordinate (z1);
\path (1.25,3) coordinate (z2);
\path (1.25,-1) coordinate (z3);
\path (5.75,-1) coordinate (z4);
\path (5.75,3) coordinate (z5);
\path (8.5,3) coordinate (z6);
\path (8.5,-2.25) coordinate (z7);
%
\path (-0.75,0.5) coordinate (x01);
\path (-0.75,1.5) coordinate (x02);
\path (-0.25,1) coordinate (x0);
\path (-0.25,1.1) coordinate (x0p);
\path (0.75,1) coordinate (x);
\path (2,0) coordinate (u1);
\path (3,0) coordinate (u2);
\path (4,0) coordinate (u3);
\path (5,0) coordinate (u4);
\path (2,2) coordinate (v1);
\path (3,2) coordinate (v2);
\path (4,2) coordinate (v3);
\path (5,2) coordinate (v4);
\path (6.25,1) coordinate (y);
\path (6.26,1) coordinate (yp);
\path (7.25,1) coordinate (y0);
\path (7.75,0.5) coordinate (y01);
\path (7.75,1.5) coordinate (y02);
\draw (x01)--(x0)--(x02);
\draw (y01)--(y0)--(y02);
\draw (x0)--(x)--(u1)--(u2)--(u3)--(u4)--(y)--(y0);
\draw (x)--(v1)--(v2)--(v3)--(v4)--(y);
\draw (u1)--(v2);
\draw (u2)--(v1);
\draw (u3)--(v4);
\draw (u4)--(v3);
\draw[densely dashed]  (z0)--(z1)--(z2)--(z3)--(z4)--(z5)--(z6)--(z7)--(z0);
\draw (x0) [fill=white] circle (\vr);
\draw (x) [fill=white] circle (\vr);
\draw (y0) [fill=white] circle (\vr);
\draw (y) [fill=white] circle (\vr);
\draw (u1) [fill=white] circle (\vr);
\draw (u2) [fill=white] circle (\vr);
\draw (u3) [fill=black] circle (\vr);
\draw (u4) [fill=white] circle (\vr);
\draw (v1) [fill=white] circle (\vr);
\draw (v2) [fill=black] circle (\vr);
\draw (v3) [fill=white] circle (\vr);
\draw (v4) [fill=white] circle (\vr);
\draw[anchor = south] (v1) node {{\small $x_2$}};
\draw[anchor = south] (v2) node {{\small $v_3$}};
\draw[anchor = south] (v3) node {{\small $u_3$}};
\draw[anchor = south] (v4) node {{\small $y_2$}};
\draw[anchor = north] (u1) node {{\small $x_1$}};
\draw[anchor = north] (u2) node {{\small $v_2$}};
\draw[anchor = north] (u3) node {{\small $u_2$}};
\draw[anchor = north] (u4) node {{\small $y_3$}};
\draw[anchor = south] (x) node {{\small $v_1$}};
\draw[anchor = south] (x0p) node {{\small $x$}};
\draw[anchor = south] (yp) node {{\small $u_1$}};
\draw[anchor = south] (y0) node {{\small $y$}};
\draw [style=dashed,rounded corners] (1.5,-0.75) rectangle (5.5,3.5);
\draw (3.5,3.15) node {{\small $X'$}};
\draw (0,-0.75) node {{\small $G'$}};
\draw (3.6,-3.25) node {{\small (a)}};
%
\path (13,-1.5) coordinate (z0);
\path (13,2.95) coordinate (z1);
\path (11,2.95) coordinate (z2);
\path (11,-2.25) coordinate (z3);
\path (21,-2.25) coordinate (z4);
\path (21,2.95) coordinate (z5);
\path (18.15,2.95) coordinate (z6);
\path (18.15,-1.5) coordinate (z7);
%

\path (12.5,1) coordinate (x);
\path (12.4,1) coordinate (xp);
\path (12.5,-1) coordinate (x0);
\path (12,-1.75) coordinate (x01);
\path (13,-1.75) coordinate (x02);

\path (14,0) coordinate (u1);
\path (15,0) coordinate (u2);
\path (17,0) coordinate (u3);
\path (17.3,0) coordinate (u3p);
\path (18.75,0) coordinate (u4);
\path (14,2) coordinate (v1);
\path (15,2) coordinate (v2);
\path (17,2) coordinate (v3);
\path (18.75,2) coordinate (v4);
\path (16,0) coordinate (x3);
\path (16,-0.85) coordinate (w);
\path (15.45,-1.75) coordinate (w1);
\path (16.55,-1.75) coordinate (w2);

\path (19.75,1) coordinate (y);
\path (19.76,1) coordinate (yp);
\path (19.75,-1) coordinate (y0);
\path (19.15,-1.75) coordinate (y01);
\path (20.35,-1.75) coordinate (y02);
\draw (x0)--(x)--(u1)--(u2);
\draw (u3)--(u4)--(y);
\draw (x)--(v1)--(v2)--(v3)--(v4)--(y);
\draw (u1)--(v2);
\draw (u2)--(v1);
\draw (u3)--(v4);
\draw (u4)--(v3);
\draw (u2)--(x3)--(u3);
\draw (w)--(x3);
\draw (w1)--(w)--(w2);
\draw (x0)--(x);
\draw (x01)--(x0)--(x02);
\draw (y0)--(y);
\draw (y01)--(y0)--(y02);
\draw[densely dashed]  (z0)--(z1)--(z2)--(z3)--(z4)--(z5)--(z6)--(z7)--(z0);
\draw [style=dashed,rounded corners] (13.5,-1.25) rectangle (17.65,3.5);
\draw (15.75,3.15) node {{\small $X^*$}};
\draw (x0) [fill=white] circle (\vr);
\draw (x) [fill=white] circle (\vr);
\draw (x3) [fill=black] circle (\vr);
\draw (w) [fill=white] circle (\vr);
\draw (y0) [fill=white] circle (\vr);
\draw (y) [fill=white] circle (\vr);
\draw (u1) [fill=white] circle (\vr);
\draw (u2) [fill=white] circle (\vr);
\draw (u3) [fill=white] circle (\vr);
\draw (u4) [fill=white] circle (\vr);
\draw (v1) [fill=white] circle (\vr);
\draw (v2) [fill=black] circle (\vr);
\draw (v3) [fill=white] circle (\vr);
\draw (v4) [fill=white] circle (\vr);
\draw[anchor = south] (v1) node {{\small $x_2$}};
\draw[anchor = south] (v2) node {{\small $v_3$}};
\draw[anchor = south] (v3) node {{\small $u_3$}};
\draw[anchor = south] (v4) node {{\small $y_2$}};
\draw[anchor = north] (u1) node {{\small $x_1$}};
\draw[anchor = north] (u2) node {{\small $v_2$}};
\draw[anchor = north] (u3p) node {{\small $u_2$}};
\draw[anchor = north] (u4) node {{\small $y_1$}};
\draw[anchor = south] (xp) node {{\small $v_1$}};
\draw[anchor = east] (x0) node {{\small $x$}};
\draw[anchor = south] (x3) node {{\small $w$}};
\draw[anchor = west] (y) node {{\small $u_1$}};
\draw[anchor = west] (w) node {{\small $z$}};
\draw[anchor = west] (y0) node {{\small $y$}};
%
\draw (11.75,2.25) node {{\small $G^*$}};
\draw (16,-3.25) node {{\small (b)}};
%
\end{tikzpicture}
\end{center}
\begin{center}
\vskip - 0.5 cm
\caption{Subgraphs of $G$ in the proof of Claim~\ref{claim:no-adj-K23}}
\label{fig:adj-K23-two-edges}
\end{center}
\end{figure}

By Claim~\ref{claim:no-adj-K23.1}, $F_1$ and $F_2$ are joined by exactly one edge. Renaming vertices if necessary, we may assume that $u_3v_3$ is the edge joining $F_1$ and $F_2$.

\begin{subclaim}
\label{claim:no-adj-K23.2}
Neither $v_1$ nor $v_2$ has a common neighbor with $u_1$ or $u_2$.
\end{subclaim}
\proof
Suppose that $v_1$ or $v_2$ has a common neighbor with $u_1$ or $u_2$. By symmetry, we may assume that $u_2$ and $v_2$ have a common neighbor, say~$w$. By Claim~\ref{claim:K23-config}, the vertex $w$ is adjacent to neither $u_1$ nor $v_1$. Let $z$ be the neighbor of $w$ different from~$u_2$ and $v_2$. Let $x$ be the neighbor of $v_1$ not in $X$, and let $y$ be the neighbor of $u_1$ not in $Y$. Possibly, $x = y$. By Claim~\ref{claim:K23-config}, the set $\{u_3,x,w\}$ is an independent set and the $\{v_3,y,w\}$ is an independent set. Therefore since $w$ is adjacent to~$z$, we note that $x \ne z$ and $y \ne z$, although possibly $x = y$.

We now let $X^* = (V(F_1) \setminus \{v_1\}) \cup \{u_2,u_3,w,z\}$ and consider the graph $G^* = G - X^*$. Every $i$-set of $G^*$ can be extended to an ID-set of $G$ by adding to it the vertices $v_3$ and~$w$ (indicated by the shaded vertices in Figure~\ref{fig:adj-K23-two-edges}(b)), and so $i(G) \le i(G^*) + 2$. Let $z_1$ and $z_2$ be the two neighbors of $z$ different from~$w$. Every vertex of $G^*$ has degree~$3$ except for three vertices of degree~$1$ (namely, $v_1$, $y_1$ and $y_2$) and two vertices of degree~$2$ (namely, $z_1$ and~$z_2$). Since $y_1$ and $y_2$ have a common neighbor in $G^*$ (namely $u_1$), we note that the only possible component of $G^*$ that belongs to the family $\cB$ is a $B_1$-component containing $v_1$, $z_1$ and $z_2$. However this would imply that $z_1$ and $z_2$ are vertices of degree~$2$ in a copy of $K_{2,3}$ that have a common neighbor, namely~$z$, that does not belong to this copy of $K_{2,3}$, which contradicts Claim~\ref{claim:K23-config}. Hence, $b(G^*) = 0$, and so, by Claim~\ref{structural-weight}, we have $\Theta(G^*) = 0$. Hence, $\Omega(G^*) = \w(G^*) \le \w(G) - \w_G(X^*) + 8 = \Omega(G) - 24 + 8 = \Omega(G) - 16$. Thus, $8i(G) \le 8(i(G^*) + 2) \le \Omega(G^*) + 16 \le \Omega(G)$, a contradiction.~\smallqed

\medskip
We now return to the proof of Claim~\ref{claim:no-adj-K23}. By Claim~\ref{claim:no-adj-K23.2}, neither $v_1$ nor $v_2$ has a common neighbor with $u_1$ or $u_2$. As before, let $x$ be the neighbor of $v_1$ not in $X$, and let $y$ be the neighbor of $u_1$ not in $Y$. Let $w$ be the neighbor of $v_2$ not in $X$, and let $z$ be the neighbor of $u_2$ not in $Y$. By our earlier observations, the vertices $x, y, w, z$ are distinct and do not belong to $V(F_1) \cup V(F_2)$. By Claim~\ref{claim:K23-config}, $xw \notin E(G)$ and $yz \notin E(G)$.

\begin{figure}[htb]
\begin{center}
\begin{tikzpicture}[scale=.75,style=thick,x=0.8cm,y=0.8cm]
\def\vr{2.5pt} 
\path (13.35,3.25) coordinate (a0);
\path (13.35,-1) coordinate (a1);
\path (16.35,-1) coordinate (a2);
\path (16.35,0.7) coordinate (a3);
\path (17.8,0.7) coordinate (a4);
\path (17.8,-0.75) coordinate (a5);
\path (19,-0.75) coordinate (a6);
\path (19,3.25) coordinate (a7);
\path (13.05,-1.25) coordinate (z0);
\path (13.05,2.95) coordinate (z1);
\path (11.5,2.95) coordinate (z2);
\path (11.5,-2.25) coordinate (z3);
\path (21,-2.25) coordinate (z4);
\path (21,2.95) coordinate (z5);
\path (19.25,2.95) coordinate (z6);
\path (19.25,-1) coordinate (z7);
\path (17.65,-1) coordinate (z8);
\path (17.65,0.5) coordinate (z9);
\path (16.65,0.5) coordinate (z10);
\path (16.65,-1.25) coordinate (z11);
%

\path (12.5,1) coordinate (x);
\path (12.4,1) coordinate (xp);
\path (12.5,-1.25) coordinate (x0);
\path (12,-1.9) coordinate (x01);
\path (13,-1.9) coordinate (x02);

\path (14,0) coordinate (u1);
\path (15,0) coordinate (u2);
\path (17,0) coordinate (u3);
\path (17.3,0) coordinate (u3p);
\path (18.25,0) coordinate (u4);
\path (14,2) coordinate (v1);
\path (15,2) coordinate (v2);
\path (17,2) coordinate (v3);
\path (18.25,2) coordinate (v4);
\path (15,-0.75) coordinate (x3);
\path (14.45,-1.5) coordinate (x31);
\path (15.55,-1.5) coordinate (x32);
\path (17,-1.25) coordinate (w);
\path (16.45,-1.9) coordinate (w1);
\path (17.55,-1.9) coordinate (w2);

\path (19.75,1) coordinate (y);
\path (19.76,1) coordinate (yp);
\path (19.75,-1.25) coordinate (y0);
\path (19.15,-1.9) coordinate (y01);
\path (20.35,-1.9) coordinate (y02);
\draw (x0)--(x)--(u1)--(u2);
\draw (u3)--(u4)--(y);
\draw (x)--(v1)--(v2)--(v3)--(v4)--(y);
\draw (u1)--(v2);
\draw (u2)--(v1);
\draw (u3)--(v4);
\draw (u4)--(v3);
\draw (u2)--(x3);
\draw (w1)--(w)--(w2);
\draw (x0)--(x);
\draw (x01)--(x0)--(x02);
\draw (y0)--(y);
\draw (y01)--(y0)--(y02);
\draw (u3)--(w);
\draw (x31)--(x3)--(x32);
\draw[densely dashed]  (z0)--(z1)--(z2)--(z3)--(z4)--(z5)--(z6)--(z7)--(z8)--(z9)--(z10)--(z11)--(z0);
\draw[densely dashed]  (a0)--(a1)--(a2)--(a3)--(a4)--(a5)--(a6)--(a7)--(a0);
%
%
\draw (x0) [fill=white] circle (\vr);
\draw (x) [fill=white] circle (\vr);
\draw (x3) [fill=white] circle (\vr);
\draw (w) [fill=white] circle (\vr);
\draw (y0) [fill=white] circle (\vr);
\draw (y) [fill=white] circle (\vr);
\draw (u1) [fill=white] circle (\vr);
\draw (u2) [fill=black] circle (\vr);
\draw (u3) [fill=white] circle (\vr);
\draw (u4) [fill=white] circle (\vr);
\draw (v1) [fill=white] circle (\vr);
\draw (v2) [fill=white] circle (\vr);
\draw (v3) [fill=black] circle (\vr);
\draw (v4) [fill=white] circle (\vr);
\draw[anchor = south] (v1) node {{\small $x_2$}};
\draw[anchor = south] (v2) node {{\small $v_3$}};
\draw[anchor = south] (v3) node {{\small $u_3$}};
\draw[anchor = south] (v4) node {{\small $y_2$}};
\draw[anchor = north] (u1) node {{\small $x_1$}};
\draw[anchor = west] (u2) node {{\small $v_2$}};
\draw[anchor = north] (u3p) node {{\small $u_2$}};
\draw[anchor = north] (u4) node {{\small $y_1$}};
\draw[anchor = south] (xp) node {{\small $v_1$}};
\draw[anchor = east] (x0) node {{\small $x$}};
\draw[anchor = west] (x3) node {{\small $w$}};
\draw[anchor = west] (y) node {{\small $u_1$}};
\draw[anchor = west] (w) node {{\small $z$}};
\draw[anchor = west] (y0) node {{\small $y$}};
\draw (12.25,2.25) node {{\small $G^*$}};
\draw (16,2.75) node {{\small $X^*$}};
\end{tikzpicture}
\end{center}
\begin{center}
\vskip - 0.5 cm
\caption{A subgraph of $G$ in the proof of Claim~\ref{claim:no-adj-K23}}
\label{fig:adj-K23-two-edges-B}
\end{center}
\end{figure}
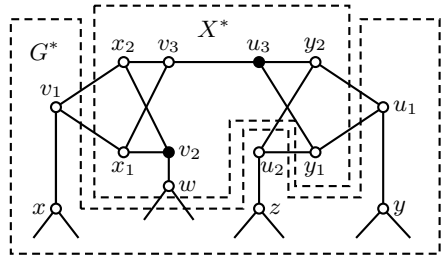

We now let $X^* = \{v_2,v_3,u_3,x_1,x_2,y_1,y_2,w\}$ and we consider the graph $G^* = G - X^*$. Every $i$-set of $G^*$ can be extended to an ID-set of $G$ by adding to it the vertices $v_2$ and~$u_3$ (indicated by the shaded vertices in Figure~\ref{fig:adj-K23-two-edges-B}), and so $i(G) \le i(G^*) + 2$. Let $w_1$ and $w_2$ be the two neighbors of $w$ different from~$v_2$. Every vertex of $G^*$ has degree~$3$ except for three vertices of degree~$1$ (namely, $v_1$, $u_1$ and $u_2$) and two vertices of degree~$2$ (namely, $w_1$ and~$w_2$). Analogous arguments as before show that there is no $B_1$-component containing the vertices $w_1$ and~$w_2$, for otherwise such a component (of order~$6$) contains $w_1$ and $w_2$ as vertices of degree~$2$ in a copy of $K_{2,3}$ that have a common neighbor, namely~$w$, that does not belong to this copy of $K_{2,3}$, which contradicts Claim~\ref{claim:K23-config}. We therefore infer that $b(G^*) = 0$, and so, by Claim~\ref{structural-weight}, we have $\Theta(G^*) = 0$. Hence, $\Omega(G^*) = \w(G^*) \le \w(G) - \w_G(X^*) + 8 = \Omega(G) - 24 + 8 = \Omega(G) - 16$. Thus, $8i(G) \le 8(i(G^*) + 2) \le \Omega(G^*) + 16 = \Omega(G)$, a contradiction. This completes the proof of Claim~\ref{claim:no-adj-K23}.~\smallqed

By Claim~\ref{claim:no-adj-K23}, if the graph $G$ contains two vertex disjoint copies of $K_{2,3}$, then there is no edge joining these two copies of $K_{2,3}$. This yields the following property of the structural weight for an arbitrary induced subgraph $G'$ of~$G$.

\begin{figure}[htb]
\begin{center}
\begin{tikzpicture}[scale=.75,style=thick,x=0.75cm,y=0.75cm]
\def\vr{2.5pt} 
\path (0,1) coordinate (a1);
\path (1,0) coordinate (a2);
\path (1,1) coordinate (a3);
\path (1,2) coordinate (a4);
\path (2,1) coordinate (a5);
\path (1,3) coordinate (u);
\path (5,1) coordinate (b1);
\path (6,0) coordinate (b2);
\path (6,1) coordinate (b3);
\path (6,2) coordinate (b4);
\path (7,1) coordinate (b5);
\path (7.5,3) coordinate (v);
\path (8,1) coordinate (c1);
\path (9,0) coordinate (c2);
\path (9,1) coordinate (c3);
\path (9,2) coordinate (c4);
\path (10,1) coordinate (c5);
\path (13,1) coordinate (d1);
\path (14,0) coordinate (d2);
\path (14,1) coordinate (d3);
\path (14,2) coordinate (d4);
\path (15,1) coordinate (d5);
\path (16,1) coordinate (e1);
\path (17,0) coordinate (e2);
\path (17,1) coordinate (e3);
\path (17,2) coordinate (e4);
\path (18,1) coordinate (e5);
\path (19,1) coordinate (f1);
\path (20,0) coordinate (f2);
\path (20,1) coordinate (f3);
\path (20,2) coordinate (f4);
\path (21,1) coordinate (f5);
\path (17,3) coordinate (w);
\draw (a1)--(a2)--(a5)--(a4)--(a1);
\draw (a1)--(a3)--(a5);
\draw (b1)--(b2)--(b5)--(b4)--(b1);
\draw (b1)--(b3)--(b5);
\draw (c1)--(c2)--(c5)--(c4)--(c1);
\draw (c1)--(c3)--(c5);
\draw (d1)--(d2)--(d5)--(d4)--(d1);
\draw (d1)--(d3)--(d5);
\draw (e1)--(e2)--(e5)--(e4)--(e1);
\draw (e1)--(e3)--(e5);
\draw (f1)--(f2)--(f5)--(f4)--(f1);
\draw (f1)--(f3)--(f5);
\draw (u)--(a4);
\draw (b4)--(v)--(c4);
\draw (d4)--(w)--(e4);
\draw (f4)--(w);
\draw (a1) [fill=white] circle (\vr);
\draw (a2) [fill=white] circle (\vr);
\draw (a3) [fill=white] circle (\vr);
\draw (a4) [fill=white] circle (\vr);
\draw (a5) [fill=white] circle (\vr);
\draw (b1) [fill=white] circle (\vr);
\draw (b2) [fill=white] circle (\vr);
\draw (b3) [fill=white] circle (\vr);
\draw (b4) [fill=white] circle (\vr);
\draw (b5) [fill=white] circle (\vr);
\draw (c1) [fill=white] circle (\vr);
\draw (c2) [fill=white] circle (\vr);
\draw (c3) [fill=white] circle (\vr);
\draw (c4) [fill=white] circle (\vr);
\draw (c5) [fill=white] circle (\vr);
\draw (d1) [fill=white] circle (\vr);
\draw (d2) [fill=white] circle (\vr);
\draw (d3) [fill=white] circle (\vr);
\draw (d4) [fill=white] circle (\vr);
\draw (d5) [fill=white] circle (\vr);
\draw (e1) [fill=white] circle (\vr);
\draw (e2) [fill=white] circle (\vr);
\draw (e3) [fill=white] circle (\vr);
\draw (e4) [fill=white] circle (\vr);
\draw (e5) [fill=white] circle (\vr);
\draw (f1) [fill=white] circle (\vr);
\draw (f2) [fill=white] circle (\vr);
\draw (f3) [fill=white] circle (\vr);
\draw (f4) [fill=white] circle (\vr);
\draw (f5) [fill=white] circle (\vr);
\draw (u) [fill=white] circle (\vr);
\draw (v) [fill=white] circle (\vr);
\draw (w) [fill=white] circle (\vr);
\draw (1,-1) node {{\small (a) $B \in \cB_1$}};
\draw (7.5,-1) node {{\small (b) $B \in \cB_2$}};
\draw (17,-1) node {{\small (c) $B \in \cB_3$}};
\end{tikzpicture}
\end{center}
\begin{center}
\vskip - 0.5 cm
\caption{The three possible bad components of $G'$ in Claim~\ref{claim:bad-comp-structure}}
\label{fig:bad-comp-structure}
\end{center}
\end{figure}
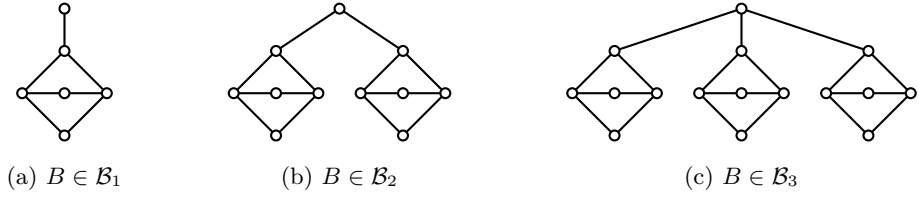

\begin{claim}
\label{claim:bad-comp-structure}
If $G'$ is an arbitrary induced subgraph of $G$ and if $B$ is a component of $G'$ in the family~$\cB$, then $B \in \cB_i$ and $B$ contains exactly $i$ copies of $K_{2,3}$ for some $i \in [3]$ as illustrated in Figure~\ref{fig:bad-comp-structure}.
\end{claim}

By Claim~\ref{claim:bad-comp-structure}, if $G'$ is an arbitrary induced subgraph of $G$ and if $B$ is a component of $G'$ that belongs to~$\cB$, then either $B \in \cB_1$ and $B$ contains one vertex of degree~$1$ and two vertices of degree~$2$ (see Figure~\ref{fig:bad-comp-structure}(a)) or $B \in \cB_2$ and $B$ contains no vertex of degree~$1$ and five vertices of degree~$2$ (see Figure~\ref{fig:bad-comp-structure}(b)) or $B \in \cB_3$ and $B$ contains no vertex of degree~$1$ and six vertices of degree~$2$ (see Figure~\ref{fig:bad-comp-structure}(c)).

\begin{figure}[htb]
\begin{center}
\begin{tikzpicture}[scale=.75,style=thick,x=0.8cm,y=0.8cm]
\def\vr{2.5pt} 
\path (10.5,0.75) coordinate (q21);
\path (10.5,1.75) coordinate (q22);
\path (10.5,2.25) coordinate (p21);
\path (10.5,3.25) coordinate (p22);
\path (10.5,-0.75) coordinate (y21);
\path (10.5,0.25) coordinate (y22);
\path (12.5,-0.25) coordinate (v1);
\path (12.6,-0.25) coordinate (v1p);
\path (12.5,1.25) coordinate (v2);
\path (12.6,1.25) coordinate (v2p);
\path (12.5,2.75) coordinate (v3);
\path (15,-0.25) coordinate (u1);
\path (15,1.25) coordinate (u2);
\path (15,2.75) coordinate (u3);
\path (16.5,0.625) coordinate (w1);
\path (16.5,1.875) coordinate (w2);
\draw (v1)--(u1);
\draw (v2)--(u2);
\draw (v3)--(u3);
\draw (w1)--(u1)--(w2);
\draw (w1)--(u2)--(w2);
\draw (w1)--(u3)--(w2);
\draw (y21)--(v1)--(y22);
\draw (p21)--(v3)--(p22);
\draw (q21)--(v2)--(q22);
\draw (v1) [fill=white] circle (\vr);
\draw (v2) [fill=white] circle (\vr);
\draw (v3) [fill=white] circle (\vr);
\draw (u1) [fill=white] circle (\vr);
\draw (u2) [fill=white] circle (\vr);
\draw (u3) [fill=white] circle (\vr);
\draw (w1) [fill=white] circle (\vr);
\draw (w2) [fill=white] circle (\vr);
%
%
\draw[anchor = north] (v1p) node {{\small $u_3$}};
\draw[anchor = north] (v2p) node {{\small $u_2$}};
\draw[anchor = south] (v3) node {{\small $u_1$}};
\draw[anchor = north] (u1) node {{\small $v_3$}};
\draw[anchor = north] (u2) node {{\small $v_2$}};
\draw[anchor = south] (u3) node {{\small $v_1$}};
\draw[anchor = west] (w1) node {{\small $x_1$}};
\draw[anchor = west] (w2) node {{\small $x_2$}};
\draw [style=dashed,rounded corners] (11.85,-1) rectangle (13,4.25);
\draw (12.5,3.75) node {{\small $C$}};
\draw [style=dashed,rounded corners] (14.35,-1) rectangle (15.5,4.25);
\draw (15,3.75) node {{\small $B$}};
\draw [style=dashed,rounded corners] (16,-1) rectangle (17.5,4.25);
\draw (16.5,3.75) node {{\small $A$}};
\draw [style=dashed,rounded corners] (13.85,-1.25) rectangle (18.5,4.5);
\draw (18,2) node {{\small $F$}};
\draw [style=dashed,rounded corners] (11.5,-1.75) rectangle (19,5.25);
\draw (12.75,4.75) node {{\small set $X''$}};
\draw [style=dashed,rounded corners] (9.25,-1.75) rectangle (11,4.25);
\draw (9.85,1.25) node {{\small $G''$}};
\end{tikzpicture}
\end{center}
\begin{center}
\vskip - 0.25 cm
\caption{A subgraph of $G$}
\label{fig:cubic-structure-1}
\end{center}
\end{figure}
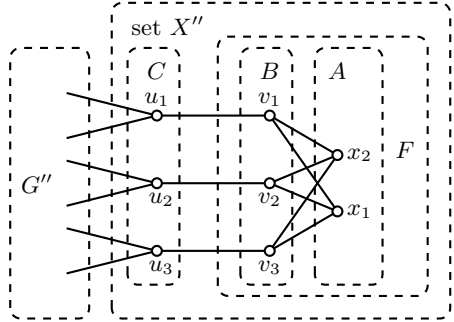

By Claim~\ref{no-triangle}, the graph $G$ contains no triangle.
Let $F$ be a $K_{2,3}$-subgraph in $G$ with partite sets $A = \{x_1,x_2\}$ and $B = \{v_1,v_2,v_3\}$. Since $G$ is triangle-free, the set $B$ is an independent set in $G$. Let $u_i$ be the neighbor of $v_i$ different from $x_1$ and $x_2$ for $i \in [3]$. By Claim~\ref{claim:K23-config}, the vertices $u_1$, $u_2$ and $u_3$ are distinct and the set $C = \{u_1,u_2,u_3\}$ is an independent set. We now let $X'' = A \cup B \cup C$ and consider the graph $G'' = G - X''$. Thus the graph illustrated in Figure~\ref{fig:cubic-structure-1} is an induced subgraph of $G$.

\begin{claim}
\label{set-no-isolates}
The graph $G''$ is isolate-free.
\end{claim}
\proof Suppose, to the contrary, that $G''$ contains at least one isolated vertex. If $G''$ contains two isolated vertices, then the graph $G$ is determined and is the graph $G_{10.2}$ shown in Figure~\ref{fig:G7-fig9b}(a). As noted earlier, in this case $8i(G) < \w(G)$, a contradiction. Hence, $G''$ contains exactly one isolated vertex, say~$w$. Let $X' = X'' \cup \{w\}$ and let $G' = G - X'$. Every $i$-set of $G'$ can be extended to an ID-set of $G$ by adding to it the vertices in the set $\{w,x_1,x_2\}$ (as indicated by the shaded vertices in Figure~\ref{fig:G7-fig9b}(b)), and so $i(G) \le i(G') + 3$. Since $G'$ has either one vertex of degree~$1$ and one vertex of degree~$2$ or three vertices of degree~$2$, we infer that $b(G') = 0$, and so by Claim~\ref{structural-weight}, we have $\Theta(G') = 0$. Hence, $\Omega(G') = \w(G') \le \w(G) - \w_G(X') + 3 = \Omega(G) - 27 + 3 = \Omega(G) - 24$. Thus, $8i(G) \le 8(i(G') + 3) \le \Omega(G') + 24 = \Omega(G)$, a contradiction.~\smallqed

\medskip
By Claim~\ref{set-no-isolates}, the graph $G''$ is isolate-free.

\begin{figure}[htb]
\begin{center}
\begin{tikzpicture}[scale=.75,style=thick,x=0.8cm,y=0.8cm]
\def\vr{2.5pt} 
\path (4,-0.25) coordinate (p1);
\path (4,1.25) coordinate (p2);
\path (4,2.75) coordinate (p3);
\path (2.75,1.25) coordinate (p);
\path (5.5,-0.25) coordinate (v1);
\path (5.6,-0.25) coordinate (v1p);
%
\path (5.5,1.25) coordinate (v2);
\path (5.6,1.25) coordinate (v2p);
\path (5.5,2.75) coordinate (v3);
\path (7,-0.25) coordinate (u1);
\path (7,1.25) coordinate (u2);
\path (7,2.75) coordinate (u3);
\path (8.5,0.5) coordinate (w1);
\path (8.5,2) coordinate (w2);
\draw (v1)--(u1);
\draw (v2)--(u2);
\draw (v3)--(u3);
\draw (w1)--(u1)--(w2);
\draw (w1)--(u2)--(w2);
\draw (w1)--(u3)--(w2);

\draw (v1)--(p1)--(v3);
\draw (v1)--(p2)--(v2);
\draw (v2)--(p3)--(v3);
\draw (p1)--(p)--(p2);
\draw (p3)--(p);
\draw (v1) [fill=white] circle (\vr);
\draw (v2) [fill=white] circle (\vr);
\draw (v3) [fill=white] circle (\vr);
\draw (u1) [fill=black] circle (\vr);
\draw (u2) [fill=black] circle (\vr);
\draw (u3) [fill=black] circle (\vr);
\draw (p) [fill=black] circle (\vr);
\draw (p1) [fill=white] circle (\vr);
\draw (p2) [fill=white] circle (\vr);
\draw (p3) [fill=white] circle (\vr);
\draw (w1) [fill=white] circle (\vr);
\draw (w2) [fill=white] circle (\vr);
%
%
\draw[anchor = north] (v1p) node {{\small $u_3$}};
\draw[anchor = north] (v2p) node {{\small $u_2$}};
\draw[anchor = south] (v3) node {{\small $u_1$}};
\draw[anchor = north] (u1) node {{\small $v_3$}};
\draw[anchor = north] (u2) node {{\small $v_2$}};
\draw[anchor = south] (u3) node {{\small $v_1$}};
\draw[anchor = north] (p1) node {{\small $w_3$}};
\draw[anchor = south] (p2) node {{\small $w_2$}};
\draw[anchor = south] (p3) node {{\small $w_1$}};
\draw[anchor = east] (p) node {{\small $w$}};
\draw[anchor = west] (w1) node {{\small $x_1$}};
\draw[anchor = west] (w2) node {{\small $x_2$}};
\draw (6.25,-1.5) node {{\small (a) $G_{12}$}};
%
%
\path (14,-0.25) coordinate (p1);
\path (14,1.25) coordinate (p2);
\path (14,2.75) coordinate (p3);
\path (12.75,-0.25) coordinate (p11);
\path (12.75,1.25) coordinate (p21);
\path (12.75,2.75) coordinate (p31);
%
%
\path (15.5,-0.25) coordinate (v1);
\path (15.6,-0.25) coordinate (v1p);
%
\path (15.5,1.25) coordinate (v2);
\path (15.6,1.25) coordinate (v2p);
\path (15.5,2.75) coordinate (v3);
\path (17,-0.25) coordinate (u1);
\path (17,1.25) coordinate (u2);
\path (17,2.75) coordinate (u3);
\path (18.5,0.5) coordinate (w1);
\path (18.5,2) coordinate (w2);
\draw (v1)--(u1);
\draw (v2)--(u2);
\draw (v3)--(u3);
\draw (w1)--(u1)--(w2);
\draw (w1)--(u2)--(w2);
\draw (w1)--(u3)--(w2);
\draw (v1)--(p1)--(v3);
\draw (v1)--(p2)--(v2);
\draw (v2)--(p3)--(v3);
\draw (p1)--(p11);
\draw (p2)--(p21);
\draw (p3)--(p31);
\draw (v1) [fill=black] circle (\vr);
\draw (v2) [fill=white] circle (\vr);
\draw (v3) [fill=black] circle (\vr);
\draw (u1) [fill=white] circle (\vr);
\draw (u2) [fill=black] circle (\vr);
\draw (u3) [fill=white] circle (\vr);
%
\draw (p1) [fill=white] circle (\vr);
\draw (p2) [fill=white] circle (\vr);
\draw (p3) [fill=white] circle (\vr);
\draw (w1) [fill=white] circle (\vr);
\draw (w2) [fill=white] circle (\vr);
%
%
\draw[anchor = north] (v1p) node {{\small $u_3$}};
\draw[anchor = north] (v2p) node {{\small $u_2$}};
\draw[anchor = south] (v3) node {{\small $u_1$}};
\draw[anchor = north] (u1) node {{\small $v_3$}};
\draw[anchor = north] (u2) node {{\small $v_2$}};
\draw[anchor = south] (u3) node {{\small $v_1$}};
\draw[anchor = north] (p1) node {{\small $w_3$}};
\draw[anchor = south] (p2) node {{\small $w_2$}};
\draw[anchor = south] (p3) node {{\small $w_1$}};
%
\draw[anchor = west] (w1) node {{\small $x_1$}};
\draw[anchor = west] (w2) node {{\small $x_2$}};
\draw (16,-1.5) node {{\small (b)}};
\draw [style=dashed,rounded corners] (11.5,-0.8) rectangle (13.25,3.4);
\draw (12.25,1.25) node {{\small $G'$}};
\path (25.1,0.5) coordinate (p2);
\path (25.1,2) coordinate (p3);
\path (23.5,0.5) coordinate (p21);
\path (23.5,2) coordinate (p31);
\path (22.75,1.25) coordinate (p);
\path (23.5,-0.25) coordinate (v11);
\path (25.5,-0.25) coordinate (v1);
\path (25.6,-0.25) coordinate (v1p);
%
\path (25.5,1.25) coordinate (v2);
\path (25.6,1.25) coordinate (v2p);
\path (25.5,2.75) coordinate (v3);
\path (23.5,2.75) coordinate (v31);
\path (27,-0.25) coordinate (u1);
\path (27,1.25) coordinate (u2);
\path (27,2.75) coordinate (u3);
\path (28.5,0.5) coordinate (w1);
\path (28.5,2) coordinate (w2);
\draw (v1)--(u1);
\draw (v2)--(u2);
\draw (v3)--(u3);
\draw (w1)--(u1)--(w2);
\draw (w1)--(u2)--(w2);
\draw (w1)--(u3)--(w2);
\draw (v1)--(v11);
\draw (v3)--(v31);
\draw (p2)--(p21);
\draw (p3)--(p31);
\draw (v1)--(p2)--(v2);
\draw (v2)--(p3)--(v3);
%
\draw (v1) [fill=white] circle (\vr);
\draw (v2) [fill=black] circle (\vr);
\draw (v3) [fill=white] circle (\vr);
\draw (u1) [fill=black] circle (\vr);
\draw (u2) [fill=white] circle (\vr);
\draw (u3) [fill=black] circle (\vr);
%
\draw (p2) [fill=white] circle (\vr);
\draw (p3) [fill=white] circle (\vr);
\draw (w1) [fill=white] circle (\vr);
\draw (w2) [fill=white] circle (\vr);
%
%
\draw[anchor = north] (v1) node {{\small $u_3$}};
\draw[anchor = east] (v2) node {{\small $u_2$}};
\draw[anchor = south] (v3) node {{\small $u_1$}};
\draw[anchor = north] (u1) node {{\small $v_3$}};
\draw[anchor = north] (u2) node {{\small $v_2$}};
\draw[anchor = south] (u3) node {{\small $v_1$}};
%
\draw[anchor = west] (p2) node {{\small $w_2$}};
\draw[anchor = west] (p3) node {{\small $w_1$}};
%
\draw[anchor = west] (w1) node {{\small $x_1$}};
\draw[anchor = west] (w2) node {{\small $x_2$}};
\draw (26.25,-1.5) node {{\small (c)}};
\draw [style=dashed,rounded corners] (22.25,-0.8) rectangle (24.25,3.4);
\draw (23.25,1.25) node {{\small $G'$}};
\end{tikzpicture}
\end{center}
\begin{center}
\vskip - 0.65 cm
\caption{Subgraphs of $G$ in the proof of Claim~\ref{claim:at-most-one-deg1}}
\label{fig:G7-fig10}
\end{center}
\end{figure}
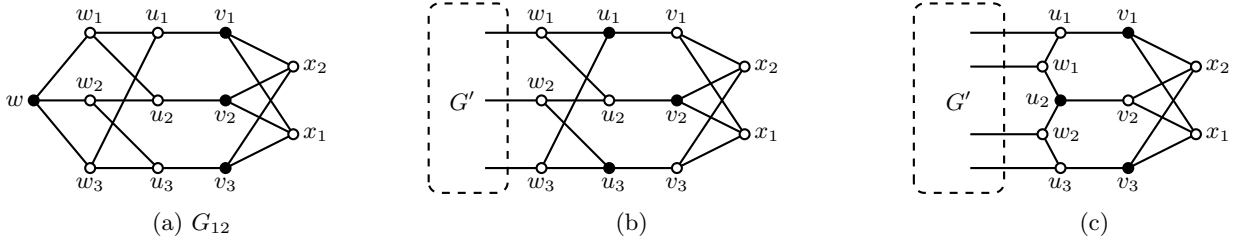

\begin{claim}
\label{claim:at-most-one-deg1}
The graph $G''$ contains at most one vertex of degree~$1$.
\end{claim}
\proof Suppose, to the contrary, that $G''$ contains at least two vertices of degree~$1$. Let $D$ be the set of vertices of degree~$1$ in $G''$. We note that $|D| \le 3$. Suppose firstly that $|D| = 3$. Let $D = \{w_1,w_2,w_3\}$. Let $w_i$ be the common neighbor of $u_i$ and $u_{i+1}$ for $i \in [3]$ where addition is taken modulo~$3$. If $w_1, w_2, w_3$ have a common neighbor, say~$w$, then the graph $G$ is determined and is the graph $G_{12}$ shown in Figure~\ref{fig:G7-fig10}(a). In this case, $i(G) = 4$ (an $i$-set is indicated by the shaded vertices in Figure~\ref{fig:G7-fig10}(a)) and $\w(G) = 36$, and so $8i(G) < \w(G)$, a contradiction. Hence, $w_1, w_2, w_3$ do not have a common neighbor.

We now let $X' = X'' \cup D$ and consider the graph $G' = G'' - D = G - X'$. Since $w_1, w_2, w_3$ have no common neighbor, the graph $G'$ is an isolate-free graph, and so $\w(G') = (\w(G) - 33) + 3 = \w(G) - 30$. Every $i$-set of $G'$ can be extended to an ID-set of $G$ by adding to it the vertices in the set $\{u_1,u_3,v_2\}$ (as indicated by the shaded vertices in Figure~\ref{fig:G7-fig10}(b)), and so $i(G) \le i(G') + 3$. Since $G'$ has either one vertex of degree~$1$ and one vertex of degree~$2$ or three vertices of degree~$2$, we infer by Claim~\ref{claim:bad-comp-structure} that $b(G') = 0$. Thus, by Claim~\ref{structural-weight}, we have $\Theta(G') = 0$, and so $\Omega(G') = \w(G') = \w(G) - 30 = \Omega(G) - 30$. Thus, $8i(G) \le 8(i(G') + 3) \le \Omega(G') + 24 < \Omega(G)$, a contradiction.

Hence, $|D| = 2$, and so $G''$ contains two vertices of degree~$1$, say $w_1$ and~$w_2$. By symmetry, we have two possibilities. First, let $w_i$ be a common neighbor of $u_i$ and $u_{i+1}$ for $i \in [2]$. Thus the graph shown in Figure~\ref{fig:G7-fig10}(c) is a subgraph of~$G$. We now let $X' = X'' \cup \{w_1,w_2\}$ and consider the graph $G' = G'' - D = G - X'$. Since the graph $G$ is triangle-free, the graph $G'$ is isolate-free, and so $\w(G') = (\w(G) - 30) + 4 = \w(G) - 26$. Every $i$-set of $G'$ can be extended to an ID-set of $G$ by adding to it the vertices in the set $\{u_2,v_1,v_3\}$ (as indicated by the shaded vertices in Figure~\ref{fig:G7-fig10}(c)), and so $i(G) \le i(G') + 3$. Since $G'$ has either two vertices of degree~$1$ and no vertex of degree~$2$ or one vertex of degree~$1$ and two vertices of degree~$2$ or four vertices of degree~$2$, we infer that $b(G') \le 1$, and so by Claim~\ref{structural-weight}, we have $\Theta(G') \le 2$. Hence, $\Omega(G') = \w(G') + \Theta(G') \le (\w(G) - 26) + 2 = \w(G) - 24$. Thus, $8i(G) \le 8(i(G') + 3) \le \w(G') + 24 = \w(G)$, a contradiction.

Second, let $u_1$ and $u_2$ have two common neighbors $w_1$ and $w_2$. Suppose that $u_3$, $w_1$ and $w_2$ have a common neighbor $z$. Let $X^*=X''\cup\{w_1,w_2,z\}$. Let $G^*=G-X^*$. By using Claims~\ref{claim:no-bad-in-Ge} and~\ref{structural-weight}, we have $\Omega(G^*) = \w(G^*) = (\w(G)-33)+1 = \Omega(G) + 32$, since there is one $X^*$-exit edge. Every $i$-set of $G^*$ can be extended to an ID-set of $G$ by adding to it the vertices in the set $\{z,v_1,v_2,v_3\}$.  Hence, $8i(G)\le 8(i(G^*)+4)\le \Omega(G^*)+32=\Omega(G)$, a contradiction. Thus,  $u_3$, $w_1$ and $w_2$ do not have a common neighbor. In this case,  let $X' = X'' \cup \{w_1,w_2\}$ and consider the graph $G' = G - X'$. Since $u_3$, $w_1$ and $w_2$ do not have a common neighbor, the graph $G'$ is isolate-free, and so $\w(G') = (\w(G) - 30) + 4 = \w(G) - 26$. Every $i$-set of $G'$ can be extended to an ID-set of $G$ by adding to it the vertices in the set $\{u_2,v_1,v_3\}$ (as indicated by the shaded vertices in Figure~\ref{fig:G7-fig10}(c)), and so $i(G) \le i(G') + 3$. Since $G'$ has either two vertices of degree~$1$ and no vertex of degree~$2$ or one vertex of degree~$1$ and two vertices of degree~$2$ or four vertices of degree~$2$, we infer that $b(G') \le 1$, and so by Claim~\ref{structural-weight}, we have $\Theta(G') \le 2$. Hence, $\Omega(G') = \w(G') + \Theta(G') \le (\w(G) - 26) + 2 = \w(G) - 24$. Thus, $8i(G) \le 8(i(G') + 3) \le \w(G') + 24 = \w(G)$, a contradiction.~\smallqed

\medskip
By Claim~\ref{claim:at-most-one-deg1}, the graph $G''$ contains at most one vertex of degree~$1$. We now define
\[
G_{i,j} = G - (N_G[u_i] \cup N_G[v_j])
\]
where $i, j \in [3]$ and $i \ne j$.

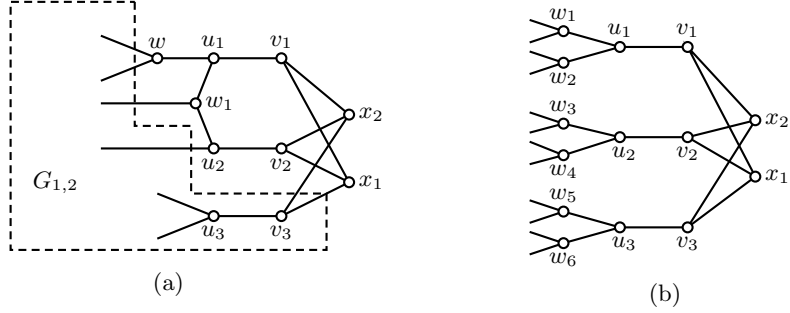
\begin{figure}[htb]
\begin{center}
\begin{tikzpicture}[scale=.75,style=thick,x=0.8cm,y=0.8cm]
\def\vr{2.5pt} 
\path (3.75,4.5) coordinate (z0);
\path (3.75,1.75) coordinate (z1);
\path (5,1.75) coordinate (z2);
\path (5,0.25) coordinate (z3);
\path (8,0.25) coordinate (z4);
\path (8,-1) coordinate (z5);
\path (1,-1) coordinate (z6);
\path (1,4.5) coordinate (z7);
\path (5.1,2.25) coordinate (p3);
\path (3.5,0.5) coordinate (p21);
\path (3,2.25) coordinate (p31);
\path (2.75,1.25) coordinate (p);
\path (4.25,-0.75) coordinate (v11);
\path (4.25,0.25) coordinate (v12);
\path (5.5,-0.25) coordinate (v1);
\path (5.6,-0.25) coordinate (v1p);
\path (3,1.25) coordinate (v21);
\path (5.5,1.25) coordinate (v2);
\path (5.6,1.25) coordinate (v2p);
\path (5.5,3.25) coordinate (v3);
\path (4.25,3.25) coordinate (r);
\path (3,2.75) coordinate (r1);
\path (3,3.75) coordinate (r2);
\path (7,-0.25) coordinate (u1);
\path (7,1.25) coordinate (u2);
\path (7,3.25) coordinate (u3);
\path (8.5,0.5) coordinate (w1);
\path (8.5,2) coordinate (w2);
\draw (v1)--(u1);
\draw (v2)--(u2);
\draw (v3)--(u3);
\draw (w1)--(u1)--(w2);
\draw (w1)--(u2)--(w2);
\draw (w1)--(u3)--(w2);
\draw (r)--(v3);
\draw (v11)--(v1)--(v12);
\draw (v2)--(v21);
\draw (p3)--(p31);
\draw (v2)--(p3)--(v3);
\draw (r1)--(r)--(r2);
\draw (v1) [fill=white] circle (\vr);
\draw (v2) [fill=white] circle (\vr);
\draw (v3) [fill=white] circle (\vr);
\draw (u1) [fill=white] circle (\vr);
\draw (u2) [fill=white] circle (\vr);
\draw (u3) [fill=white] circle (\vr);
\draw (r) [fill=white] circle (\vr);
\draw (p3) [fill=white] circle (\vr);
\draw (w1) [fill=white] circle (\vr);
\draw (w2) [fill=white] circle (\vr);
\draw[densely dashed]  (z0)--(z1)--(z2)--(z3)--(z4)--(z5)--(z6)--(z7)--(z0);
\draw (2,0.5) node {{\small $G_{1,2}$}};
%
%
\draw[anchor = north] (v1) node {{\small $u_3$}};
\draw[anchor = north] (v2) node {{\small $u_2$}};
\draw[anchor = south] (v3) node {{\small $u_1$}};
\draw[anchor = north] (u1) node {{\small $v_3$}};
\draw[anchor = north] (u2) node {{\small $v_2$}};
\draw[anchor = south] (u3) node {{\small $v_1$}};
\draw[anchor = west] (p3) node {{\small $w_1$}};
\draw[anchor = west] (w1) node {{\small $x_1$}};
\draw[anchor = west] (w2) node {{\small $x_2$}};
\draw[anchor = south] (r) node {{\small $w$}};
\draw (4.5,-1.75) node {{\small (a)}};
%
\path (16,3.5) coordinate (u3);
\path (14.5,3.5) coordinate (v3);
\path (13.25,3.15) coordinate (p21);
\path (13.25,3.85) coordinate (p22);
\path (12.5,2.9) coordinate (a1);
\path (12.5,3.4) coordinate (a2);
\path (12.5,3.6) coordinate (b1);
\path (12.5,4.1) coordinate (b2);
\path (12.5,0.9) coordinate (c1);
\path (12.5,1.4) coordinate (c2);
\path (12.5,1.6) coordinate (d1);
\path (12.5,2.05) coordinate (d2);
\path (12.5,-1.1) coordinate (e1);
\path (12.5,-0.6) coordinate (e2);
\path (12.5,-0.4) coordinate (f1);
\path (12.5,0.1) coordinate (f2);
\path (13.25,1.1) coordinate (q21);
\path (13.25,1.8) coordinate (q22);
\path (13.25,-0.85) coordinate (y21);
\path (13.25,-0.15) coordinate (y22);
\path (14.5,-0.5) coordinate (v1);
\path (14.6,-0.5) coordinate (v1p);
\path (14.5,1.5) coordinate (v2);
\path (14.6,1.5) coordinate (v2p);
\path (16,-0.5) coordinate (u1);
\path (16,1.5) coordinate (u2);
\path (17.5,0.625) coordinate (w1);
\path (17.5,1.875) coordinate (w2);
\draw (v1)--(u1);
\draw (v2)--(u2);
\draw (v3)--(u3);
\draw (w1)--(u1)--(w2);
\draw (w1)--(u2)--(w2);
\draw (w1)--(u3)--(w2);
\draw (y21)--(v1)--(y22);
\draw (p21)--(v3)--(p22);
\draw (q21)--(v2)--(q22);
\draw (a1)--(p21)--(a2);
\draw (b1)--(p22)--(b2);
\draw (c1)--(q21)--(c2);
\draw (d1)--(q22)--(d2);
\draw (e1)--(y21)--(e2);
\draw (f1)--(y22)--(f2);
%
%
%
\draw (v1) [fill=white] circle (\vr);
\draw (v2) [fill=white] circle (\vr);
\draw (v3) [fill=white] circle (\vr);
\draw (u1) [fill=white] circle (\vr);
\draw (u2) [fill=white] circle (\vr);
\draw (u3) [fill=white] circle (\vr);
\draw (w1) [fill=white] circle (\vr);
\draw (w2) [fill=white] circle (\vr);
\draw (p21) [fill=white] circle (\vr);
\draw (p22) [fill=white] circle (\vr);
\draw (y21) [fill=white] circle (\vr);
\draw (y22) [fill=white] circle (\vr);
\draw (q21) [fill=white] circle (\vr);
\draw (q22) [fill=white] circle (\vr);
%
%
\draw[anchor = west] (w1) node {{\small $x_1$}};
\draw[anchor = west] (w2) node {{\small $x_2$}};
\draw[anchor = north] (v1p) node {{\small $u_3$}};
\draw[anchor = north] (v2p) node {{\small $u_2$}};
\draw[anchor = south] (v3) node {{\small $u_1$}};
\draw[anchor = north] (u1) node {{\small $v_3$}};
\draw[anchor = north] (u2) node {{\small $v_2$}};
\draw[anchor = south] (u3) node {{\small $v_1$}};
\draw[anchor = north] (p21) node {{\small $w_2$}};
\draw[anchor = south] (p22) node {{\small $w_1$}};
\draw[anchor = north] (q21) node {{\small $w_4$}};
\draw[anchor = south] (q22) node {{\small $w_3$}};
\draw[anchor = north] (y21) node {{\small $w_6$}};
\draw[anchor = south] (y22) node {{\small $w_5$}};
\draw (15.5,-2) node {{\small (b)}};
\end{tikzpicture}
\end{center}
\begin{center}
\vskip - 0.65 cm
\caption{Subgraphs of $G$ in the proof of Claim~\ref{claim:at-most-one-deg1-b}}
\label{fig:G7-fig10-b}
\end{center}
\end{figure}

\begin{claim}
\label{claim:at-most-one-deg1-b}
At least one of the subgraphs $G_{i,j}$ where $i, j \in [3]$ and $i \ne j$ contains no isolated vertex.
\end{claim}
\proof Suppose firstly that $G''$ contains a vertex of degree~$1$, say~$w_1$. By symmetry, we may assume that $w_1$ is a common neighbor of $u_1$ and $u_2$. Let $w$ be the neighbor of $u_1$ different from~$v_1$ and~$w_1$. Thus the graph shown in Figure~\ref{fig:G7-fig10-b}(a) is a subgraph of~$G$, where the graph $G_{1,2}$ is the subgraph in the dashed box. If $G_{1,2}$ contains an isolated vertex~$v$, then $N_G(v) = \{w,w_1,u_2\}$. However, then, $G[\{v,u_2,w_1\}] = K_3$, contradicting Claim~\ref{no-triangle}. Hence, $G_{1,2}$ contains no isolated vertex, and desired result of the claim follows.

Suppose next that $G''$ contains no vertex of degree~$1$. Thus, every vertex in $G''$ is adjacent to at most one vertex that belongs to the set $C = \{u_1,u_2,u_3\}$, implying that $N_G[u_i] \cap N_G[u_j] = \emptyset$ for $i, j \in [3]$ and $i \ne j$. Hence the graph shown in Figure~\ref{fig:G7-fig10-b}(b) is a subgraph of $G$. In particular, we note that the vertices $w_1,w_2,\ldots,w_6$ are distinct. We show that at least one of the subgraphs $G_{i,j}$ where $i, j \in [3]$ and $i \ne j$ contains no isolated vertex. Suppose, to the contrary, that $G_{i,j}$ contains an isolated vertex for all $i, j \in [3]$ and $i \ne j$.

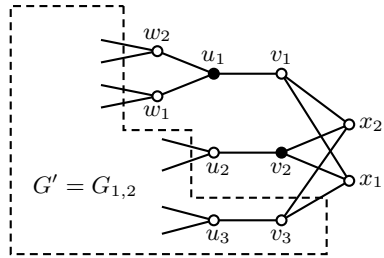
\begin{figure}[htb]
\begin{center}
\begin{tikzpicture}[scale=.75,style=thick,x=0.8cm,y=0.8cm]
\def\vr{2.5pt} 
\path (10.5,4.5) coordinate (z0);
\path (10.5,1.75) coordinate (z1);
\path (12,1.75) coordinate (z2);
\path (12,0.25) coordinate (z3);
\path (15,0.25) coordinate (z4);
\path (15,-1) coordinate (z5);
\path (8,-1) coordinate (z6);
\path (8,4.5) coordinate (z7);
\path (14,3) coordinate (u3);
\path (12.5,3) coordinate (v3);
\path (11.25,2.5) coordinate (p21);
\path (11.25,3.5) coordinate (p22);
\path (10,2.25) coordinate (a1);
\path (10,2.75) coordinate (a2);
\path (10,3.25) coordinate (b1);
\path (10,3.75) coordinate (b2);
\path (11.35,0.85) coordinate (q21);
\path (11.35,1.55) coordinate (q22);
\path (11.35,-0.55) coordinate (y21);
\path (11.35,0.05) coordinate (y22);
\path (12.5,-0.25) coordinate (x1);
\path (12.6,-0.25) coordinate (x1p);
%
\path (12.5,1.25) coordinate (v2);
\path (12.6,1.25) coordinate (v2p);
\path (14,-0.25) coordinate (u1);
\path (14,1.25) coordinate (u2);
\path (15.5,0.625) coordinate (w1);
\path (15.5,1.875) coordinate (w2);
\draw (x1)--(u1);
\draw (v2)--(u2);
\draw (v3)--(u3);
\draw (w1)--(u1)--(w2);
\draw (w1)--(u2)--(w2);
\draw (w1)--(u3)--(w2);
\draw (y21)--(x1)--(y22);
\draw (p21)--(v3)--(p22);
\draw (q21)--(v2)--(q22);
\draw (a1)--(p21)--(a2);
\draw (b1)--(p22)--(b2);
\draw[densely dashed]  (z0)--(z1)--(z2)--(z3)--(z4)--(z5)--(z6)--(z7)--(z0);
\draw (9.65,0.5) node {{\small $G' = G_{1,2}$}};
\draw (x1) [fill=white] circle (\vr);
\draw (v2) [fill=white] circle (\vr);
\draw (v3) [fill=black] circle (\vr);
\draw (u1) [fill=white] circle (\vr);
\draw (u2) [fill=black] circle (\vr);
\draw (u3) [fill=white] circle (\vr);
\draw (w1) [fill=white] circle (\vr);
\draw (w2) [fill=white] circle (\vr);
\draw (p21) [fill=white] circle (\vr);
\draw (p22) [fill=white] circle (\vr);
%
%
\draw[anchor = north] (x1p) node {{\small $u_3$}};
\draw[anchor = north] (v2p) node {{\small $u_2$}};
\draw[anchor = south] (v3) node {{\small $u_1$}};
\draw[anchor = north] (u1) node {{\small $v_3$}};
\draw[anchor = north] (u2) node {{\small $v_2$}};
\draw[anchor = south] (u3) node {{\small $v_1$}};
\draw[anchor = west] (w1) node {{\small $x_1$}};
\draw[anchor = west] (w2) node {{\small $x_2$}};
\draw[anchor = north] (p21) node {{\small $w_1$}};
\draw[anchor = south] (p22) node {{\small $w_2$}};
%
%
%
%
\end{tikzpicture}
\end{center}
\begin{center}
\caption{A subgraph of $G$ in the proof of Claim~\ref{claim:at-most-one-deg1-b}}
\label{fig:G12}
\end{center}
\end{figure}

We consider the graph $G_{1,2}$ (illustrated by the subgraph inside the dashed region in Figure~\ref{fig:G12}). By supposition, $G_{1,2}$ contains an isolated vertex. The only possible such isolated vertex is a neighbor of $u_2$ that is adjacent to both~$w_1$ and~$w_2$. Renaming $w_3$ and $w_4$ if necessary, we may assume that $w_3$ is isolated in $G_{1,2}$, that is, $N_G(w_3) = \{u_2,w_1,w_2\}$. We next consider the graph $G_{1,3}$. By supposition, $G_{1,3}$ contains an isolated vertex. The only possible such isolated vertex is a neighbor of $u_3$ that is adjacent to both~$w_1$ and~$w_2$. Renaming $w_5$ and $w_6$ if necessary, we may assume that $w_5$ is isolated in $G_{1,3}$, that is, $N_G(w_5) = \{u_3,w_1,w_2\}$. In particular, we note that $C \colon w_1w_3w_2w_5w_1$ is a $4$-cycle in $G$. However, this implies that the graph $G_{2,1}$ does not contain an isolated vertex. Indeed, in $G_{2,1}$, vertices $w_1,w_2$ and $u_3$ have degree $1$, the neighbors of $w_4$ in $G_{2,1}$ have degree $2$ in $G_{2,1}$, whereas other vertices of $G_{2,1}$ have degree $3$. This is a contradiction, proving the claim.~\smallqed

By Claim~\ref{claim:at-most-one-deg1-b}, at least one of the subgraphs $G_{i,j}$ where $i, j \in [3]$ and $i \ne j$ contains no isolated vertex. Renaming vertices if necessary, we may assume that $G_{1,2}$ is isolate-free. For notational convenience, let $G' = G_{1,2}$ and let $X_{1,2} = \{u_1,u_2,v_1,v_2,x_1,x_2,w_1,w_2\}$, and so $G' = G - X_{1,2}$. Every $i$-set of $G'$ can be extended to an ID-set of $G$ by adding to it the vertices $u_1$ and $v_2$ (as indicated by the shaded vertices in Figure~\ref{fig:G12}), and so $i(G) \le i(G') + 2$. If $b(G') =  0$, then by Claim~\ref{structural-weight}, $\Theta(G') = 0$, and so $\Omega(G') = \w(G') = \w(G) - 24 + 8 = \w(G) - 16 = \Omega(G) - 16$. In this case, $8i(G) \le 8(i(G') + 2) \le \w(G') + 16 = \Omega(G)$, a contradiction. Hence, $b(G') \ge 1$. Let $B'$ be a bad component in $G'$, and so $B' \in \cB$. We note that $\deg_{G'}(v_3) = 1$, and so $G'$ contains at least one vertex of degree~$1$.

\begin{figure}[htb]
\begin{center}
\begin{tikzpicture}[scale=.75,style=thick,x=0.8cm,y=0.8cm]
\def\vr{2.5pt} 
\path (7,-1) coordinate (d1);
\path (6,-0.5) coordinate (d2);
\path (7,-0.5) coordinate (d3);
\path (9,0) coordinate (d31);
\path (8,-0.5) coordinate (d4);
\path (9,-0.5) coordinate (d41);
\path (7,0) coordinate (d5);
\path (7,1) coordinate (e1);
\path (6,1.5) coordinate (e2);
\path (7,1.5) coordinate (e3);
\path (9,2) coordinate (e31);
\path (8,1.5) coordinate (e4);
\path (9,1.5) coordinate (e41);
\path (7,2) coordinate (e5);
\path (7,3) coordinate (f1);
\path (6,3.5) coordinate (f2);
\path (7,3.5) coordinate (f3);
\path (9,4) coordinate (f31);
\path (8,3.5) coordinate (f4);
\path (9,3.5) coordinate (f41);
\path (7,4) coordinate (f5);
\path (4.5,1.5) coordinate (w);
\path (10,-0.5) coordinate (z3);
\path (15,-0.5) coordinate (z4);
\path (15,-1.75) coordinate (z5);
\path (10,-1.75) coordinate (z6);
\path (14,3) coordinate (u3);
\path (12.5,3) coordinate (v3);
\path (11.25,2.5) coordinate (p21);
\path (11.25,3.5) coordinate (p22);
\path (10,2.25) coordinate (a1);
\path (10,2.75) coordinate (a2);
\path (10,3.25) coordinate (b1);
\path (10,3.75) coordinate (b2);
\path (10,1) coordinate (q21);
\path (10,1.5) coordinate (q22);
\path (11.35,-1.25) coordinate (y21);
\path (11.35,-0.75) coordinate (y22);
\path (12.5,-1) coordinate (x1);
\path (12.6,-1) coordinate (x1p);
%
\path (11.25,1.25) coordinate (v2);
\path (11.3,1.25) coordinate (v2p);
\path (14,-1) coordinate (u1);
\path (14,1.25) coordinate (u2);
\path (15.5,0.625) coordinate (w1);
\path (15.5,1.875) coordinate (w2);
\draw (x1)--(u1);
\draw (v2)--(u2);
\draw (v3)--(u3);
\draw (w1)--(u1)--(w2);
\draw (w1)--(u2)--(w2);
\draw (w1)--(u3)--(w2);
\draw (y21)--(x1)--(y22);
\draw (p21)--(v3)--(p22);
\draw (q21)--(v2)--(q22);
\draw (a1)--(p21)--(a2);
\draw (b1)--(p22)--(b2);
\draw [style=dashed,rounded corners] (3,-1.75) rectangle (8.5,4.5);
\draw (3.65,-1) node {{\small $B'$}};
\draw [style=dashed,rounded corners] (9.5,0) rectangle (16.5,4.5);
\draw (15.5,4) node {{\small $X_{1,2}$}};
\draw [style=dashed,rounded corners] (9.5,-1.75) rectangle (15.5,-0.5);
%
\draw (10.5,-1) node {{\small $G^*$}};
\draw (d1)--(d2)--(d5)--(d4)--(d1);
\draw (d1)--(d3)--(d5);
\draw (e1)--(e2)--(e5)--(e4)--(e1);
\draw (e1)--(e3)--(e5);
\draw (f1)--(f2)--(f5)--(f4)--(f1);
\draw (f1)--(f3)--(f5);
\draw (d2)--(w)--(e2);
\draw (f2)--(w);
\draw (d4)--(d41);
\draw (e4)--(e41);
\draw (f4)--(f41);
\draw (d3) to[out=45,in=180, distance=0.5cm ] (d31);
\draw (e3) to[out=45,in=180, distance=0.5cm ] (e31);
\draw (f3) to[out=45,in=180, distance=0.5cm ] (f31);
\draw (x1) [fill=white] circle (\vr);
\draw (v2) [fill=black] circle (\vr);
\draw (v3) [fill=white] circle (\vr);
\draw (u1) [fill=white] circle (\vr);
\draw (u2) [fill=white] circle (\vr);
\draw (u3) [fill=black] circle (\vr);
\draw (w1) [fill=white] circle (\vr);
\draw (w2) [fill=white] circle (\vr);
\draw (p21) [fill=black] circle (\vr);
\draw (p22) [fill=black] circle (\vr);
\draw (d1) [fill=white] circle (\vr);
\draw (d2) [fill=black] circle (\vr);
\draw (d3) [fill=white] circle (\vr);
\draw (d4) [fill=white] circle (\vr);
\draw (d5) [fill=white] circle (\vr);
\draw (e1) [fill=white] circle (\vr);
\draw (e2) [fill=black] circle (\vr);
\draw (e3) [fill=white] circle (\vr);
\draw (e4) [fill=white] circle (\vr);
\draw (e5) [fill=white] circle (\vr);
\draw (f1) [fill=white] circle (\vr);
\draw (f2) [fill=black] circle (\vr);
\draw (f3) [fill=white] circle (\vr);
\draw (f4) [fill=white] circle (\vr);
\draw (f5) [fill=white] circle (\vr);
\draw (w) [fill=white] circle (\vr);

%
\draw[anchor = north] (x1p) node {{\small $u_3$}};
\draw[anchor = north] (v2p) node {{\small $u_2$}};
\draw[anchor = south] (v3) node {{\small $u_1$}};
\draw[anchor = north] (u1) node {{\small $v_3$}};
\draw[anchor = north] (u2) node {{\small $v_2$}};
\draw[anchor = south] (u3) node {{\small $v_1$}};
\draw[anchor = west] (w1) node {{\small $x_1$}};
\draw[anchor = west] (w2) node {{\small $x_2$}};
\draw[anchor = north] (p21) node {{\small $w_1$}};
\draw[anchor = south] (p22) node {{\small $w_2$}};
\draw[anchor = east] (w) node {{\small $v$}};
\end{tikzpicture}
\end{center}
\begin{center}
\vskip - 0.25 cm
\caption{A subgraph of $G$ in the proof of Claim~\ref{claim:no-B3}}
\label{fig:G12b}
\end{center}
\end{figure}
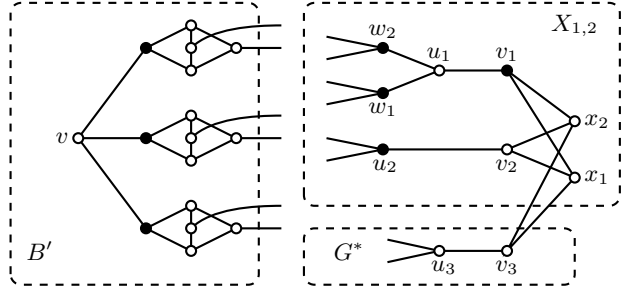

\newpage
\begin{claim}
\label{claim:no-B3}
$B' \notin \cB_3$.
\end{claim}
\proof Suppose, to the contrary, that $B' \in \cB_3$, and so $B'$ is as illustrated in Figure~\ref{fig:bad-comp-structure}(c). In this case, the six $X_{1,2}$-exit edges from $w_1$, $w_2$ and $u_2$ are all incident with vertices of degree~$2$ in $B'$. Let $v$ denote the vertex of degree~$3$ in $B'$ that does not belong to a copy of $K_{2,3}$. We now let $X^* = V(B') \cup X_{1,2}$ and we let $G^* = G - X^*$. We note that $G^*$ contains one vertex of degree~$1$, namely $v_3$, and all other vertices of $G^*$ have degree~$3$. Hence, $b(G^*) = 0$, and so, by Claim~\ref{structural-weight}, $\Theta(G^*) = 0$. Thus, $\Omega(G^*) = \w(G^*) = \w(G) - 24 \times 3 + 2 = \w(G) - 70$. Every $i$-set of $G^*$ can be extended to an ID-set of $G$ by adding to it the seven vertices in the set $N(v) \cup \{u_2,v_1,w_1,w_2\}$  (as indicated by the shaded vertices in Figure~\ref{fig:G12b}), and so $i(G) \le i(G^*) + 7$. Thus, $8i(G) \le 8(i(G^*) + 7) \le \Omega(G^*) + 56 < \Omega(G) $, a contradiction.~\smallqed

\begin{claim}
\label{claim:no-B2}
$B' \notin \cB_2$.
\end{claim}
\proof Suppose, to the contrary, that $B' \in \cB_2$, and so $B'$ is as illustrated in Figure~\ref{fig:bad-comp-structure}(b). In this case, five $X_{1,2}$-exit edges from $w_1$, $w_2$ and $u_2$ are all incident with vertices of degree~$2$ in $B'$. Let $v$ denote the vertex of degree~$2$ in $B'$ that does not belong to a copy of $K_{2,3}$. We now let $X^* = V(B') \cup X_{1,2}$ and we let $G^* = G - X^*$. We note that $G^*$ contains one vertex of degree~$1$, one vertex of degree~$2$, and all other vertices of $G^*$ have degree~$3$. Hence, $b(G^*) = 0$, and so, by Claim~\ref{structural-weight}, we have $\Theta(G^*) = 0$. Thus, $\Omega(G^*) = \w(G^*) = \w(G) - 19 \times 3 + 3 = \w(G) - 54 = \Omega(G) - 54$. We show that every $i$-set of $G^*$ can be extended to an ID-set of $G$ by adding to it at most six vertices. We will use the property given in Claim~\ref{claim:K23-config} that every $K_{2,3}$-subgraph of $G$ belongs to a $G_{8.2}$-configuration that is an induced subgraph in $G$, where $G_{8.2}$ is the graph in Figure~\ref{fig:K23-config}(b).

\begin{figure}[htb]
\begin{center}
\begin{tikzpicture}[scale=.75,style=thick,x=0.75cm,y=0.75cm]
\def\vr{2.5pt} 
%
\path (6,1) coordinate (d1);
\path (5,1.5) coordinate (d2);
\path (6,1.5) coordinate (d3);
\path (8,2) coordinate (d31);
\path (7,1.5) coordinate (d4);
\path (8,1.5) coordinate (d41);
\path (6,2) coordinate (d5);
\path (6,3) coordinate (e1);
\path (5,3.5) coordinate (e2);
\path (6,3.5) coordinate (e3);
\path (8,4) coordinate (e31);
\path (7,3.5) coordinate (e4);
\path (8,3.5) coordinate (e41);
\path (6,4) coordinate (e5);
\path (3.5,2.5) coordinate (w);

\path (7,-0.25) coordinate (z3);
\path (14,-0.25) coordinate (z4);
\path (14,-1.5) coordinate (z5);
\path (7,-1.5) coordinate (z6);
\path (13,3) coordinate (u3);
\path (11.5,3) coordinate (v3);
\path (10.25,2.5) coordinate (p21);
\path (8,-0.75) coordinate (z1);
\path (7.25,-0.5) coordinate (z11);
\path (7.25,-1) coordinate (z12);
\path (10.25,3.5) coordinate (p22);
\path (9,2.25) coordinate (a1);
\path (9,2.75) coordinate (a2);
\path (9,3.25) coordinate (b1);
\path (9,3.75) coordinate (b2);
\path (9,1) coordinate (q21);
\path (9,1.5) coordinate (q22);
\path (10.35,-1) coordinate (y21);
\path (10.35,-0.5) coordinate (y22);
\path (11.5,-0.75) coordinate (x1);
\path (11.6,-0.75) coordinate (x1p);
%
\path (10.25,1.25) coordinate (v2);
\path (10.3,1.25) coordinate (v2p);
\path (13,-0.75) coordinate (u1);
\path (13,1.25) coordinate (u2);
\path (14.5,0.625) coordinate (w1);
\path (14.5,1.875) coordinate (w2);
\draw (x1)--(u1);
\draw (v2)--(u2);
\draw (v3)--(u3);
\draw (w1)--(u1)--(w2);
\draw (w1)--(u2)--(w2);
\draw (w1)--(u3)--(w2);
\draw (y21)--(x1)--(y22);
\draw (p21)--(v3)--(p22);
\draw (z11)--(z1)--(z12);
\draw [style=dashed,rounded corners] (2,0.5) rectangle (7.5,4.5);
\draw (2.65,4) node {{\small $B'$}};
\draw [style=dashed,rounded corners] (8.5,0.25) rectangle (15.5,4.5);
\draw (14.5,4) node {{\small $X_{1,2}$}};
\draw[densely dashed]  (z3)--(z4)--(z5)--(z6)--(z3);
\draw (9.5,-0.75) node {{\small $G^*$}};
\draw (d1)--(d2)--(d5)--(d4)--(d1);
\draw (d1)--(d3)--(d5);
\draw (e1)--(e2)--(e5)--(e4)--(e1);
\draw (e1)--(e3)--(e5);
\draw (d2)--(w)--(e2);
%
%
\draw (z1) to[out=90,in=180, distance=1.5cm ] (p21);
\draw (w) to[out=90,in=90, distance=2cm ] (p22);
\draw (v2) to[out=90,in=0, distance=0.5cm ] (d4);
\draw (v2) to[out=90,in=0, distance=0.5cm ] (e4);
\draw (d3) to[out=45,in=135, distance=0.65cm ] (p21);
\draw (e3) to[out=45,in=135, distance=0.65cm ] (p22);
%
\draw (x1) [fill=white] circle (\vr);
\draw (z1) [fill=white] circle (\vr);
\draw (v2) [fill=black] circle (\vr);
\draw (v3) [fill=white] circle (\vr);
\draw (u1) [fill=white] circle (\vr);
\draw (u2) [fill=white] circle (\vr);
\draw (u3) [fill=black] circle (\vr);
\draw (w1) [fill=white] circle (\vr);
\draw (w2) [fill=white] circle (\vr);
\draw (p21) [fill=white] circle (\vr);
\draw (p22) [fill=white] circle (\vr);
\draw (d1) [fill=white] circle (\vr);
\draw (d2) [fill=white] circle (\vr);
\draw (d3) [fill=black] circle (\vr);
\draw (d4) [fill=white] circle (\vr);
\draw (d5) [fill=white] circle (\vr);
\draw (e1) [fill=white] circle (\vr);
\draw (e2) [fill=white] circle (\vr);
\draw (e3) [fill=black] circle (\vr);
\draw (e4) [fill=white] circle (\vr);
\draw (e5) [fill=white] circle (\vr);
\draw (w) [fill=black] circle (\vr);
%
\draw[anchor = north] (x1p) node {{\small $u_3$}};
\draw[anchor = north] (v2p) node {{\small $u_2$}};
\draw[anchor = south] (v3) node {{\small $u_1$}};
\draw[anchor = north] (u1) node {{\small $v_3$}};
\draw[anchor = north] (u2) node {{\small $v_2$}};
\draw[anchor = south] (u3) node {{\small $v_1$}};
\draw[anchor = west] (w1) node {{\small $x_1$}};
\draw[anchor = west] (w2) node {{\small $x_2$}};
\draw[anchor = north] (p21) node {{\small $w_1$}};
\draw[anchor = north] (p22) node {{\small $w_2$}};
\draw[anchor = east] (w) node {{\small $v$}};
0%
\draw (10,-2.5) node {{\small (a)}};
%
%
\path (20,1) coordinate (d1);
\path (19,1.5) coordinate (d2);
\path (20,1.5) coordinate (d3);
\path (22,2) coordinate (d31);
\path (21,1.5) coordinate (d4);
\path (22,1.5) coordinate (d41);
\path (20,2) coordinate (d5);
\path (20,3) coordinate (e1);
\path (19,3.5) coordinate (e2);
\path (20,3.5) coordinate (e3);
\path (22,4) coordinate (e31);
\path (21,3.5) coordinate (e4);
\path (22,3.5) coordinate (e41);
\path (20,4) coordinate (e5);
\path (17.5,2.5) coordinate (w);

\path (21,-0.25) coordinate (z3);
\path (28,-0.25) coordinate (z4);
\path (28,-1.5) coordinate (z5);
\path (21,-1.5) coordinate (z6);
\path (27,3) coordinate (u3);
\path (25.5,3) coordinate (v3);
\path (24.25,2.5) coordinate (p21);
\path (22,-0.75) coordinate (z1);
\path (21.25,-0.5) coordinate (z11);
\path (21.25,-1) coordinate (z12);
\path (24.25,3.5) coordinate (p22);
\path (23,2.25) coordinate (a1);
\path (23,2.75) coordinate (a2);
\path (23,3.25) coordinate (b1);
\path (23,3.75) coordinate (b2);
\path (23,1) coordinate (q21);
\path (23,1.5) coordinate (q22);
\path (24.35,-1) coordinate (y21);
\path (24.35,-0.5) coordinate (y22);
\path (25.5,-0.75) coordinate (x1);
\path (25.6,-0.75) coordinate (x1p);
%
\path (24.25,1.25) coordinate (v2);
\path (24.4,1.25) coordinate (v2p);
\path (27,-0.75) coordinate (u1);
\path (27,1.25) coordinate (u2);
\path (28.5,0.625) coordinate (w1);
\path (28.5,1.875) coordinate (w2);
\draw (x1)--(u1);
\draw (v2)--(u2);
\draw (v3)--(u3);
\draw (w1)--(u1)--(w2);
\draw (w1)--(u2)--(w2);
\draw (w1)--(u3)--(w2);
\draw (y21)--(x1)--(y22);
\draw (p21)--(v3)--(p22);
\draw (z11)--(z1)--(z12);
%
%
%
\draw[densely dashed]  (z3)--(z4)--(z5)--(z6)--(z3);
\draw (23.5,-0.75) node {{\small $G^*$}};
\draw [style=dashed,rounded corners] (22.5,0.25) rectangle (29.5,4.5);
\draw (28.5,4) node {{\small $X_{1,2}$}};
\draw [style=dashed,rounded corners] (16.5,0.5) rectangle (21.5,4.5);
\draw (17,4) node {{\small $B'$}};
\draw (d1)--(d2)--(d5)--(d4)--(d1);
\draw (d1)--(d3)--(d5);
\draw (e1)--(e2)--(e5)--(e4)--(e1);
\draw (e1)--(e3)--(e5);
\draw (d2)--(w)--(e2);
%
%
\draw (z1) to[out=90,in=200, distance=1.5cm ] (p21);
\draw (w) to[out=270,in=240, distance=1.5cm ] (v2);
\draw (p21) to[out=140,in=90, distance=0.5cm ] (d4);
\draw (v2) to[out=90,in=0, distance=0.5cm ] (e4);
\draw (d3) to[out=45,in=180, distance=0.65cm ] (p22);
\draw (e3) to[out=45,in=130, distance=0.65cm ] (p22);
%
\draw (x1) [fill=white] circle (\vr);
\draw (z1) [fill=white] circle (\vr);
\draw (v2) [fill=black] circle (\vr);
\draw (v3) [fill=white] circle (\vr);
\draw (u1) [fill=white] circle (\vr);
\draw (u2) [fill=white] circle (\vr);
\draw (u3) [fill=black] circle (\vr);
\draw (w1) [fill=white] circle (\vr);
\draw (w2) [fill=white] circle (\vr);
\draw (p21) [fill=white] circle (\vr);
\draw (p22) [fill=black] circle (\vr);
\draw (d1) [fill=white] circle (\vr);
\draw (d2) [fill=black] circle (\vr);
\draw (d3) [fill=white] circle (\vr);
\draw (d4) [fill=black] circle (\vr);
\draw (d5) [fill=white] circle (\vr);
\draw (e1) [fill=white] circle (\vr);
\draw (e2) [fill=black] circle (\vr);
\draw (e3) [fill=white] circle (\vr);
\draw (e4) [fill=white] circle (\vr);
\draw (e5) [fill=white] circle (\vr);
\draw (w) [fill=white] circle (\vr);

%

%
\draw[anchor = north] (x1p) node {{\small $u_3$}};
\draw[anchor = north] (v2p) node {{\small $u_2$}};
\draw[anchor = south] (v3) node {{\small $u_1$}};
\draw[anchor = north] (u1) node {{\small $v_3$}};
\draw[anchor = north] (u2) node {{\small $v_2$}};
\draw[anchor = south] (u3) node {{\small $v_1$}};
\draw[anchor = west] (w1) node {{\small $x_1$}};
\draw[anchor = west] (w2) node {{\small $x_2$}};
\draw[anchor = north] (p21) node {{\small $w_1$}};
\draw[anchor = north] (p22) node {{\small $w_2$}};
\draw[anchor = east] (w) node {{\small $v$}};
%
%
%
%
\draw (24,-2.5) node {{\small (b)}};
%
\end{tikzpicture}
\end{center}
\begin{center}
\vskip - 0.65 cm
\caption{Subgraphs of $G$ in the proof of Claim~\ref{claim:no-B2}}
\label{fig:G12c}
\end{center}
\end{figure}
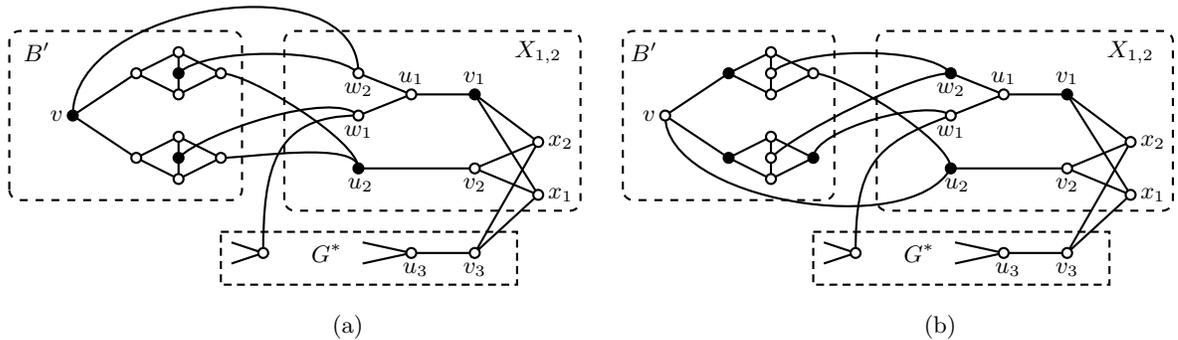

Suppose firstly that $u_2$ is adjacent to two vertices in $B'$. If $u_2$ is not adjacent to~$v$, then we let $I'$ consist of the three vertices of degree~$2$ in $B'$ that are not adjacent to~$u_2$, together with the vertices~$v_1$ and~$u_2$ (as indicated by the shaded vertices in Figure~\ref{fig:G12c}(a)). If $u_2$ is adjacent to~$v$, then renaming $w_1$ and $w_2$ if necessary, we may assume that $w_2$ is adjacent to two vertices in $B'$. In this case, we let $I'$ consist of the two neighbors of $v$ in $B'$, the vertex of degree~$2$ in $B'$ adjacent to~$w_1$, together with the vertices in the set $\{u_2,v_1,w_2\}$ (as indicated by the shaded vertices in Figure~\ref{fig:G12c}(b)). In both cases, every $i$-set of $G^*$ can be extended to an ID-set of $G$ by adding to the vertices in the set~$I'$, and so $i(G) \le i(G^*) + |I'| \le i(G^*) + 6$.

\begin{figure}[htb]
\begin{center}
\begin{tikzpicture}[scale=.75,style=thick,x=0.75cm,y=0.75cm]
\def\vr{2.5pt} 
%
\path (6,1) coordinate (d1);
\path (5,1.5) coordinate (d2);
\path (6,1.5) coordinate (d3);
\path (8,2) coordinate (d31);
\path (7,1.5) coordinate (d4);
\path (8,1.5) coordinate (d41);
\path (6,2) coordinate (d5);
\path (6,3) coordinate (e1);
\path (5,3.5) coordinate (e2);
\path (6,3.5) coordinate (e3);
\path (8,4) coordinate (e31);
\path (7,3.5) coordinate (e4);
\path (8,3.5) coordinate (e41);
\path (6,4) coordinate (e5);
\path (3.5,2.5) coordinate (w);

\path (7,-0.25) coordinate (z3);
\path (14,-0.25) coordinate (z4);
\path (14,-1.5) coordinate (z5);
\path (7,-1.5) coordinate (z6);
\path (13,3) coordinate (u3);
\path (11.5,3) coordinate (v3);
\path (10.25,2.5) coordinate (p21);
\path (8,-0.75) coordinate (z1);
\path (7.25,-0.5) coordinate (z11);
\path (7.25,-1) coordinate (z12);
\path (10.25,3.5) coordinate (p22);
\path (9,2.25) coordinate (a1);
\path (9,2.75) coordinate (a2);
\path (9,3.25) coordinate (b1);
\path (9,3.75) coordinate (b2);
\path (9,1) coordinate (q21);
\path (9,1.5) coordinate (q22);
\path (10.35,-1) coordinate (y21);
\path (10.35,-0.5) coordinate (y22);
\path (11.5,-0.75) coordinate (x1);
\path (11.6,-0.75) coordinate (x1p);
%
\path (10.25,1.25) coordinate (v2);
\path (10.3,1.25) coordinate (v2p);
\path (13,-0.75) coordinate (u1);
\path (13,1.25) coordinate (u2);
\path (14.5,0.625) coordinate (w1);
\path (14.5,1.875) coordinate (w2);
\draw (x1)--(u1);
\draw (v2)--(u2);
\draw (v3)--(u3);
\draw (w1)--(u1)--(w2);
\draw (w1)--(u2)--(w2);
\draw (w1)--(u3)--(w2);
\draw (y21)--(x1)--(y22);
\draw (p21)--(v3)--(p22);
\draw (z11)--(z1)--(z12);
%
%
%
\draw[densely dashed]  (z3)--(z4)--(z5)--(z6)--(z3);
%
%
\draw (d1)--(d2)--(d5)--(d4)--(d1);
\draw (d1)--(d3)--(d5);
\draw (e1)--(e2)--(e5)--(e4)--(e1);
\draw (e1)--(e3)--(e5);
\draw (d2)--(w)--(e2);
%
%
\draw (z1) to[out=70,in=180, distance=0.75cm ] (v2);
\draw (w) to[out=90,in=90, distance=2cm ] (p22);
\draw (e4) to[out=0,in=90, distance=0.5cm ] (p21);
\draw (v2) to[out=145,in=45, distance=0.5cm ] (d4);
\draw (d3) to[out=45,in=135, distance=0.65cm ] (p21);
\draw (e3) to[out=45,in=135, distance=0.65cm ] (p22);
%
\draw (x1) [fill=white] circle (\vr);
\draw (z1) [fill=white] circle (\vr);
\draw (v2) [fill=white] circle (\vr);
\draw (v3) [fill=white] circle (\vr);
\draw (u1) [fill=white] circle (\vr);
\draw (u2) [fill=black] circle (\vr);
\draw (u3) [fill=black] circle (\vr);
\draw (w1) [fill=white] circle (\vr);
\draw (w2) [fill=white] circle (\vr);
\draw (p21) [fill=black] circle (\vr);
\draw (p22) [fill=white] circle (\vr);
\draw (d1) [fill=white] circle (\vr);
\draw (d2) [fill=white] circle (\vr);
\draw (d3) [fill=white] circle (\vr);
\draw (d4) [fill=black] circle (\vr);
\draw (d5) [fill=white] circle (\vr);
\draw (e1) [fill=white] circle (\vr);
\draw (e2) [fill=white] circle (\vr);
\draw (e3) [fill=black] circle (\vr);
\draw (e4) [fill=white] circle (\vr);
\draw (e5) [fill=white] circle (\vr);
\draw (w) [fill=black] circle (\vr);

%
\draw[anchor = north] (x1p) node {{\small $u_3$}};
\draw[anchor = north] (v2p) node {{\small $u_2$}};
\draw[anchor = south] (v3) node {{\small $u_1$}};
\draw[anchor = north] (u1) node {{\small $v_3$}};
\draw[anchor = north] (u2) node {{\small $v_2$}};
\draw[anchor = south] (u3) node {{\small $v_1$}};
\draw[anchor = west] (w1) node {{\small $x_1$}};
\draw[anchor = west] (w2) node {{\small $x_2$}};
\draw[anchor = north] (p21) node {{\small $w_1$}};
\draw[anchor = north] (p22) node {{\small $w_2$}};
\draw[anchor = east] (w) node {{\small $v$}};
%
\draw [style=dashed,rounded corners] (2,0.5) rectangle (7.5,4.5);
\draw (2.65,4) node {{\small $B'$}};
\draw [style=dashed,rounded corners] (8.5,0.25) rectangle (15.5,4.5);
\draw (14.5,4) node {{\small $X_{1,2}$}};
\draw[densely dashed]  (z3)--(z4)--(z5)--(z6)--(z3);
\draw (9.5,-0.75) node {{\small $G^*$}};
%
\draw (10,-2.5) node {{\small (a)}};
%
%
\path (20,1) coordinate (d1);
\path (19,1.5) coordinate (d2);
\path (20,1.5) coordinate (d3);
\path (22,2) coordinate (d31);
\path (21,1.5) coordinate (d4);
\path (22,1.5) coordinate (d41);
\path (20,2) coordinate (d5);
\path (20,3) coordinate (e1);
\path (19,3.5) coordinate (e2);
\path (20,3.5) coordinate (e3);
\path (22,4) coordinate (e31);
\path (21,3.5) coordinate (e4);
\path (22,3.5) coordinate (e41);
\path (20,4) coordinate (e5);
\path (17.5,2.5) coordinate (w);

\path (21,-0.25) coordinate (z3);
\path (28,-0.25) coordinate (z4);
\path (28,-1.5) coordinate (z5);
\path (21,-1.5) coordinate (z6);
\path (27,3) coordinate (u3);
\path (25.5,3) coordinate (v3);
\path (24.25,2.5) coordinate (p21);
\path (22,-0.75) coordinate (z1);
\path (21.25,-0.5) coordinate (z11);
\path (21.25,-1) coordinate (z12);
\path (24.25,3.5) coordinate (p22);
\path (23,2.25) coordinate (a1);
\path (23,2.75) coordinate (a2);
\path (23,3.25) coordinate (b1);
\path (23,3.75) coordinate (b2);
\path (23,1) coordinate (q21);
\path (23,1.5) coordinate (q22);
\path (24.35,-1) coordinate (y21);
\path (24.35,-0.5) coordinate (y22);
\path (25.5,-0.75) coordinate (x1);
\path (25.6,-0.75) coordinate (x1p);
%
\path (24.25,1.25) coordinate (v2);
\path (24.4,1.25) coordinate (v2p);
\path (27,-0.75) coordinate (u1);
\path (27,1.25) coordinate (u2);
\path (28.5,0.625) coordinate (w1);
\path (28.5,1.875) coordinate (w2);
\draw (x1)--(u1);
\draw (v2)--(u2);
\draw (v3)--(u3);
\draw (w1)--(u1)--(w2);
\draw (w1)--(u2)--(w2);
\draw (w1)--(u3)--(w2);
\draw (y21)--(x1)--(y22);
\draw (p21)--(v3)--(p22);
\draw (z11)--(z1)--(z12);
%
\draw[densely dashed]  (z3)--(z4)--(z5)--(z6)--(z3);
\draw (23.5,-0.75) node {{\small $G^*$}};
\draw [style=dashed,rounded corners] (22.5,0.25) rectangle (29.5,4.5);
\draw (28.5,4) node {{\small $X_{1,2}$}};
\draw [style=dashed,rounded corners] (16.5,0.5) rectangle (21.5,4.5);
\draw (17,4) node {{\small $B'$}};
%
%
%
\draw (d1)--(d2)--(d5)--(d4)--(d1);
\draw (d1)--(d3)--(d5);
\draw (e1)--(e2)--(e5)--(e4)--(e1);
\draw (e1)--(e3)--(e5);
\draw (d2)--(w)--(e2);
%
%
\draw (e4) to[out=0,in=90, distance=0.5cm ] (p21);
\draw (z1) to[out=0,in=270, distance=0.75cm ] (v2);
\draw (w) to[out=270,in=240, distance=1.5cm ] (v2);
\draw (p21) to[out=140,in=90, distance=0.5cm ] (d4);
\draw (d3) to[out=45,in=180, distance=0.65cm ] (p22);
\draw (e3) to[out=45,in=130, distance=0.65cm ] (p22);
%
\draw (x1) [fill=white] circle (\vr);
\draw (z1) [fill=white] circle (\vr);
\draw (v2) [fill=white] circle (\vr);
\draw (v3) [fill=white] circle (\vr);
\draw (u1) [fill=white] circle (\vr);
\draw (u2) [fill=black] circle (\vr);
\draw (u3) [fill=black] circle (\vr);
\draw (w1) [fill=white] circle (\vr);
\draw (w2) [fill=white] circle (\vr);
\draw (p21) [fill=black] circle (\vr);
\draw (p22) [fill=black] circle (\vr);
\draw (d1) [fill=white] circle (\vr);
\draw (d2) [fill=black] circle (\vr);
\draw (d3) [fill=white] circle (\vr);
\draw (d4) [fill=white] circle (\vr);
\draw (d5) [fill=white] circle (\vr);
\draw (e1) [fill=white] circle (\vr);
\draw (e2) [fill=black] circle (\vr);
\draw (e3) [fill=white] circle (\vr);
\draw (e4) [fill=white] circle (\vr);
\draw (e5) [fill=white] circle (\vr);
\draw (w) [fill=white] circle (\vr);

%
\draw[anchor = north] (x1p) node {{\small $u_3$}};
\draw[anchor = east] (v2) node {{\small $u_2$}};
\draw[anchor = south] (v3) node {{\small $u_1$}};
\draw[anchor = north] (u1) node {{\small $v_3$}};
\draw[anchor = north] (u2) node {{\small $v_2$}};
\draw[anchor = south] (u3) node {{\small $v_1$}};
\draw[anchor = west] (w1) node {{\small $x_1$}};
\draw[anchor = west] (w2) node {{\small $x_2$}};
\draw[anchor = north] (p21) node {{\small $w_1$}};
\draw[anchor = south] (p22) node {{\small $w_2$}};
\draw[anchor = east] (w) node {{\small $v$}};
%
%
%
%
\draw (24,-2.5) node {{\small (b)}};
%
\end{tikzpicture}
\end{center}
\begin{center}
\vskip - 0.65 cm
\caption{Subgraphs of $G$ in the proof of Claim~\ref{claim:no-B2}}
\label{fig:G12d}
\end{center}
\end{figure}
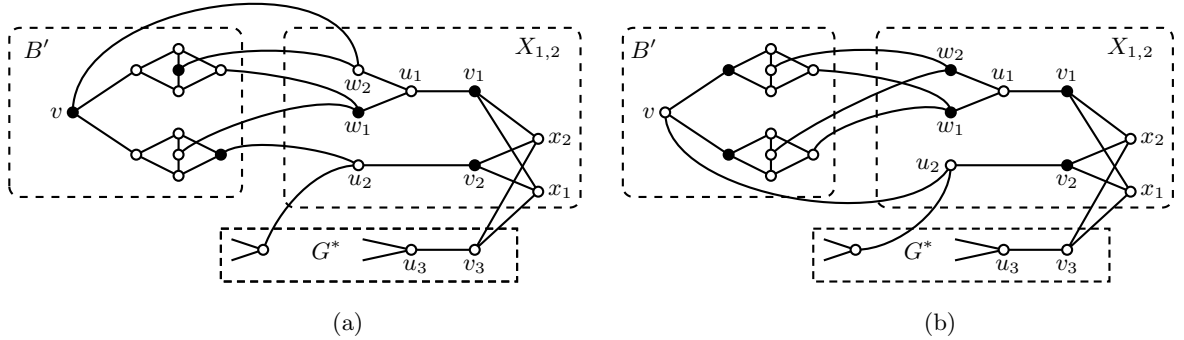

Suppose next that $u_2$ is adjacent to exactly one vertex in $B'$. Suppose that $u_2$ is not adjacent to~$v$.  Renaming $w_1$ and $w_2$ if necessary, we may assume that $w_2$ is adjacent to $v$. We now let $I'$ consist of the vertex~$v$, the vertices of degree~$2 $ in $B'$ adjacent to~$w_2$ and~$u_2$, and the vertices in the set $\{v_1,v_2,w_1\}$ (as indicated by the shaded vertices in Figure~\ref{fig:G12d}(a)). If $u_2$ is adjacent to~$v$, then we let $I'$ consist of the two neighbors of $v$ in $B'$ together with the vertices in the set $\{v_1,v_2,w_1,w_2\}$ (as indicated by the shaded vertices in Figure~\ref{fig:G12d}(b)). In both cases, every $i$-set of $G^*$ can be extended to an ID-set of $G$ by adding to the vertices in the set~$I'$, and so $i(G) \le i(G^*) + |I'| = i(G^*) + 6$.

As observed earlier, $\Omega(G^*) = \Omega(G) - 54$. We have shown that every $i$-set of $G^*$ can be extended to an ID-set of $G$ by adding to it at most six vertices from~$X^*$, and so $i(G) \le i(G^*) + 6$. Thus, $8i(G) \le 8(i(G^*) + 6) \le \Omega(G^*) + 48 < \Omega(G)$, a contradiction.~\smallqed

\medskip
By Claims~\ref{claim:no-B3} and~\ref{claim:no-B2}, we have $B' \in \cB_1$. Suppose that $v_3 \in V(B')$. Since $v_3$ has degree~$1$ in $B'$, by Claim~\ref{claim:bad-comp-structure} such a component $B' \in \cB_1$ is as illustrated in Figure~\ref{fig:bad-comp-structure}(a). In this case, the graph $G$ contains two vertex disjoint copies of $K_{2,3}$ joined by the edge $u_3v_3$, contradicting Claim~\ref{claim:no-adj-K23}. Hence, $v_3 \notin V(B')$.

Let $y_1$ and $y_2$ be the two vertices of degree~$3$ in $B'$ that are not adjacent, and let $N(y_1) = N(y_2) = \{a,b,c\}$. Further, let $a'$ be the vertex of degree~$1$ in $B'$, where $aa'$ is an edge. Let $b'$ and $c'$ be the neighbors of $b$ and $c$, respectively, that do not belong to~$B'$. By Claim~\ref{claim:K23-config}, the subgraph of $G$ induced by the set $V(B') \cup \{b',c'\}$ is a $G_{8.2}$-configuration.
In particular, we note that the set $\{a',b',c'\}$ is an independent set. Since $v_3 \notin V(B')$, we note that $a' \ne v_3$. The two neighbors of $a'$ that are not in $B'$ belong to the set $\{u_2,w_1,w_2\}$. Further, the vertices $b'$ and $c'$ belong to the set $\{u_2,w_1,w_2\}$, that is, $\{b',c'\} \subset \{u_2,w_1,w_2\}$. Since the vertex~$a'$ is adjacent to two vertices that belong to the set $\{u_2,w_1,w_2\}$, the vertex $a'$ is adjacent to at least one of~$b'$ and $c'$, contradicting our earlier observation that $\{a',b',c'\}$ is an independent set. This final contradiction completes our proof of Theorem~\ref{thm:main2}.~\QED

\section*{Acknowledgments}

The first and the second author were supported by the Slovenian Research and Innovation agency (grants P1-0297, N1-0285, and N1-0431). Research of the third author was supported in part by the University of Johannesburg.

\end{document}